\documentclass{conm-p-l}

\usepackage{amsmath,amsthm,amsfonts,amssymb}
\usepackage{graphicx}

\usepackage{amssymb}
\usepackage{psfig}
\usepackage{amsmath}
\usepackage{amssymb}

\newtheorem{lm}{Lemma}[section]
\newtheorem{prop}[lm]{Proposition}
\newtheorem{theorem}[lm]{Theorem}

\newtheorem{cy}[lm]{Corollary}
\newtheorem{Claim}{Claim}

\theoremstyle{definition}
\newtheorem{df}[lm]{Definition}

\newtheorem{rk}[lm]{Remark}
\newtheorem{notation}[lm]{Notation}

\newcommand{\beq}{\begin{equation}}
\newcommand{\eeq}{\end{equation}}
\newcommand{\be}{\begin{enumerate}}
\newcommand{\ee}{\end{enumerate}}
\newcommand{\bp}{\begin{proof}}
\newcommand{\ep}{\end{proof}}
\newcommand{\bi}{\begin{itemize}}
\newcommand{\ei}{\end{itemize}}
\newcommand{\bea}{\begin{eqnarray*}}
\newcommand{\eea}{\end{eqnarray*}}
\newcommand{\bml}{\begin{multline*}}
\newcommand{\eml}{\end{multline*}}

\usepackage{amssymb}
\usepackage{psfig}
\usepackage{amsmath}
\usepackage{amssymb}

\edef\storedcatcodeat{\the\catcode`\@} \catcode`\@=11
\def\ig#1#2#3{{\setlength\arraycolsep{0pt}\def\arraystretch{0.8}
\setbox250\vbox{\strut}
\begin{array}[t]{|lcr|}
\multicolumn3{|c|}{\: #1 \:} \\ \hline \scriptstyle \: #2 &
\scriptstyle\hspace{1em} & \scriptstyle #3 \:
\end{array}}}
\def\igl#1#2{{\setlength\arraycolsep{0pt}\def\arraystretch{0.8}
\setbox250\vbox{\strut}
\begin{array}[t]{|l}
\multicolumn1{|r}{\: #1} \\[-\ht250]
\rlap{\vbox{\hrule width0.5em}} \strut \\[-\dp250]
\scriptstyle \: #2 \end{array}}}
\def\igr#1#2{{\setlength\arraycolsep{0pt}\def\arraystretch{0.8}
\setbox250\vbox{\strut}
\begin{array}[t]{r|}
\multicolumn1{l|}{\: #1} \\[-\ht250]
\strut \llap{\vbox{\hrule width0.5em}} \\[-\dp250]
\scriptstyle #2 \: \end{array}}}
\catcode`\@=\storedcatcodeat\relax

\begin{document}
\setcounter{page}{213}
\title[Algebraic Geometry over Free Groups]{Algebraic Geometry over Free Groups: Lifting Solutions into Generic Points}
\author[O. Kharalmpovich]{Olga Kharlampovich}
\address{Department of Mathematics and Statistics, McGill University, Montreal, QC, Canada, H3A2K6}

\email{olga@math.mcgill.ca}
\thanks{The first author was supported
 by a NSERC Grant.}

\author[A.~G.~Myasnikov]{Alexei Myasnikov}
\address{Department of Mathematics and Statistics, McGill University, Montreal, QC, Canada, H3A2K6}
\email{alexeim@att.net}
\thanks{The second author was supported by a NSERC Grant and by NSF
GrantDMS-9970618}

\subjclass{Primary 20F10; Secondary 03C05}
\date{ May 21, 2004}

\pagestyle{myheadings} \markright{{\sl $\bullet$ Liftings, version
10 $\bullet$ 05.15.04}}

\date{Version 13, May 21, 2004}

\keywords{Free group, quadratic equation, lifting}

\begin{abstract}

In this paper we prove  Implicit Function Theorems (IFT) for
algebraic varieties defined by  regular quadratic equations and,
more generally,  regular NTQ systems over free groups.  In the
model theoretic language these results state the existence of very
simple Skolem functions for particular $\forall\exists$-formulas
over free groups. We construct these functions effectively.  In
non-effective form IFT  first appeared in \cite{Imp}.
 From algebraic geometry view-point IFT can be  described as lifting
solutions of equations into generic points of algebraic varieties.

Moreover, we show that the converse is also true, i.e., IFT holds
only for algebraic varieties defined by regular NTQ systems. This
implies that if a finitely generated group $H$ is
$\forall\exists$-equivalent to a free non-abelian group then $H$
is  isomorphic to the coordinate group of a regular NTQ system.

\end{abstract}

\maketitle

\tableofcontents

\section*{Introduction}

 The classical algebraic geometry  is one of the main tools to deal
with   polynomial equations over fields. To study solutions of
equations in free groups one needs a similar theory over groups.
Recently  basics of algebraic geometry over groups were developed
in a series of papers \cite{BMR1,KMNull,KMIrc}. This provides the
necessary topological machinery to transcribe  geometric notions
into the language of pure group theory. In this  paper, following
\cite{BMR1} and \cite{KMNull},  we  freely  use the standard
algebraic geometric  notions such as algebraic sets, the Zariski
topology, Noetherian domains, irreducible   varieties, radicals
and coordinate groups to organize an approach to finding a
solution of Tarski's problems in \cite{KM7}. Our goal here is to
prove several  variations of so-called implicit function theorem
 (IFT) for free groups.
 The basic version of IFT was announced at the Model Theory
conference  at MSRI in 1998 \cite{KM4,MSRI}. In \cite{KM5} we used
the basic version of implicit function theorem to solve the genus
problem for quadratic non-orientable equations, and showed also
that the abelianization  of the cartesian power of infinitely many
copies of a free non-abelian group has 2-torsion. The preprint
\cite{Imp} contains proofs of several variations of IFT
 in terms of liftings.

In a  sense some   formulations of  IFT can be viewed as analogs
of the corresponding results from analysis, hence the name. To
demonstrate this we start with  a very basic version of the
implicit function theorem which holds for regular quadratic
equations.

Let $G$ be a group generated by $A$, $F(X)$ be a free group with
basis $X = \{x_1, x_2, \dots, x_n\}$, $G[X] = G \ast F(X)$ be a
free product of $G$ and $F(X)$. If $S \subset G[X]$ then the
expression $S = 1$ is called {\em a system of equations} over $G$.
A solution of the system $S = 1$ over $G$ can be described as a
$G$-homomorphism $\phi : G[X] \longrightarrow G$ such that
$\phi(S) = 1$. By $V_G(S)$ we denote the set of all solutions in
$G$ of the system $ S = 1$, it is called the {\em algebraic set
defined by} $S$. This algebraic set $V_G(S)$ uniquely corresponds
to the radical $R(S)$:
$$ R(S) = \{ T(x) \in G[X] \ \mid \ \forall A\in G^n (S(A) = 1
\rightarrow T(A) = 1 \}. $$  The quotient group
$$G_{R(S)}=G[X]/R(S)$$ is the {\em coordinate group} of the
algebraic set  $V(S).$  Every solution of $S(X) = 1$ in $G$ can be
described as a $G$-homomorphism $G_{R(S)} \rightarrow G$.

Recall that a standard quadratic equation $S(X) = 1$ over group
$G$ is an equation in one of the following forms (below $d,c_i$
are nontrivial elements from $G$):
\begin{equation}\label{eq:st1}
\prod_{i=1}^{n}[x_i,y_i] = 1, \ \ \ n > 0;
\end{equation}
\begin{equation}\label{eq:st2}
\prod_{i=1}^{n}[x_i,y_i] \prod_{i=1}^{m}z_i^{-1}c_iz_i d = 1,\ \ \
n,m \geqslant 0, m+n \geqslant 1 ;
\end{equation}
\begin{equation}\label{eq:st3}
\prod_{i=1}^{n}x_i^2 = 1, \ \ \ n > 0;
\end{equation}
\begin{equation}\label{eq:st4}
\prod_{i=1}^{n}x_i^2 \prod_{i=1}^{m}z_i^{-1}c_iz_i d = 1, \ \ \
n,m \geqslant 0, n+m \geqslant 1.
\end{equation}

Equations (\ref{eq:st1}), (\ref{eq:st2}) are called {\em
orientable} and  equations (\ref{eq:st3}), (\ref{eq:st4}) are
called {\em non-orientable}. The numbers  $n$ and $n+m$ are called
{\em genus} and {\em atomic rank}  of $S(X) = 1$.  Put
 $$\kappa(S) = |X| + \varepsilon(S),$$
  where $\varepsilon(S) = 1$  if
the coefficient $d$ occurs in $S$, and  $\varepsilon(S) = 0$
otherwise. A standard quadratic equation $S(X) = 1$ is  {\em
regular} if $\kappa(S) \geqslant 4$ and there is a non-commutative
solution of $S(X) = 1$ in $G$ (see \cite{JSJ} for details), or it
is an equation of the type $[x,y]d = 1$. Notice, that if $S(X)=1$
has a solution in $G$,  $\kappa(S) \geqslant 4$, and $n
> 0$ in the orientable case ($n > 1$ in the non-orientable case),
then the equation $S = 1$ has a non-commutative solution, hence
regular.

\medskip
{\bf Basic Form of IFT.} {\it  Let $S(X) = 1$ be a regular
standard quadratic equation over a non-abelian free group $F$ and
let $T(X,Y) = 1$ be an equation over $F$, $|X| = m, \ |Y| = n.$
Suppose that for any solution $U \in V_F(S)$ there exists a tuple
of elements $W \in F^n$ such that $T(U,W) = 1.$ Then there exists
a tuple of words $P = (p_1(X), \ldots, p_n(X))$, with constants
from $F$, such that $T(U,P(U)) = 1$ for any $U \in V_F(S)$.
Moreover, one can fund a tuple $P$ as above effectively. }

\medskip

We define a Zariski topology on $G^n$ by taking algebraic sets in
$G^n$ as a sub-basis for the closed sets of this topology. If $G$
is a  non-abelian fully residually free group (for every finite
set of non-trivial elements in $G$ there exists a homomorphism
from $G$ to a free group such that the images of these elements
are non-trivial), then the closed sets in the Zariski topology
over $G$ are precisely the algebraic sets.

The Basic Form of IFT implies that locally (in terms of Zariski
topology in $F^n$), i.e., in the neighborhood defined by  the
equation $S(X) = 1$, the implicit functions $y_1,\dots ,y_m$ can
be expressed as explicit words in variables $x_1, \ldots, x_n$ and
constants from $F$, say $Y = P(X)$.  This   allows one to
eliminate a quantifier from the following formula (if it holds in
a free group $F$)
$$ \Phi = \forall X \exists Y (S(X) = 1 \ \ \rightarrow \ \ T(X,Y) = 1).$$
  Indeed, in this event the sentence $\Phi$ is equivalent in $F$ to the following
one:
$$
\Psi = \forall X (S(X) = 1 \ \ \rightarrow \ \ T(X, P(X)) = 1).
$$
 From the point of view of model theory Theorem A  states the
existence of very simple Skolem functions for particular
$\forall\exists$-formulas over free groups. Observe, that Theorem
A reinforces the results of \cite{Imp} by giving the corresponding
explicit Skolem functions effectively.

 From algebraic geometry view-point the implicit
function theorem tells one that (in the notations above)  $T(X,Y)
= 1$ has a solution at a generic point of the equation $S(X) = 1$.
Indeed, since the coordinate group $F_{R(S)}$ of the equation
$S(X) = 1$ is discriminated by the free group $F$ the equation
$T(X,Y) =1$ has a solution in the group $F_{R(S)}$ (where elements
 from $X$ are viewed as constants). This shows the Theorem A can be stated
 in  the following form.

\medskip
{\bf Theorem $A^\prime$.} {\it  Let $S(X) = 1$ be a regular
standard quadratic equation over a non-abelian free group $F$ and
let $T(X,Y) = 1$ be an equation over $F$, $|X| = m, \ |Y| = n.$
Suppose that for any solution $U \in V_F(S)$ there exists a tuple
of elements $W \in F^n$ such that $T(U,W) = 1.$ Then the equation
$T(X,Y) = 1$ has a solution in the group $F_{R(S)}$ {\rm (}where
elements from $X$ are viewed as constants from $F_{R(S)}${\rm )}. }

\medskip
This approach allows one to generalize  the results above by
replacing the equation $T(X,Y) = 1$ by an  arbitrary system of
equations and inequalities or even by an  arbitrary boolean
formula. Notice, that such  generalizations in the form of Theorem
A are impossible. To this end we need to introduce a few
definitions.

Let $S(X) = 1$ be a system of  equations over a group $G$ which
has a solution in $G$.  We say that a system of equations $T(X,Y)
= 1$ is {\em compatible}  with $S(X) = 1$ over $G$ if for every
solution $U$ of $S(X) = 1$ in $G$ the equation $T(U,Y) = 1$ also
has a solution in $G$. More generally,   a formula $\Phi(X,Y)$ in
the language $L_A$ is {\em compatible}  with $S(X) = 1$ over $G$,
if for every solution $\bar{a}$ of $S(X) = 1$ in $G$ there exists
a tuple $\bar{b}$ over $G$ such that the formula $\Phi(\bar a,
\bar b )$ is true in $G$, i.e., the algebraic set $V_G(S)$ is a
projection of the truth  set of the formula $\Phi(X,Y) \ \wedge \
(S(X) = 1). $

Suppose now that a  formula $\Phi(X,Y)$ is compatible with $S(X)=
1 $ over $G$. We say that $\Phi(X,Y)$ {\em admits a lift to a
generic point} of $S = 1$ over $G$ (or shortly {\em $S$-lift} over
$G$), if the formula $\exists Y \Phi(X^\mu,Y)$ is true in
$G_{R(S)}$ (here $Y$ are variables and $X^\mu$ are constants from
$G_{R(S)}$).  Finally, an equation $T(X,Y) = 1$, which is
compatible with $S(X) = 1$,   admits a {\em complete $S$-lift} if
every formula $T(X,Y) = 1 \ \& \ W(X,Y) \neq 1$, which is
compatible with $S(X) = 1$  over $G$,  admits an $S$-lift. We say
that the lift (complete lift) is {\em effective} if there is an
algorithm to decide for any equation $T(X,Y)=1$ (any formula
$T(X,Y) = 1 \ \& \ W(X,Y) \neq 1$) whether $T(X,Y)=1$ (the formula
$T(X,Y) = 1 \ \& \ W(X,Y) \neq 1$) admits an $S$-lift, and if it
does, to construct a solution in $G_{R(S)}.$

Now the Implicit Function Theorem (IFT) for regular quadratic
equations can be stated in the following general form. This is the
main technical result of the paper, we prove it in Sections 3--6.

\medskip
{\bf Theorem A.} {\it  Let $S(X,A)=1$ be a regular standard
quadratic equation over $F(A)$. Every equation $T(X,Y,A) = 1$
compatible with $S(X,A) = 1$ admits an effective  complete
$S$-lift. }

\medskip
Furthermore,  the IFT still holds if one replaces $S(X) = 1$ by an
arbitrary  system of a certain type, namely, by a {\em regular
NTQ} system (see \cite{JSJ} for details). To explain  this we need
to introduce a few definitions.

 Let $G$ be a group with a generating set $A$. A system of equations $S
= 1$ is called {\em triangular quasi-quadratic} (shortly, TQ) if
it can be partitioned into the following subsystems
$$
\begin{array}{rrr}
S_1(X_1, X_2, \ldots, X_n,A) & = &1 \\
S_2(X_2, \ldots, X_n,A) &  = & 1\\
&\vdots  &  \\
S_n(X_n,A) &  = & 1
\end{array}
$$
 where for each $i$ one of the following holds:
\begin{enumerate}
\item [1)] $S_i$ is quadratic  in variables $X_i$;
 \item [2)] $S_i= \{[y,z]=1, [y,u]=1 \mid y, z \in X_i\}$ where $u$ is a
group word in $X_{i+1} \cup  \cdots \cup X_n \cup A$ such that its
canonical image  in $G_{i+1}$ is not a proper power. In this case
we say that $S_i=1$ corresponds to an extension of a centralizer;
 \item [3)] $S_i= \{\,[y,z]=1 \mid y, z \in X_i\,\}$;
 \item [4)] $S_i$ is the empty equation.
  \end{enumerate}

Define $G_{i}=G_{R(S_{i}, \ldots, S_n)}$ for $i = 1, \ldots, n$
and put $G_{n+1}=G.$ The  TQ system $S = 1$ is called {\em
non-degenerate} (shortly, NTQ) if each system  $S_i=1$, where
$X_{i+1}, \ldots, X_n$ are viewed as the corresponding constants
from $G_{i+1}$ (under the canonical maps $X_j \rightarrow
G_{i+1}$, $j = i+1, \ldots, n$, has a solution in $G_{i+1}$. The
coordinate group of an NTQ system is called an {\em NTQ group}.

 An NTQ system $S = 1$ is called {\em regular} if for each $i$ the
 system $S_i = 1$ is either of the type 1) or 4), and in the former
 case the quadratic equation $S_i$ is in standard form and regular.

In Section \ref{se:7.4} we prove  IFT for  {\em regular NTQ}
systems.

\medskip
{\bf Theorem B.} {\it  Let $U(X,A)=1$ be a regular NTQ-system.
Every equation $V(X,Y,A)=1$ compatible with $U=1$ admits a
complete effective $U$-lift.}

Notice, that by definition we allow empty  equations in  regular
NTQ systems.  In the case when the  whole system $U = 1$ is empty
 there exists a very strong generalization of the basic implicit
 function theorem due to Merzljakov \cite{Merz}.

\medskip
{\bf Merzljakov's Theorem.} {\it If
 $$F\models \forall X_1\exists Y_1\cdots\forall X_k\exists Y_k
(S(X,Y,A)=1),$$ where $X=X_1\cup\cdots\cup X_k,
Y=Y_1\cup\cdots\cup Y_k$,  then there exist words {\rm (}with
constants from $F$\/{\rm )} $q_1(X_1),\dots , q_k(X_1,\dots ,X_k)
\in F[X]$, such that
$$F[X]\models S(X_1, q_1(X_1),\dots ,X_k, q_k(X_1,\dots
,X_k,A))=1,$$
 i.e., the  equation
  $$S(X_1,Y_1, \dots, X_k,Y_k,A) = 1$$
  {\rm (}in variables  $Y${\rm )} has a solution $Y_i = q_i(X_1, \dots, X_i,A)$ in the free group $F[X]$,
or equivalently,
$$ F \models \forall X_1 \ldots \forall X_n (S(X_1, q_1(X_1,A),
  \ldots, X_k,q_k(X_1,\dots
,X_k,A)) = 1).$$ }

  \medskip
In \cite{Imp} we gave a short proof of Merzljakov's  theorem based
on generalized equations.  The key idea of all known proofs of
this result is to consider a set of {\em Merzljakov's words} as
values of variables from $X_i=\{x_{i1},\dots ,x_{ik_i}\}$:
$$x_{ij}=ba^{m_{ij1}}ba^{m_{ij2}}b\cdots ba^{m_{ijn_{ij}}}b,$$
where $a,b$ are two different generators of $F= F(A)$.  If
$S(X,Y,A)=1$ has a solution for any Merzljakov' words as values of
variables from $X$, then it has a solution of the type
$Y_i=q_i(X_1,\dots ,X_i)$, $i=1,\dots ,k$.

Unfortunately, Merzljakov's words are not, in general, solutions
of a regular quadratic equation $S(X) = 1$ over $F$. In this case,
one needs to find sufficiently many solutions of $S(X) = 1$ over
$F$ with sufficiently complex periodic structure of subwords. To
this end we consider analogs of Merzljakov's words in the group of
automorphisms of $F[X]$ that fix the standard quadratic word
$S(X)$ and the corresponding set of solutions of $S(X) = 1$ in
$F$. In Sections  \ref{se:7.2.5} and \ref{se:7.2} we study in
detail the periodic structure of these solutions. This is the most
technically demanding part of the paper.

There are two more important generalizations of the implicit
function theorem, one -- for arbitrary NTQ-systems, and another --
 for arbitrary systems. We need a few more definitions  to explain
 this. Let $U(X_1, \ldots, X_n,A) = 1$ be an NTQ-system:
$$
\begin{array}{rr}
S_1(X_1, X_2, \ldots, X_n,A) & = 1 \\
S_2(X_2, \ldots, X_n,A) &  = 1\\
\vdots  &  \\
S_n(X_n,A) &  = 1
\end{array}
$$
and $G_{i}=G_{R(S_{i}, \ldots, S_n)}$, $G_{n+1} = F(A)$.

A $G_{i+1}$-automorphism $\sigma$ of $G_i$ is called a {\em
canonical automorphism }
 if the following holds:
 \begin{enumerate}
  \item [1)] if $S_i$ is quadratic  in variables $X_i$ then $\sigma$ is
  induced by a $G_{i+1}$-automorphism of the  group
$G_{i+1}[X_i]$ which fixes $S_i$;
 \item [2)] if $S_i= \{[y,z]=1, [y,u]=1 \mid y, z \in X_i\}$ where $u$ is a
group word in $X_{i+1} \cup  \cdots \cup X_n \cup A$,  then $G_i=
G_{i+1} \ast_{u = u} Ab(X_i\cup \{u\})$, where $Ab(X_i\cup \{u\})$
is a free abelian group with basis $X_i\cup \{u\}$, and in this
event $\sigma$ extends an automorphism of $Ab(X_i\cup \{u\})$
(which fixes $u$);
 \item [3)] If $S_i= \{[y,z]=1 \mid y, z \in X_i\}$ then $G_i =G_{i+1} \ast  Ab(X_i)$,  and in this event $\sigma$ extends an
automorphism of $Ab(X_i)$;
 \item [4)] If $S_i$ is the empty equation then $G_i = G_{i+1}[X_i]$, and in this
 case $\sigma$ is just the identity automorphism of $G_i$.
\end{enumerate}

Let $\pi_i$   be a fixed $G_{i+1}[Y_{i}]$-homomorphism
$$\pi_i : G_i[Y_i]  \rightarrow G_{i+1}[Y_{i+1}],$$
where $\emptyset = Y_1 \subseteq Y_2 \subseteq \ldots \subseteq
Y_n \subseteq Y_{n+1}$ is an ascending chain of finite sets of
parameters, and $G_{n+1} = F(A)$.  Since the system $U = 1$ is
non-degenerate such homomorphisms $\pi_i$ exist. We assume also
that if $S_i(X_i) = 1$ is a standard quadratic equation (the case
1) above) which has a non-commutative solution in $G_{i+1}$, then
$X^{\pi_i}$ is a non-commutative  solution  of $S_i(X_i) = 1$ in
$G_{i+1}[Y_{i+1}].$

A {\em fundamental sequence} (or a {\em fundamental set}) of
solutions of the system $U(X_1,\dots ,X_n,A)=1$ in $F(A)$ with
respect to the  fixed homomorphisms $\pi_1, \ldots, \pi_n$ is a
set of all solutions of $U = 1$ in $F(A)$ of the form
$$\sigma _1\pi _1\cdots \sigma _n\pi _n\tau ,$$
where  $\sigma_i$ is $Y_i$-automorphism of $G_i[Y_i]$ induced by a
canonical automorphism of $G_i$, and  $\tau$ is an
$F(A)$-homomorphism $\tau: F(A \cup Y_{n+1}) \rightarrow F(A)$.
 Solutions
from a given fundamental set of $U$ are called {\em fundamental }
solutions.

\medskip
{\bf Theorem C (Parametrization theorem).} {\it  Let $U(X,A) =1$
be an NTQ-system and $V_{\rm fund}(U)$  a fundamental set of
solutions of $U = 1$ in $F = F(A)$. If a formula
 $$\Phi=\forall X (U(X) = 1 \rightarrow \exists Y (W(X, Y,A)=1  \wedge W_1(X, Y,A)\not =1)$$
 is true in $F$ then one can effectively find finitely many NTQ
 systems $U_1 = 1, \ldots, U_k = 1$ and embeddings $\theta_i:  F_{R(U)}  \rightarrow F_{R(U_i)}$
 such that the formula
 $$\exists Y (W(X^{\theta_{i}}, Y,A)=1 \wedge W_1(X^{\theta _{i}}, Y,A)\not =1)$$
is true in each group $F_{R({ U}_{i})}$.  Furthermore,  for every
solution $\phi: F_{R(U)}  \rightarrow F$ of $U = 1$ from $V_{\rm
fund}(U)$  there exists $i \in \{1, \ldots, k\}$   and a
fundamental solution $\psi: F_{R({ U}_{i})} \rightarrow F$ such
that $\phi = \theta_i \psi$. }

\medskip
As a corollary of this theorem and results from \cite[Section
11]{JSJ},  we obtain the following result.

\medskip
{\bf Theorem D.} {\it Let $S(X)=1$ be an arbitrary  system of
equations over $F$. If a formula
$$ \Phi = \forall X \exists Y (S(X) = 1 \ \ \rightarrow \ \ (W(X, Y,A)=1  \wedge W_1(X, Y,A)\not =1))$$
is  true in $F$  then one can effectively find finitely many NTQ
 systems $U_1 = 1, \ldots, U_k = 1$ and $F$-homomorphisms $\theta_i:  F_{R(S)}  \rightarrow F_{R(U_i)}$
 such that the formula
 $$\exists Y (W(X^{\theta_{i}}, Y,A)=1 \wedge W_1(X^{\theta _{i}}, Y,A)\not =1)$$
is true in each group $F_{R({ U}_{i})}$.  Furthermore,  for every
solution $\phi: F_{R(S)}  \rightarrow F$ of $S = 1$    there
exists $i \in \{1, \ldots, k\}$   and a fundamental solution
$\psi: F_{R({ U}_{i})} \rightarrow F$ such that $\phi = \theta_i
\psi$. }

In Section \ref{se:eefg} we show that the converse of Theorem B
holds. Namely, we prove the following theorem.

\medskip
{\bf Theorem E.} {\it Let $F$ be a free non-abelian group and
$S(X) = 1$ a consistent system of equations over $F$. Then the
following conditions are equivalent:

\begin{enumerate}
\item   The system $S(X) = 1$ is rationally equivalent to a
regular NTQ system.
 \item  Every equation $T(X,Y)=1$ which is
compatible with $S(X)=1$ over $F$ admits an $S$-lift.
 \item   Every
equation $T(X,Y)=1$ which is compatible with $S(X)=1$ over $F$
admits a complete $S$-lift.
\end{enumerate}
}

\medskip
Theorem E immediately implies  the following remarkable property
of regular NTQ systems. Denote by  $L_A$ the  first-order group
theory language with constants from the free group $F(A)$. If
$\Phi$ is  a set of first order sentences of the language $L_A$
then two groups $G$ and $H$ are called {\em $\Phi$-equivalent} if
they satisfy precisely the same sentences from the set $\Phi$. In
this event we write $G\equiv_\Phi H$. In particular,
$G\equiv_{\forall \exists} H$ ($G\equiv_{ \exists \forall} H$)
means that $G$ and $H$ satisfy precisely the same $\forall
\exists$-sentences ($\exists \forall$-sentences). We have shown in
\cite{KMIrc} that for a finitely generated group $G$ if
$G\equiv_{\forall \exists} H$ then $G$ is torsion-free hyperbolic
and fully residually free. Now we improve on this result.

\medskip
{\bf Theorem F.} {\it Let $G$ be a finitely generated group. If\/
$G$ is $\forall \exists$-equivalent to a free non-abelian group
$F$ then $G$ is isomorphic to the coordinate group $F_{R(S)}$ of a
regular NTQ system $S = 1$ over $F$. }

\medskip
Notice, that we prove in the consequent paper  \cite{KM7} that the
converse is also true, moreover, it holds in the strongest
possible form. Namely, the coordinate group $F_{R(S)}$ of a
regular NTQ system $S = 1$ over $F$ is elementary equivalent to a
free non-abelian group $F$. Combining this result with Theorem E
one obtains a complete  algebraic characterization
 of finitely generated groups
which are elementary equivalent to a free non-abelian group.
Similar characterization in different terms is given in
\cite{Sela6}.

\section{Scheme of the proof}\label{se:scheme}
We sketch here the proof of Theorem A  for the orientable
quadratic equation.
\begin{equation}\label{eq:st2..}
\prod_{i=1}^{n}[x_i,y_i] \prod_{i=1}^{m}z_i^{-1}c_iz_i c = 1,\ \ \
n \geqslant 1, m+n \geqslant 1, c\neq 1 .
\end{equation}

We begin with the definition of compatibility. Let $X,Y$ be
families of variables
\begin{df}
Let $S(X) = 1$ be a system of  equations over a group $G$ which
has a solution in $G$.  We say that a system of equations  $T(X,U)
= 1$ is compatible with $S(X) = 1$ over $G$ if for every solution
$B$ of $S(X) = 1$ in $G$ the equation $T(B,U) = 1$ also has a
solution in $G$.
\end{df}

Let $F=F(A)$ be a free group with alphabet $A$. Denote by $S(X)=1$
 equation (\ref{eq:st2..}), where $X=\{x_1,y_1,\dots ,
x_n,y_n,z_1,\dots ,z_m\}$, and suppose that an equation $T(X,U) 1$
is compatible with $S(X)=1$.

STEP 1. The following statement can be obtained using the
Elimination process similar to Makanin-Razborov's process
described in \cite{JSJ}.

One can effectively find a finite disjunction of systems $\Pi
(M,X)$ of graphic equations (without cancellation) in variables
$M,X$ with the following properties.
\bi
\item[1)] Each equation  in the system $\Pi (M,X)$ has form $x\circeq \mu
_{i_1}\circ \cdots\circ \mu _{i_k}$, where $x\in X$, $\mu _i\in
M$, ``$\circeq$ '' stands for graphic equality and ``$\circ$''
means multiplication without cancellation. A solution of such a
graphic equation is a tuple of reduced words $x^{\alpha},\mu
_{i_1}^{\alpha},\dots ,\mu_{i_k}^{\alpha}$ in $F$ such that
$x^{\alpha}$ is graphically equal to the product $\mu
_{i_1}^{\alpha}\circ \cdots\circ \mu _{i_k}^{\alpha}.$

\item[2)] For every solution $B$ of $S(X)=1$ written in reduced form
there exists a graphic solution $B,D$ of one of the systems $\Pi
(M,X)$ in this disjunction.

\item[3)] Let $U=\{u_1,\dots ,u_k\}.$ For every system $Q(X,M)$ one can
find words \newline $f_1(M),\dots ,f_k(M)$ such that for every
solution $B,D$ (not necessary graphic) of the system $Q(X,M)$ in
$F$ one has $T(X, f_1(D),\dots ,f_k(D))=1$.
\ei

Such system of graphic equations $\Pi (M,X)$ is called in Section
\ref{se:cut} a ``cut equation'' (see Definition \ref{df:cut} and
Theorem \ref{th:cut}.) Indeed, variables $X$ are ``cut'' into
pieces. We can think about the cut equation as  a system of
intervals labelled by solutions of $S(X)=1$ that are cut into
several parts corresponding to variables in $M$.

STEP 2. Now we construct a  discriminating family of solutions of
$S(X)=1$ (see the definition in  \cite[Section 1.4]{JSJ}) which
later will be called a {\em generic} family. Consider a group
$F[X]=F*F(X)$ and construct a particular sequence ${\Gamma}
=(\gamma _1,\dots ,\gamma _K)$ of $F$-automorphisms of $F[X]$ that
fix the quadratic word $S(X)$. This is done in Section
\ref{se:7.2.5}. These automorphisms have the property that any two
neighbors in the sequence do not commute and it is in some sense
maximal with this property. For any natural number $j$ define
$\gamma _{j}=\gamma _{r}$, where $r$ is the remainder when $j$ is
divided by $K$.

 For example, for the equation $[x,y]=[a,b]$ we can
take
$$\gamma _1:x\rightarrow x,\ y\rightarrow xy;$$ $$\gamma _2:
x\rightarrow yx,\ y\rightarrow y,$$ in this case $K=2$ and
$$\gamma _{2s-1}=\gamma
_1, \gamma _{2s}=\gamma _2.$$ These automorphisms are, actually,
Dehn twists. Notice that
$$\gamma _1^q : x\rightarrow x,\ y\rightarrow x^qy;\ \ \gamma
_2^q: x\rightarrow y^qx,\ y\rightarrow y,$$
 therefore big
powers of automorphisms produce big powers of elements. Let $L$ be
a multiple of $K$. Define
$$\phi _{L,p}=\gamma _{L}^{p_L}\gamma _{L-1}^{p_{L-1}}\cdots \gamma
_1^{p_1},$$ where $p=(p_1,\dots ,p_L).$ Now we take a suitable
(with small cancellation, in general position) solution of
$S(X)=1$. Denote $F_{{\rm Rad}(S)}=F*F[X]/{\rm ncl}(S).$ This
solution is a homomorphism $\beta: F_{{\rm Rad}(S)}\rightarrow F$.
The family of mappings $$\Psi _L =\{\psi_{L,p}=\phi _{L,p}\beta ,\
p\in P\},$$ where $L$ is large and $P$ is an infinite set of
$L$-tuples of large natural numbers, is a family of solutions of
$S(X)=1$. It is very important that this is a discriminating
family.

For example, take for the equation $[x,y]=[a,b]$ $ x^{\beta}=a,
y^{\beta}=b$, then for $L=4$ we have
\begin{equation}\label{eq:ex}
x=(((a^{p_1}b)^{p_2}a)^{p_3}a^{p_1}b )^{p_4}(a^{p_1}b)^{p_2}a,\ \
\  y=((a^{p_1}b)^{p_2}a)^{p_3}a^{p_1}b.\end{equation}

The word $((a^{p_1}b)^{p_2}a)^{p_3}a^{p_1}b $ is called a period
in rank 4. Notice that the period of rank $4$ is, actually,
$y^{\psi _{3,p}}.$

Since the family of cut equations is finite, some infinite set of
solutions $X^{\Psi _L}$ satisfies the same cut equation $\Pi
(M,X)$. Therefore, it is enough to consider one of the cut
equations $\Pi (M,X)$.

In the example (\ref{eq:ex}) there is no cancellation between $a$
and $b$ and, therefore, it does not matter whether we label
intervals of the cut equation by $X^{\psi _{L,p}}$ or by $X^{\phi
_{L,p}}.$ In Section \ref{se:7.2} we show how to choose a solution
$\beta$ with relatively small cancellation, so that we can forget
about the cancellation and label the intervals of $\Pi (M,X)$ by
$X^{\phi _{L,p}}$.

STEP 3. We can see now that for different $L$-tuples $p$ all
values of
 $X^{\phi _{L,p}}$ (in $F[X]$)  have similar periodic
structure and must be ``cut '' the same way into pieces $\mu\in
M$. Therefore big powers are similarly distributed between pieces
$\mu\in M$. In Section \ref{se:7.3} we introduce the notion of
{\em complexity} of a cut equation.

Let $\Pi (M,X)$ be a cut equation. For a positive integer $n$ by
$k_n(\Pi)$ we denote the number of equations (intervals) $x\circeq
\mu _{i_1}\circ \cdots\circ \mu _{i_{n}}$ that have  right hand
side of length $n$. The following sequence of integers
$$Comp (\Pi)=\left(k_2(\Pi),k_3(\Pi ),\dots ,k_{{\rm length} (\Pi )}(\Pi)\right)$$
is called the {\em complexity} of $\Pi$.

We well-order complexities of cut equations in the (right)
shortlex order:
 if $\Pi$ and $\Pi^\prime$ are two cut equations then
  ${\rm Comp}(\Pi) <  {\rm Comp}(\Pi^\prime)$ if and only if ${\rm length}(\Pi) <
 {\rm length}(\Pi^\prime)$ or ${\rm length}(\Pi) = {\rm length}(\Pi^\prime)$ and there exists $1 \leqslant i\leqslant {\rm length}(\Pi)$ such
 that
  $k_j(\Pi) = k_j(\Pi^\prime)$
 for all $j > i$ but $k_i(\Pi) < k_i(\Pi^\prime)$.

Observe that equations of the form $x\circeq \mu _i$ have no input
into the complexity of a cut equation. In particular, cut
equations that have all graphic equations of length one have the
minimal possible complexity among equations of a given length. We
will write ${\rm Comp}(\Pi) = { 0}$ in the case when $k_i(\Pi) 0$
for every $i = 2, \dots, {\rm length}(\Pi)$.

We introduce the process of transformations of the cut equation
$\Pi (M,X)$. This process consists in ``cutting out'' big powers
of largest periods from the interval and replacing one interval
labelled by $X^{\phi _{i,p}}$ by several intervals labelled by
$X^{\phi _{i-1, p}}$. After such a transformation the left sides
of the graphic equalities in the cut equation correspond to values
$X^{\phi _{i-1,p}}$ (or very short words in $X^{\phi _{i-1,p}}$)
and the complexity either decreases  or stabilizes during several
steps of the process. Suppose ${\rm Comp}(\Pi) = { 0}$ after $t$
transformations, so that each graphic equality
 has form $x^{\phi _{L-t,p}}\circeq\mu$ or $y^{\phi _{L-t,p}}\circeq\nu$.
 Therefore, $x^{\psi _{L-t,p}}\circeq\mu$ or $y^{\psi _{L-t,p}}\circeq\nu$ for
 a discriminating family of solutions $\Psi _{L-t,P}$. By the properties
of discriminating families, $\mu =x,\ \nu=y$ in the group $F_{{\rm
Rad}(S)}$. Substituting $\mu$ and $\nu$ into words $f_1,\dots
,f_k$ we obtain a solution $U$ of the equation $T(X,U)=1$ given by
a formula in $x,y$ in $F_{{\rm Rad}(S)}.$

In a general case, when the length of the right hand side of the
cut equation does not decrease during several steps in the process
of transformations, the situation is, certainly, a bit more
complicated. In this case one can show that in each graphic
equation all the variables $\mu _i$ except the first and the last
one are very short and can be taken almost arbitrary, and the
other variables can be expressed in terms of them and $X^{\Psi
_{L-t,P}}$.

\section{Elementary properties of
liftings}\label{se:lift}
 In this section we discuss some basic properties  of liftings of equations
and inequalities into generic points.

Let $G$ be a group and let $S(X) = 1$ be a system of equations
over $G$. Recall that by $G_{S}$ we denote the quotient group
$G[X]/{\rm ncl}(S)$, where ${\rm ncl}(S)$ is the normal closure of
$S$ in $G[X]$. In particular, $G_{R(S)} = G[X]/R(S)$ is the
coordinate group defined by $S(X) = 1$. The radical $R(S)$ can be
described as follows. Consider a set of $G$-homomorphisms
 $$\Phi_{G,S} = \{\phi \in {\rm Hom}_G(G[S],G)  \mid
\phi(S) = 1\}.$$ Then  \[R(S) =
\left\{\begin{array}{ll}\bigcap_{\phi \in \Phi_{G,S}} \ker \phi &
\mbox{if $\Phi_{G,S} \neq \emptyset$}\\ G[X] & \mbox{otherwise}
\end{array} \right. \]

  Now we put these definitions in a more general framework.
  Let  $H$ and $K$ be $G$-groups   and $M \subset H$. Put
  $$\Phi_{K,M} = \{\phi \in {\rm Hom}_G(H,K)  \mid
\phi(M) = 1\}.$$ Then the following subgroup
  is termed the {\it  $G$-radical of $M$
  with respect to $K$:}
 \[{\rm Rad}_K(M) = \left\{\begin{array}{ll}\bigcap_{\phi \in
\Phi_{K,M}} \ker \phi , & \mbox{if $\Phi_{K,M} \neq \emptyset$,}\\
G[X] & \mbox{otherwise.} \end{array} \right. \] Sometimes, to
emphasize that $M$ is a subset of $H$, we write ${\rm
Rad}_K(M,H)$. Clearly, if $K = G$, then $R(S) = {\rm
Rad}_G(S,G[X])$.

Let $$H_K^* = H/{\rm Rad}_K(1).$$ Then $H_K^*$ is either a
$G$-group or trivial.  If $H_K^* \neq 1$, then it is $G$-separated
by $K$. In the case $K = G$ we omit $K$ in the notation above and
simply write $H^*$. Notice that
$$(H/{\rm ncl}(M))^*_K \simeq H/{\rm Rad}_K(M),$$ in particular, $(G_S)^* G_{R(S)}$.

\begin{lm}
Let $\alpha: H_1 \rightarrow H_2$ be a $G$-homomorphism and
suppose $$\Phi = \{\phi:H_2 \rightarrow K \}$$ be a separating
family of\/ $G$-homomorphisms. Then
$$ \ker \alpha = \bigcap \{\ker(\alpha  \phi) \ \mid \ \phi \in \Phi \} $$
\end{lm}

\bp Suppose $h \in H_1$ and $h \not \in \ker(\alpha).$ Then
$\alpha(h) \neq 1$ in $H_2$. Hence there exists $\phi \in \Phi$
such that $\phi(\alpha(h)) \neq 1$. This shows that $ \ker \alpha
\supset \bigcap \{\ker(\alpha \circ \phi) \ \mid \ \phi \in \Phi
\} $. The other inclusion is obvious. \ep

\begin{lm}
\label{le:6.2} Let $H_1$, $H_2$, and $K$  be $G$-groups.
\begin{enumerate}
\item  Let  $\alpha: H_1 \rightarrow H_2$ be a $G$-homomorphism
and let $H_2$  be
 $G$-separated by $K$. If $M \subset \ker \alpha$,  then
 ${\rm Rad}_K(M) \subseteq \ker \alpha$.
\item Every $G$-homomorphism $\phi: H_1 \rightarrow H_2$ gives
rise to a unique homomorphism $$\phi^*:(H_1)_K^* \rightarrow
(H_2)_K^*$$ such that $\phi\eta_2   = \eta_1\phi^*   $, where
$\eta_i: H_i \rightarrow H_i^*$ is the canonical epimorphism.
\end{enumerate}
\end{lm}

\begin{proof}
 (1) We have
 \bea
 {\rm Rad}_K(M,H_1) &=& \bigcap \{\ker \phi \mid \phi: H_1
\rightarrow_G K  \ \wedge \ \phi(M) = 1\}\\
&\subseteq& \bigcap
\{\ker(\alpha  \beta) \mid \beta : H_2 \rightarrow_G K\}\\
&=& \ker \alpha. \eea

(2) Let $\alpha:H_1 \rightarrow (H_2)_K^*$ be the composition of
the following homomorphisms $$ H_1 \stackrel{\phi}{\rightarrow}
H_2 \stackrel{\eta_2}{\rightarrow} (H_2)_K^*. $$ Then by assertion
1 ${\rm Rad}_K(1,H_1) \subseteq \ker \alpha$, therefore $\alpha$
induces the canonical $G$-homomorphism $\phi^*:(H_1)_K^*
\rightarrow (H_2)_K^*$. \end{proof}

\begin{lm}\
\label{le:6.3}
\begin{enumerate}
\item The canonical map $\lambda:G \rightarrow G_S$ is an
embedding $\Longleftrightarrow$ $S(X) = 1$ has a solution in some
$G$-group $H$. \item The canonical map $\mu:G \rightarrow
G_{R(S)}$ is an embedding $\Longleftrightarrow$ $S(X) = 1$ has a
solution in some $G$-group $H$ which is $G$-separated by $G$.
\end{enumerate}
\end{lm}

\bp
  (1) If $S(x_1, \ldots,x_m) = 1$ has a solution
$(h_1,\dots,h_m)$ in some $G$-group $H$, then the $G$-homomorphism
$x_i \rightarrow h_i, \ (i =1,\dots,m)$ from $G[x_1,\dots,x_m]$
into $H$ induces a homomorphism  $\phi: G_S \rightarrow H$. Since
$H$ is a $G$-group all non-trivial elements from $G$ are also
non-trivial in the factor-group $G_S$, therefore $\lambda:G
\rightarrow G_S$ is an embedding. The converse is obvious.

(2) Let $S(x_1, \ldots,x_m) = 1$ have a solution $(h_1,\dots,h_m)$
in some $G$-group $H$ which is $G$-separated by $G$. Then there
exists the canonical $G$-homomorphism $\alpha: G_S \rightarrow H$
defined as in the proof of the first assertion. Hence $R(S)
\subseteq \ker \alpha$ by Lemma \ref{le:6.2},  and $\alpha$
induces a homomorphism from $G_{R(S)}$ into $H$, which is monic on
$G$. Therefore $G$ embeds into $G_{R(S)}$. The converse is
obvious. \ep

Now we apply Lemma \ref{le:6.2} to coordinate groups of nonempty
algebraic sets.

\begin{lm}
\label{le:6.4} Let subsets $S$ and $T$ from $ G[X]$ define
non-empty algebraic sets in a group $G$. Then every
$G$-homomorphism $\phi:G_S \rightarrow G_T$ gives rise to a
$G$-homomorphism $\phi^*: G_{R(S)} \rightarrow G_{R(T)}$.
\end{lm}

\bp The result follows from Lemma \ref{le:6.2}  and Lemma
\ref{le:6.3}. \ep

Now we are in a position to give the following

Recall that for a consistent  system of  equations $S(X) = 1$ over
a group $G$, a system of equations $T(X,Y) = 1$ is compatible with
$S(X) = 1$ over $G$ if for every solution $U$ of $S(X) = 1$ in $G$
the equation $T(U,Y) = 1$ also has a solution in $G$, i.e., the
algebraic set $V_G(S)$ is a projection of the algebraic set
$V_G(S\cup T).$

The next proposition describes compatibility of two equations in
terms of their coordinate groups.

\begin{prop}
\label{prop:6.1} Let $S(X) = 1$ be a system of  equations over a
group $G$ which has a solution in $G$. Then $T(X,Y) = 1$ is
compatible with $S(X) = 1$ over $G$ if and only if $G_{R(S)}$ is
canonically embedded into $G_{R(S\cup T)}$, and every
$G$-homomorphism $\alpha: G_{R(S)} \rightarrow G$ extends to a
$G$-homomorphisms $\alpha^\prime: G_{R(S\cup T)} \rightarrow G$.
\end{prop}

\bp Suppose first that $T(X,Y) = 1$ is compatible with $S(X) = 1$
over $G$ and suppose that $V_G(S) \neq \emptyset.$ The identity
map $X \rightarrow X$ gives rise to a $G$-homomorphism $$ \lambda:
G_S \longrightarrow G_{S\cup T} $$ (notice that both $G_S$ and
$G_{S\cup T}$ are $G$-groups by Lemma \ref{le:6.3}), which by
Lemma \ref{le:6.4} induces a $G$-homomorphism $$ \lambda^*:
G_{R(S)} \longrightarrow G_{R(S\cup T)}. $$ We claim that
$\lambda^*$ is an embedding. To show this  we need to prove first
the  statement about the extensions of homomorphisms. Let
$\alpha:G_{R(S)} \rightarrow G$ be  an arbitrary $G$-homomorphism.
It follows that $ \alpha(X)$ is a solution of $S(X) = 1$ in $G$.
Since $T(X,Y) = 1$ is compatible with $S(X) = 1$ over $G$, there
exists a solution, say $\beta(Y)$,  of  $T(\alpha(X),Y) = 1$ in
$G$. The  map $$ X \rightarrow \alpha(X), Y \rightarrow \beta(Y)
$$ gives rise to  a $G$-homomorphism $ G[X,Y] \rightarrow G$,
which induces a $G$-homomorphism $\phi: G_{S\cup T} \rightarrow
G$. By Lemma \ref{le:6.4} $\phi$ induces a $G$-homomorphism $$
\phi^* : G_{R(S\cup T)} \longrightarrow G.$$ Clearly, $\phi^*$
makes the following diagram to commute.
\bigskip
\begin{center}
\begin{picture}(100,100)(0,0)
\put(0,100){$G_{R(S)}$} \put(100,100){$G_{R(S\cup T)}$}
\put(0,0){$G$} \put(25,103){\vector(1,0){70}}
\put(5,93){\vector(0,-1){78}} \put(95,95){\vector(-1,-1){80}}
\put(-10,50){$\alpha$} \put(55,45){$\phi^*$}
\put(50,108){$\lambda^*$}
\end{picture}
\end{center}
Now to prove  that $\lambda^*$ is an embedding, observe that
$G_{R(S)}$ is $G$-separated by $G$. Therefore for every
non-trivial $h \in G_{R(S)}$ there exists a $G$-homomorphism
$\alpha: G_{R(S)} \rightarrow G$ such that $\alpha(h) \neq 1$. But
then $\phi^*(\lambda^*(h)) \neq 1$ and consequently $h \not \in
\ker \lambda^*$. The converse statement is obvious. \ep

\smallskip
Let $S(X) = 1$ be a system of  equations over $G$ and suppose
$V_G(S) \neq \emptyset$. The canonical embedding $X \rightarrow
G[X]$ induces the canonical map $$\mu: X \rightarrow G_{R(S)}.$$
We are ready to formulate the main definition.
\begin{df} Let $S(X) = 1$ be a system of  equations over $G$ with
$V_G(S) \neq \emptyset$ and let $\mu :X \rightarrow G_{R(S)}$ be
the canonical map.  Let a system  $T(X,Y) = 1$ be compatible with
$S(X) = 1 $ over $G$. We say that $T(X,Y) = 1$ admits a lift to a
generic point of $S = 1$ over $G$ (or, shortly, $S$-lift over $G$)
if $T(X^\mu,Y) = 1$ has a solution in $G_{R(S)}$ (here $Y$ are
variables and $X^\mu$ are constants from $G_{R(S)}$).
\end{df}
\begin{lm}
\label{le:oger} Let $T(X,Y) = 1$ be compatible with $S(X) = 1 $
over $G$. If $T(X,Y) = 1$ admits an $S$-lift, then the identity
map $Y \rightarrow Y$ gives rise to  a canonical
$G_{R(S)}$-epimorphism from $G_{R(S\cup T)}$ onto the coordinate
group of $T(X^\mu,Y) = 1$ over $G_{R(S)}$:
 $$  \psi^* : G_{R(S\cup T)} \rightarrow
G_{R(S)}[Y]/{\rm Rad}_{G_{R(S)}}(T(X^\mu,Y)).$$ Moreover, every
solution $U$
 of  $T(X^\mu,Y) = 1$ in $G_{R(S)}$ gives rise to a $G_{R(S)}
 $-homomorphism $\phi_U:  G_{R(S\cup T)} \rightarrow G_{R(S)}$, where
$\phi_U(Y) = U$.
\end{lm}

\bp
 Observe that the following chain of isomorphisms hold:
 \bea
G_{R(S\cup T)} &\simeq_G& G[X][Y]/{\rm Rad}_G(S\cup T)\\
& \simeq_G &
G[X][Y]/{\rm Rad}_G({\rm Rad}_G(S,G[X])\cup T)\\
& \simeq_G &
\left(G[X][Y]/{\rm ncl}({\rm Rad}_G(S,G[X]) \cup T)\right)^*\\
&\simeq_G& \left(G_{R(S)}[Y]/{\rm ncl}(T(X^\mu,Y))\right)^*. \eea
 Denote by $\overline{G_{R(S)}}$ the canonical image of $G_{R(S)}$ in
 $(G_{R(S)}[Y]/{\rm ncl}(T(X^\mu,Y)))^* .$

 Since ${\rm Rad}_{G_{R(S)}}(T(X^\mu,Y))$ is a normal subgroup in
$G_{R(S)}[Y]$\linebreak containing
 $T(X^\mu,Y)$  there exists  a canonical $G$-epimorphism
 $$ \psi:  G_{R(S)}[Y]/{\rm ncl}(T(X^\mu,Y)) \rightarrow
G_{R(S)}[Y]/{\rm Rad}_{G_{R(S)}}(T(X^\mu,Y)).$$
 By Lemma \ref{le:6.2} the homomorphism $\psi$  gives rise to a
canonical
  $G$-homomorphism
   $$ \psi^*:  (G_{R(S)}[Y]/{\rm ncl}(T(X^\mu,Y)))^* \rightarrow
(G_{R(S)}[Y]/{\rm Rad}_{G_{R(S)}}(T(X^\mu,Y)))^*.$$
  Notice that the group $G_{R(S)}[Y]/{\rm Rad}_{G_{R(S)}}(T(X^\mu,Y))$ is the
  coordinate group of the system $T(X^\mu,Y)= 1$ over  $G_{R(S)}$  and
  this system has a solution in  $G_{R(S)}$. Therefore this group is a
  $G_{R(S)}$-group and  it is $G_{R(S)}$-separated by
  $G_{R(S)}$. Now since $G_{R(S)}$ is the coordinate group of $S(X) = 1$
over $G$ and
  this system has a solution in $G$, we see that $G_{R(S)}$ is $G$-separated by $G$.
  It follows that the group $G_{R(S)}[Y]/{\rm Rad}_{G_{R(S)}}(T(X^\mu,Y))$  is
$G$-separated
   by $G$. Therefore
  $$G_{R(S)}[Y]/{\rm Rad}_{G_{R(S)}}(T(X^\mu,Y)) =(G_{R(S)}[Y]/{\rm Rad}_{G_{R(S)}}(T(X^\mu,Y)))^*.$$
  Now we can see that
  $$ \psi^* : G_{R(S\cup T)} \rightarrow
G_{R(S)}[Y]/{\rm Rad}_{G_{R(S)}}(T(X^\mu,Y))$$
  is a $G$-homomorphism which maps the subgroup $\overline{G_{R(S)}}$
from
  $G_{R(S \cup T)}$   onto the subgroup $G_{R(S)}$ in
  $G_{R(S)}[Y]/{\rm Rad}_{G_{R(S)}}(T(X^\mu,Y))$.

  This shows  that $\overline{G_{R(S)}} \simeq_G G_{R(S)}$ and
$\psi^*$ is a  $G_{R(S)}$-homomorphism. If $U$ is a solution of
$T(X^\mu,Y) = 1$ in  $G_{R(S)}$, then there exists a
$G_{R(S)}$-homomorphism
 $$\phi_U:
G_{R(S)}[Y]/{\rm Rad}_{G_{R(S)}}(T(X^\mu,Y)) \rightarrow
G_{R(S)}.$$
  such that $\phi_U(Y) = U$.  Obviously, composition of $\phi_U$ and
$\psi^*$
  gives a  $G_{R(S)}$-homomorphism from $G_{R(S\cup T)}$ into
$G_{R(S)}$,
 as desired.
 \ep

The next result characterizes lifts in terms of the coordinate
groups of the corresponding equations.

\begin{prop}
\label{pr:Slift}
 Let $S(X) = 1$ be an equation over $G$ which has a solution in
$G$. Then for an arbitrary  equation  $T(X,Y) = 1$ over $G$  the
following conditions are equivalent:
\begin{enumerate}
\item  $T(X,Y) = 1$ is compatible with $S(X) = 1$ and $T(X,Y) = 1$
admits $S$-lift over $G$; \item  $G_{R(S)}$ is a retract of
$G_{R(S,T)}$, i.e., $G_{R(S)}$ is a subgroup of $G_{R(S,T)}$ and
there exists a $G_{R(S)}$-homomorphism $G_{R(S,T)} \rightarrow
G_{R(S)}.$
\end{enumerate}
\end{prop}

\begin{proof} (1) $\Longrightarrow$ (2). By Proposition
\ref{prop:6.1} $G_{R(S)}$ is a subgroup of $G_{R(S,T)}$. Moreover,
$T(X^\mu,Y) = 1$ has a solution in $G_{R(S)}$, so by Lemma
\ref{le:oger} there exists a $G_{R(S)}$- homomorphism $G_{R(S,T)}
\rightarrow G_{R(S)}$, i.e., $G_{R(S)}$ is a retract of
$G_{R(S,T)}$.

(2) $\Longrightarrow$ (1). If $\phi: G_{R(S,T)} \rightarrow
G_{R(S)}$ is a retract then every $G$-homomorphism $\alpha:
G_{R(S)}  \rightarrow G$ extends to a $G$-homomorphism $\phi\alpha
:  G_{R(S,T)} \rightarrow G$. It follows from
 Proposition \ref{prop:6.1} that $T(X,Y) = 1$ is compatible with $S(X)=1$ and
 $\phi$ gives a solution of $T(X^\mu,Y) = 1$ in $G_{R(S)}$, as desired.
\end{proof}

Denote by $\mathcal C$ (respectively ${\mathcal C}^\ast$) the
class of all finite systems $S(X) = 1$ over $F$ such that every
equation
 $T(X,Y) = 1$ compatible with $S = 1$ admits an $S$-lift (complete
 $S$-lift).

The following result shows that the classes $\mathcal C$ and
${\mathcal C}^\ast $ are closed under rational equivalence.

\begin{lm} \label{le:rat-equiv} Let systems $S = 1$ and $U = 1$
be  rationally equivalent. Then:
 \begin{enumerate}
  \item  If $U = 1$ is in $\mathcal C$ then $S = 1$ is $\mathcal
  C$;
   \item  If $U = 1$ is in ${\mathcal C}^\ast $ then $S = 1$ is ${\mathcal C}^\ast
   $.
   \end{enumerate}
\end{lm}

\begin{proof} We prove (2), a similar argument proves (1).
 Suppose that a system $S(X)
= 1$ is rationally equivalent to a system $U(Z) = 1$ from
${\mathcal C}^\ast $. Then (see \cite{BMR1}) their coordinate
groups $F_{R(S)}$ and $F_{R(U)}$ are $F$-isomorphic. Let $\phi:
F_{R(S)} \rightarrow F_{R(U)}$ be an $F$-isomorphism. Then
 $X^{\phi} = P(Z)$ for some word mapping $P$.
Suppose now that a formula
 $$T(X,Y)=1 \wedge W(X,Y)\neq 1 $$
   is compatible with $S(X)=1$ over $F$. One needs to show that this formula
    admits an $S$-lift. Notice that
$$T(P(Z),Y) = 1 \wedge W(P(Z),Y)\neq 1 $$
 is compatible with $U(Z)= 1$, hence it admits a $U$-lift.  So there exists an element, say
$D(Z) \in F_{R(U)}$,  such that in $F_{R(U)}$ the following holds
 $$T(P(Z),D(Z)) = 1 \wedge W(P(Z),D(Z))\neq 1. $$
   Now
$$1 = T(P(Z),D(Z))^{\phi^{-1}} =T(P(Z)^{\phi^{-1}},D(Z^{\phi^{-1}})) = T(X,D(Z^{\phi^{-1}}))$$
 and
  $$ 1 \neq W(P(Z),D(Z))^{\phi^{-1}} = W(X,D(Z^{\phi^{-1}}))$$
so
 $$T(P(Z),Y) = 1 \wedge W(P(Z),Y)\neq 1 $$
    admits a complete $S$-lift, as required.
    \end{proof}

\section{Cut equations}\label{se:cut}

We refer to \cite{JSJ} for the notion of a generalized equation.
In the proof of the implicit function theorems it will be
convenient to use a modification of the notion of a generalized
equation.  The following definition provides a framework for such
a modification.

\begin{df}\label{df:cut}
A cut equation $\Pi = ({\mathcal E}, M, X, f_M, f_X)$ consists of
a set of intervals $\mathcal E$, a set of variables $M$, a set of
parameters $X$, and two  labeling functions $$f_X: {\mathcal E}
\rightarrow F[X], \ \ \ f_M:  {\mathcal E} \rightarrow F[M] .$$
For an interval $\sigma \in {\mathcal E}$ the image  $f_M(\sigma)
= f_M(\sigma)(M)$ is a reduced word in variables $M^{\pm 1}$ and
constants from $F$, we call it a {\it partition}  of
$f_X(\sigma)$.
\end{df}

Sometimes we write $\Pi = ({\mathcal E}, f_M, f_X)$ omitting $M$
and $X$.

\begin{df}
A solution of a cut equation $\Pi = ({\mathcal E}, f_M, f_X)$ with
respect to an $F$-homomorphism $\beta : F[X] \rightarrow F$  is an
$F$-homomorphism $\alpha : F[M] \rightarrow F$ such that: 1) for
every $\mu \in M$ $\alpha(\mu)$ is a reduced non-empty word; 2)
for every reduced word $f_M(\sigma)(M)\ (\sigma \in {\mathcal E})$
the replacement $m \rightarrow \alpha(m) \ (m \in M) $ results in
a word $f_M(\sigma)(\alpha(M))$ which is  a  reduced word as
written and  such that $f_M(\sigma)(\alpha(M))$ is graphically
equal to the reduced form of $\beta(f_X(\sigma))$; in particular,
the following diagram is commutative.
\begin{center}
\begin{picture}(100,100)(0,0)
\put(50,100){$\mathcal E$} \put(0,50){$F(X)$} \put(100,50){$F(M)$}
\put(50,0){$F$} \put(47,97){\vector(-1,-1){30}}
\put(53,97){\vector(1,-1){30}} \put(3,47){\vector(1,-1){30}}
\put(97,47){\vector(-1,-1){30}} \put(16,73){$f_X$}
\put(77,73){$f_M$} \put(16,23){$\beta$} \put(77,23){$\alpha$}
\end{picture}
\end{center}
\end{df}

If $\alpha: F[M] \rightarrow F$ is a solution of a cut equation
$\Pi = ({\mathcal E}, f_M, f_X)$ with respect to an
$F$-homomorphism $\beta : F[X] \rightarrow F$, then we write
$(\Pi, \beta,\alpha)$ and refer to $\alpha$ as a \emph{solution
of} $\Pi$ {\it modulo} $\beta$. In this event,   for a given
$\sigma \in {\mathcal E}$ we say that $f_M(\sigma)(\alpha(M))$ is
a {\it partition} of $\beta(f_X(\sigma))$. Sometimes we also
consider homomorphisms $\alpha:F[M] \rightarrow F$, for which the
diagram above is still commutative, but cancellation may occur in
the words $f_M(\sigma)(\alpha(M))$. In this event we refer to
$\alpha$ as a {\em group} solution of $\Pi$ with respect to
$\beta$.

\begin{lm}
\label{le:cut}
 For a generalized equation $\Omega(H)$ one can effectively construct a cut
equation
 $\Pi_{\Omega} = ({\mathcal E}, f_X, f_M)$ such that  the following
conditions hold:
\begin{enumerate}
 \item $ X$ is a partition of the whole interval $[1,\rho_{\Omega}]$ into
disjoint
  closed subintervals;

 \item  $M$ contains the set of variables $H$;

 \item for any solution $U = (u_1, \ldots, u_\rho)$ of $\Omega$  the cut
equation
 $\Pi_{\Omega}$ has a solution
 $\alpha$ modulo the canonical homomorphism $$\beta_U: F(X) \rightarrow
F$$
 {\rm (}$\beta_U(x) = u_i u_{i+1} \cdots u_j$ where $i,j$ are,
correspondingly,
  the left and the right  end-points of the interval $x${\rm )};

 \item for any solution $(\beta,\alpha)$ of the cut equation $\Pi_{\Omega}$
the
 restriction of $\alpha$ on $H$ gives a solution of the generalized
equation
 $\Omega$.
 \end{enumerate}
 \end{lm}

\bp
 We begin with defining the sets $X$ and $M$. Recall that  a
closed interval of $\Omega$ is a union of closed sections of
$\Omega$.  Let $X$ be an arbitrary partition of the whole interval
$[1,\rho_{\Omega}]$ into closed subintervals (i.e., any two
intervals in $X$ are disjoint and the union of $X$ is the whole
interval $[1,\rho_{\Omega}]$).

 Let $B$ be a set of representatives of dual bases of $\Omega$, i.e.,
for every base
  $\mu$  of $\Omega$ either $\mu$ or $\Delta(\mu)$ belongs to $B$, but
not both.
 Put $M = H \cup B$.

Now let $\sigma \in X$.  We denote  by $B_{\sigma}$  the set of
all bases over $\sigma$ and by  $H_\sigma$ the set of all items in
$\sigma$. Put $S_\sigma = B_{\sigma} \cup H_\sigma.$  For $e \in
S_\sigma$ let $s(e)$ be the interval $[i,j]$, where $i < j$ are
the endpoints of $e$.  A sequence  $P = (e_1, \ldots,e_k)$ of
elements from $S_\sigma$  is called a {\it partition} of $\sigma$
if $s(e_1) \cup \cdots \cup s(e_k) = \sigma$ and $s(e_i) \cap
s(e_j) = \emptyset$ for $i \neq j$.  Let ${\rm Part}_\sigma$ be
the set of all partitions of $\sigma$. Now put
 $${\mathcal E} = \{P \mid P \in {\rm Part}_\sigma, \sigma \in X\}.$$
Then for every $P \in {\mathcal E}$
 there exists one and only one $\sigma \in X$ such that $P \in
{\rm Part}_\sigma$.
 Denote this $\sigma$ by $f_X(P)$. The map $f_X: P \rightarrow f_X(P)$
is a
 well-defined function from ${\mathcal E}$ into $F(X)$.

 Each partition $P = (e_1, \ldots,e_k) \in {\rm Part}_\sigma$ gives rise to a
word
  $w_P(M) = w_1 \ldots w_k$  as follows.  If $e_i \in H_\sigma$ then
$w_i = e_i$.
  If $e_i = \mu \in B_\sigma$ then $w_i = \mu^{\varepsilon(\mu)}$.
  If $e_i = \mu$ and $\Delta(\mu) \in B_\sigma$ then
  $w_i = \Delta(\mu)^{\varepsilon(\mu)}$. The map $f_M(P) = w_P(M)$ is a
  well-defined function from ${\mathcal E}$ into $F(M)$.

  Now set $\Pi_{\Omega} = ({\mathcal E}, f_X, f_M)$.  It is not hard to see
from the
  construction that the cut equation $\Pi_\Omega$ satisfies all the
required properties.
  Indeed, (1) and (2) follow directly from the construction.

  To verify (3), let's consider a
  solution $U = (u_1, \ldots, u_{\rho_\Omega})$ of $\Omega$. To define
  corresponding functions $\beta_U$ and $\alpha$, observe that the
function $s(e)$
  (see above) is defined for every $e \in X \cup M$.
   Now  for $\sigma \in X$ put
 $\beta_U(\sigma) = u_i\ldots u_j$, where $s(\sigma) = [i,j]$, and for $m
\in M$
 put $\alpha(m) =  u_i\ldots u_j$, where $s(m) = [i,j]$.  Clearly,
$\alpha$ is a
 solution of $\Pi_\Omega$ modulo $\beta$.

 To verify (4) observe that if $\alpha$ is a solution of $\Pi_\Omega$ modulo
$\beta$,
  then the restriction of $\alpha$ onto the subset $H \subset M$ gives a
solution
 of the generalized equation $\Omega$. This follows from the
construction
 of the words $w_p$ and the fact that the words $w_p(\alpha(M))$ are
reduced as
 written (see definition of a solution of a cut equation). Indeed, if a
base
 $\mu$ occurs in a partition $P \in {\mathcal E}$, then there is a partition
 $P^\prime \in {\mathcal E}$  which is obtained from $P$ by replacing $\mu$
by the
 sequence $h_i\ldots h_j$. Since there is no cancellation in words
$w_P(\alpha(M))$ and
 $w_{P^\prime}(\alpha(M))$, this implies that
 $\alpha(\mu)^{\varepsilon(\mu)} = \alpha(h_i\ldots h_j)$. This shows
that
 $\alpha_H$ is a solution of $\Omega$.
 \ep

\begin{theorem}\label{th:cut}
Let $S(X,Y,A)) = 1$ be a system of  equations over $F = F(A)$.
Then one can effectively construct a finite set of cut equations
 $${\mathcal CE}(S) = \{\Pi_i \mid  \Pi_i =({\mathcal E}_i, f_{X_i},  f_{M_i}),
i = 1 \ldots,k \}$$
  and a finite set of tuples of words $\{Q_i(M_i) \mid i = 1,
\dots,k\}$  such that:
\begin{enumerate}
\item for every equation $\Pi_i =({\mathcal E}_i, f_{X_i},
f_{M_i}) \in {\mathcal CE}(S)$,
 one has $X_i = X$ and   $f_{X_i}({\mathcal E}_i) \subset  X^{\pm 1}$;

\item  for  any  solution $(U,V)$ of $S(X,Y,A) = 1$ in $F(A)$,
there exists a number $i$
 and a
 tuple of words $P_{i,V}$ such that the  cut
equation $\Pi_i \in {\mathcal CE}(S)$  has a solution $\alpha: M_i
\rightarrow F$ with respect to the $F$-homomorphism  $\beta_U:
F[X] \rightarrow F$ which is induced by the map $ X \rightarrow
U$. Moreover, $U = Q_i(\alpha(M_i))$, the word $Q_i(\alpha(M_i))$
is reduced as written,  and  $V = P_{i,V}(\alpha(M_i))$;

\item  for any $\Pi_i \in {\mathcal CE}(S)$ there exists a tuple
of words $P_{i,V}$ such that for  any solution (group solution)
$(\beta, \alpha)$ of  $\Pi_i$ the pair $(U,V),$ where $U =
Q_i(\alpha(M_i))$ and $V = P_{i,V}(\alpha(M_i)),$ is a solution of
$S(X,Y) = 1$ in $F$.
\end{enumerate}
\end{theorem}

\bp
 Let $S(X,Y) = 1$ be a system of  equations over a free group $F$.
In  \cite[Section 4.4]{JSJ}    we have constructed a set of
initial parameterized generalized equations  ${\mathcal GE}_{\rm
par}(S) = \{\Omega_1, \ldots, \Omega_r\}$ for $S(X,Y) = 1$ with
respect to the set of parameters $X$.
 For each $\Omega \in {\mathcal GE}_{\rm par}(S)$ in \cite[Section 8]{JSJ} we
constructed  the finite tree $T_{\rm sol}(\Omega)$ with respect to
parameters $X$. Observe that parametric part
$[j_{v_0},\rho_{v_0}]$ in the root equation $\Omega =
\Omega_{v_0}$ of the tree $T_{\rm sol}(\Omega)$ is  partitioned
into a disjoint union of closed sections  corresponding to
$X$-bases and constant bases (this follows from the construction
of the initial equations in the set  ${\mathcal GE}_{\rm
par}(S)$). We label every closed section $\sigma$ corresponding to
a variable  $x \in X^{\pm 1}$ by $x$, and every constant section
corresponding to a constant $a$ by $a$.  Due to our construction
of the tree $T_{\rm sol}(\Omega)$ moving along a branch $B$ from
the initial vertex $v_0$ to a terminal vertex $v$, we transfer all
the bases from the active and non-active parts into parametric
parts until, eventually,  in $\Omega_v$ the whole interval
consists of the parametric part.
 Observe also that, moving along $B$ in the parametric part,  we neither
introduce
  new closed sections nor  delete any. All we do is we split
(sometimes) an item in a closed  parametric section into two new
ones. In any event we keep the same label of the section.

Now for a terminal vertex $v$ in $T_{\rm sol}(\Omega)$  we
construct a cut equation $\Pi^\prime_v = ({\mathcal E}_v, f_{X_v},
f_{M_v})$ as in Lemma \ref{le:cut} taking the set of all  closed
sections of $\Omega_v$ as the partition $X_v$.  The set of cut
equations
 $$ {\mathcal CE}^\prime(S) = \{\Pi^\prime_v \mid \Omega \in {\mathcal
GE}_{\rm par}(S), v
 \in VTerm(T_{\rm sol}(\Omega))\}$$
satisfies all the requirements of the theorem except $X_v$ might
not be equal to $X$. To satisfy this condition we adjust slightly
the equations $\Pi_v^\prime$.

To do this, we denote by  $l:X_v \rightarrow X^{\pm 1} \cup A^{\pm
1}$ the labelling function on the set of closed sections of
$\Omega_v$. Put $\Pi_v = ({\mathcal E}_v, f_{X}, f_{M_v})$ where
$f_X$ is the composition of $f_{X_v}$ and $l$. The set of  cut
equations
 $$ {\mathcal CE}(S) = \{\Pi_v \mid  \Omega \in {\mathcal GE}_{\rm par}(S), v
 \in VTerm(T_{\rm sol}(\Omega))\}$$
satisfies all the conditions of the theorem. This follows from
\cite[Theorem 8.1]{JSJ},  and from Lemma \ref{le:cut}. Indeed, to
satisfy  3) one can take the words $P_{i,V}$ that correspond to a
minimal solution of $\Pi_i$, i.e., the words $P_{i,V}$ can be
obtained from a given particular way to transfer all bases from
$Y$-part onto $X$-part.

\ep

The next result shows that for every cut equation $\Pi$ one can
effectively and canonically  associate a generalized equation
$\Omega_{\Pi}$.

   For  every  cut equation $\Pi= ({\mathcal E}, X, M, f_X, f_M)$  one
  can  canonically associate a generalized equation $\Omega_{\Pi}(M,X)$
as follows.
  Consider the following word

$$ V = f_X(\sigma_1)f_M(\sigma_1)  \cdots f_X(\sigma_k)f_M(\sigma_k). $$
 Now we are going to mimic the construction of the generalized equation
in \cite[Lemma 4.6]{JSJ}.
  The set of boundaries $BD$ of $\Omega_{\Pi}$ consists of
positive integers
 $1, \ldots, |V| +1$. The set of bases $BS$ is union of the following
sets:

a)  every letter $\mu$ in  the word $V$. Letters  $X^{\pm 1} \cup
M^{\pm 1}$ are variable bases, for every two different occurrences
$\mu^{\varepsilon_1}, \mu^{\varepsilon_2}$ of  a letter  $\mu \in
X^{\pm 1} \cup M^{\pm 1}$ in $V$ we say that these bases  are dual
and they have the same orientation if $\varepsilon_1\varepsilon_2
= 1$, and different orientation otherwise. Each occurrence of a
letter  $a \in A^{\pm 1}$ provides a constant base with the label
$a$.  Endpoints of these bases correspond to their positions in
the word $V$ \cite[Lemma 4.6]{JSJ}.

 b)   every pair of subwords $f_X(\sigma_i), f_M(\sigma_i)$ provides a
pair of dual
 bases $\lambda_i, \Delta(\lambda_i)$, the base $\lambda_i$ is located
above the
 subword $f_X(\sigma_i)$, and  $\Delta(\lambda_i)$ is located above
 $f_M(\sigma_i)$ (this defines the endpoints of the bases).

 Informally, one can visualize the generalized equation $\Omega_{\Pi}$
as
 follows.   Let ${\mathcal E} = \{\sigma_1, \ldots, \sigma_k\}$ and let
  ${\mathcal E}^\prime = \{ \sigma^\prime \mid \sigma \in {\mathcal E}\}$
   be another disjoint copy of the set ${\mathcal E}$. Locate intervals from
${\mathcal
   E} \cup {\mathcal E}^\prime$ on a segment $I$ of a straight line from
left to the
   right in the following order $\sigma_1, \sigma_1^\prime,
   \ldots, \sigma_k, \sigma_k^\prime$; then  put bases over $I$
according to the word
   $V$.
 The next result summarizes the discussion above.
 \begin{lm}
  \label{le:5.4.3}
 For  every  cut equation $\Pi= ({\mathcal E}, X, M, f_x, f_M)$,  one
  can  canonically associate a generalized equation $\Omega_{\Pi}(M,X)$
such that
   if $\alpha_\beta : F[M] \rightarrow F$ is a solution of the cut
 equation $\Pi$, then the maps $\alpha: F[M] \rightarrow F$ and $\beta:
F[X]
 \rightarrow F$ give rise to a solution of the group equation (not
generalized!)
 $\Omega_{\Pi}^* = 1$ in such a way that  for every $\sigma \in {\mathcal E}$
 $f_M(\sigma)(\alpha(M))$ is a reduced word which is graphically equal
to
 $\beta(f_X(\sigma)(X))$, and vice versa.
  \end{lm}

\section{Basic automorphisms of orientable quadratic equations}
\label{se:7.2.5}

In this section, for a finitely generated fully residually free
group $G$ we introduce some particular $G$-automorphisms of a free
$G$-group $G[X]$ which fix a given standard orientable quadratic
word with coefficients in $G$. Then we describe some cancellation
properties of  these automorphisms.

 Let $G$ be a group and let $S(X) = 1$ be a regular standard orientable
quadratic  equation over $G:$
\begin{equation}
\label{5}
 \prod_{i=1}^{m}z_i^{-1}c_iz_i\prod_{i=1}^{n}[x_i,y_i] d^{-1} = 1,
 \end{equation}
where $c_i, d$ are non-trivial constants from $G$,  and
  $$X = \{x_i, y_i, z_j \mid i = 1, \dots, n, j = 1, \dots, m\}$$
is the set of  variables.  Sometimes we omit $X$ and write simply
$S = 1$. Denote by
$$C_S =  \{c_1, \ldots, c_m, d\}$$
 the set of constants which occur in the equation $S = 1$.

 Below we  define a {\em basic sequence}
$$\Gamma = (\gamma_1, \gamma_2, \ldots, \gamma_{K(m,n)})$$
 of $G$-automorphisms of the free $G$-group $G[X]$,  each of which
 fixes  the element
$$S_0 = \prod_{i=1}^{m}z_i^{-1}c_iz_i\prod_{i=1}^{n}[x_i,y_i] \in G[X].$$  We
assume that each $\gamma \in \Gamma$    acts identically on all
the generators from $X$ that are not mentioned in the description
of $\gamma$.

\medskip \noindent
Let $m \geqslant 1, n = 0$. In this case $K(m,0) = m-1.$ Put

$ \gamma _{i} \ \ \ : \ z_i\rightarrow
z_i(c_i^{z_i}c_{i+1}^{z_{i+1}}), \ \ \ z_{i+1}\rightarrow
z_{i+1}(c_i^{z_i}c_{i+1}^{z_{i+1}})$, \ \ \ for $i=1,\dots ,m-1$.

\medskip \noindent
Let $m = 0$, $n\geqslant 1$. In this case $K(0,n) = 4n-1.$ Put

\bea \gamma _{4i-3} &:& y_i\rightarrow x_iy_i,  \ \ \ \mbox{ for }
i=1,\dots ,n;\\
\gamma _{4i-2} &:& x_i\rightarrow y_ix_i,  \ \ \ \mbox{ for }
i=1,\dots ,n;\\
\gamma _{4i-1} &:&  y_i\rightarrow x_iy_i,  \ \ \ \mbox{ for }
i=1,\dots ,n;\\
 \gamma _{4i}  &:& x_i\rightarrow
(y_ix_{i+1}^{-1})^{-1}x_i,\ \ \ y_i\rightarrow
y_i^{y_ix_{i+1}^{-1}},\ \ \  x_{i+1}\rightarrow
x_{i+1}^{y_ix_{i+1}^{-1}},\\
&&   y_{i+1}\rightarrow (y_ix_{i+1}^{-1})^{-1}y_{i+1}\ ,
  \mbox{ for }  i=1,\dots ,n-1.
\eea

 \noindent
 Let $m \geqslant 1$, $n\geqslant 1$. In this case $K(m,n) = m + 4n-1.$ Put

\bea \gamma _{i}  &:& z_i\rightarrow
z_i(c_i^{z_i}c_{i+1}^{z_{i+1}}), \ \ \ z_{i+1}\rightarrow
z_{i+1}(c_i^{z_i}c_{i+1}^{z_{i+1}}), \ \ \ \mbox{ for } i=1,\dots ,m-1;\\
\gamma _{m} &:& z_m\rightarrow z_m(c_m^{z_m}x_1^{-1}),\ \ \
x_1\rightarrow x_1^{c_m^{z_m}x_1^{-1}},\ \ \  y_1\rightarrow
(c_m^{z_m}x_1^{-1})^{-1}y_1;\\ \gamma _{m + 4i-3}&:&
y_i\rightarrow x_iy_i,  \ \ \ \mbox{ for } i=1,\dots ,n;\\
\gamma _{m + 4i-2}&:& x_i\rightarrow y_ix_i,  \ \ \ \mbox{ for }
i=1,\dots ,n;\\
\gamma _{m + 4i-1}&:& y_i\rightarrow x_iy_i,  \ \ \  \mbox{ for }
i=1,\dots ,n;\\
\gamma _{m + 4i} &:&  x_i\rightarrow (y_ix_{i+1}^{-1})^{-1}x_i,\ \
\ y_i\rightarrow y_i^{y_ix_{i+1}^{-1}},\ \ \  x_{i+1}\rightarrow
x_{i+1}^{y_ix_{i+1}^{-1}},\\  && y_{i+1}\rightarrow
(y_ix_{i+1}^{-1})^{-1}y_{i+1}, \mbox{  for } i=1,\dots ,n-1. \eea

It is easy to check that each $\gamma \in \Gamma$ fixes the word
$S_0$ as well as
 the word $S$. This shows that  $\gamma$ induces a $G$-automorphism on the group
$G_S = G[X]/{\rm ncl}(S)$. We denote  the induced automorphism
again by $\gamma$, so $\Gamma \subset Aut_G(G_S)$.  Since $S = 1$
is regular, $G_S = G_{R(S)}$. It follows
 that composition of any product of automorphisms from
$\Gamma$ and  a particular  solution $\beta$ of $S = 1$  is again
a solution of $S = 1$.

Observe, that in the case $m \neq 0, n\neq 0$  the basic sequence
of automorphisms $\Gamma$
 contains the basic automorphisms from the other two cases. This allows
us, without loss of
 generality, to formulate some of the  results below only for the case $K(m,n)= m + 4n - 1$.
 Obvious adjustments provide the proper argument in the other cases.
 From now on we  order elements of the set $X$
 in the following  way
 $$z_1<\ldots <z_m<x_1<y_1<\ldots <x_n<y_n.$$
For a word $w \in F(X)$ we denote by $v(w)$ the {\em leading}
variable (the highest variable with respect to the order
introduced above) that occurs in $w$. For $v = v(w)$ denote by
$j(v)$ the following number

 \[ j(v) = \left\{\begin{array}{ll}
                  m+4i, & \mbox{if $v = x_i$ or $v = y_i$ and $i < n$,}\\
                  m+4i-1,  & \mbox{if $v = x_i$ or $v = y_i$ and $i =                  n$,}\\
                  i, &   \mbox{if $v = z_i$ and $n\neq 0$,} \\
                  m-1, &  \mbox{if $v = z_m$, n= 0.}
                 \end{array}
                   \right. \]

The following lemma describes the action of powers of basic
automorphisms from
 $\Gamma$ on $X$. The proof is obvious, and we omit it.

\begin{lm}
\label{le:7.1.29} Let $\Gamma = (\gamma_1, \ldots,
\gamma_{m+4n-1})$ be the basic sequence of automorphisms and $p$
be a  positive integer. Then the following hold: \bea \gamma
_{i}^p \ \ \ &:& \ z_i\rightarrow
z_i(c_i^{z_i}c_{i+1}^{z_{i+1}})^p, \ \ \ z_{i+1}\rightarrow
z_{i+1}(c_i^{z_i}c_{i+1}^{z_{i+1}})^p, \\
&& \qquad\qquad  \mbox{for}\;\; i=1,\dots
,m-1;\\
\gamma _{m}^p \ \  &:&\ z_m\rightarrow z_m(c_m^{z_m}x_1^{-1})^p,\
\ \ x_1\rightarrow x_1^{(c_m^{z_m}x_1^{-1})^p},\ \ \
y_1\rightarrow
(c_m^{z_m}x_1^{-1})^{-p}y_1;\\
\gamma _{m+4i-3}^p &:&\ y_i\rightarrow x_i^py_i,  \ \ \
\mbox{for}\;\;
i=1,\dots ,n;\\
\gamma _{m+4i-2}^p &:&\ x_i\rightarrow y_i^px_i,  \ \ \
\mbox{for}\;\;
i=1,\dots ,n;\\
\gamma _{m+4i-1}^p &:&\ y_i\rightarrow x_i^py_i,  \ \ \
\mbox{for}\;\;
i=1,\dots ,n;\\
\gamma _{m+4i}^p \ \ \ &:&\ x_i\rightarrow
(y_ix_{i+1}^{-1})^{-p}x_i,\ y_i\rightarrow
y_i^{(y_ix_{i+1}^{-1})^p},\\
 &&\ x_{i+1}\rightarrow x_{i+1}^{(y_ix_{i+1}^{-1})^p},\ \ \
y_{i+1}\rightarrow
 (y_ix_{i+1}^{-1})^{-p}y_{i+1},\\
&& \qquad\qquad \mbox{for}\;\; i=1,\dots ,n-1. \eea
 \end{lm}

The $p$-powers of elements that occur in Lemma  \ref{le:7.1.29}
play an important part in what follows, so we describe them in  a
separate definition.
\begin{df}
\label{de:A(gamma)} Let $\Gamma = (\gamma_1, \ldots,
\gamma_{m+4n-1})$ be the basic sequence of automorphism for $S=
1$. For every $\gamma \in \Gamma$ we define the leading term
$A(\gamma)$ as follows:

\smallskip $A(\gamma _{i}) = c_i^{z_i}c_{i+1}^{z_{i+1}}$, \ \ \ for
$i=1,\dots ,m-1$;

\smallskip
 $A(\gamma _{m}) = c_m^{z_m}x_1^{-1};$

\smallskip
 $A(\gamma _{m+4i-3}) = x_i, $ \ \ \ for $i=1,\dots ,n$;

\smallskip $A(\gamma _{m+4i-2}) =  y_i, $ \ \ \ for $i=1,\dots ,n$;

\smallskip $A( \gamma _{m+4i-1}) =  x_i, $ \ \ \ for $i=1,\dots ,n$;

\smallskip $A(\gamma _{m+4i}) = y_ix_{i+1}^{-1} $,  for $i=1,\dots
,n-1$.

\end{df}


Now we introduce  vector notations for automorphisms of particular
type.

 Let ${\mathbb{N}}$ be the set of all positive integers and  ${\mathbb{N}}^k$
 the set of all $k$-tuples of elements from ${\mathbb{N}}$.
 For $s \in {\mathbb{N}}$ and $p \in {\mathbb{N}}^k$ we say that the tuple $p$ is
{\it $s$-large} if every coordinate of $p$  is greater then  $s$.
Similarly, a subset $P \subset  {\mathbb{N}}^k$ is {\em $s$-large}
if every tuple in $P$ is $s$-large.   We say that the set $P$ is
{\em unbounded} if for any $s \in  {\mathbb{N}}$ there exists an
$s$-large tuple in $P$.

Let $\delta = (\delta_1, \ldots, \delta_k)$ be a sequence of
$G$-automorphisms of the group $G[X]$, and $p = (p_1, \ldots,p_k)
\in {\mathbb{N}}^k$.  Then by $\delta^p$ we denote the following
automorphism of $G[X]$: $$\delta^p = \delta_1^{p_1} \cdots
\delta_k^{p_k}.$$

\begin{notation}
  Let $\Gamma = (\gamma_1, \dots, \gamma_K)$ be the
 basic sequence of automorphisms for $S = 1$. Denote by
$\Gamma_{\infty}$ the
 infinite periodic sequence with period $\Gamma$, i.e.,
 $\Gamma_{\infty} = \{\,\gamma_i\,\}_{i \geqslant 1}$ with $\gamma_{i+K}= \gamma_i$.
 For $j \in {\mathbb{N}}$ denote by $\Gamma_j$ the initial segment of
 $\Gamma_{\infty}$ of length $j$. Then for a given $j$ and  $p \in
{\mathbb{N}}^j$ put
 $$ \phi_{j,p} =\stackrel{\leftarrow}{\Gamma}_j^{\stackrel{\leftarrow}{p}} =\gamma_j^{p_j}\gamma_{j-1}^{p_{j-1}} \cdots \gamma_1^{p_1}.$$ Sometimes
we omit $p$ from $ \phi_{j,p}$ and write simply $\phi_j$.
\end{notation}

{\bf Agreement.} {\it From now on we fix an arbitrary positive
multiple $L$ of the number $K = K(m,n)$,  a $2$-large tuple $p \in
{\mathbb{N}}^L$, and the automorphism $\phi = \phi_{L,p}$ {\rm
(}as well as all the automorphism $\phi_j$, $j \leqslant L${\rm
)}. }

\begin{df}
The leading term $A_j=A(\phi _j)$ of the automorphism $\phi_j$ is
defined to be the
 cyclically  reduced form of the word

 \[ \left\{\begin{array}{ll}
                  A(\gamma_j)^{\phi_{j- 1}}, & \mbox{if $j \neq  m+4i-1 +sK$ \ for \ any \
                  $i = 1, \ldots,n, s \in {\mathbb{N}};$},\\
                  y_i^{-\phi _{j-2}} A(\gamma _{j})^{\phi _{j-1}}y_i^
                  {\phi _{j-2}},
                   & \mbox{if $j = m+4i-1 +sK$ \ for \ some \
                   $i = 1, \ldots,n, s \in {\mathbb{N}}.$ }
                 \end{array}
                   \right. \]
\end{df}

\begin{lm}
\label{le:Apower} For every $j \leqslant L$ the element $A_j$ is
not a proper power in $G[X]$. \end{lm}

\begin{proof}
 It is easy to check that $A(\gamma_s)$ from  Definition  \ref{de:A(gamma)}
 is not a proper power for $s = 1, \dots, K.$ Since  $A_j)$ is
the  image of some $A(\gamma_s)$ under an automorphism of $G[X]$
it is not a proper power in $G[X]$. \end{proof}

For words $w, u, v \in G[X]$,  the notation $$\ig{w}{u}{v}$$ \ means
that $w = u \circ w^\prime \circ v$ for some $w^\prime \in G[X]$,
where the length of elements and reduced form defined as in the
free product $G*\left<X\right>$. Similarly, notations $\igl{w}{u}$
\  and   $ \igr{w}{v}$ \ mean that $w = u \circ w^\prime$ and $w =
w^\prime \circ v$. Sometimes we write $\ig{w}{u}{\ast}$ \ or
$\ig{w}{\ast}{v}$ \ when the corresponding words are irrelevant.


 If $n$ is a positive integer and $w \in G[X]$, then by $Sub_n(w)$ we denote the
set of all $n$-subwords of $w$, i.e.,
$$Sub_n(w) = \{u \ \mid \ |u| = n \ and \  w = w_1 \circ u \circ w_2 \ for \ some \
w_1, w_2 \in G[X] \}.$$ Similarly, by $SubC_n(w)$ we denote all
$n$-subwords of the {\it cyclic} word $w$. More generally, if $W
\subseteq G[X]$, then
$$Sub_n(W) = \bigcup_{w \in W} Sub_n(w), \ \ \ SubC_n(W) = \bigcup_{w \in W}
SubC_n(w).$$ Obviously, the set $Sub_i(w)$ ($SubC_i(w)$) can be
effectively reconstructed from $Sub_n(w)$  ($SubC_n(w)$) for $i
\leqslant n$.

 In the following  series of lemmas we write  down explicit
expressions for images of elements of $X$ under the automorphism
 $$\phi_K = \gamma_K^{p_K} \cdots \gamma_1^{p_1}, \ \ \ K = K(m,n).$$
 These lemmas are very easy and straightforward, though tiresome in terms of
 notations. They provide  basic data needed to  prove the implicit function theorem.
 All elements that occur in the lemmas below can be viewed as
 elements (words) from the free group $F(X \cup C_S)$. In particular, the notations $\circ$,
 $\raisebox{1ex}{\ig{w}{u}{v}}$ \ , and $Sub_n(W)$ correspond to the  standard length function on $F(X \cup C_S)$.
Furthermore, until the end of this section we assume that the
elements $c_1, \dots, c_m$ are {\em pairwise different}.

\begin{lm}
\label{le:7.1.zforms}  Let $m \neq 0$,  $K = K(m,n)$, $p = (p_1,
\ldots, p_K)$ be a 3-large tuple, and
 $$\phi_K = \gamma_K^{p_K} \cdots \gamma_1^{p_1}.$$
 The following statements hold.
\begin{enumerate}
\item  [(1)] All automorphisms from $\Gamma$, except for
$\gamma_{i-1}, \gamma_{i}$ {\rm (}if  defined\/{\rm )},  fix
$z_i$, $i=1,\dots, m$. It follows that
 $$z_i^{\phi_K}= \ldots = z_i^{\phi_{i}} \ $$
for $i = 1, \dots,
 m-1$.

\item [(2)]
 Let ${\tilde z}_i = z_i^{\phi_{i-1}}$ \  ($i = 2, \ldots, m$), $ \tilde z_1=z_1.$ Then
$$ {\tilde z}_i  = \ig{ z_i \circ ( c_{i-1}^{{\tilde z}_{i-1}} \circ
c_i^{{z_i}})^{p_{i-1}} }{z_iz_{i-1}^{-1}}{c_iz_i}$$
for  $i = 2, \dots, m$.

\item  [(3)] The reduced forms of the leading terms of the
corresponding automorphisms are listed below:
\bea A_1 &=& \ig{c_1^{z_1} \circ c_2^{z_2}}{z_1^{-1}c_1}{c_2z_2}\ ,\\
&&\qquad\qquad
 \  \ A_2=A_1^{-p_1}c_2^{z_2}A_1^{p_1}c_3^{z_3}, {\rm (}m\geqslant 2{\rm )}\\
SubC_3(A_1) &=& \{z_1^{-1}c_1z_1, \ c_1z_1z_2^{-1}, \
z_1z_2^{-1}c_2, \ z_2^{-1}c_2z_2, \ c_2z_2z_1^{-1}, \
z_2z_1^{-1}c_1\};\\[1ex]
A_{i} &=& \ig{A_{i-1}^{-p_{i-1}}}{z_i^{-1}c_i^{-1}}{c_{i-1}z_{i-1}}\
 c_i^{z_i} \ \ig{A_{i-1}^{p_{i-1}}}{z_{i-1}^{-1}c_{i-1}^{-1}}{c_iz_i}
   \ig{c_{i+1}^{z_{i+1}}}{z_{i+1}^{-1}}{c_{i+1}z_{i+1}},\\
   &&\qquad\qquad i=3,\dots ,m-1,\\
SubC_3(A_{i}) &=& SubC_3(A_{i-1})^{\pm 1}\\
&&\qquad \cup
\{c_{i-1}z_{i-1}z_i^{-1}, \ z_{i-1}z_i^{-1}c_i, \ z_i^{-1}c_iz_i,
\ c_iz_iz_{i-1}^{-1}, z_iz_{i-1}^{-1}c_{i-1}^{-1},\\
&&\qquad \
c_iz_iz_{i+1}^{-1}, \ z_iz_{i+1}^{-1}c_{i+1}, \
z_{i+1}^{-1}c_{i+1}z_{i+1}, \ c_{i+1}z_{i+1}z_i^{-1}, \
z_{i+1}z_i^{-1}c_i^{-1} \};\\
A_2 & =&A_1^{-p_1}c_2^{z_2}A_1^{p_1}x_1^{-1}  (m=2);\\
[2ex] A_{m} &=&
\ig{A_{m-1}^{-p_{m-1}}}{z_m^{-1}c_m^{-1}}{c_{m-1}z_{m-1}} \
c_m^{z_m} \
\ig{A_{m-1}^{p_{m-1}}}{z_{m-1}^{-1}c_{m-1}^{-1}}{c_mz_m} \
x_1^{-1}\\
&&\qquad\qquad  (n \neq 0, m>2),\\
SubC_3(A_{m}) &=& SubC_3(A_{m-1})^{\pm 1}\\
&&\qquad \cup
\{c_{m-1}z_{m-1}z_m^{-1}, \ z_{m-1}z_m^{-1}c_m, \
z_{m-1}^{-1}c_mz_m,\\
&&\qquad\qquad \ c_mz_mz_{m-1}^{-1}, \ c_mz_mx_1^{-1},
z_mx_1^{-1}z_m^{-1}, \ x_1^{-1}z_m^{-1}c_m^{-1} \}.
\eea

\item   [(4)] The reduced forms of $z_i^{\phi_{i-1}},
z_i^{\phi_{i}}$ are listed below:

\bea
 z_1^{\phi_K} &=& z_1^{\phi _{1}} =c_1  \ig{z_1 c_2^{z_2}}{z_1z_2^{-1}}{c_2z_2}
\ig{A_1^{p_1-1}}{z_1^{-1}c_1}{c_2z_2}\ \;\; (m\geq 2)\ ,\\
SubC_3(z_1^{\phi_K}) &=& \{ \ c_1z_1z_2^{-1}, \ z_1z_2^{-1}c_2, \
z_2^{-1}c_2z_2, \ c_2z_2z_1^{-1}, \ z_2z_1^{-1}c_1, \
z_1^{-1}c_1z_1\};\\
z_i^{\phi_{i-1}} &=& z_i
\ig{A_{i-1}^{p_{i-1}}}{z_{i-1}^{-1}c_{i-1}^{-1}}{c_iz_i} \ ,\\
z_{i}^{\phi _{K}} &=&   z_{i}^{\phi _{i}}   = c_i z_i \
\ig{A_{i-1}^{p_{i-1}}}{z_{i-1}^{-1}c_{i-1}^{-1}}{c_iz_i} \
c_{i+1}^{z_{i+1}}
\ig{A_{i}^{p_{i}-1}}{z_i^{-1}c_i^{-1}}{c_{i+1}z_{i+1}}  \\
&&\qquad ( i=3,\dots ,m-1),\\
Sub_3(z_i^{\phi_K}) &=& SubC_3(A_{i-1}) \cup SubC_3(A_{i}) \cup \{
c_iz_iz_{i-1}^{-1}, \ z_iz_{i-1}^{-1}c_{i-1}^{-1},\\
 && \ c_iz_iz_{i+1}^{-1}, \ z_iz_{i+1}^{-1}c_{i+1}, \ z_{i+1}^{-1}c_{i+1}z_{i+1}, \
c_{i+1}z_{i+1}z_i^{-1}, \ z_{i+1}z_i^{-1}c_i^{-1} \} \ ;\\
z_m^{\phi_K} &=& z_m^{\phi_{m-1}} =  z_m
\ig{A_{m-1}^{p_{m-1}}}{z_{m-1}^{-1}c_{m-1}^{-1}}{c_mz_m}\ , \ \ \
(n = 0),\\
Sub_3(z_m^{\phi_K}) &=& SubC_3(A_{m-1}) \cup \{
c_mz_mz_{m-1}^{-1}, \ z_mz_{m-1}^{-1}c_{m-1}^{-1} \} {(\mbox{when}
\ n = 0)}\ ; \eea \bea z_m^{\phi_K} &=& c_mz_m^{\phi_{m}} =  c_m
z_{m} \ig{A_{m-1}^{p_{m-1}}}{z_{m-1}^{-1}c_{m-1}^{-1}}{c_mz_m}\
x_1^{-1} \ \ig{A_{m}^{p_{m}-1}}{z_m^{-1}c_m^{-1}}{z_mx_1^{-1}}\ \
\ (n
\neq 0),\\
Sub_3(z_m^{\phi_K}) &=& Sub_3(z_m^{\phi_K})_{(when \ n = 0)} \cup
\{c_mz_mx_1^{-1}, \ z_mx_1^{-1}z_m^{-1}, \
x_1^{-1}z_m^{-1}c_m^{-1} \}. \eea

\item [(5)] The elements $z_i^{\phi _K}$ have the following
properties:
$$z_i^{\phi _K}=c_iz_i\hat z_i \ \ \ (i=1,\dots ,m-1),$$ where $\hat
z_i$ is a word in the alphabet $\{ c_1^{z_1}, \ldots,
c_{i+1}^{z_{i+1}}, \}$ which begins with $c_{i-1}^{-z_{i-1}}$, if
$i\neq 1$, and with $c_2^{z_2}$, if $i=1$;

\medskip\noindent
$z_m^{\phi _K}=z_m\hat z_m  \ \ \ (n=0)$, where $\hat z_m$ is a
word in the alphabet $\left\{ c_1^{z_1}, \ldots,
c_{m}^{z_{m}}\right \};$

\medskip\noindent
$z_m^{\phi _K}=c_mz_m\hat z_m  \ \ \ (n \neq 0)$, where $\hat z_m$
is a word in the alphabet $$\{ c_1^{z_1}, \ldots, c_{m}^{z_{m}},
x_1 \};$$
 Moreover, if $m\geq 3$, the word $(c_m^{z_m})^{\pm 1}$ occurs in  $z_i^{\phi _K} \ \ \ (i = m-1,m)$
 only  as a part of the subword $\left(\prod _{i=1}^m c_i^{z_i}\right)^{\pm 1}.$
\end{enumerate}
\end{lm}

\begin{proof}
 (1) is obvious. We prove (2) by induction. For $i\geqslant 2$,
$${\tilde z}_i = z_i^{\phi_{i-1}} = z_i^{
 \gamma_{i-1}^{p_{i-1}}\phi_{i-2
}}.$$ Therefore
$${\tilde z}_i = z_i ( c_{i-1}^{{\tilde z}_{i-1}}
c_i^{{z_i}})^{p_{i-1}} = z_i \circ ( c_{i-1}^{{\tilde z}_{i-1}}
\circ c_i^{{z_i}})^{p_{i-1}}, $$
and the claim follows by induction.

 Now we prove (3) and (4) simultaneously.  By the straightforward verification one has:
\medskip
$A_1 = \ig{c_1^{z_1} \circ c_2^{z_2}}{z_1^{-1}}{z_2}$;

\medskip
$z_1^{\phi_1} = z_1^{\gamma _1^{p_1}} =
z_1(c_1^{z_1}c_2^{z_2})^{p_1} = \ig{c_1 \circ z_1 \circ c_2^{z_2}
\circ A_1^{p_1-1}}{c_1}{z_2}$.

Denote by cycred $(w)$ the cyclically reduced form of $w$.

\medskip
$A_i =$ cycred $( \left (c_i^{z_i}c_{i+1}^{z_{i+1}} \right
)^{\phi_{i-1}})= \ig {c_i^{{\tilde z}_i} \circ
c_{i+1}^{z_{i+1}}}{z_i^{-1}}{z_{i+1}}\ (i\leq m-1).$

\medskip \noindent
Observe that in the notation above

\medskip
${\tilde z}_i =  z_iA_{i-1}^{p_{i-1}} \ (i\geq 2).$

\medskip \noindent
This shows that we can rewrite  $A(\phi_{i})$ as follows:

\medskip
$A_i =  A_{i-1}^{-p_{i-1}}\circ c_i^{z_i}\circ
  A_{i-1}^{p_{i-1}} \circ
  c_{i+1}^{z_{i+1}}$,

\medskip \noindent
 beginning with $z_i^{-1}$ and ending with $z_{i+1}$ \ \ ( $i=2,\dots ,m-1$);

\medskip
$A_{m} =$ cycred $(c_m^{\tilde z_m}x_1^{-1}) = c_m^{\tilde
z_m}x_1^{-1} = A_{m-1}^{-p_{m-1}} \circ c_m^{z_m}\circ
A_{m-1}^{p_{m-1}} \circ x_1^{-1}\ (m\geq 2).$

\medskip \noindent
beginning with $z_m^{-1}$ and ending with $x_1^{-1}$  \ \ ( $n
\neq 0$).

$$z_i^{\phi_{i-1}} = \left (z_i
(c_{i-1}^{z_{i-1}}c_i^{z_i})^{p_{i-1}} \right )^{\phi_{i-2}} = z_i
(c_{i-1}^{{\tilde z}_{i-1}}c_i^{z_i})^{p_{i-1}} = z_i \circ
A_{i-1}^{p_{i-1}},$$ beginning with $z_i$ and ending with $z_i$;
\bea z_i^{\phi_{i}} &=& \left
(z_i(c_i^{z_i}c_{i+1}^{z_{i+1}})^{p_{i}} \right )^{\phi_{i-1}}
= {\tilde z}_i(c_i^{{\tilde z}_i}c_{i+1}^{z_{i+1}})^{p_{i}}\\
&=& c_i \circ {\tilde z}_i \circ c_{i+1}^{z_{i+1}} \circ
(c_i^{{\tilde
z}_i}c_{i+1}^{z_{i+1}})^{p_{i}-1}\\
&=& c_i \circ z_i \circ A_{i-1}^{p_{i-1}} \circ c_{i+1}^{z_{i+1}}
\circ A_{i}^{p_{i}-1}, \eea beginning with $c_i$ and ending with
$z_{i+1}$ \ \  ($i = 2, \ldots, m-1$); \bea z_m^{\phi_{m}} &=&
\left (z_m(c_m^{z_m}x_1^{-1})^{p_{m}} \right )^{\phi_{m-1}}
=  {\tilde z}_m(c_m^{{\tilde z}_m}x_1^{-1})^{p_{m}}\\
&=& c_m {\tilde z}_m x_1^{-1} (c_m^{{\tilde z}_m}x_1^{-1})^{p_{m}-1}\\
&=& c_m \circ z_{m} \circ A_{m-1}^{p_{m-1}}\circ x_1^{-1} \circ
A_{m}^{p_{m}-1}\ \  \ (n \neq 0), \eea beginning with $c_m$ and
ending with $x_1^{-1}.$  This proves the lemma. \end{proof}

In the following two lemmas  we describe the  reduced expressions
of the elements $x_1^{\phi_K}$ and $y_1^{\phi_K}$.
\begin{lm}
\label{le:7.1.x1formsm0}
 Let $m = 0$,  $K = 4n-1$, $p = (p_1, \ldots, p_K)$ be a $3$-large tuple, and
 $$\phi_K = \gamma_K^{p_K} \cdots \gamma_1^{p_1}.$$

\begin{enumerate}
\item  [(1)] All automorphisms from $\Gamma $, except for  $\gamma
_{2}, \gamma _{4} $, fix $x_1$,  and all automorphisms from
$\Gamma $, except  for $\gamma _{1}, \gamma_3, \gamma _{4} $,  fix
$y_1$. It follows that

 \medskip
$ x_1^{\phi_K} = x_1^{\phi_4}, \ y_1^{\phi_K} = y_1^{\phi_4} \ \ \
(n \geqslant 2).$

\item [(2)] Below we list the reduced forms of the leading terms
of the corresponding automorphisms  (the words on the right are
reduced as written)

\medskip
$A_1 = x_1$;

\medskip
$A_2 =  x_1^{p_1}y_1 = A_1^{p_1}\circ y_1 $ ;

\medskip
$A_3 = \ig{A_2^{p_2-1}}{x_1^2}{x_1y_1} \ x_1^{p_1+1}  y_1$,

\medskip
$SubC_3(A_3) = SubC_3(A_2) = \{x_1^3, \  x_1^2y_1, \ x_1y_1x_1,
y_1x_1^2 \} \ ;$

\medskip
$A_4 =  \ig{\left(\ig{A_2^{p_2}}{x_1^2}{x_1y_1} \
x_1\right)^{p_3}}{x_1^2}{y_1x_1} \ig{A_2}{x_1^2}{x_1y_1}\
x_2^{-1}\ \ (n\geq 2),$

\medskip
$SubC_3(A_4) = SubC_3(A_2) \cup \{x_1y_1x_2^{-1}, \
y_1x_2^{-1}x_1, \ x_2^{-1}x_1^2\}.$

 \item [(3)] Below we list reduced forms of $x_1^{\phi _j}, y_1^{\phi _j}$ for $j = 1,
 \ldots, 4$:

\medskip
$x_1^{\phi_1} = x_1$;

\medskip
$y_1^{\phi_1} = x_1^{p_1} y_1; $

\medskip
$x_1^{\phi_2} =  \ig{A_2^{p_2}}{x_1^2}{x_1y_1}\  x_1;$

\medskip
$y_1^{\phi_2} = x_1^{p_1}y_1;$

\medskip
$x_1^{\phi_3} =  x_1^{\phi_2} = \ig{A_2^{p_2}}{x_1^2}{x_1y_1}\
x_1$;

\medskip
$Sub_3(x_1^{\phi_K})_{(when \ n = 1)} = SubC_3(A_2)$;

\medskip
$y_1^{\phi_3} =   \ig{(\ig{A_2^{p_2}}{x_1^2}{x_1y_1} \
x_1)^{p_3}}{x_1^2}{x_1y_1} \ x_1^{p_1}y_1 ;$

\medskip
$Sub_3(y_1^{\phi_K})_{(when \ n = 1)} = SubC_3(A_2);$

 \medskip
 $x_1^{\phi_4} = x_1^{\phi_K} = \ig{A_4^{-(p_4-1)}}{x_2y_1^{-1}}{x_1^{-2}} \  x_2
 \ \ig{A_2^{-1}}{y_1^{-1} x_1^{-1}}{x_1^{-2}}  \ig{(x_1^{-1} \
 \ig{A_2^{-p_2}}{y_1^{-1}x_1^{-1}}{x_1^{-2}})^{p_3-1}}{x_1^{-1}y_1^{-1}x_1^{-1}}{x_1^{-2}}$
 \ \  $(n \geqslant 2)$,

 \medskip
 $Sub_3(x_1^{\phi_K}) = SubC_3(A_4)^{-1} \cup SubC_3(A_2)^{-1} \cup
 \{x_1^{-2}x_2, \ x_1^{-1}x_2y_1^{-1}, $ $$ x_2y_1^{-1}x_1^{-1}, \
 x_1^{-3},
 \ x_1^{-2}y_1^{-1}, x_1^{-1}y_1^{-1}x_1^{-1} \}\ \ \ (n\geqslant 2);$$

\medskip
$y_1^{\phi_4} =  \ig{A_4^{-(p_4-1)}}{x_2y_1^{-1}}{x_1^{-1}} \ x_2
\ \ig{A_4^{p_4}}{x_1^2}{y_1x_2^{-1}}\ \ \ (n\geqslant 2), $

\medskip
$Sub_3(y_1^{\phi_K}) = SubC_3(A_4)^{\pm 1} \cup \{x_1^{-2}x_2, \
x_1^{-1}x_2x_1, \ x_2x_1^2 \} \ \ \ (n\geqslant 2).$
\end{enumerate}
 \end{lm}
\begin{proof} (1) follows directly from definitions.

To show (2) observe that

\medskip
$A_1 = A(\gamma_1) = x_1$;

\medskip
$x_1^{\phi_1} = x_1$;

\medskip
$y_1^{\phi_1} = x_1^{p_1}y_1 = A_1^{p_1} \circ y_1.$

\medskip\noindent
Then \bea A_2 &=& {\rm cycred}(A(\gamma_2)^{\phi_1}) = {\rm
cycred}(y_1^{\phi_1}) = x_1^{p_1} \circ
y_1 = A_1^{p_1} \circ y_1;\\
x_1^{\phi_2} &=& (x_1^{\gamma_2^{p_2}})^{\gamma_1^{p_1}}=
(y_1^{p_2}x_1)^{\gamma_1^{p_1}}
= (x_1^{p_1}y_1)^{p_2} x_1 =A_2^{p_2} \circ x_1;\\
y_1^{\phi_2} &=& (y_1^{\gamma_2^{p_2}})^{\gamma_1^{p_1}}=
y_1^{\gamma_1^{p_1}} = x_1^{p_1}y_1 = A_2.\eea Now \bea A_3
&=&{\rm cycred}(y_1^{-\phi_1} A(\gamma_3)^{\phi_2}y_1^{\phi_1})=
{\rm cycred}((x_1^{p_1}y_1)^{-1} x_1^{\phi_2} (x_1^{p_1}y_1))\\
&=&{\rm cycred}(
(x_1^{p_1}y_1)^{-1}(x_1^{p_1}y_1)^{p_2}x_1(x_1^{p_1}y_1))\\
& =&   (x_1^{p_1}y_1)^{p_2-1} x_1^{p_1+1}y_1
 =  A_2^{p_2-1}\circ A_1^{p_1+1} \circ y_1. \eea

\medskip\noindent
It follows that

\medskip
$x_1^{\phi_3} = (x_1^{\gamma_3^{p_3}})^{\phi_2} = x_1^{\phi_2}$;

\medskip
$y_1^{\phi_3} = (y_1^{\gamma_3^{p_3}})^{\phi_2}
=(x_1^{p_3}y_1)^{\phi_2} = (x_1^{\phi_2})^{p_3}y_1^{\phi_2}
=(A_2^{p_2} \circ x_1)^{p_3}\circ
 A_2.$

 \medskip\noindent
 Hence
\[
A_4 = {\rm cycred}(A(\gamma_4)^{\phi_3})
={\rm cycred}(
y_1^{\phi_3}x_2^{-\phi_3})  = (A_2^{p_2} \circ x_1)^{p_3}\circ
 A_2\circ x_2^{-1}.
 \]
 Finally:
\bea x_1^{\phi_4} &=& (x_1^{\gamma_4^{p_4}})^{\phi_3}  =\left
((y_1x_2^{-1})^{-p_4}x_1 \right )^{\phi_3}\\ &=&
 \left ((y_1x_2^{-1})^{\phi_3})\right )^{-p_4}x_1^{\phi_3}= A_4^{-p_4}A_2^{p_2} \circ x_1\\
& = &  A_4^{-(p_4-1)} \circ x_2
 \circ A_2^{-1} \circ (x_1^{-1} \circ
A_2^{-p_2})^{p_3-1}\\
y_1^{\phi_4} &=& (y_1^{\gamma_4^{p_4}})^{\phi_3} =(y_1^{(y_1x_2^{-1})^{p_4}})^{\phi_3}\\
& = &  \left ((y_1x_2^{-1})^{\phi_3}\right ) ^{-p_4} y_1^{\phi_3}
\left
((y_1x_2^{-1})^{\phi_3}\right )^{p_4}\\
& = & A_4^{-p_4} y_1^{\phi_3}  A_4^{p_4}=
A_4^{-(p_4-1)} A_4^{-1} y_1^{\phi_3} A_4^{p_4}\\
&=&
 A_4^{-(p_4-1)}\circ x_2 \circ
A_4^{p_4} . \eea This proves the lemma. \end{proof}

\begin{lm}
\label{le:7.1.x1formsmneq0}
 Let $m \neq 0$, $n \neq 0$,  $K = m + 4n-1$, $p = (p_1, \ldots, p_K)$ be a $3$-large tuple, and
 $$\phi_K = \gamma_K^{p_K} \cdots \gamma_1^{p_1}.$$

\begin{enumerate}
\item [(1)] All automorphisms from $\Gamma $ except  for
$\gamma_{m}, \gamma _{m+2}, \gamma _{m+4} $ fix $x_1$; and all
automorphisms from $\Gamma $ except for $\gamma_{m}, \gamma
_{m+1}, \gamma_{m+3}, \gamma _{m+4} $ fix $y_1$.  It follows that

\medskip
$x_1^{\phi_K} = x_1^{\phi_{m+4}}, \  y_1^{\phi_K}
=y_1^{\phi_{m+4}} \ \ \ (n \geqslant 2).$

\item [(2)] Below we list the reduced forms of the leading terms
of the corresponding automorphisms  (the words on the right are
reduced as written)
\bea
A_{m+1} &=&  x_1,\\
A_{m+2} &=& y_1^{\phi_{m+1}}\\
&&
=\ig{A_{m}^{-p_{m}}}{x_1z_m^{-1}}{c_mz_m} \ x_1^{p_{m+1}} y_1,\\
SubC_3(A_{m+2}) &=& SubC_3(A_{m})^{-1}\\
&&\qquad \cup \{c_mz_mx_1, \
z_mx_1^2, \ x_1^3, \ x_1^2y_1, \ x_1y_1x_1, \ y_1x_1z_m^{-1} \};\\
A_{m+3} &=& \ig{A_{m+2}^{p_{m+2}-1}}{x_1z_m^{-1}}{x_1y_1}
 \ig{A_{m}^{-p_{m}}}{x_1z_m^{-1}}{c_mz_m}\  x_1^{p_{m+1}+1}
  y_1,\\
SubC_3(A_{m+3}) &=& SubC_3(A_{m+2});\\
A_{m+4} &=& \ig{A_{m}^{-p_{m}}}{x_1z_m^{-1}}{c_mz_m} \\
&&\qquad \circ \left
(x_1^{p_{m+1}} y_1 \
\ig{A_{m+2}^{p_{m+2}-1}}{x_1z_m^{-1}}{x_1y_1}\ig{A_{m}^{-p_{m}}}{x_1z_m^{-1}}{c_mz_m}\
x_1\right )^{p_{m+3}}
 x_1^{p_{m+1}} y_1  x_2^{-1}\\
 && \qquad\qquad (n \geqslant 2),\\
SubC_3(A_{m+4}) &=& SubC_3(A_{m+2}) \cup
\{x_1y_1x_2^{-1}, \
 y_1x_2^{-1}x_1, \ x_2^{-1}x_1z_m^{-1} \}\\
 &&\qquad\qquad (n
\geqslant 2).\eea

 \item [(3)]Below we list reduced forms of $x_1^{\phi _j}, y_1^{\phi _j}$ for
$j=m,\dots ,m+4$ and their expressions via the leading terms:

$x_1^{\phi _{ m}}=  A_{m}^{-p_{m}}\circ x_1 \circ A _{m}^{p_{m}}$,
\medskip

$y_1^{\phi _{m}}=  A_{m}^{-p_{m}} \circ y_1$,
\medskip
$x_1^{\phi _{m+1}}=x_1^{\phi _{m}}$,

\medskip
$y_1^{\phi _{ m+1}}= A _{m}^{-p_{m}} \circ x_1^{p_{m+1}} \circ
y_1$,

\medskip
$x_1^{\phi _{ m+2}}= {x_1^{\phi_K}}_{(when \ n = 1)}
=\ig{A_{m+2}^{p_{m+2}}}{x_1z_m^{-1}}{x_1y_1}\ig{A_{m}^{-p_{m}}}{x_1z_m^{-1}}{c_mz_m}
\ x_1 \ \ig{A_{m}^{p_{m}}}{z_m^{-1}c_m^{-1}}{z_mx_1^{-1}}\ $,

\medskip\noindent
$Sub_3(x_1^{\phi_K})_{(when \ n = 1)} = SubC_3(A_{m+2}) \cup
SubC_3(A_{m}) \cup \{z_mx_1z_m^{-1}, \ x_1z_m^{-1}c_m^{-1} \};$

\medskip
$y_1^{\phi _{m+2}}=y_1^{\phi _{m+1}}$,

\medskip
$x_1^{\phi _{m+3}}=x_1^{\phi _{m+2}}$,
\medskip
\begin{multline*} y_1^{\phi _{ m+3}}= {y_1^{\phi_K}}_{(when \ n
= 1)} =\ig{A_{m}^{-p_{m}}}{x_1z_m^{-1}}{c_mz_m}  \\ \left
(x_1^{p_{m+1}}y_1  \ig{A_{m+2}^{p_{m+2}-1}}{x_1z_m^{-1}}{x_1y_1}
\ig{A_{m})^{-p_{m}}}{x_1z_m^{-1}}{c_mz_m} \ x_1\right )^{p_{m+3}}
x_1^{p_{m+1}} y_1.\end{multline*}

\medskip
\vspace{.2cm} $Sub_3(y_1^{\phi_K})=_{(when \ n = 1)}
SubC_3(A_{m+2});$

\medskip
\begin{multline*} x_1^{\phi_{m+4}} = {x_1^{\phi_K}}_{(when \ n \geqslant 2)}= \ig{A_{m+4}^{-p_{m+4}+1}}{x_2y_1^{-1}}{z_mx_1^{-1}} \
x_2y_1^{-1}x_1^{-p_{m+1}} \circ \\ \left( x_1^{-1} \
\ig{A_{m}^{p_{m}}}{z_m^{-1}c_m^{-1}}{z_mx_1^{-1}}
\ig{A_{m+2}^{-p_{m+2}}}{y_1^{-1}x_1^{-1}}{z_mx_1^{-1}}\ y_1^{-1}
x_1^{-p_{m+1}}\right )^{p_{m+3}-1}
\ig{A_{m}^{p_{m}}}{z_m^{-1}c_m^{-1}}{z_mx_1^{-1}} \ \ \ (n
\geqslant 2),\end{multline*}

\medskip\vspace{.2cm}
$Sub_3(x_1^{\phi_K}) = SubC_3(A_{m+2})^{-1} \cup \{z_mx_1^{-1}x_2,
\ x_1^{-1}x_2y_1^{-1}, \ x_2y_1^{-1}x_1^{-1} \} ;$

\vspace{.2cm}
 $y_1^{\phi _{m+4}} = {y_1^{\phi_K}}_{(when \ n \geqslant 2)} =\ig{A_{m+4}^{-(p_{m+4}-1)}}{x_2y_1^{-1}}{z_mx_1^{-1}} \ x_2 \
\ig{A_{m+4}^{p_{m+4}}}{x_1z_m^{-1}}{y_1x_2^{-1}} \ \ \ (n\geqslant
2)$,

\medskip
$Sub_3(y_1^{\phi_K}) = SubC_3(A_{m+4})^{\pm 1} \cup
\{z_mx_1^{-1}x_2, \ x_1^{-1}x_2x_1, x_2x_1z_m^{-1}\}$ \ \newline
$(n \geqslant 2).$
\end{enumerate}
 \end{lm}

\begin{proof}
 Statement (1) follows immediately from definitions of
automorphisms of $\Gamma$.

\medskip \noindent
We prove formulas in the second and third statements
simultaneously:
$$x_1^{\phi _{m}}= \left (x_1^{(c_m^{{
z}_m}x_1^{-1})^{p_{m}}} \right )^{\phi_{m-1}}= x_1^{A_{m}^{p_{m}}}
=  A _{m}^{-p_{m}}\circ x_1 \circ A _{m}^{p_{m}},$$ beginning with
$x_1$ and ending with $x_1^{-1}$.
$$y_1^{\phi _{m}}= \left ((c_m^{{
z}_m}x_1^{-1})^{-p_{m}}y_1\right )^{\phi_{m-1}}= A_{m}^{-p_{m}}
\circ y_1,$$ beginning with $x_1$ and ending with $y_1$. Now
 $A_{m+1}$ is the cyclically reduced form of
 $ A(\gamma_{m+1})^{\phi_{m}} = x_1^{\phi_{m}} = A
_{m}^{-p_{m}}\circ x_1 \circ A _{m}^{p_{m}}.$

\medskip\noindent
$A _{m+1}=x_1.$

\bea x_1^{\phi _{m+1}} &=& x_1^{\phi _{m}},\\
y_1^{\phi _{ m+1}} &=& \left (y_1^{\gamma_{m+1}^{p_{m+1}}}\right
)^{\phi_{m}} = (x_1^{p_{m+1}}y_1)^{\phi_{m}}\\
&=& (x_1^{\phi_{m}})^{p_{m+1}}y_1^{\phi_{m}}\\
 &=& A_{m}^{-p_{m}} \circ x_1^{p_{m+1}} \circ y_1, \eea beginning
with $x_1$ and ending with $y_1$, moreover, the element that
cancels in reducing

$A_{m+1}^{p_{m+1}} A_{m}^{-p_{m}} y_1$ is equal to
$A_{m}^{p_{m}}.$

\medskip
$A_{m+2} = {\rm cycred}(A(\gamma_{m+2})^{\phi_{m+1}}) = {\rm
cycred}(y_1^{\phi_{m+1}})= A_{m}^{-p_{m}} \circ x_1^{p_{m+1}}
\circ y_1$,

\medskip\noindent
 beginning with $x_1$ and ending with $y_1$.
\bea x_1^{\phi _{ m+2}} &=& \left
(x_1^{\gamma_{m+2}^{p_{m+2}}}\right
)^{\phi_{m+1}}\\
&=&  (y_1^{\phi _{ m+1}})^{p_{m+2}} x_1^{\phi _{
m+1}}\\
&=&  A _{m+2}^{p_{m+2}} \circ  A _{m}^{-p_{m}}\circ x_1
\circ A _{m}^{p_{m}}\\
& =& A_{m}^{-p_{m}} \circ \left (x_1^{p_{m+1}} \circ y_1 \circ A
_{m+2}^{p_{m+2}-1}  \circ  A _{m})^{-p_{m}}\circ x_1\right ) \circ
A _{m}^{p_{m}}, \eea beginning with $x_1$ and ending with
$x_1^{-1}$;
\bea y_1^{\phi _{m+2}} & =& y_1^{\phi _{m+1}}.\\
A_{m+3} &=&  y_1^{-\phi_{m+1}} x_1^{\phi_{m+2}} y_1^{\phi_{m+1}}\\
&  =& A _{m+2}^{p_{m+2}-1}\circ A _{m}^{-p_{m}}\circ
x_1^{p_{m+1}+1}
  \circ y_1,
\eea
beginning with $x_1$ and ending with $y_1$;
\bea
x_1^{\phi _{m+3}} &=& x_1^{\phi _{m+2}},\\
\ && \\
y_1^{\phi _{ m+3}} &=& (x_1^{\phi _{ m+2}})^{p_{m+3}} y_1^{\phi _{
m+1}}\\
&=& A_{m}^{-p_{m}} \circ \left (x_1^{p_{m+1}} \circ y_1 \circ A
_{m+2}^{p_{m+2}-1}  \circ  A _{m}^{-p_{m}}\circ
x_1\right )^{p_{m+3}} \circ  x_1^{p_{m+1}} \circ y_1,\\
\eea
 beginning with $x_1$ and ending with $y_1$. Finally,
$$A_{m+4} = {\rm
cycred}(A(\gamma_{m+4})^{\phi_{m+3}})= {\rm
cycred}((y_1x_2^{-1})^{\phi_{m+3}}) = y_1^{\phi_{m+3}}\circ
x_2^{-1},$$ beginning with $x_1$ and ending with $x_2^{-1}$; \bea
x_1^{\phi_{m+4}} &=& \left ((y_1x_2^{-1})^{-p_{m+4}}x_1\right
)^{\phi_{m+3}}\\
&=& \left ((x_2 y_1^{-\phi_{m+3}}\right
)^{p_{m+4}}x_1^{\phi_{m+3}}\\
&=&   \left ((x_2
y_1^{-\phi_{m+1}}(x_1^{\phi_{m+2}})^{-p_{m+3}}\right
)^{p_{m+4}}x_1^{\phi_{m+2}}\\
&=& (x_2y_1^{-\phi
_{m+3}})^{p_{m+4}-1}\circ x_2\circ y_1^{-1} \circ x_1^{-p_{m+1}}
\\ && \qquad \circ \left( x_1^{-1} \circ A_{m}^{p_{m}} \circ
A_{m+2}^{-p_{m+2}} \circ y_1^{-1} \circ x_1^{-p_{m+1}}\right
)^{p_{m+3}-1}\circ A_{m}^{p_{m}},\eea

\medskip\noindent
beginning with $x_2$ and ending with $x_1^{-1}$, moreover,  the
element that is cancelled out is $x_1^{\phi _{m+2}}$. Similarly,
\bea
y_1^{\phi _{m+4}} &=& (x_2y_1^{-\phi _{m+3}})^{p_{m+4}}y_1^{\phi
_{m+3}}(y_1^{\phi _{m+3}}x_2^{-1})^{p_{m+4}}\\
&=& (x_2y_1^{-\phi _{m+3}})^{p_{m+4}-1} \circ x_2 \circ (y_1^{\phi
_{m+3}}x_2^{-1})^{p_{m+4}}\\
&=& A _{m+4}^{-(p_{m+4}-1)} \circ x_2 \circ A _{m+4}^{p_{m+4}},
\eea
  beginning with $x_2$ and ending with $x_2^{-1}$,
moreover, the element that is cancelled out is $y_1^{\phi
_{m+3}}$.

This proves the lemma. \end{proof}

In the following lemmas  we describe the  reduced expressions of
the elements $x_i^{\phi_j}$ and $y_i^{\phi_j}$.
\begin{lm}
\label{le:7.1.xiforms}
 Let $n \geqslant 2$, $K = K(m,n)$, $p = (p_1, \dots,p_K)$ be a 3-large tuple,  and
 $$\phi_K = \gamma_K^{p_K} \ldots \gamma_1^{p_1}.$$
 Then for any $i$, $n \geqslant i \geqslant 2$,  the following holds:

\begin{enumerate}
\item  [(1)] All automorphisms from $\Gamma $, except for
$\gamma_{m+4(i-1)}, \gamma_{m+4i-2},
 \gamma_{m+4i}$ fix $x_i$, and all automorphisms from $\Gamma $, except for $\gamma_{m+4(i-1)}$,
 $\gamma_{m+4i-3}$,
$\gamma_{m+4i-1}$, $\gamma_{m+4i}$ fix $y_i$. It follows that

\medskip
$x_i^{\phi_K} = x_i^{\phi_{K-1}} = \ldots = x_i^{\phi_{m+4i}}, $

\medskip
$y_i^{\phi_K} = y_i^{\phi_{K-1}} = \ldots = y_i^{\phi_{m+4i}}. $

\item  [(2)] Let $\tilde  y_{i}=y_{i}^{\phi _{ m+4i-1}}$. Then
$${\tilde y}_{i} = \ig{{\tilde y}_{i}}{x_iy_{i-1}^{-1}}{x_iy_i}\ $$
where (for $i = 1$) we assume that  $y_0 = x_1^{-1}$ for $m = 0$,
and $y_0 = z_m$ for $m \neq 0$;

 \item [(3)] Below we list the reduced forms of the leading terms
of the corresponding automorphisms.   Put $q_j = p_{m+4(i-1)+j}$
for $j = 0, \ldots, 4$. In the formulas below we assume that $y_0
= x_1^{-1}$ for $m = 0$, and $y_0 = z_m$ for $m \neq 0$.
\bea
A_{m+4i-4}  &=& \ig{{\tilde y}_{i-1}\circ
x_i^{-1}}{x_{i-1}y_{i-2}^{-1}}{x_{i-1}y_{i-1}x_i^{-1}} ,\\
SubC_3(A_{m+4i-4}) &=& Sub_3({\tilde y}_{i-1})\\
&&\qquad \cup
\{x_{i-1}y_{i-1}x_i^{-1}, \ y_{i-1}x_i^{-1}x_{i-1}, \
x_i^{-1}x_{i-1}y_{i-2}^{-1} \} ;\\
A_{m+4i-3} &=&  x_i; \\
A_{m+4i-2}
&=&\ig{A_{m+4i-4}^{-q_0}}{x_iy_{i-1}^{-1}}{y_{i-2}x_{i-1}^{-1}}
 \  x_i^{q_1} y_i,\\
SubC_3(A_{m+4i-2}) &=& SubC_3(A_{m+4i-4})\\
&&\qquad \cup
\{y_{i-2}x_{i-1}^{-1}x_i, \ x_{i-1}^{-1}x_i^2, \ x_i^2y_i, \
x_iy_ix_i, \ y_ix_iy_{i-1}^{-1} x_i^3\} ;\\
A_{m+4i-1} &=& \ig{A_{m+4i-2}^{q_2-1}}{x_iy_{i-1}^{-1}}{x_iy_i}
\ig{A_{m+4i-4}^{-q_0}}{x_iy_{i-1}^{-1}}{y_{i-2}x_{i-1}^{-1}} \
 x_i^{q_1+1}y_i, \\
SubC_3(A_{m+4i-1}) &=& SubC_3(A_{m+4i-2}).
\eea

\item [(4)] Below we list the reduced forms of elements
$x_i^{\phi_{m+4(i-1)+j}}, y_i^{\phi_{m+4(i-1)+j}}$
 for $j = 0, \ldots, 4.$ Again, in the formulas below we assume that $y_0
= x_1^{-1}$ for $m = 0$, and $y_0 = z_m$ for $m \neq 0$.

\medskip
$x_i^{\phi _{ m+4i-4}} = A _{m+4i-4}^{-q_0}\circ x_i\circ A
_{m+4i-4}^{q_0}$,

\medskip
 $y_i^{\phi _{ m+4i-4}} = A _{m+4i-4}^{-q_0}\circ y_i$,

\medskip
$x_i^{\phi _{m+4i-3}}=x_i^{\phi_{m+4i-4}}$,

\medskip
$y_i^{\phi _{ m+4i-3}}= A _{m+4i-4}^{-q_0}\circ x_i^{q_1} \circ
y_i$,

\medskip
$x_i^{\phi _{
m+4i-2}}=\ig{A_{m+4i-2}^{q_2}}{x_iy_{i-1}^{-1}}{x_iy_i}
\ig{A_{m+4i-4}^{-q_0}}{x_iy_{i-1}^{-1}}{y_{i-2}x_{i-1}^{-1}}\ x_i\
\ig{A_{m+4i-4}^{q_0}}{x_{i-1}y_{i-2}^{-1}}{y_{i-1}x_i^{-1}}\ $,

\medskip
$y_i^{\phi_{m+4i-2}} = y_i^{\phi_{m+4i-3}}, $

\medskip
$x_i^{\phi_{m+4i-1}} = x_i^{\phi_{m+4i-2}} =_{(when \ i = n)}
x_i^{\phi_K}$,

\medskip
$Sub_3(x_i^{\phi_K}) =_{(when \ i = n)} \ SubC_3(A_{m+4i-2})\ \cup
\ SubC_3(A_{m+4i -4})^{\pm 1}\  \cup\ $ $$
\{y_{i-2}x_{i-1}^{-1}x_i, \ x_{i-1}^{-1}x_ix_{i-1}, \
x_ix_{i-1}y_{i-2}^{-1} \} ;$$

\medskip
$y_i^{\phi_{m+4i-1}}= {\tilde y}_i =_{(when \ i = n)} y_i^{\phi_K}
=$\newline $
\ig{A_{m+4i-4}^{-q_0}}{x_iy_{i-1}^{-1}}{y_{i-2}x_{i-1}^{-1}} \left
(x_i^{q_1} y_i \ \ig{A_{m+4i-2}^{q_2-1}}{x_iy_{i-1}^{-1}}{x_iy_i}
\ig{A_{m+4i-4}^{-q_0}}{x_iy_{i-1}^{-1}}{y_{i-2}x_{i-1}^{-1}} \ x_i
\right )^{q_3} \ x_i^{q_1} y_i $,

\medskip\begin{multline*}
Sub_3({\tilde y}_i) = SubC_3(A_{m+4i-2}) \cup
SubC_3(A_{m+4i-4})^{- 1} \cup \{y_{i-2}x_{i-1}^{-1}x_i, \
x_{i-1}^{-1}x_i^2,\\ \ x_i^3, \ x_iy_ix_i, \ y_ix_iy_{i-1}^{-1}, \
\ x_i^2y_i \}\end{multline*}

\medskip \begin{multline*}
x_i^{\phi_{m+4i}} =_{(when \ i \neq n)} x_i^{\phi_K}=
\ig{A_{m+4i}^{-q_4+1}}{x_{i+1}y_i^{-1}}{y_{i-1}x_i^{-1}}  \
x_{i+1} \circ y_i^{-1}x_i^{-q_1} \circ \\
\circ  \left ( x_i^{-1} \
\ig{A_{m+4i-4}^{q_0}}{x_{i-1}y_{i-2}^{-1}}{y_{i-1}x_i^{-1}}
\ig{A_{m+4i-2}^{-q_2+1}}{y_{i}^{-1}x_i^{-1}}{y_{i-1}x_i^{-1}} \
y_i^{-1}x_i^{-q_1} \right )^{q_3-1}\
\ig{A_{m+4i-4}^{q_0}}{x_{i-1}y_{i-2}^{-1}}{y_{i-1}x_i^{-1}} \
,\end{multline*}

\medskip\begin{multline*}
Sub_3(x_i^{\phi_K}) = SubC_3(A_{m+4i})^{-1} \cup
SubC_3(A_{m+4i-2})^{-1} \cup SubC_3(A_{m+4i-4}) \\
\cup \{y_{i-1}x_i^{-1}x_{i+1}, \ x_i^{-1}x_{i+1}y_i^{-1}, \
x_{i+1}y_i^{-1}x_i^{-1}, \ y_i^{-1}x_i^{-2}, \ x_i^{-3}, \
x_i^{-2}x_{i-1},\\ \ x_i^{-1}x_{i-1}y_{i-2}^{-1},\
y_{i-1}x_i^{-1}x_{i-1}, \ y_{i-1}x_i^{-1}y_i^{-1}, \
x_i^{-1}y_i^{-1}x_i^{-1} \};\end{multline*}

\medskip
\begin{multline*} y_i^{\phi_{m+4i}} = y_i^{\phi_K} =\ig{A_{m+4i}^{-q_4+1}}{x_{i+1}y_i^{-1}}{y_{i-1}x_i^{-1}}  \
x_{i+1} \ \ig{{\tilde y}_i}{x_iy_{i-1}^{-1}}{x_iy_i} \
x_{i+1}^{-1} \
\ig{A_{m+4i}^{q_4-1}}{x_iy_{i-1}^{-1}}{y_ix_{i+1}^{-1}}\
,\end{multline*}

\medskip
\begin{multline*}Sub_3(y_i^{\phi_K}) = SubC_3(A_{m+4i})^{\pm 1} \cup Sub_3({\tilde
y}_i) \cup \{y_{i-1}x_i^{-1}x_{i+1},\\ \ x_i^{-1}x_{i+1}x_i, \
x_{i+1}x_iy_{i-1}^{-1}, \ x_iy_ix_{i+1}^{-1}, \
y_ix_{i+1}^{-1}x_i, \
x_{i+1}^{-1}x_iy_{i-1}^{-1}\}.\end{multline*}
\end{enumerate}
\end{lm}
\begin{proof} Statement (1) is obvious. We prove  statement (2) by
induction on $i \geqslant 2$. Notice that by Lemmas
\ref{le:7.1.x1formsm0} and \ref{le:7.1.x1formsmneq0}
  ${\tilde y}_1 =  y_1^{\phi_{m+3}}$ begins with $x_1$  and ends  with $y_1$.
   Now let $i \geqslant  2$. Then
\bea
{\tilde y}_i  &=& y_i^{\phi_{m+4i-1}}\\
&=&
(x_i^{q_3}y_i)^{\phi_{m+4i-2}}\\
&=& \left ((y_i^{q_2}x_i)^{q_3}y_i
\right )^{\phi_{m+4i-3}}\\
& =&   \left ( \left ((x_i^{q_1}y_i)^{q_2}x_i\right
)^{q_3}x_i^{q_1}y_i\right )^{\phi_{m+4i-4}}.
\eea

 Before we  continue, and to avoid huge formulas,  we compute separately
$x_i^{\phi_{m+4i-4}}$ and $y_i^{\phi_{m+4i-4}}$:
\bea
x_i^{\phi_{m+4i-4}} &=& \left (x_i^{(y_{i-1}x_i^{-1})^{q_0}}\right
)^{\phi_{m+4(i-1)-1}}\\
&=& x_i^{({\tilde y}_{i-1}x_i^{-1})^{q_0}}\\
&=& (x_i {\tilde y}_{i-1}^{-1})^{q_0} \circ x_i \circ ({\tilde
y}_{i-1}x_i^{-1})^{q_0},
\eea
  by induction (by Lemmas \ref{le:7.1.x1formsm0} and
\ref{le:7.1.x1formsmneq0}  in the case $i = 2$) beginning with
$x_iy_{i-1}^{-1}$ and ending with $y_{i-1}x_i^{-1}$.

\bea
y_i^{\phi_{m+4i-4}} &=&  \left ((y_{i-1}x_i^{-1})^{-q_0} y_i \right
)^{\phi_{m+4(i-1)-1}}\\
&=&  ({\tilde y}_{i-1}x_i^{-1})^{-q_0}y_i\\
&=& (x_i \circ {\tilde y}_{i-1}^{-1})^{q_0} \circ y_i,
\eea
beginning with $x_iy_{i-1}^{-1}$ and ending with
$x_{i-1}^{-1}y_i$. It follows that
\bea
(x_i^{q_1}y_i )^{\phi_{m+4i-4}} &=&  (x_i {\tilde
y}_{i-1}^{-1})^{q_0} \circ x_i^{q_1}  \circ ({\tilde
y}_{i-1}x_i^{-1})^{q_0}(x_i \circ {\tilde y}_{i-1}^{-1})^{q_0}
\circ y_i\\
&=& (x_i {\tilde y}_{i-1}^{-1})^{q_0}  \circ x_i^{q_1}  \circ y_i,
\eea
beginning with $x_iy_{i-1}^{-1}$ and ending with $x_iy_i$.  Now
looking at the formula
$${\tilde y}_i  = \left ( \left ((x_i^{q_1}y_i)^{q_2}x_i\right
)^{q_3}x_i^{q_1}y_i\right )^{\phi_{m+4i-4}}$$
it is obvious  that ${\tilde y}_i $ begins with $x_iy_{i-1}^{-1}$
and ends with $x_iy_i$, as required.

Now we prove statements (3) and (4) simultaneously.

\medskip
$A_{m+4i-4} = {\rm cycred}((y_{i-1}x_i^{-1})^{\phi_{m+4(i-1)-1}})
={\tilde y}_{i-1}\circ x_i^{-1}$,
\medskip \noindent
beginning with $x_{i-1}$ and ending with  $x_i^{-1}$. As we have
observed in proving (2)

\medskip
$x_i^{\phi_{m+4i-4}} = (x_i {\tilde y}_{i-1}^{-1})^{q_0} \circ x_i
\circ ({\tilde y}_{i-1}x_i^{-1})^{q_0} = $
\medskip
$A_{m+4i-4}^{-q_0}\circ x_i \circ A_{m+4i-4}^{q_0}, $

\medskip \noindent
beginning with $x_i$  and ending with $x_i^{-1}$.
$$y_i^{\phi_{m+4i-4}} =   (x_i \circ {\tilde y}_{i-1}^{-1})^{q_0}
\circ y_i  = A_{m+4i-4}^{-q_0} \circ y_i,$$ beginning with $x_i$
and ending with $y_i$. Now

\medskip
$A_{m+4i-3} ={\rm cycred}( x_i^{\phi_{m+4i-4}}) = x_i$,
\medskip \noindent
beginning with $x_i$  and ending with $x_i$.
\bea
x_i^{\phi_{m+4i-3}} &=& x_i^{\phi_{m+4i-4}},\\
y_i^{\phi_{m+4i-3}} &=&  (x_i^{q_1}y_i)^{\phi_{m+4i-4}}\\
& =& A_{m+4i-4}^{-q_0}\circ x_i^{q_1} \circ
A_{m+4i-4}^{q_0} A_{m+4i-4}^{-q_0} \circ y_i\\
& = &  A_{m+4i-4}^{-q_0}\circ x_i^{q_1}  \circ y_i, \eea
 beginning with $x_i$ and ending with $y_i$.
Now \bea
A_{m+4i-2} &=& y_i^{\phi_{m+4i-3}},\\
x_i^{\phi_{m+4i-2}} &=&  (y_i^{q_2}x_i)^{\phi_{m+4i-3}}\\
&=&  A_{m+4i-2}^{q_2}\circ A_{m+4i-4}^{-q_0}\circ x_i \circ
A_{m+4i-4}^{q_0}, \eea beginning with $x_i$ and ending with
$x_i^{-1}$. It is also convenient to rewrite $x_i^{\phi_{m+4i-2}}$
(by rewriting the subword $A_{m+4i-2}$) to show its cyclically
reduced form: \bea x_i^{\phi_{m+4i-2}} &=& A_{m+4i-4}^{-q_0}\circ
\left ( x_i^{q_1} \circ y_i  \circ A_{m+4i-2}^{q_2-1}\circ
A_{m+4i-4}^{-q_0}\circ x_i \right )\\
&& \qquad\qquad \circ
A_{m+4i-4}^{q_0}.\\
y_i^{\phi_{m+4i-2}} &=& y_i^{\phi_{m+4i-3}}.
\eea

Now we can write down the next set of formulas: \bea A_{m+4i-1}
&=&{\rm cycred}(
y_i^{-\phi_{m+4i-3}}x_i^{\phi_{m+4i-2}}y_i^{\phi_{m+4i-3}})\\
& =& {\rm cycred}(A_{m+4i-2}^{-1} A_{m+4i-2}^{q_2}
A_{m+4i-4}^{-q_0}\\
&& \qquad\qquad  x_i A_{m+4i-4}^{q_0}
A_{m+4i-2}\\
&=& A_{m+4i-2}^{q_2-1}\circ  A_{m+4i-4}^{-q_0}  \circ x_i^{q_1+1}
\circ y_i, \eea beginning with $x_i$ and ending with $y_i$,

\medskip\noindent
$x_i^{\phi_{m+4i-1}} = x_i^{\phi_{m+4i-2}}$,
\medskip
$ y_i^{\phi_{m+4i-1}} = {\tilde y}_i = (x_i^{q_3}y_i)^{
\phi_{m+4i-2}}  = (x_i^{ \phi_{m+4i-2}})^{q_3}y_i^{ \phi_{m+4i-2}}
= $

\medskip \noindent
substituting the cyclic decomposition of $x_i^{ \phi_{m+4i-2}}$
from above one has

\medskip
$= A_{m+4i-4}^{-q_0} \circ \left (x_i^{q_1} \circ y_i \circ
A_{m+4i-2}^{q_2-1} \circ  A_{m+4i-4}^{-q_0} \circ x_i \right
)^{q_3} \circ  x_i^{q_1} \circ y_i .$

\medskip \noindent
beginning with $x_i$ and ending with $y_i$.

\medskip \noindent
Finally
$$A_{m+4i} = (y_ix_{i+1}^{-1})^{\phi_{m+4i-1}} = {\tilde y}_i
\circ x_{i+1}^{-1},$$
beginning with $x_i$ and ending with  $x_{i+1}^{-1}$.

\bea
x_i^{\phi_{m+4i}} &=& \left ((y_ix_{i+1}^{-1})^{-q_4}x_i\right
)^{\phi_{m+4i-1}}\\
 &=& ({\tilde
y}_ix_{i+1}^{-1})^{-q_4}x_i^{\phi_{m+4i-1}}\\
& =& A_{m+4i}^{-q_4+1}x_{i+1}{\tilde y}_i^{-1}
x_i^{\phi_{m+4i-1}}\\
&=& A_{m+4i}^{-q_4+1} \circ x_{i+1} \circ  \left (
(x_i^{\phi_{m+4i-2}})^{q_3-1} y_i^{\phi_{m+4i-2}} \right )^{-1}.
\eea Observe that computations similar to that for $
y_i^{\phi_{m+4i-1}}$ show that
\begin{multline*}\left ( (x_i^{\phi_{m+4i-2}})^{q_3-1} y_i^{\phi_{m+4i-2}} \right
)^{-1} = \\ \left (A_{m+4i-4}^{-q_0} \circ \left (x_i^{q_1} \circ
y_i \circ A_{m+4i-2}^{q_2-1} \circ A_{m+4i-4}^{-q_0} \circ x_i
\right )^{q_3-1} \circ x_i^{q_1} \circ y_i\right )^{-1}.
\end{multline*}

\noindent Therefore
\begin{multline*}x_i^{\phi_{m+4i}} = A_{m+4i}^{-q_4+1} \circ x_{i+1} \circ
\\
\circ  \left (A_{m+4i-4}^{-q_0} \circ \left (x_i^{q_1} \circ y_i
\circ A_{m+4i-2}^{q_2-1} \circ A_{m+4i-4}^{-q_0} \circ x_i \right
)^{q_3-1} \circ x_i^{q_1} \circ y_i\right )^{-1} ,\end{multline*}

\medskip \noindent
beginning with $x_{i+1}$ and ending with  $x_{i}^{-1}$.

\bea y_i^{\phi_{m+4i}} &=& \left ( y_i^{(y_ix_{i+1}^{-1})^{q_4}}\right
)^{\phi_{m+4i-1}}\\
& =& (x_{i+1}{\tilde y}_i^{-1})^{q_4} {\tilde y}_i
({\tilde y}_ix_{i+1}^{-1})^{q_4}\\
&=& A_{m+4i}^{-q_4+1} \circ x_{i+1} \circ {\tilde y}_i \circ
x_{i+1}^{-1} \circ  A_{m+4i}^{q_4-1}, \eea beginning with
$x_{i+1}$ and ending with  $x_{i+1}^{-1}$. This finishes the proof
of the lemma. \end{proof}

\begin{lm}
\label{le:7.1.2words}  Let $m>2, n=0$, $K = K(m,n)$, $p = (p_1,
\ldots,p_K)$ be a $3$-large tuple,  $\phi_K = \gamma_K^{p_K}
\cdots \gamma_1^{p_1}$,  and $X^{\pm \phi_K}  = \{x^{\phi_K} \mid
x \in X^{\pm 1} \}.$
 Then the following holds:

\medskip
\begin{enumerate}
 \item [(1)] $
 Sub_2(X^{\pm \phi_K})  = \left\{
 \begin{array}{ll}    c_jz_j,\ z_j^{-1}c_j & (1 \leqslant j \leqslant m), \\
z_jz_{j+1}^{-1}    & (1 \leqslant j \leqslant m-1), \\
z_mx_1^{-1},\ z_mx_1 &  (if \ m \neq 0, n \neq 0),  \\
x_i^2, \ x_iy_i,\ y_ix_i &  (1 \leqslant i \leqslant n), \\
x_{i+1}y_{i}^{-1},  \ x_i^{-1}x_{i+1},\ x_{i+1}x_i & (1 \leqslant
i \leqslant n-1) \ \end{array} \right\}^{\pm1}$

\medskip \noindent moreover, the word $z_j^{-1}c_j$, as well as $c_jz_j$,
 occurs only as a part of the subword $(z_j^{-1}c_jz_j)^{\pm 1}$
 in $x^{\phi_K}$ $( x \in X^{\pm 1})$;

\item [(2)] $ Sub_3(X^{\pm \phi_K})  =$ \begin{multline*} \left\{
\begin{array}{ll}
z_{j}^{-1}c_{j}z_{j}, & (1 \leqslant j \leqslant m), \\
 c_jz_jz_{j+1}^{-1}, \ z_jz_{j+1}^{-1}c_{j+1}^{-1}, \
z_jz_{j+1}^{-1}c_{j+1},  &  (1 \leqslant j \leqslant m-1),  \\
y_1x_1^2, \ & (m = 0, n = 1), \\
 x_2^{-1}x_1^2, \ x_2x_1^2,  & (m = 0, n \geqslant 2) \\
 c_m^{-1}z_mx_1, & (m= 1, n\neq 0)\\
   c_mz_mx_1^{-1}, \ c_mz_mx_1, \ z_mx_1^{-1}z_m^{-1}, \
z_mx_1^2, \  z_mx_1^{-1}y_1^{-1},  & (\ m
\neq 0, n \neq 0), \\
   z_mx_1^{-1}x_2, \ z_mx_1^{-1}x_2^{-1},  &  (m \neq 0, n \geqslant 2),  \\
  c_1^{-1}z_1z_2^{-1}, & (m \geqslant 2),  \\
     x_i^3, \ x_i^2y_i,\ \ x_iy_ix_i,  &  (1 \leqslant i \leqslant n),  \\
  x_i^{-1}x_{i+1}x_i,  \ y_ix_{i+1}^{-1}x_i, \
x_iy_ix_{i+1}^{-1},  &  (1 \leqslant i \leqslant n-1),  \\
 x_{i-1}^{-1}x_i^2,  \ y_ix_{i}y_{i-1}^{-1},  &  (2 \leqslant i \leqslant n), \\
 y_{i-2}x_{i-1}^{-1}x_i^{-1}, \ y_{i-2}x_{i-1}^{-1}x_i &  (3 \leqslant i \leqslant n).\   \end{array} \right\}^{\pm 1}.\end{multline*}

\item [(3)] for any $2$-letter word $uv \in  Sub_2(X^{\pm \phi_K})
$
  one has
  $$Sub_2(u^{\phi_K}v^{\phi_K}) \subseteq Sub_2(X^{\pm \phi_K}),
   \ \ \ Sub_3(u^{\phi_K}v^{\phi_K}) \subseteq Sub_3(X^{\pm
   \phi_K}).$$\end{enumerate}
\end{lm}
 \begin{proof} (1) and (2) follow by  straightforward inspection of the reduced
forms of elements $x^{\phi_K}$ in Lemmas \ref{le:7.1.zforms},
\ref{le:7.1.x1formsm0},  \ref{le:7.1.x1formsmneq0}, and
\ref{le:7.1.xiforms}.

To prove (3) it suffices for every word $uv \in Sub_2(X^{\pm
\phi_K}) $ to write down  the product $u^{\phi_K}v^{\phi_K}$
(using formulas from the lemmas mentioned above),  then make all
possible cancellations and check whether   3-subwords  of the
resulting word all lie in $Sub_3(X^{\pm \phi_K})$. Now we do the
checking one by one for all possible 2-words from $Sub_2(X^{\pm
\phi_K}) $.

\begin{enumerate}
\item [1)] For $uv  \in \{  c_jz_j, \   z_j^{-1}c_j\}$ the
checking is obvious and we omit it.
 \item [2)] Let $uv =z_jz_{j+1}^{-1}$. Then there are three cases to consider:
\begin{enumerate}

\item [2.a)]  Let $j \leqslant m-2$, then
$$(z_jz_{j+1}^{-1})^{\phi_K}=
\ig{z_j^{\phi_K}}{*}{c_{j+1}z_{j+1}}\ig{z_{j+1}^{-\phi_K}}{z_{j+2}^{-1}c_{j+2}^{-1}}{*}
,$$
in this case there is no cancellation in $u^{\phi_K}v^{\phi_K}$.
All 3-subwords of $u^{\phi_K}$ and $v^{\phi_K}$ are obviously in
$Sub_3(X^{\pm \phi_K})$. So one needs only to check the new
3-subwords which arise ``in between"  $u^{\phi_K}$ and $v^{\phi_K}$
(below we will check only subwords of this type). These subwords
are $c_{j+1}z_{j+1}z_{j+2}^{-1}$ and
$z_{j+1}z_{j+2}^{-1}c_{j+2}^{-1}$  which both lie in $Sub_3(X^{\pm
\phi_K})$.

\item [2.b)] Let $j = m-1$ and $n \neq 0$. Then
$$(z_{m-1}z_{m}^{-1})^{\phi_K}=
\ig{z_{m-1}^{\phi_K}}{*}{c_mz_m}\ig{z_{m}^{-\phi_K}}{x_1z_m^{-1}}
{*},$$
 again, there is no cancellation in this case and the words ``in
between" are $c_mz_mx_1$ and $z_mx_1z_m^{-1}$, which are in
$Sub_3(X^{\pm \phi_K})$.

\item [2.c)] Let $j = m-1$ and $n = 0$. Then ( below we put
$\cdot$ at the place where the corresponding   initial segment of
$u^{\phi_K}$ and the corresponding terminal segment of
$v^{\phi_K}$ meet)
\bea
(z_{m-1}z_{m}^{-1})^{\phi_K} &=& z_{m-1}^{\phi_K} \cdot
z_{m}^{-\phi_K}\\
&=& c_{m-1}z_{m-1}A_{m-4}^{p_{m-4}}c_m^{z_m}A_{m-1}^{p_{m-1}-1} \cdot
A_{m-1}^{-p_{m-1}}z_m^{-1}\\
&& \mbox{(cancelling $A_{m-1}^{p_{m-1}-1}$ and substituting for}\\
&& \mbox{
$A_{m-1}^{-1}$ its expression via the leading terms)}\\
&=& c_{m-1}z_{m-1}A_{m-4}^{p_{m-4}}c_m^{z_m} \cdot (c_m^{-z_m}
A_{m-4}^{-p_{m-4}} c_{m-1}^{-z_{m-1}}A_{m-4}^{p_{m-4}})z_m^{-1}\\
&=&
z_{m-1}\ \ig{A_{m-4}^{p_{m-4}}}{z_{m-2}^{-1}}{*} \ z_m^{-1}.
\eea
Here $z_{m-1}^{\phi _K}$ is completely cancelled.

\item [3.a)] Let $n=1.$ Then
\bea (z_mx_1^{-1})^{\phi _K} &=& c_m z_mA_{m-1}^{p_{m-1}}x_1^{-1}A_{m}^{p_{m}-1}\cdot
A_{m}^{-p_{m}}x_1^{-1}A_{m}^{p_{m}}A_{m+2}^{p_{m+2}}\\
&=& c_mz_mA_{m-1}^{p_{m-1}}x_1^{-1}\cdot
x_1A_{m-1}^{-p_{m-1}}c_m^{-z_m}A_{m-1}^{p_{m-1}}x_1^{-1}A_{m}^{p_{m}}A_{m+2}^{p_{m+2}}\\
&=& \ig
{z_mA_{m-1}^{p_{m-1}}x_1^{-1}A_{m}^{p_{m}}A_{m+2}^{p_{m+2}}}{z_mz_{m-1}^{-1}}{},
\eea
 and $z_m^{\phi _K}$ is completely cancelled.

\item [3.b)] Let $n>1$. Then \begin{multline*}(z_mx_1^{-1})^{\phi
_K}=c_mz_mA_{m-1}^{p_{m-1}}x_1^{-1}A_{m}^{p_{m}-1}\\
A_{m}^{-p_{m}}(x_1^{-1}A_{m}^{p_{m}}A_{m+2}^{-p_{m+2}+1}y_1^{-1}x_1^{-p_{m+1}})^{-p_{m+3}+1}x_1^{p_{m+1}}y_1x_2^{-1}
A_{m+4}^{p_{m+4}-1}\\ =c_mz_mA_{m-1}^{p_{m-1}}x_1^{-1}
A_{m}^{-1}(x_1^{-1}A_{m}^{p_{m}}A_{m+2}^{-p_{m+2}+1}y_1^{-1}x_1^{-p_{m+1}})^{-p_{m+3}+1}
x_1^{p_{m+1}}y_1x_2^{-1}A_{m+4}^{p_{m+4}-1}\\
=c_mz_mA_{m-1}^{p_{m-1}}x_1^{-1}\cdot
x_1A_{m-1}^{-p_{m-1}}c_m^{-z_m}A_{m-1}^{p_{m-1}}
(x_1^{-1}A_{m}^{p_{m}}A_{m+2}^{-p_{m+2}+1}y_1^{-1}x_1^{-p_{m+1}})^{-p_{m+3}+1}\\
x_1^{p_{m+1}}y_1x_2^{-1}A_{m+4}^{p_{m+4}-1}\\ =\ig{z_mA_{m-1}^{p_{m-1}}}{z_mz_{m-1}^{-1}c_{m-1}^{-1}}{}
\ ,\end{multline*} and $z_m^{\phi _K}$ is completely cancelled.

\item [4.a)]  Let $n=1$. Then \bea (z_mx_1)^{\phi _K} &=&
z_mA_{m-1}^{p_{m-1}}x_1^{-1}A_{m}^{p_{m}-1}\cdot
A_{m+2}^{p_{m+2}}A_{m}^{-p_{m}}x_1A_{m}^{p_{m}}\\
&=& \ig{z_mA_{m-1}^{p_{m-1}}**}{z_mz_{m-1}^{-1}c_{m-1}^{-1}}{*}\
,\eea and $z_m^{\phi _K}$ is completely cancelled.

\item [4.b)] Let $n>1.$ Then $$(z_mx_1)^{\phi _K}=\ig{z_m^{\phi
_K}}{*}{z_mx_1^{-1}}\ig{x_1^{\phi _K}}{x_2y_1^{-1}}{*}\ .$$

\item [5.a)] Let $n=1.$ Then \bea x_1^{2\phi _K} &=&
A_{m+2}^{p_{m+2}}A_{m}^{-p_{m}}x_1A_{m}^{p_{m}}\cdot
A_{m+2}^{p_{m+2}}A_{m}^{-p_{m}}x_1A_{m}^{p_{m}}\\
&=& A_{m+2}^{p_{m+2}}A_{m}^{-p_{m}}x_1A_{m}^{p_{m}}\cdot
(A_{m}^{-p_{m}}x_1^{p_{m+1}}y_1)A_{m+2}^{p_{m+2}-1}A_{m}^{-p_{m}}x_1A_{m}^{p_{m}}\\
&=& A_{m+2}^{p_{m+2}}\ig{A_{m}^{-p_{m}}x_1}{*}{z_mx_1}\cdot
x_1^{p_{m+1}}y_1**\ .\eea

\item [5.b)] Let $n>1$. Then $$x_1^{2\phi _K}=\ig{x_1^{\phi
_K}}{}{z_mx_1^{-1}}\ig{x_1^{\phi _K}}{x_2y_1^{-1}}{}\ .$$

\item [6.a)] Let $1<i<n.$
Then
\bea
x_i^{2\phi
_K}
&=& A_{m+4i}^{-q_4+1}x_{i+1}y_i^{-1}x_i^{-q_1}(x^{-1}A_{m+4i-4}^{q_0}A_{m+4i-2}^{-q_2+1}y_i^{-1}x_i^{-q_1})^{q_3-1}\\
&&\qquad \ig{A_{m+4i-4}^{q_0}}{} {y_{i-1}x_i^{-1}} \cdot
\ig{A_{m+4i}^{-q_4+1}}{x_{i+1}y_i^{-1}}{}**.\eea

\item [6.b)] \bea
x_{n}^{2\phi
_K} &=& A_{m+4n-2}^{q_2}A_{m+4n-4}^{-q_0}x_nA_{m+4n-4}^{q_0}\\
&&\qquad \cdot
A_{m+4n-2}^{q_2}A_{m+4n-4}^{-q_0}x_nA_{m+4n-4}^{q_0}\\
&=&
A_{m+4n-2}^{q_2}A_{m+4n-4}^{-q_0}x_nA_{m+4n-4}^{q_0}\\
&&\qquad \cdot
A_{m+4n-4}^{-q_0}x_n^{q_1}y_n
A_{m+4n-2}^{q_2-1}A_{m+4n-4}^{-q_0}x_nA_{m+4n-4}^{q_0}\\
&=&
A_{m+4n-2}^{q_2}\ig{A_{m+4n-4}^{-q_0}x_n}{}{x_{n-1}x_n}\cdot
x_n^{q_1}**.\eea

\item [7.a)] If $n=1$. Then $(x_1y_1)^{\phi
_K}=A_{m+2}^{p_{m+2}}A_{m}^{-p_{m}}x_1\cdot x_1^{p_{m+1}}**.$

\medskip
\item  [7.b)] If $n>1$. Then $(x_1y_1)^{\phi _K}=\ig{x_1^{\phi
_K}}{}{z_{m}x_1^{-1}}\ig{y_1^{\phi _K}}{x_2y_1^{-1}}{}\ .$

\medskip
\item [7.c)] If $1<i<n.$ Then $$(x_iy_i)^{\phi _K}=\ig{x_i^{\phi
_K}}{}{y_{i-1}x_i^{-1}}\ig{y_i^{\phi _K}}{x_{i+1}y_i^{-1}}{}\ .$$

\medskip
\item [7.d)] $(x_ny_n)^{\phi _K}=\ig{x_n^{\phi
_K}}{}{x_{n-1}^{-1}x_n}\ig{y_n^{\phi _K}}{x_{n}^2}{}\ .$

\item [8a)] If $n=1$. Then $$(y_1x_1)^{\phi _K}=\ig{y_1^{\phi
_K}}{}{x_1y_1}\ig{x_1^{\phi _K}}{x_1z_m^{-1}}{}\ .$$

\item [8.b)] If $n>1$. Then \bea
(y_1x_1)^{\phi
_K} &=& A_{m+4}^{-p_{m+4}+1}x_2A_{m+4}^{p_{m+4}}\cdot
A_{m+4}^{-p_{m+4}+1}x_2y_1^{-1}x_1^{-p_{m+1}}\circ **\\
&=& A_{m+4}^{-p_{m+4}+1}x_2A_{m+4}\cdot
x_2y_1^{-1}x_1^{-p_{m+1}}\circ **\\
&=& A_{m+4}^{-p_{m+4}+1}x_2A_{m}^{-p_{m}}
(x_1^{p_{m+1}}y_1A_{m+2}^{p_{m+2}-1}A_{m}^{-p_{m}}x_1)^{p_{m+3}}
x_1^{p_{m+1}}y_1x_2^{-1}\\
&&\qquad x_2y_1^{-1}x_1^{-p_{m+1}}( )^{p_{m+3}-1}A_{m}^{p_{m}}\\
&=& A_{m+4}^{-p_{m+4}+1}x_2A_{m}^{-p_{m}}
(x_1^{p_{m+1}}y_1A_{m+2}^{p_{m+2}-1}\ig{A_{m}^{-p_{m}}x_1}{}
{z_mx_1}\ig{A_{m}}{z_m^{-1}c_m^{-1}}{}\ .
\eea

\item[8.c)] $(y_nx_n)^{\phi _K}=\ig{y_n^{\phi
_K}}{}{x_ny_n}\ig{x_n^{\phi _K}}{x_ny_{n-1}}{}\ .$

\item [9.a)] If $n=2$, then
$$
 (x_2y_1^{-1})^{\phi _K}
 = A_{m+6}^{q_2}A_{m+4}^{-1}.$$

\item [9.b)] If $n>2$,
 $1<i<n$. Then
\bea
(x_iy_{i-1}^{-1})^{\phi_K}
&=& \ig{A_{m+4i}^{-q_4+1}}{x_{i+1}y_i^{-1}}{y_{i-1}x_i^{-1}}\
x_{i+1} \circ y_i^{-1}x_i^{-q_1}\\
&&\qquad \circ  \left ( x_i^{-1} \
\ig{A_{m+4i-4}^{q_0}}{x_{i-1}y_{i-2}^{-1}}{y_{i-1}x_i^{-1}}
\ig{A_{m+4i-2}^{-q_2+1}}{x_{i-1}y_{i-2}^{-1}}{y_{i-1}x_i^{-1}}\
y_i^{-1}x_i^{-q_1} \right )^{q_3-1}\\
&&\qquad\qquad
\ig{A_{m+4i-4}^{q_0}}{x_{i-1}y_{i-2}^{-1}}{y_{i-1}x_i^{-1}} \cdot
A_{m+4i-4}^{-q_0+1}\circ x_i\circ \tilde y_{i-1}\circ x_i^{-1}
\ig{A_{m+4i-4}^{q_0-1}}{x_{i-1}y_{i-2}^{-1}}{y_{i-1}x_i^{-1}}\\
&=& \ig{A_{m+4i}^{-q_4+1}}{x_{i+1}y_i^{-1}}{y_{i-1}x_i^{-1}} \
x_{i+1} \circ y_i^{-1}x_i^{-q_1}\\
&&\qquad  \circ   \left ( x_i^{-1} \
\ig{A_{m+4i-4}^{q_0}}{x_{i-1}y_{i-2}^{-1}}{y_{i-1}x_i^{-1}}
\ig{A_{m+4i-2}^{-q_2+1}}{x_{i-1}y_{i-2}^{-1}}{y_{i-1}x_i^{-1}} \
y_i^{-1}x_i^{-q_1} \right )^{q_3-1}\\
&&\qquad\qquad
x_i^{-1}\ig{A_{m+4i-4}^{q_0-1}}{x_{i-1}y_{i-2}^{-1}}{y_{i-1}x_i^{-1}}\
.\eea

\item [9.c)] $(x_ny_{n-1}^{-1})^{\phi
_K}=\ig{A_{m+4n-2}^{q_2}}{}{x_ny_n}\ig{A_{m+4n-4}}{x_ny_{n-1}^{-1}}{}\
\ .$

\item [10.a)] Let $n=2$, then \bea & &(x_1^{-1}x_2)^{\phi _K}
\\&=& A_{m}^{-p_{m}}(x_1^{p_{m+1}}y_1A_{m+2}^{p_{m+2}-1}
A_{m}^{-p_{m}}x_1)^{p_{m+3}-1}x_1^{p_{m+1}}y_1x_2^{-1}A_{m+4}^{p_{m+4}-1}
A_{m+6}^{p_{m+6}}A_{m+4}^{-p_{m+4}}x_2A_{m+4}^{p_{m+4}}\\
&=&A_{m}^{-p_{m}}(x_1^{p_{m+1}}y_1A_{m+2}^{p_{m+2}-1}
A_{m}^{-p_{m}}x_1)^{p_{m+3}-1}x_1^{p_{m+1}}y_1x_2^{-1}A_{m+4}^{p_{m+4}-1}\\
&&
 \qquad (A_{m+4}^{-p_{m+4}}x_2^{p_{m+5}}y_2)^{p_{m+6}}A_{m+4}^{-p_{m+4}}x_2A_{m+4}^{p_{m+4}}\\
&=& A_{m}^{-p_{m}}(x_1^{p_{m+1}}y_1A_{m+2}^{p_{m+2}-1}
A_{m}^{-p_{m}}x_1)^{p_{m+3}-1}x_1^{p_{m+1}}y_1x_2^{-1}\\
&&\qquad \cdot  A_{m+4}^{-1}x_2^{p_m+5}y_2(A_{m+4}^{-p_{m+4}}
x_2^{p_{m+5}}y_2)^{p_{m+6}-1}A_{m+4}^{-p_{m+4}}x_2A_{m+4}^{p_{m+4}}\\
&=& \ig{A_{m}^{-p_{m}}}{}{c_mz_m}\\
&& \qquad
\ig{x_1^{-1}A_{m}^{p_{m}}}{x_1^{-1}z_m^{-1}}{}A_{m+2}^{-p_{m+2}+1}y_1^{-1}x_1^{-p_{m+1}}A_{m}^{p_{m}}x_2^{p_{m+5}}
y_2(A_{m+4}^{-p_{m+4}}x_2^{p_{m+5}}y_2)^{p_{m+6}-1}\\
&&\qquad\qquad
A_{m+4}^{-p_{m+4}}x_2A_{m+4}^{p_{m+4}}.
\eea

\item [10.b)] If $1<i<n-1$, then $$(x_i^{-1}x_{i+1})^{\phi
_K}=\ig{x_i^{-\phi _K}}{}{y_{i}x_{i+1}^{-1}}\ig{x_{i+1}^{\phi
_K}}{x_{i+2}y_{i+1}^{-1}}{}.$$

\item [10.c)] Similarly to 10.a) we get $$(x_{n-1}^{-1}x_n)^{\phi
_K}=\ig{A_{2n+4n-8}^{-p_{m+4n-8}}}{}{y_{n-3}x_{n-2}^{-1}}\ \cdot\
\ig{x_{n-1}^{-1}A_{m+4n-8}^{p_{m+4n-8}}}{x_{n-1}^{-1}x_{n-2}}{}A_{m+4n-6}^{p_{m+4n-6}+1}**.$$

\item [11.a)] If $1<i<n-1$, then \bea &&(x_{i+1}x_i)^{\phi _K}\\
&= &
A_{m+4i+4}^{-q_8+1}x_{i+2}y_{i+1}^{-1}x_{i+1}^{-q_5}\left(x_{i+1}^{-1}A_{m+4i}^{q_4}A_{m_4i+2}^{-q_6+1}y_{i+1}^{-1}x_{i+1}^{-q_5}\right
)^{q_7-1}A_{m+4i}^{q_4}\\
&& \qquad A_{m+4i}^{-q_4+1}x_{i+1}y_{i}^{-1}x_{i}^{-q_1}\left(x_{i}^{-1}A_{m+4i-4}^{q_0}A_{m_4i-2}^{-q_2+1}y_{i}^{-1}x_{i}^{-q_1}\right
)^{q_3-1}A_{m+4i-4}^{q_0}\\
& =& A_{m+4i+4}^{-q_8+1}x_{i+2}y_{i+1}^{-1}x_{i+1}^{-q_5}\left(x_{i+1}^{-1}A_{m+4i}^{q_4}A_{m_4i+2}^{-q_6+1}y_{i+1}^{-1}x_{i+1}^{-q_5}\right
)^{q_7-1}A_{m+4i}\\
&&\qquad x_{i+1}y_{i}^{-1}x_{i}^{-q_1}\left(x_{i}^{-1}A_{m+4i-4}^{q_0}A_{m_4i-2}^{-q_2+1}y_{i}^{-1}x_{i}^{-q_1}\right
)^{q_3-1}A_{m+4i-4}^{q_0}\\
&=& A_{m+4i+4}^{-q_8+1}x_{i+2}y_{i+1}^{-1}x_{i+1}^{-q_5}\left(x_{i+1}^{-1}A_{m+4i}^{q_4}A_{m_4i+2}^{-q_6+1}y_{i+1}^{-1}x_{i+1}^{-q_5}\right
)^{q_7-1}\\
&&\qquad A_{m+4i-4}^{-q_0}x_i^{q_1}y_iA_{m+4i-2}^{q_2-1}\ig{A_{m+4i-4}^{-q_0}x_i}{}{x_{i-1}^{-1}x_i}\ig
{A_{m+4i-4}^{q_0}}{x_{i-1}y_{i-2}^{-1}}{}.
\eea

\item [11.b)] If $n>2$, then $$(x_2x_1)^{\phi _K}= \ **
\ig{A_{m}^{-q_0}x_1}{}{z_mx_1}\ig{A_{m}^{q_0}}{z_m^{-1}c_m^{-1}}{}.$$

\item [11.c)] \bea (x_nx_{n-1})^{\phi _K} &=&
A_{m+4n-2}^{q_6}A_{m+4n-4}^{-q_4}x_nA_{m+4n-4}^{q_4}\cdot
A_{m+4n-4}^{-q_4+1}x_ny_{n-1}^{-1}x_{n-1}^{-q_1}\\
&& \qquad
(x_{n-1}^{-1}A_{m+4n-8}^{q_0}A_{m+4n-6}^{-q_2+1}y_{n-1}^{-1}x_{n-1}^{-q_1}
 )^{q_3-1}A_{m+4n-8}^{q_0}\\
 &=&
\ **\ig{A_{m+4n-8}^{-q_0}x_{n-1}}{}{x_{n-2}^{-1}x_{n-1}}\ \cdot
 \ig{A_{m+4n-8}^{q_0}}{x_{n-2}y_{n-3}^{-1}}{}.\eea

 \item [11.d)] Similarly, if $n=2$, then $$(x_2x_1)^{\phi _K}=\
 **\ig{A_{m}^{-p_{m}}x_1}{}{z_mx_1}\ig
 {A_{m}^{p_{m}}}{z_m^{-1}c_m^{-1}}{}.$$

\end{enumerate}

\end{enumerate}

This proves the lemma. \end{proof}

\begin{lm} \label{le:7.1.2words.new}  Let $m >2,\  n=0, $ $K =K(m,0)$. Let  $p = (p_1, \ldots,p_K)$ be a 3-large tuple,  $\phi_K
= \gamma_K^{p_K} \ldots \gamma_1^{p_1}$,  and $X^{\pm \phi_K}=
\{x^{\phi_K} \mid x \in X^{\pm 1} \}.$ Denote the element
$$c_1^{z_1}\cdots
 c_m^{z_m}\in F(X\cup C_S)$$ by a new letter $d$.
Then the following holds:

\bi \item[(1)]  Every element from  $X^{\phi _K}$ can be uniquely
presented as a reduced product of elements and their inverses from
the set
$$X\cup \{c_1,\dots ,c_{m-1}, d\}$$
  Moreover:
 \begin{itemize}
 \item  all elements $z_i^{\phi _K},i\neq m$ have the form $z_i^{\phi _K}=c_iz_i\hat z_i$, where
 $\hat z_i$ is a words in the alphabet
 $\{c_1^{z_1}, \ldots, c_{m-1}^{z_{m-1}},
 d\}$,  \item $z_m^{\phi _K}=z_{m}\hat z_m$, where $\hat z_m$ is a
 word in the alphabet $\{c_1^{z_1},
  \ldots, c_{m-1}^{z_{m-1}},
 d\}$.
\end{itemize}

 When
viewing elements from  $X^{\phi _K}$ as elements in $$F(X\cup
\{c_1,\dots ,c_{m-1},d\}),$$ the following holds:

\item[(2)] \noindent
 $
 Sub_2(X^{\pm \phi_K})  = \left\{
 \begin{array}{ll}    c_jz_j\ & (1 \leqslant j \leqslant m), \\
  z_j^{-1}c_j, \ z_jz_{j+1}^{-1}    & (1 \leqslant j \leqslant m-1), \\
z_2d,\  dz_{m-1}^{-1} &
 \end{array} \right\}^{\pm1}$

\medskip \noindent Moreover:
\begin{itemize}
 \item the word $z_mz_{m-1}^{-1}$ occurs  only in the beginning
of $z_m^{\phi _K}$ as a part of the subword \newline$
z_mz_{m-1}^{-1}c_{m-1}^{-1}z_{m-1}$
 \item  the words $z_2d,\
dz_{m-1}^{-1}$ occur  only as  parts of  subwords
$$(c_1^{z_1}c_2^{z_2})^2dz_{m-1}^{-1}c_{m-1}^{-1}z_{m-1}c_{m-1} $$
and $(c_1^{z_1}c_2^{z_2})^2d$. \end{itemize}

\item[(3)] \noindent $ Sub_3(X^{\pm \phi_K})  = \left\{
\begin{array}{ll}

 z_{j}^{-1}c_{j}z_{j},\ c_jz_jz_{j+1}^{-1}, \ z_jz_{j+1}^{-1}c_{j+1}^{-1}, \
  &  (1 \leqslant j \leqslant m-1),  \\
z_jz_{j+1}^{-1}c_{j+1}  &  (1 \leqslant j \leqslant m-2),\\
c_2z_2d,\ z_2dz_{m-1}^{-1},\  dz_{m-1}^{-1}c_{m-1}^{-1},\ &
\end{array} \right\}^{\pm 1}.$
\ei
\end{lm}

\begin{proof} The lemma follows from Lemmas \ref {le:7.1.zforms}
and \ref{le:7.1.2words}
 by replacing all the products $c_1^{z_1}\ldots
c_m^{z_m}$
 in subwords of  $X^{\pm \phi_K}$
 by the letter $d$.
\end{proof}

\begin{notation} Let $m \neq 0,$ and if $m=1$, then $n\neq 1; \ K = K(m,n),
p=(p_1,\dots ,p_K)$ be a 3-large tuple, and  $\phi_K = \gamma
_K^{p_K} \ldots \gamma _1^{p_1}.$ Let ${\mathcal W}$ be the set of
words in $F(X \cup C_S)$ with the following properties:
\begin{enumerate} \item If $v\in W$ then $ Sub_3(v) \subseteq
Sub_3(X^{\pm \phi _K}), Sub_2(v) \subseteq Sub_2(X^{\pm \phi
_K});$ \item Every subword $x_i^{\pm 2}$ of $v\in W$ is contained
in a subword $x_i^{\pm 3};$ \item Every subword $c_1^{\pm z_1}$ of
$v\in W$ is contained in $(c_1^{z_1}c_2^{z_2})^{\pm 3}$ when
$m\geq 2$ or in $(c_1^{z_1}x_1^{-1})^{\pm  3}$ when $m=1$; \item
Every subword $c_m^{\pm z_m}\ (m\geq 3)$ is contained in
$\left(\prod _{i=1}^{m} c_i^{z_i}\right)^{\pm 1}.$ \item  every
subword $c_2^{\pm z_2}$ of $v\in W$ is contained either in
$(c_1^{z_1}c_2^{z_2})^{\pm 3}$ or as the central occurrence of
$c_2^{\pm z_2}$ in $(c_2^{-z_2}c_1^{-z_1})^{ 3}c_2^{\pm z_2}
(c_1^{z_1}c_2^{z_2})^{3}$ or in
$(c_1z_1c_2^{z_2}(c_1^{z_1}c_2^{z_2})^{3})^{\pm 1}$.
\end{enumerate}
\end{notation}

\begin{df}
 The following words are called {\em elementary periods}:
$$x_i, \  \ \ c_1^{z_1}c_2^{z_2} \ (if \ m \geq 2), \ \ c_1^{z_1}x_1^{-1} \  (if \
m=1).$$
  We call the squares (cubes) of elementary periods or
 their inverses elementary squares (cubes).
\end{df}

\begin{notation}
 Denote by $Y$ the following set of words
\begin{itemize}
 \item  [1)] if $n \neq 0$ then
 $Y = \{x_i,y_i,c_j^{z_j} \mid  i=1,\dots ,n,\ j=1,\dots ,m\}.$
 \item   [2)] if $n=0$ then
 $Y = \{c_1^{z_1}, \ldots,   c_{m-1}^{z_{m-1}},\ d\}.$
\end{itemize}
\end{notation}

\begin{notation}
 \begin{itemize}
 \item  [1)] Denote by  ${\mathcal W}_{\Gamma}$  the set of
all subwords of words in $\mathcal W$.
 \item   [2)] Denote by  $\bar{\mathcal W}_{\Gamma}$  the set of all  words $v \in {\mathcal
W}_{\Gamma}$ that are freely reduced forms of products of elements
from $Y^{\pm 1}$. In this case we say that these elements $v$ are
(group) words in the alphabet $Y$.
\end{itemize}
\end{notation}

\begin{lm}
\label{le:xyu} Let $v\in {\mathcal W}_{\Gamma}$. Then  the
following holds:
\begin{enumerate}
\item  [(1)] If $v$ begins and ends with an elementary square but
not an elementary cube, then $v$ belongs to the following set:

$$
 \left\{
\begin{array}{ll}
x_{i-2}^2y_{i-2}x_{i-1}^{-1}x_ix_{i-1}y_{i-2}^{-1}x_{i-2}^{-2}, \
x_i^2y_ix_iy_{i-1}^{-1}x_{i-1}^{-2}, &   m \geqslant 2, n\neq 0\\
x_{i-2}^2y_{i-2}x_{i-1}^{-1}x_i^2,\
x_{i-2}^2y_{i-2}x_{i-1}^{-1}x_iy_{i-1}^{-1}x_{i-1}^{-2}, \ & \\
 x_1^2y_1x_1c_m^{\pm z_m}CD,\ D_1C_1c_m^{z_m}x_1c_m^{-z_m}C_2D_2, & \\

 D_1^{-1}C_1c_m^{z_m}x_1^{-1}x_2x_1c_m^{-z_m}C_2D_2 & \\
D_1^{-1}C_1 c_m^{z_m}x_1^2, \ x_1^2y_1x_2^{-1}x_1c_m^{-z_m}C_2D_2,
x_2^{-2}x_1c_m^{- z_m}C_3D_3,&
\\
D_1^{-1}CD_2 & \\
     &   \\
(c_1^{z_1}c_2^{z_2})^2dc_{m-1}^{-z_{m-1}}\ldots
(c_2^{-z_2}c_1^{-z_1})^2,\ z_mc_{m-1}^{-z_{m-1}}\ldots
(c_2^{-z_2}c_1^{-z_1})^2,& m\geqslant 3, n=0\\
 D_1^{-1}CD_2&   \\
\prod _{i=m-1}^1c_i^{-z_i}(c_2^{-z_2}c_1^{-z_1})^2
\\
     &   \\
x_1^2y_1(x_1c_1^{-z_1})^2, \
(c_1^{z_1}x_1^{-1})^{2}x_2(x_1c_1^{-z_1})^2,\ & m = 1, n\geq 2 \\
(x_1c_1^{-z_1})^{2}x_1^2,\ x_1^2y_1x_2^{-1}(x_1c_1^{-z_1})^2,
x_2^{-2}(x_1c_1^{-z_1})^2,
& \\

 & \\
x_{i-2}^2y_{i-2}x_{i-1}^{-1}x_ix_{i-1}y_{i-2}^{-1}x_{i-2}^{-2}, \
x_i^2y_ix_iy_{i-1}^{-1}x_{i-1}^{-2},  & m = 0, \ n > 1 \\
x_{i-2}^2y_{i-2}x_{i-1}^{-1}x_iy_{i-1}^{-1}x_{i-1}^{-2}, \
x_1^2y_1x_2^{-1}x_1^2,\ x_1^{-2}x_2^{-1}x_1^2, & \\
 x_1^2y_1x_1, \ x_2^2y_2x_2 &
  \end{array} \right\}^{\pm 1},$$
where $C_k$ is an arbitrary product of the type $\prod
_jc_{i_j}^{\pm z_{i_j}}$ with $i_{j+1}=i_j\pm 1$,
$$D_k =(c_1^{z_1}c_2^{z_2})^{\pm 2}.$$

\item  [(2)] If $v$ does not contain two elementary squares and
begins (ends) with an elementary square, or contains  no
elementary squares, then $v$ is a subword of one of the words
above.
\end{enumerate}
\end{lm}

\begin{proof}
 Straightforward verification using the description of the set
$Sub_3(X^{\pm \phi_K})$ from Lemma \ref{le:7.1.2words}.
\end{proof}

\begin{df}
Let $Y$ be an alphabet and $E$  a set of words of length at least
2 in $Y$. We say that an occurrence of a word $w \in Y \cup E$ in
a word $v$ is {\em maximal} relative to $E$ if it is not contained
 in any other (distinct from $w$)  occurrence of a word from $E$ in $v$.
We say that a set of words $W$ in the alphabet $Y$ admits {\em
Unique Factorization Property (UF)} with respect to $E$ if  every
word $w \in W$ can be uniquely presented as a product
  $$w = u_1 \ldots u_k$$
  where $u_i$ are maximal occurrences of words from $Y \cup E$. In
  this event the decomposition above is called {\em irreducible}.
\end{df}

\begin{lm}
\label{le:UF} Let $Y$ be an alphabet and $E$  a set of words of
length at least 2 in $Y$. If a set of words $W$ in the alphabet
$Y$ satisfies the following condition:

\begin{itemize}
 \item if $w_1w_2w_3$ is a subword of a word from $W$ and
  $w_1w_2, w_2w_3 \in E$ then $w_1w_2w_3 \in E$
   then $W$ admits (UF) with respect to $E$.
\end{itemize}
\end{lm}

\begin{df}
\label{de:nielsen} Let $Y$ be an alphabet, $E$  a set of words of
length at least 2 in $Y$ and $W$ a set of words in $Y$ which
admits (UF) relative to $E$.  An automorphism $\phi \in Aut F(Y)$
satisfies the Nielsen property with respect to $W$ with exceptions
$E$ if for any word $z \in Y \cup E$ there exists a decomposition
 \begin{equation}
\label{eq:LMR} z^\phi  = L_z \circ M_z \circ R_z,
 \end{equation}
  for some words $L_z, M_z, R_z   \in F(Y)$
  such that for any
  $u_1,u_2  \in Y \cup E$ with $u_1u_2 \in Sub(W) \smallsetminus
  E$ the words $L_{u_1} \circ M_{u_1}$ and $M_{u_2} \circ R_{u_2}$
occur as written  in the reduced form of $u_1^\phi u_2^\phi$.
\end{df}

\begin{lm}
 \label{le: Nielsen-FU}
 Let $W$ be a set of words in the alphabet $Y$ which admits (UF) with respect to a set of words $E$.
 If an automorphism $\phi
\in Aut F(Y)$ satisfies the Nielsen property with respect to $W$
with exceptions $E$ then for every $w \in W$ if $w = u_1 \ldots
u_k$ is the irreducible decomposition of $w$ then the words
$M_{u_i}$ occur as written (uncancelled) in the reduced form of
$w^{\phi}$.
\end{lm}
  {\it Proof.} follows directly from definitions.

It is easy to show that if an automorphism $\phi$ satisfies the
Nielsen property with respect to $W$ and $E$ as above, then for
each word  $z \in Y \cup E$  there exists a unique decomposition
(\ref{eq:LMR}) with maximal length of $M_z$. In this event we call
$M_z$  the {\em middle} of $z^\phi$ (with respect to $\phi$).

Set
\bea
 T(m,1) &=& \left\{c_s^{z_s} (s=1,\dots ,m),
\prod _{i=1}^{m}c_i^{z_i}x_1\prod _{i=m}^1c_i^{-z_i}\right\}^{\pm 1}, m\neq 1,\\
T(m,2) &=& T(m,1)\\
&& \quad \cup \left\{ \
\prod_{i=1}^{m}c_i^{z_i}x_1^{-1}x_2x_1\prod _{i=m}^{1}c_i^{-z_i},\
y_1x_2^{-1}x_1\prod _{i=m}^{1}c_i^{-z_i}, \ \prod
_{i=1}^{m}c_i^{z_i}x_1^{-1}y_1^{-1}\right\}^{\pm 1},
\eea
if $n \geqslant 3$ then put
$$T(m,n) = T(m,1) \cup \left\{ \ \prod
_{i=1}^{m}c_i^{z_i}x_1^{-1}x_2^{-1},\ \prod
_{i=1}^{m}c_i^{z_i}x_1^{-1}y_1^{-1} \right\}^{\pm 1} \cup T_1(m,n),
$$
where
\bea
T_1(m,n) &=& \{ y_{n-2}x_{n-1}^{-1}x_nx_{n-1}y_{n-2}^{-1},\
y_{r-2}x_{r-1}^{-1}x_r^{-1},\ y_{r-1}x_{r}^{-1}y_{r}^{-1},\\
&& \qquad\qquad
 y_{n-1}x_n^{-1}x_{n-1}y_{n-2}^{-1} \ \  (n > r \geqslant 2) \}^{\pm 1}.
\eea

\medskip
Now, let
 $E(m,n) = \bigcup_{i \geqslant 2} Sub_i(T(m,n))\cap \bar {\mathcal W}_{\Gamma}.$

\begin{lm}
\label{main} Let $m \neq 0, n \neq 0, K = K(m,n), p=(p_1,\dots
,p_K)$ be a 3-large tuple.  Then the following holds:
\begin{enumerate}
 \item [(1)] Let  $w\in E(m,n)$, $v = v(w)$ be the leading
variable of $w$, and $j = j(v)$ (see notations at the beginning of
Section \ref{se:7.2.5}). Then the period $A _j^{p_j-1}$ occurs in
$w^{\phi_K}$ and  each occurrence of $A _j^2$ in $w^{\phi _j}$ is
contained in some occurrence of $A _j^{p_j-1}.$ Moreover, no
square $A _k^2$ occurs in $w$ for $k > j$.

 \item [(2)]  The automorphism $\phi_K$ satisfies
the Nielsen property with respect to $\bar{\mathcal W}_{\Gamma}$
with exceptions $E(m,n)$. Moreover, the following conditions hold:
  \begin{enumerate}
\item $M_{x_j} = A_{m+4r-8}^{-p_{m+4r-8}+1}x_{r-1}$, for $j \neq
n$.
 \item $M_{x_n} = x_n^{q_1}  \circ y_n  \circ A_{m+4n-2}^{q_2-1}\circ
A_{m+4n-4}^{-q_0}\circ x_n $
 \item $M_{y_j} = y_j^{\phi_K}$, for $j < n$.
 \item $M_{y_n} = \left (x_n^{q_1} y_n \
\ig{A_{m+4n-2}^{q_2-1}}{x_ny_{n-1}^{-1}}{x_ny_n}
\ig{A_{m+4n-4}^{-q_0}}{x_ny_{n-1}^{-1}}{y_{n-2}x_{n-1}^{-1}} \ x_n
\right )^{q_3} \ x_n^{q_1} y_n$.
  \item   $M_w=w^{\phi_K}$ for any
$w \in E(m,n)$ except for the following words:
 \begin{itemize}
 \item  $w_1=y_{r-2}x_{r-1}^{-1}x_r^{-1}, 3\leq r\leq n-1$,
  $w_2=y_{r-1}x_{r}^{-1}y_{r}^{-1}, 2\leq r\leq n-1$,
 \item $w_3=y_{n-2}x_{n-1}^{-1}x_n$,
 $w_4=y_{n-2}x_{n-1}^{-1}x_ny_{n-1}^{-1},$ $w_5=y_{n-2}x_{n-1}^{-1}x_nx_{n-1}^{-1}y_{n-2}^{-1},$
 $w_6=y_{n-2}x_{n-1}^{-1}x_nx_{n-1},$
 $w_7=y_{n-2}x_{n-1}^{-1}x_n^{-1}$, $w_8=y_{n-1}x_n^{-1},$ $w_9=x_{n-1}^{-1}x_n,$
$w_{10}=x_{n-1}^{-1}x_ny_{n-1}^{-1},$
 $w_{11}=x_{n-1}^{-1}x_nx_{n-1}y_{n-2}^{-1}.$
\end{itemize}
  \item  The only letter that may occur in a word from ${\mathcal W}_{\Gamma}$ to the left of a subword
 $w\in\{w_1,\ldots ,w_8\}$ ending with $y_i$ ($i=r-1, r-2,n-1,n-2,\ i\geq 1$) is $x_{i}$
 the maximal number $j$ such that $L_{w}$
 contains $A_j^{p_j-1}$ is $j=m+4i-2$, and $R_{w_1}=R_{w_2}=1$,

 \end{enumerate}

\end{enumerate}

\end{lm}

\begin{proof}   We first exhibit the formulas for $u^{\phi _K}$,
where $u\in \bigcup_{i \geqslant 2} Sub_i(T_1(m,n)).$

(1.a) Let $i<n$. Then
\bea
(x_iy_{i-1}^{-1})^{\phi
_{m+4i}} &=& (x_iy_{i-1}^{-1})^{\phi
_K}\\
&=& \ig{A_{m+4i}^{-q_4+1}}{x_{i+1}y_i^{-1}}{y_{i-1}x_i^{-1}} \
x_{i+1} \circ y_i^{-1}x_i^{-q_1}\\
&& \qquad \circ  \left ( x_i^{-1} \
\ig{A_{m+4i-4}^{q_0}}{x_{i-1}y_{i-2}^{-1}}{y_{i-1}x_i^{-1}}
\ig{A_{m+4i-2}^{-q_2+1}}{x_{i-1}y_{i-2}^{-1}}{y_{i-1}x_i^{-1}} \
y_i^{-1}x_i^{-q_1} \right )^{q_3-1}\\[1ex]
&&  \qquad\qquad \ig{A_{m+4i-4}^{q_0}}{x_{i-1}y_{i-2}^{-1}}{y_{i-1}x_i^{-1}}\\[1ex]
&&  \qquad\qquad\qquad \cdot A_{m+4i-4}^{-q_0+1}\circ x_i\circ \tilde y_{i-1}\circ
x_i^{-1} \ig{A_{m+4i-4}^{q_0-1}}{x_{i-1}y_{i-2}^{-1}}{y_{i-1}x_i^{-1}}\\
&=& \ig{A_{m+4i}^{-q_4+1}}{x_{i+1}y_i^{-1}}{y_{i-1}x_i^{-1}} \
x_{i+1} \circ y_i^{-1}x_i^{-q_1}\\
&&  \qquad \circ  \left ( x_i^{-1} \
\ig{A_{m+4i-4}^{q_0}}{x_{i-1}y_{i-2}^{-1}}{y_{i-1}x_i^{-1}}
\ig{A_{m+4i-2}^{-q_2+1}}{x_{i-1}y_{i-2}^{-1}}{y_{i-1}x_i^{-1}} \
y_i^{-1}x_i^{-q_1} \right )^{q_3-1}\\
&&\qquad\qquad \cdot
x_i^{-1}\ig{A_{m+4i-4}^{q_0-1}}{x_{i-1}y_{i-2}^{-1}}{y_{i-1}x_i^{-1}}\
.
\eea

\vspace{.2cm} (1.b) Let $i=n$.  Then \bea
(x_ny_{n-1}^{-1})^{\phi
_{m+4n-1}} &=& (x_ny_{n-1}^{-1})^{\phi
_K}\\
&=& \ig{A_{m+4n-2}^{q_2}}{x_ny_{n-1}^{-1}}{x_ny_n}
\ig{A_{m+4n-4}^{-1}}{x_ny_{n-1}^{-1}}{y_{n-2}x_{n-1}^{-1}}\ .
\eea
Here $y_{n-1}^{-\phi _K}$ is completely cancelled.

\vspace{.2cm} (2.a) Let $i<n-1.$ Then
\bea
(x_{i+1}x_iy_{i-1}^{-1})^{\phi _K} &=& (x_{i+1}x_iy_{i-1}^{-1})^{\phi
_{m+4i+4}}\\
&=& A_{m+4i+4}^{-q_8+1}\circ x_{i+2}\circ y_{i+1}^{-1}\circ
x_{i+1}^{-q_5}\\
&& \qquad \circ
\left (x_{i+1}^{-1}\circ A_{m+4i}^{q_4}\circ
A_{m+4i+2}^{-q_6+1}\circ y_{i+1}^{-1}x_{i+1}^{-q_5}\right
)^{q_7-1}A_{m+4i-4}^{-q_0}\\
&&\qquad\qquad \circ x_i^{q_1}y_i\circ
A_{m+4i-2}^{q_2-1}\circ A_{m+4i-4}^{-1}.
\eea
Here $(x_iy_{i-1}^{-1})^{\phi _{m+4i+4}}$ was completely
cancelled.

\vspace{.2cm} (2.b) Similarly, $(x_iy_{i-1}^{-1})^{\phi
_{m+4i+3}}$ is completely cancelled in
$(x_{i+1}x_iy_{i-1}^{-1})^{\phi _{m+4i+3}}$ and
$$(x_{i+1}x_iy_{i-1}^{-1})^{\phi _{m+4i+3}}= A_{m+4i+2}^{q_6}\circ
A_{m+4i}^{-q_4}\circ x_{i+1}\circ
A_{m+4i-4}^{-q_0}A_{m+4i-2}^{q_2-1}\circ A_{m+4i-4}^{-1}.$$

 (2.c) \bea
 (x_n^{-1}x_{n-1}y_{n-2}^{-1})^{\phi
_{m+4n-1}} &=& A_{m+4n-4}^{-q_4}\circ x_n^{-1}\circ
A_{m+4n-4}^{q_4}\circ  A_{m+4n-2}^{-q_6+1}\circ y_n^{-1}\circ
x_n^{-q_5}\\
&& \qquad \circ A_{m+4n-8}^{-q_0}\circ x_{n-1}^{q_1}\circ
y_{n-1}\circ A_{m+4n-6}^{q_2-1}\circ A_{m+4n-8}^{-1} , \eea
and $(x_{n-1}y_{n-2}^{-1})^{\phi _{m+4n-1}}$ is completely
cancelled.

\vspace{.2cm}
\medskip
(3.a) $$(y_ix_iy_{i-1}^{-1})^{\phi _{m+4i}}= A_{m+4i}^{-q_4+1}\circ
x_{i+1}\circ A_{m+4i-4}^{-q_0}\circ x_i^{q_1}\circ y_i\circ
A_{m+4i-2}^{q_2-1}\circ A_{m+4i-4}^{-1},$$
 and
$(x_iy_{i-1}^{-1})^{\phi _{m+4i}}$ is completely cancelled.

\vspace{.2cm}
\medskip
(3.b) $(y_nx_ny_{n-1}^{-1})^{\phi _{K}}=y_n^{\phi _{K}}\circ
(x_ny_{n-1}^{-1})^{\phi _{K}}.$

\vspace{.2cm}
\medskip
(3.c) \bea (y_{n-1}x_{n}^{-1}x_{n-1}y_{n-2}^{-1})^{\phi
_K}& = &
    A_{m+4n-4}\circ A_{m+4n-2}^{-q_6+1}\circ y_n^{-1}\circ
x_n^{-q_5}\\
    && \qquad \circ A_{m+4n-8}^{-q_0}\circ x_{n-1}^{q_1}\circ
y_{n-1}\circ A_{m+4n-6}^{q_2-1}\circ A_{m+4n-8}^{-1} ,
\eea \noindent
and $y_{n-1}^{\phi _K}$ and $(x_{n-1}y_{n-2}^{-1})^{\phi _K}$ are
completely cancelled.

(4.a) Let $n\geqslant 2$. \bea (x_1c_m^{-z_m})^{\phi _{m+4i}} &=&
(x_1c_m^{-z_m})^{\phi
_K}\\
&=& \ig{A_{m+4}^{-q_4+1}}{x_{1}y_1^{-1}}{c_m^{z_m}x_1^{-1}} \
x_{2} \circ y_1^{-1}x_1^{-q_1} \\
&& \qquad \circ  \left ( x_1^{-1}\circ A_{m}^{q_0}\circ
A_{m+2}^{-q_2+1}\circ y_1^{-1}\circ x_1^{-q_1}\right
)^{q_3-1}\\
&& \qquad\qquad \circ A_{m}^{q_0} \cdot A_{m}^{-q_0}\circ x_1^{-1}\circ
A_{m}^{q_0-1}\\
&=& A_{m+4}^{-q_4+1}\circ x_2\circ y_1^{-1}\circ x_1^{-q_1}\\
&& \qquad \circ
\left ( x_1^{-1}\circ A_{m}^{q_0}\circ A_{m+2}^{-q_2+1}\circ
y_1^{-1}\circ x_1^{-q_1}\right )^{q_3-1}\circ x_1^{-1}\circ
A_{m}^{q_0-1}.
\eea

Let $n=1$. \bea (x_1z_m^{-c_m})^{\phi _K} &=& A_{m}^{-p_{m}}\circ
x_1^{p_{m+1}}\circ
y_1\circ A_{m+2}^{p_{m+2}-1}\circ A_{m}^{-1},\\
(y_1x_1z_m^{-c_m})^{\phi _K} &=&y_1^{\phi _K}\circ
(x_1z_m^{-c_m})^{\phi _K}.
\eea

(4.b) $(x_1c_m^{-z_m})^{\phi _K}$ is completely cancelled in
$x_2^{\phi _K}$ and for $n>2$:
\bea
(x_2x_1c_m^{-z_m})^{\phi _K} &=&
    A_{m+8}^{-q_8+1}\circ x_3\circ y_2^{-1}\circ x_3^{-q_5}\\
    && \qquad \circ \left(x_3^{-1}\circ A_{m+4}^{q_4}\circ A_{m+6}^{-q_6+1}\circ
y_2^{-1}\circ x_3^{-q_5}\right)^{q_7-1}\\
    && \qquad\qquad\circ A_{m}^{-q_0}\circ
x_1^{q_1}\circ y_1\circ A_{m+2}^{q_2-1}\circ A_{m}^{-1} \eea
and for $n=2$:
$$(x_2x_1c_m^{-z_m})^{\phi _K}= A_{m+6}^{q_6}\circ
A_{m+4}^{-q_4}\circ x_i\circ A_{m}^{-q_0}\circ x_1^{q_1}\circ
y_1\circ A_{m+2}^{q_2-1}\circ A_{m}^{-1}.$$

 (4.c)  The cancellation between
$(x_2x_1c_m^{-z_m})^{\phi _K}$ and $c_{m-1}^{-z_{m-1}}$ is the
same as the cancellation between $A_{m}^{-1}$ and
$c_{m-1}^{-z_{m-1}^{\phi _K}}$, namely,
\bea
A_{m}^{-1}c_{m-1}^{-z_{m-1}^{\phi _K}} &=& \left ( x_1\circ
A_{m-1}^{-p_{m-1}}\circ c_m^{-z_m}\circ A_{m-1}^{p_{m-1}}\right )\\
&& \quad\left(A_{m-1}^{-p_{m-1}+1}\circ c_m^{-z_m}\circ
A_{m-4}^{-p_{m-4}}\circ c_{m-1}^{-z_{m-1}}\circ
A_{m-4}^{p_{m-4}}\circ  c_m^{z_m}\circ A_{m-1}^{p_{m-1}-1}\right
)\\
&=& x_1A_{m-1}^{-1},
\eea

 and $c_{m-1}^{-z_{m-1}^{\phi _K}}$ is
completely cancelled.

(4.d)  The cancellations between $(x_2x_1c_m^{-z_m})^{\phi _K}$
(or between $(y_1x_1c_m^{-z_m})^{\phi _K}$) and $\prod
_{i=m-1}^1c_{i}^{-z_{i}^{\phi _K}}$ are the same as the
cancellations between $A_{m}^{-1}$ and $\prod
_{i=m-1}^1c_{i}^{-z_{i}^{\phi _K}}$ namely, the product $\prod
_{i=m-1}^1c_{i}^{-z_{i}^{\phi _K}}$ is completely cancelled and
\[
A_{m}^{-1}\prod _{i=m-1}^1c_{i}^{-z_{i}^{\phi _K}}= x_1\prod
_{i=m}^1c_i^{-z_i}.
\]

 Similarly one can write expressions for $u^{\phi _K}$ for all
$u\in E(m,n).$ The first statement of the lemma  now follows from
these formulas.

 Let us verify the second statement.
Suppose $w\in E(m,n)$  is a maximal subword from $E(m,n)$ of  a
 word $u$ from ${\mathcal W}_{\Gamma}$.
 If $w$ is a subword
of a word in $T(m,n)$, then either $u$ begins with $w$ or $w$ is
the leftmost subword of  a word in $T(m,n).$ All the words in
$T_1(m,n)$ begin with some $y_j$, therefore  the only possible
letters in $u$ in front of $w$ are $x_{j}^2$.

We have $x_{j}^{\phi _K}x_{j}^{\phi _K}w^{\phi _K }=x_{j}^{\phi _K
}\circ x_{j}^{\phi _K}\circ w^{\phi _K}$ if $w$ is a two-letter
word, and $x_{j}^{\phi _K}x_{j}^{\phi _K }w^{\phi _K }=x_{j}^{\phi
_K }\circ x_{j}^{\phi _K} w^{\phi _K}$ if $w$ is more than a
two-letter word. In this last case there are some cancellations
between $x_{j}^{\phi _K}$ and  $w^{\phi _K}$, and the middle of
$x_j$ is the non-cancelled part of $x_j$ because  $x_j$ as a
letter not belonging to $E(m,n)$  appears only in $x_j^n$.

We still have to consider all letters that can appear to the right
of $w$, if $w$ is the end of some word in $T_1(m,n)$ or
$w=y_{n-1}x_n^{-1}x_{n-1}$, $w=y_{n-1}x_n^{-1}$. There are the
following possibilities: \bi \item[(i)] $w$ is an end of $
y_{n-2}x_{n-1}^{-1}x_nx_{n-1}y_{n-2}^{-1};$

\item[(ii)] $w$ is an end of $y_{r-2}x_{r-1}^{-1}x_r^{-1}, r<i$;

\item[(iii)] $w$ is an end of   $y_{n-2}x_{n-1}^{-1}y_{n-1}^{-1}$.
\ei Situation (i) is equivalent to  the situation when $w^{-1}$ is
the beginning of the word $
y_{n-2}x_{n-1}^{-1}x_nx_{n-1}y_{n-2}^{-1}$, we have considered
this case already. In the situation (ii)  the only possible word
to the right of $w$ will be left end of
$x_{r-1}y_{r-2}^{-1}x_{r-2}^{-2}$ and $w^{\phi _K }x_{r-1}^{\phi
_K}y_{r-2}^{-\phi _K}x_{r-2}^{-2\phi _K }=w^{\phi _K }\circ
x_{r-1}^{\phi _K}y_{r-2}^{-\phi _K}\circ x_{r-2}^{-2\phi _K},$ and
$w^{\phi _K }x_{r-1}^{\phi _K}=w^{\phi _K }\circ x_{r-1}^{\phi
_K}.$ In the situation (iii) the first two letters to the right of
$w$ are $x_{n-1}x_{n-1}$, and $w^{\phi _K}x_{n-1}^{\phi
_K}=w^{\phi _K}\circ x_{n-1}^{\phi _K}.$

There is no cancellation in the words $ (c_j^{z_j})^{\phi_K}\circ
(c_{j+1}^{\pm z_{j+1}})^{\phi_K},
 (c_m^{z_m})^{\phi_K}\circ  x_1^{\pm\phi_K}, \  x_1^{\phi_K}\circ x_1^{\phi_K}. $
 For  all the
other occurrences of $x_i$ in the words from ${\mathcal
W}_{\Gamma}$, namely for occurrences in $x_i^n,\ x_i^2y_i$, we
have $(x_i^2y_i)^{\phi _K}=x_i^{\phi _K}\circ x_i^{\phi _K}\circ
y_i^{\phi _k}$ for $i<n$.

In the case $n=i$, the bold subword of the word

\noindent $x_n^{\phi_{_K}} = A_{m+4n-4}^{-q_0}\circ {\bf \left (
x_n^{q_1}  \circ y_n  \circ A_{m+4n-2}^{q_2-1}\circ
A_{m+4n-4}^{-q_0}\circ x_n \right ) }\circ A_{m+4n-4}^{q_0} $

\noindent is $M_{x_n}$ for $\phi _K$, and the bold subword in the
word

\noindent $y_n^{\phi_K}= \ig{{
A_{m+4n-4}^{-q_0}}}{x_ny_{n-1}^{-1}}{y_{n-2}x_{n-1}^{-1}}
{\bf\left (x_n^{q_1} y_n \ \ig{{\bf A_{m+4n-2}^{q_2-1}}}{{\bf
x_ny_{n-1}^{-1}}}{{\bf x_ny_n}} \ig{{\bf A_{m+4n-4}^{-q_0}}}{{\bf
x_ny_{n-1}^{-1}}}{{\bf y_{n-2}x_{n-1}^{-1}}} \ x_n \right )^{q_3}
\ x_n^{q_1} y_n }$,

\noindent is $M_{y_n}$ for $\phi _K.$

\end{proof}

\begin{cy} \label{cy:middles}
  Let $m \neq 0, n \neq 0, K = K(m,n), p=(p_1,\dots
,p_K)$ be a 3-large tuple, $L = Kl$.  Then  for any $u \in X \cup
E(m,n)$ the element $M_u$ with respect to ${\phi_L}$ contains
$A_j^q$ for some $j> L-K$ and $q > p_j - 3$.
\end{cy}
\begin{proof} This follows from the formulas for $M_u$ with respect
to ${\phi_K}$ in the lemma above.\end{proof}

\begin{notation} 1) Denote by ${\mathcal W}_{\Gamma ,L}$ the least set of words
in the alphabet $Y$  that contains ${\bar {\mathcal W}}_{\Gamma
},$ is closed under taking subwords, and is $\phi _K$-invariant.

 2) Let $\bar {\mathcal W}_{\Gamma ,L}$ be union of ${\mathcal
W}_{\Gamma ,L}$ and the set of all initial subwords of $z_i^{\phi
_{Kj}}$ which are of the form $c_i^j\circ z_i \circ w,$ where
$w\in{\mathcal W}_{\Gamma ,L}.$
\end{notation}

\begin{rk}
The set $\bar{\mathcal W}_{\Gamma ,L}$ is  $\phi _K$-invariant.
\end{rk}
\begin{proof}
 Indeed,  if $c_i^jz_iw\in\bar{\mathcal W}_{\Gamma
,L}$, then  $c_i^{c_i^jz_iw}=w^{-1}\circ c_i^{z_i}\circ w\in
{\mathcal W}_{\Gamma ,L}$ and $c_i^{(c_i^jz_iw)^{\phi _K}}=w^{-\phi
_K}\circ c_i^{z_i^{\phi _K}}\circ w^{\phi _K}\in {\mathcal
W}_{\Gamma ,L}$, therefore $c_i^{j+1}{z_i^{\phi _K}}\circ w^{\phi
_K}\in \bar{\mathcal W}_{\Gamma ,L}$.

\end{proof}

\begin{notation}  Denote by $Exc$ the following set of words in the alphabet
$Y$.
 $$Exc = \{ c_1^{-z_1}c_i^{-z_i}c_{i-1}^{-z_{i-1}}, \
 c_1^{-z_1}x_1c_m^{-z_m}, \ c_1^{-z_1}x_jy_{j-1}^{-1}\}.$$
 \end{notation}

\begin{lm}
\label{main11} The following holds:

\begin{enumerate}
 \item [(1)] $Sub_{3,Y}({\mathcal W}_{\Gamma ,L})=Sub_{3,Y}(X^{\pm\phi _K})\cup Exc$.
\item [(2)] Let $v \in {\mathcal W}_{\Gamma ,L}$ be a word that
begins and ends with an elementary square and does not contain any
elementary cubes. Then either $v \in {\bar {\mathcal W}}_{\Gamma
}$ or $v = v_1v_2$   where $v_1,v_2\in {\bar {\mathcal W}}_{\Gamma
}$ and these words are exhibited below:
\begin{enumerate}
\item for $m>2,\ n\geq 2$,
 $$v_1 \in \{v_{11}=(c_1^{z_1}c_2^{z_2})^2 \prod_{i = 3}^{m} c_i^{z_i}x_1x_2x_1 \prod_{i = m}^{1}c_i^{-z_i},
 \ v_{12}= x_1^2y_1x_1\prod_{i = m}^{1} c_i^{-z_i} \}, $$

  $$v_2 \in \{  v_{2i}=c_i^{-z_i}\ldots c_3^{-z_3}(c_2^{-z_2}c_1^{-z_1})^2,  u_{2,1}=x_1c_m^{-z_m}\ldots c_3^{-z_3}
  (c_2^{-z_1}c_1^{-z_1})^2,$$
  $$u_{2,j}=x_jy_{j-1}^{-1}x_{j-1}^2 \};$$

  \item for $m=2,\ n\geq 2$,
 $$v_1 \in \{v_{11}=(c_1^{z_1}c_2^{z_2})^2x_1x_2x_1 \prod_{i = m}^{1}c_i^{-z_i},
 \ v_{12}= x_1^2y_1x_1\prod_{i = m}^{1} c_i^{-z_i} \},$$

  $$v_2 \in \{
   u_{2,1}=x_1
  (c_2^{-z_1}c_1^{-z_1})^2, \
  u_{2,j}=x_jy_{j-1}^{-1}x_{j-1}^2 \};$$

  \item for $m>2,\ n=1$,
 $v_1=
 x_1^2y_1x_1\prod_{i = m}^{1} c_i^{-z_i},$

  $$v_2 \in \{  v_{2i}=c_i^{-z_i}\ldots c_3^{-z_3}(c_2^{-z_2}c_1^{-z_1})^2, \ u_{2,1}=x_1c_m^{-z_m}\ldots c_3^{-z_3}
  (c_2^{-z_1}c_1^{-z_1})^2 \};$$

 \item for $m=2,\ n=1$,\
 $v_1=x_1^2y_1x_1\prod_{i = m}^{1} c_i^{-z_i}, $
 $v_2=x_1
  (c_2^{-z_1}c_1^{-z_1})^2 ;$
  \item for $m=1,\ n\geq 2$,
 $$v_1 \in \{v_{11}=(c_1^{z_1}x_1^{-1})^2x_2x_1 c_1^{-z_1},
 \ v_{12}= x_1^2y_1x_1c_1^{-z_1} \},
\  v_2=x_jy_{j-1}^{-1}x_{j-1}^2 .$$
\end{enumerate}
\end{enumerate}
\end{lm}

\begin{proof} Let $T=Kl.$ We will consider only the case $m\geq 2,\ n\geq
2$. We will prove the statement of the lemma by induction on $l$.
If $l=1$, then $T=K$ and the statement is true. Suppose now that
$$Sub _{3,Y}(\bar {\mathcal W}_{\Gamma}^{\phi _{T-K}})=Sub
_{3,Y}(\bar {\mathcal W}_{\Gamma})\cup Exc.$$ Formulas in the
beginning of the proof of Lemma \ref{main} show that  $$Sub
_{3,Y}(E(m,n)^{\pm\phi _K})\subseteq Sub _{3,Y}(\bar {\mathcal
W}_{\Gamma}).$$ By the second statement the automorphism $\phi _K$
satisfies the Nielsen property with exceptions $E(m,n)$. Let us
verify that new 3-letter subwords do not occur "between" $u^{\phi
_K}$ for $u\in T_1(m,n)$ and the power of the corresponding $x_i$
to the left and right of it. All the cases are similar to the
following:

$$(x_nx_{n-1}y_{n-2}^{-1})^{\phi _K}\cdot x_{n-2}^{\phi _K}\ldots
\ig {A_{m+4n-10}^{-q+1}}{*}{y_{n-3}x_{n-2}^{-1}}\cdot
x_{n-1}^{-1}\ig{A_{m+4n-8}^{q_0-1}}{x_{n-2}}{*}\ .$$

Words $$ (v_1v_2)^{\phi _K}$$ produce the subwords from $Exc$.
Indeed, $[(x_2x_1\prod _{i=m}^1c_i^{-z_i})]^{\phi _{Kj}}$ ends
with $v_{12}$ and $v_{12}^{\phi _K}$ ends with $v_{12}.$
Similarly, $v_{2,j}^{\phi _K}$ begins with $v_{2,j+1}$ for $j<m$
and with $u_{2,1}$ for $j=m$. And $u_{2,j}^{\phi _K}$ begins with
$u_{2,j+1}$ for $j<n$ and with $u_{2,j}$ for $j=n$.

This and the second part of Lemma \ref{le:7.1.2words} finish the
proof.\end{proof}

 Let  $W \in G[X]$.
 We say that a word  $U \in G[X]$
 {\it occurs } in $W$ if $W = W_1 \circ U \circ W_2$ for some $W_1, W_2 \in G[X]$.
An occurrence of $U^q$ in $W$ is called {\it maximal} with respect
to a property $P$ of words if $U^q$ is not a part of any
occurrence of $U^r$ with $q < r$ and which satisfies $P$. We say
that an occurrence of $U^q$ in $W$ is {\it stable} if $q \geqslant
1$ and $W = W_1 \circ U U^q U \circ W_2$ (it follows that $U$ is
cyclically reduced). Maximal stable occurrences $U^q$ will play an
important part in what follows. If $(U^{-1})^q$
 is a stable occurrence of $U^{-1}$ in $W$ then, sometimes, we say
 that  $U^{-q}$ is a stable occurrence of $U$ in $W$. Two given occurrences
$U^q$ and $U^p$ in a word $W$ are {\it disjoint} if they do not
have a common letter as subwords of  $W$. Observe that if integers
$p$ and $q$ have different signs then any two occurrences of $A^q$
and $A^p$ are disjoint. Also, any two different maximal stable
occurrences of powers of $U$ are disjoint. To explain the main
property of stable occurrences of powers of $U$, we need the
following definition.
 We say that a given occurrence of $U^q$ {\em occurs correctly} in a given occurrence
  of $U^p$ if  $|q |\leqslant |p|$ and for these
occurrences $U^q$ and $U^p$ one has $U^p = U^{p_1} \circ U^q \circ
U^{p_1}$. We say, that two
 given non-disjoint occurrences of $U^q, U^p$ {\it overlap correctly} in $W$ if
 their common subword occurs correctly in each of them.

A cyclically reduced word $A$ from $G[X]$ which is not a proper
power and does not belong to $G$ is called {\it a period}.

\begin{lm}
\label{le:stable-overlap} Let $A$ be a period in $G[X]$ and $W \in
G[X]$. Then any two stable occurrences of powers of $A$ in $W$ are
either disjoint or they overlap correctly. \end{lm}
  \begin{proof} Let $A^q$, $A^p$ ($q
\leqslant p$) be two non-disjoint stable occurrences of powers of
$A$ in $W$. If they overlap  incorrectly then $A^2 = u \circ A
\circ v$ for some elements $u,v \in G[X]$. This implies that $A =
u \circ v = v \circ u$ and hence $u$ and $v$ are (non-trivial)
powers of some element in $G[X]$. Since $A$ is  not a proper power
it follows that $u = 1$ or $v = 1$ - contradiction. This shows
that $A^q$ and $A^p$ overlap correctly.  \end{proof}

Let $W \in G[X]$ and ${\mathcal O} = {\mathcal O}(W,A) =\{A^{q_1},
\ldots, A^{q_k}\}$ be a set of pair-wise disjoint stable
occurrences of powers of a period $A$ in $W$ (listed according to
their appearance  in $W$ from the left to the right). Then
${\mathcal O}$ induces an ${\mathcal O}$-decomposition of $W$ of
the following form:
\begin{equation}
\label{eq:new-O-decomp} W = B_1 \circ A^{q_1} \circ \cdots \circ
B_k \circ A^{q_k} \circ B_{k+1}\end{equation}

For example, let $P$ be a property of words (or just a property of
occurrences in  $W$) such that if two powers of $A$ (two
occurrences of powers of $A$ in $W$) satisfy $P$ and overlap
correctly then their union also satisfies $P$. We refer to such
$P$ as  {\it preserving correct overlappings}. In this event, by
${\mathcal O}_P = {\mathcal O}_P(W,A)$ we denote the uniquely
defined set of all maximal stable occurrences of powers of $A$ in
$W$ which satisfy the property $P$. Notice, that occurrences in
${\mathcal O}_P$ are pair-wise disjoint by Lemma
\ref{le:stable-overlap}. Thus, if $P$ holds on every power of $A$
then ${\mathcal O}_P(W,A) = {\mathcal O}(W,A)$ contains all
maximal stable occurrences of powers of $A$ in $W$. In this case,
the decomposition (\ref{eq:new-O-decomp}) is unique and it is
called the {\it canonical (stable)} $A$-decomposition of $W$.

The following example provides another property $P$ that will be
in use later. Let $N$ be a positive integer and let $P_N$ be the
property of $A^q$ that $|q| \geqslant N$. Obviously, $P_N$
preserves correct overlappings. In this case the set ${\mathcal
O}_{P_N}$ provides the so-called {\it canonical $N$-large }
$A$-decompositions of $W$ which are also uniquely defined.

\begin{df}
 Let
$$W = B_1 \circ A^{q_1} \circ \cdots \circ B_k
\circ A^{q_k} \circ B_{k+1}$$ be the  decomposition
(\ref{eq:new-O-decomp}) of $W$ above. Then   the numbers
$$ \max_A(W) = \max \{q_i \mid i = 1, \dots, k\}, \ \ \
 \min_A(W) = \min \{q_i \mid i = 1, \dots, k\}$$ are called,
correspondingly,  the {\em upper} and the {\em lower} $A$-bounds
of $W$.
\end{df}

\begin{df}
Let $A$ be a period in $G[X]$ and $W \in G[X]$.   For a positive
integer $N$ we say that the $N$-large $A$-decomposition of $W$
 $$W = B_1 \circ
A^{q_1} \circ \cdots \circ B_k \circ A^{q_k} \circ B_{k+1}$$
 has $A$-size $(l,r)$ if $\min_A(W) \geqslant l$ and
 $\max_A(B_i) \leqslant r$ for every $i = 1, \dots, k$.
 \end{df}

 Let ${\mathcal A} = \{A_1, A_2,
\ldots, \}$ be a sequence of periods from $G[X]$.   We say that a
word $W \in G[X]$ has ${\mathcal A}$-rank $j$ (${\rm
rank}_{\mathcal A}(W) = j$) if $W$ has a stable occurrence of
$(A_j^{\pm 1})^q$ ($q \geqslant 1$) and $j$ is maximal with this
property. In this event, $A_j$ is called the {\em ${\mathcal
A}$-leading term} (or just the {\em leading term}) of $W$
(notation $LT_{{\mathcal A}}(W) = A_j$ or $LT(W) = A_j$).

We now fix an arbitrary  sequence ${\mathcal A}$ of periods  in
the group $G[X]$. For a period $A = A_j$ one can consider
canonical $A_j$-decompositions of a word $W$ and define the
corresponding $A_j$-bounds and $A_j$-size. In this case
 we, sometimes, omit $A$ in the writings and simply write $max_j(W)$ or
$min_j(W)$ instead of $max_{A_j}(W)$,  $min_{A_j}(W)$.

In the case when   ${\rm rank}_{{\mathcal A}}(W) = j$  the
canonical $A_j$-decomposition of $W$ is called the {\em canonical}
${\mathcal A}$-decomposition of $W$.

Now we turn to an analog of ${\mathcal O}$-decompositions of $W$
with respect to ``periods" which are not necessarily cyclically
reduced words. Let $U = D^{-1} \circ A \circ D$, where  $A$ is a
period. For  a set ${\mathcal O} = {\mathcal O}(W,A) = \{A^{q_1},
\ldots, A^{q_k}\}$ as above  consider the ${\mathcal
O}$-decomposition of a word $W$
 \begin{equation} \label{eq:A-dec}
 W = B_1 \circ A^{q_1} \circ \cdots \circ B_k \circ A^{q_k} \circ B_{k+1}
 \end{equation}
  Now it can be rewritten in the form:
 $$ W = (B_1D) (D^{-1} \circ A^{q_1} \circ D) \cdots (D^{-1}B_kD) (D^{-1}\circ A^{q_k} \circ D)
 (D^{-1} B_{k+1}).$$
Let $\varepsilon _i, \delta _i={\rm sgn}(q_i).$ Since every
occurrence of $A^{q_i}$ above is stable, $B_1 = {\bar B}_1 \circ
A^{\varepsilon _1}$, $B_i = (A^{\delta _{i-1}} \circ {\bar B}_i
\circ A^{\varepsilon _i})$, $B_{k+1} = A^{\delta _k} \circ {\bar
B}_{k+1}$ for suitable words ${\bar B}_i$.  This shows that the
decomposition above can be written as
$$W = ({\bar B}_1 A^{\varepsilon _1} D) (D^{-1} A^{q_1} D) \cdots (D^{-1}
A^{\delta _{i-1}} {\bar B}_i  A^{\varepsilon _i} D) \cdots (D^{-1}
A^{q_k}  D)
 (D^{-1} A^{\delta _k}{\bar B}_{k+1}) = $$
$$({\bar B}_1 D) (D^{-1}  A^{\varepsilon _1} D) (D^{-1} A^{q_1} D) \cdots (D^{-1}
A^{\delta _{i-1}}D)( D^{-1}{\bar B}_i D) (D^{-1} A^{^{\varepsilon
_i}} D) \cdots $$ $$(D^{-1} A^{q_k} D)
 (D^{-1} A^{\delta _k} D)(D^{-1} {\bar B}_{k+1})$$
  $$= ({\bar B}_1D) (U^{\varepsilon _1}) (U^{q_1}) \cdots (U^{\delta _{k-1}})( D^{-1}{\bar B}_k D)
 (U^{\varepsilon _k}) (U^{q_k}) (U^{\delta _k})(D^{-1} {\bar B}_{k+1}).$$
 Observe, that the cancellation between parentheses  in the
decomposition above does not exceed the length  $d = |D|$ of $D$.
Using notation $w = u \circ_d v$ to indicate that the cancellation
between $u$ and $v$ does not exceed the number $d$, we can rewrite
the decomposition above in the following form:
$$W = ({\bar B}_1D) \circ_{d}  U^{\varepsilon _1} \circ_{d} U^{q_1} \circ_{d}  U^{\delta _1} \circ_d \cdots
\circ_d U ^{\varepsilon _k} \circ_{d} U^{q_k} \circ_{d} U^{\delta
_k} \circ_{d} (D^{-1} {\bar B}_{k+1}),$$
 hence
 \begin{equation}
 \label{eq:U-dec}
 W = D_1 \circ_d U^{q_1} \circ_d \cdots
\circ_d D_k \circ_d U^{q_k} \circ_d D_{k+1}, \end{equation}
 where
$D_1 = {\bar B}_1D, \ D_{k+1} = D^{-1} {\bar B}_{k+1},  \ D_i=
D^{-1}{\bar B}_i D \ (2 \leqslant i \leqslant k)$, and the
occurrences $U^{q_i}$ are stable (with respect to $\circ_d$).  We
will refer to this decomposition of $W$ as  $U$-decomposition with
respect to ${\mathcal O}$ (to get a rigorous definition of
$U$-decompositions one has to replace in the definition of the
${\mathcal O}$-decomposition of $W$ the period $A$ by $U$ and
$\circ$ by $\circ_{|D|}$). In the case when an $A$-decomposition
of $W$ (with respect to ${\mathcal O}$) is unique then the
corresponding $U$-decomposition of $W$  is also unique, and in
this event  one can easily rewrite $A$-decompositions of $W$ into
$U$-decomposition and vice versa.

We  summarize the discussion above in the following lemma.

\begin{lm} \label{le:A-U}
Let $A \in G[X]$ be a period and $U = D^{-1} \circ A \circ D \in
G[X]$. Then for a word $W \in G[X]$ if

 $$W = B_1 \circ A^{q_1} \circ \cdots \circ B_k \circ A^{q_k} \circ
 B_{k+1}$$
is a stable $A$-decomposition of $W$ then
 $$W = D_1 \circ_d U^{q_1} \circ_d \cdots
\circ_d D_k \circ_d U^{q_k} \circ_d D_{k+1}$$
 is a stable $U$-decomposition of $W$, where $D_i$ are defined as in (\ref{eq:U-dec}). And vice versa.

\end{lm}

 From now on we fix the following  set of leading terms
$${\mathcal A}_{L,p}
= \{ A_j \mid j \leqslant L,  \phi = \phi_{L,p} \}$$ for a given
multiple $L$ of $K = K(m,n)$ and a given tuple $p$.

\begin{df}
Let $W \in G[X]$
 and $N$ be a positive integer.   A word of the type $A_s$ is termed
 the  $N$-large leading term $LT_N(W)$  of the word $W$ if
$A_s^q$ has a stable occurrence in $W$ for some  $q \geq N$, and
$s$ is maximal with this property. The number $s$ is called the
$N$-rank of $W$ ($s = {\rm rank}_N(W), s\geq 1$).
\end{df}

 \begin{lm}
 \label{le:Phiwords}
Let  $W \in G[X]$, $N \geq 2$,  and let $A= LT_N(W)$. Then $W$ can
be presented in the form
\begin{equation}
\label{eq:pPhiform} W = B_1 \circ A^{q_1} \circ \ldots B_k \circ
A^{q_k} \circ B_{k+1}
 \end{equation}
where $A^{q_i}$ are maximal stable occurrences,  $q_i \geq N$, and
${\rm rank}_N(B_i) < {\rm rank}_N(W)$.  This presentation is
unique and it is called   the $N$-large $A$-presentation of $W$.
\end{lm}
\begin{proof} Existence follows from the definition of the leading term
$LT_N(W)$. To prove uniqueness it is suffice to notice  that two
stable occurrences $A^q$ and $A^r$ do not intersect.  Since $A =
LT_N(W)$ is cyclically reduced and it is not a proper power it
follows that an equality $A^2 = u \circ A \circ v$ holds in
$F(X\cup C_S)$ if and only if $u = 1$ or $v = 1$. So, stable
occurrences of $A^q$ and $A^r$ are protected from overlapping by
the neighbors  of $A$ on each side of them. \end{proof}

In Lemmas  \ref{le:7.1.zforms}, \ref{le:7.1.x1formsm0},
\ref{le:7.1.x1formsmneq0}, and   \ref{le:7.1.xiforms}  we
described precisely the leading terms $A _j, j = 1, \ldots,K$ as
the cyclically reduced forms of some words $A_j$.  It is not easy
to describe $A _j$ for an arbitrary $j > K$.  So we are not going
to do it  here, instead, we chose a  compromise by introducing a
modified version of $A_j$ which is not cyclically reduced, in
general, but which  is ``more cyclically reduced" then the initial
word $A_j$.

  Let $L$ be a multiple of $K$ and  $1\leq j\leq K.$ Define
 $$A^*_{L+j}= A^*(\phi _{L+j}) = A _j^{\phi _L}.$$

\begin{lm}
\label{le:cyclicAA*}  Let $L$ be a multiple of $K$ and  $1\leq
j\leq K.$  Let $p=(p_1,\dots ,p_n)$ be $N+3$-large tuple. Then
  $A _{L+j}=cycred(
 A^*(\phi _{L+j}))$.   Moreover, if
 $$A^*(\phi _{L+j})  = R^{-1} \circ A _{L+j} \circ R$$
then  $rank _N(R)\leq L-K+j +2$ and $|R| < |A _{L+j}|.$
\end{lm}

\begin{proof} First, let $L=K.$ Consider elementary periods $x_i=A_{m+4i-3}$
and $A_1=c_1^{z_1}c_2^{z_2}$. For $i\neq n$, $x_i^{2\phi
_K}=x_i^{\phi _k}\circ x_i^{\phi _K}$.  For $i=n$,
$$A^*(\phi _{K+m+4n-3})=R^{-1}\circ A _{K+m+4n-3}\circ R,$$
where $R=A_{m+4i-4}^{p_{m+4n-4}}$, therefore ${\rm
rank}_N(R)=m+4n-4.$ For the other elementary period,
$(c_1^{z_1}c_2^{z_2})^{2\phi _K}=(c_1^{z_1}c_2^{z_2})^{\phi
_K}\circ (c_1^{z_1}c_2^{z_2})^{\phi _K}.$

Any other $A_j$ can be written in the  form $A_j=u_1\circ v_1\circ
u_2\circ v_2\circ u_3$, where $v_1, v_2$ are the first and the
last elementary squares in $A_j$, which are parts of big powers of
elementary periods.  The Nielsen property of $\phi _K$ implies
that the word $R$ for  $A^*(\phi _{K+j})$  is  the word that
cancels between $(v_2u_3)^{\phi _K}$ and $(u_1v_1)^{\phi _K}.$ It
definitely has $N$-large rank  $\leq K$, because the element
$(v_2u_3u_1v_1)^{\phi _K}$ has $N$-large rank $\leq K$. To give an
exact bound for the rank of $R$ we consider all possibilities for
$A _j$:
\begin{enumerate}
\item $A _{i}$ begins with $z_i^{-1}$ and ends with $z_{i+1}$,
$i=1,\dots ,m-1$, \item $A _{m}$ begins with $z_m^{-1}$ and ends
with $x_1^{-1}$, \item $A _{m+4i-4}$ begins with
$x_{i-1}y_{i-2}^{-1}x_{i-2}^{-2}$, if $i=3,\dots n$, and ends with
$x_{i-1}^2y_{i-1}x_{i}^{-1}$ if $i=2,\dots ,n$, If $i=2$ it begins
with $x_1\Pi _{j=m}^1c_j^{-z_j}(c_2^{-z_2}c_1^{-z_1})^2$.

\item $A _{m+4i-2}$  and $A _{m+4i-1}$ begins with
$x_iy_{i-1}^{-1}x_{i-1}^{-2}$ and ends with $x_i^2y_{i}$ if
$i=1,\dots ,n.$

\end{enumerate}

Therefore, $A _{i}^{\phi _K}$ begins with $z_{i+1}^{-1}$ and ends
with $z_{i+2}$, $i=1,\dots ,m-2$, and is cyclically reduced.

 $A _{m-1}^{\phi _K}$ begins with $z_m^{-1}$ and ends with $x_1$, and is cyclically reduced,
 $A _{m}^{\phi _K}$ begins with $z_m^{-1}$ and ends with $x_1^{-1}$ and is cyclically reduced.

 We have already considered $A_{m+4i-3}^{\phi _K}$.

Elements $A_{m+4i-4}^{\phi _K}, A_{m+4i-2}^{\phi _K},
A_{m+4i-1}^{\phi _K}$ are not cyclically reduced. By Lemma
\ref{main}, for $A^*(\phi _K+m+4i-4)$, one has
$R=(x_{i-1}y_{i-2}^{-1})^{\phi _K}$ ($rank _N(R)=m+4i-4$); for
$A^*(\phi _K+m+4i-2)$, and $A^*(\phi _K+m+4i-2)$,
$R=(x_{i}y_{i-1}^{-1})^{\phi _K}$ ($rank _N(R)=m+4i$).

This proves the statement of the Lemma for $L=K$.

We can suppose by induction that $A^*(\phi _{L-K+j})=R^{-1}\circ A
_{L-K+j}\circ R$, and $rank _N(R)\leq L-2K+j+2$. The cancellations
between $A _{L-K+j}^{\phi _K}$ and $R^{\phi _K}$ and between $A
_{L-K+j}^{\phi _K}$ and $A _{L-K+j}^{\phi _K}$ correspond to
cancellations in words $u^{\phi _K}$, where $u$ is a word in
${\mathcal W}_{\Gamma}$ between two elementary squares. These
cancellations  are in rank $\leq K$, and the statement of the
lemma follows. \end{proof}

 \begin{lm}
 \label{le:Phiwords1}
Let  $W \in F(X\cup C_S)$ and  $A = A_j= LT_N(W)$, and $A^* =
R^{-1} \circ A \circ R$. Then $W$ can be presented in the form
\begin{equation}
\label{eq:Phiform} W = B_1 \circ_d A^{* q_1} \circ_d B_2 \circ_d
\cdots \circ_d B_k \circ _d A^{* q_k} \circ_d B_{k+1}
 \end{equation}
where  $A^{* q_i}$ are maximal stable $N$-large occurrences of
$A^*$
 in $W$ and $d \leq |R|$.  This presentation  is unique and it is called
 the
canonical $N$-large $A^*$-decomposition of $W$.
\end{lm}
\begin{proof} The result follows from existence and  uniqueness of the
canonical $A$-decompositions. Indeed, if
 $$W = B_1  \circ A^{ q_1} \circ  B_2 \circ \cdots  \circ B_k \circ A^{q_k}
\circ B_{k+1} $$ is the canonical  $A$-decomposition  of $W$, then
 $$ (B_1R)(R^{-1}AR)^{ q_1}(R^{-1} B_2R) \cdots  (R^{-1}B_kR) (R^{-1}AR)^{q_k}
 (R^{-1}B_{k+1}) $$
is the canonical $A^*$-decomposition of $W$. Indeed, since every
$A^{q_i}$ is a stable occurrence, then every $B_i$ starts with $A$
(if $i \neq 1$) and ends with $A$ (if $i = k+1$). Hence $R^{-1}
B_iR= R^{-1} \circ B_i \circ R$.

Conversely if
 $$ W = B_1 A^{* q_1} B_2 \cdots B_k A^{* q_k} B_{k+1}$$
  is an $A^*$-representation of $W$ then
 $$W = (B_1R^{-1})  \circ A^{ q_1} \circ  (RB_2R^{-1}) \circ \cdots  \circ (RB_kR^{-1})
  \circ A^{q_k}
\circ (RB_{k+1}) $$ is the canonical  $A$-decomposition for $W$.
\end{proof}

In the following lemma we collect various properties of words
$x_i^{\phi_L}, y_i^{\phi_L}, z_j^{\phi_L}$ where $L = Kl$  is a
multiple of $K$.

\begin{lm}
\label{le:properties-gamma-words}
 Let $X = \{x_i, y_i, z_j \mid i = 1, \ldots, n, j = 1, \ldots, m\}$,
 let $K = K(m,n)$, and $L = Kl$ be a multiple of $K$.  Then
 for any  number $N\geq 5$ and for any $N$-large
 tuple $p \in N^L$ the following holds (below $\phi  = \phi_{L,p}$, $A_j = A_j$):

 \begin{enumerate}
 \item  [(1)] If $i < j \leq L$ then $A_j^2$ does not occur in $A_i$;
\item  [(2)] Let $i \leq K,\ j=i+L$. There are positive integers
$s$, $1 \leq j_1, \ldots, j_s \leq j$, integers $\varepsilon_1,
\ldots, \varepsilon_s$ with $|\varepsilon_t| \leq 3$, and words
$w_1, \ldots, w_{s+1} \in F(X \cup C_S)$ ( which do not depend on
the tuple $p$ and do not contain any square of leading terms) such
that the leading term $A_{j}$ ($A_j^*$) of $\phi_{j}$ has the
following form:
\begin{equation}\label{55} w_1 \circ A_{j_1}^{p_{j_1} + \varepsilon_1} \circ w_2 \circ
\cdots \circ w_s \circ A_{j_s}^{p_{j_s} + \varepsilon_s} \circ
w_{s+1},\end{equation}
 i.e., the "periodic structure" of $A_{j}$ ($A_j^*$) does not depend on the
 tuple $p$.
\item  [(3)] Let $i\leq K$, $u\in {\mathcal W}_{\Gamma,L}$ such
that
$$u=v_1 \circ A_{j_1}^{p_{j_1} + \varepsilon _1} \circ v_2 \circ \cdots
\circ v_r \circ A_{j_r}^{p_{j_r} + \varepsilon_r} \circ v_{r+1},$$
where $j_1,\dots ,j_r\leq i$, and at least one of $j_t$ is equal
to $i$, $|\varepsilon_t| \leq 1$, and words $v_1, \ldots, v_{r+1}
\in F(X \cup C_S)$ do not depend on $p$. Then
$$u^{\phi _L}=v_1^{\phi _L}A_{j_1}^{\sigma _1\phi _L}W_1^{-1}
\circ A_{j_1+L}^{(p_{j_1} + \varepsilon _1-2\sigma _1)} \circ
W_1A_{j_1}^{\sigma _1\phi _L}v_2^{\phi _L} \  \ldots \ $$ $$
v_r^{\phi _L}A_{j_r}^{\sigma _r\phi _L}W_r^{-1} \circ
A_{j_r+L}^{(p_{j_r} + \varepsilon_r-2\sigma _r)} \circ
W_rA_{j_r}^{\sigma _r\phi _L}v_{r+1}^{\phi L},$$
 where
$A_{j_t}^{\phi _L}=W_t^{-1}\circ A_{j_t+L}\circ W_t;\  \sigma _t
=1$ if $p_t$ is positive and $\sigma _t=-1$ if $p_t$ is negative.
In addition, for each $t=1,\dots ,r$ the product
$$W_tA_{j_t}^{\sigma _t\phi _L}v_{t+1}^{\phi _L}A_{j_{t+1}}^{\sigma
_{t+1}\phi _L}W_{t+1}^{-1}$$ has form {\rm (\ref{55})} with
$j_1,\dots ,j_s< i+L.$

\item  [(4)]  For any  $i\leq K$ and any $x \in X^{\pm 1}$ there
is a positive integer $s$ and there are indices $1 \leq j_1,
\ldots, j_s \leq i$, integers $\varepsilon_1, \ldots,
\varepsilon_s$ with $|\varepsilon_t| \leq 1$, and words $w_1,
\ldots, w_{s+1} \in F(X \cup C_S)$  which do not depend on the
tuple $p$  such that the element $x^{\phi_i}$ can be presented in
the following form:
$$x^{\phi_i} = w_1 \circ A_{j_1}^{p_{j_1} + \varepsilon_1}
\circ w_2 \circ \cdots \circ w_s \circ A_{j_s}^{p_{j_s} +
\varepsilon_s} \circ w_{s+1}.$$

\end{enumerate}
\end{lm}

\begin{proof} Statement (1) follows from Lemmas
\ref{le:7.1.zforms}--\ref{le:7.1.x1formsmneq0}.

Statements (2) and (3) will be proved by simultaneous induction on
$j=i+L$. Case $l=0$ corresponds to $i\leq K$. In this case
statement (2) follows from Lemmas \ref{le:7.1.zforms} -
\ref{le:7.1.x1formsmneq0} and statement (3) is simply the
assumption of the lemma.  $A_i$ has form (\ref{55}) with
$j_1,\dots ,j_s< i$ and $|\varepsilon_1|, \ldots,
|\varepsilon_s|\leq 1$.

We know that $A_{j_t}$ contains an elementary square (actually,
big power) for any $t=1,\dots s$, $A_{j_t}^{\phi
_K}=R_{j_t}^{-1}\circ A_{j_t+K}\circ R_{j_t}$, where $R_{j_t}$
does not contain big powers of $A_k$ for $k\geq j_t+2$.  Then it
follows from the second statement of Lemma \ref{main} that
\begin{multline*} A_i^{\phi _K}= w_1^{\phi _K}R_{j_1}^{-\sigma _1
}A_{j_1+K}^{\sigma _1} \circ A_{j_1+K}^{p_{j_1} +
\varepsilon_1-2\sigma _1} \circ A_{j_1+K}^{\sigma
_1}R_{j_1}^{\sigma _1} w_2^{\phi _K}\
 \ldots  \ w_s^{\phi _K}R_{j_s}^{-\sigma
_s }A_{j_s+K}^{\sigma _s}\\ \circ A_{j_s+K}^{p_{j_s} +
\varepsilon_s-2\sigma _s} \circ A_{j_s+K}^{\sigma _
s}R_{j_s}^{\sigma _s} w_{s+1}^{\phi _K},\end{multline*} where
$\sigma _t=1$ if $p_{j_t}$ is positive and $\sigma _t=-1$ if
$p_{j_t}$ is negative.

When we apply $\phi _K$,  the images of elementary big powers in
$A_{j_t}$ by Lemma \ref{main} are not touched by cancellations
between $w_{t-1}^{\phi _K}$ and $A_{j_t}^{\phi _K}$, and between
$A_{j_t}^{\phi _K}$ and $w_{t+1}^{\phi _K}$, therefore $A_i^{\phi
_L}=$ $$w_1^{\phi _L}R_{j_1}^{-\sigma _1\phi _{L-K}
}A_{j_1+K}^{\sigma _1\phi _{L-K}}W_1^{-1} \circ A_{j_1+L}^{p_{j_1}
+ \varepsilon_1-2\sigma _1} \circ W_1A_{j_1+K}^{\sigma _1\phi
_{L-K}}R_{j_1}^{\sigma _1\phi _{L-K}} w_2^{\phi _L}\
 \cdots $$ $$ \ w_s^{\phi _L}R_{j_s}^{-\sigma _s\phi _{L-K}
}A_{j_s+K}^{\sigma _s\phi _{L-K}}W_s^{-1}
 \circ
A_{j_s+L}^{p_{j_s} + \varepsilon_s-2\sigma _s} \circ
W_sA_{j_s+K}^{\sigma _ s\phi _{L-K}}R_{j_s}^{\sigma _s\phi _{L-K}}
w_{s+1}^{\phi _L},$$ where $A_{j_t+K}^{\phi _{L-K}}=W_t^{-1}\circ
A_{j_t+L}\circ W_t$,  $\sigma _t=1$ if $p_{j_t}$ is positive and
$\sigma _t=-1$ if $p_{j_t}$ is negative ( $t=1,\dots ,s,$). We can
now apply statement 3) for $i_1+Kl, i_1<i$ to elements

\noindent
 $w_1^{\phi _L}R_{j_1}^{-\sigma _1\phi _{L-K}
}A_{j_1+K}^{\sigma _1\phi _{L-K}W_1^{-1}},\dots ,\
W_sA_{j_s+K}^{\sigma _ s\phi _{L-K}}R_{j_s}^{\sigma _s\phi _{L-K}}
w_{s+1}^{\phi _L}.$

To prove statement (3) for $i+Kl$, we use it for $i_1+Kl$ and
statement (2) for $i+Kl.$

(4) Existence of such a decomposition follows from Lemmas
\ref{le:7.1.zforms}--\ref{le:7.1.x1formsmneq0}. \end{proof}

\begin{cy}\label{R} If $L$ is a multiple of $K$, then the automorphism $\phi _{L}$
satisfies  the Nielsen property with respect to $\bar {\mathcal
W}_ {\Gamma}$ with exceptions $ E(n,m)$. \end{cy}

\begin{proof} The middles  $M_x$ of elements from $X$ and from $E(m,n)$ with
respect to $\phi _K$ contain  big powers of some $A _j$, where
$j=1,\dots ,K.$ By Lemma \ref{le:properties-gamma-words}
 these big powers cannot disappear after application of $\phi
_{L-K}$. Therefore, $M_x^{\phi _{L-K}}$ contains the middle of $x$
with respect to $\phi _L.$ \end{proof}

\begin{cy}\label{co:12-new}
Let $u, v \in \bar {\mathcal W}_ {\Gamma}$. If  the canceled
subword
 in the product  $u^{\phi_{K}}v^{\phi_{K}}$ does not contain
$A_j^l$ for some $j \leq K$ and $l \in \mathbb{Z}$ then the
canceled subword in the product $u^{\phi_{K+L}}v^{\phi_{K+L}}$
does not contain the subword $A_{L+j}^l$.
\end{cy}

\begin{lm}\label{le:A-A^*} Let $W \in {\mathcal W}_{\Gamma, L}$. Suppose that
 $1\leq r\leq K$, $ L_1$ is a multiple of $K$, and $j=r+L_1.$
   Then the following conditions are equivalent:
    \begin{enumerate}
    \item [1)]    $rank_N(W) = r$ and
    $$W=D_1\circ
A_r^{q_1}\circ D_2\ldots D_k\circ A_r^{q_k}\circ D_{k+1}$$ is a
stable $5$-large  $A_r$-decomposition of $W$;

  \item [2)]
$rank_N(W^{\phi_{L_1}}) = j$  and
$$W^{\phi _{L_1}}=(D_1A_j^{\varepsilon _1})^{\phi _{L_1}}\circ _d
A_j^{*q_1-\varepsilon _1 -\delta_1}\circ _d(A_j^{\delta _1}D_2
A_j^{\varepsilon _2})^{\phi _{L_1}}\ldots $$ $$(A_j^{\delta
_{k-1}} D_kA_j^{\varepsilon _k})^{\phi _{L_1}}\circ_d
A_j^{*q_k-\varepsilon _k -\delta_k}\circ (A_j^{\delta
_k}D_{k+1})^{\phi _{L_1}}$$
  is a stable
$A_j^*$-decomposition of $W^{\phi _{L_1}}$, where $\delta _s,
\varepsilon_s \in \{0,\pm 1\}$ depending on the sign of $q_s$ and
$\beta$.
 \end{enumerate}
\end{lm}
\begin{proof} It follows from Lemmas \ref{le:Phiwords1} and
 \ref{le:properties-gamma-words}. Indeed, let $W \in  {\mathcal W}_{\Gamma, L}$ and
    $$W=D_1\circ
A_r^{q_1}\circ D_2\ldots D_k\circ A_r^{q_k}\circ D_{k+1}$$  the
canonical $N$-large $A_r$-decomposition of $W$. Then by Lemma
\ref{le:properties-gamma-words} (3)
 $$
W^{\phi _{L_1}}=(D_1A_r^{\sigma _1})^{\phi _{L_1}}w_r^{-1} \circ
A_j^{q_1-2\sigma _1}\circ w_r(A_r^{\sigma _1}D_2 A_r^{\sigma
_2})^{\phi _{L_1}}\ldots $$ $$(A_r^{\sigma _{k-1}} D_kA_r^{\sigma
_k})^{\phi _{L_1}}w_r^{-1}\circ  A_j^{q_k-2\sigma _k}\circ
w_r(A_r^{\sigma _k}D_{k+1})^{\phi _{L_1}}
$$
where $A_r^{\phi_{L_1}} = w_r^{-1} \circ A_j \circ w_r$, $\sigma_t
\in \{1,-1\}$. This implies that the canonical $A^*$-decomposition
of $W^{\phi _{L_1}}$ takes the form described in 2).

Conversely, suppose  2) is the canonical $A^*$-decomposition of
$W^{\phi _{L_1}}$, but 1) is not the canonical $A_r$-decomposition
of $W$. Then taking the canonical $A_r$-decomposition of $W$ and
applying $\phi_{L_1}$ by 1) we get another canonical decomposition
of $W^{\phi _{L_1}}$ - contradiction with uniqueness of
$A^*$-decompositions.
\end{proof}

\begin{lm}\label{le:7.1.gammawords}

  Suppose $p$ is an $(N+3)$-large tuple, $\phi _j=\phi _{jp}$. Let $L$ be a
multiple of $K$. Then:

\begin{enumerate}\item [(1)]
\begin{enumerate}
   \item $x_i^{\phi _j}$ has a canonical $N$-large
$A^*_j$-decomposition of size $(N,2)$ if  either $j\equiv
m+4(i-1)(mod \ K)$, or $j \equiv m+4i-2 (mod \ K)$, or $j \equiv
m+4i (mod \ K)$. In all other cases $rank(x_i^{\phi _j}) < j$.
    \item $ y_i^{\phi _j}$ has a canonical $N$-large $A^*_j$-decomposition
 of size $(N,2)$ if either $j\equiv m+4(i-1)(mod \ K)$, or $j \equiv m+4i-3(mod \
K)$, or $j \equiv m+4i-1(mod \ K)$, or $j \equiv  m+4i\ (mod \
K).$ In all other cases $rank(y_i^{\phi _j}) < j$.
   \item $z_i^{\phi _{j}}$ has a canonical $N$-large
   $A^*_j$-decomposition of size $(N,2)$
 if $j\equiv i\ (mod \ K)$ and either $1 \leq i \leq m-1$ or $i = m$
 and $n \neq 0$. In all other
cases $rank(z_i^{\phi _j}) < j$.
    \item if $n = 0$ then $z_m^{\phi _{j}}$ has a canonical $N$-large
    $A^*_j$-decomposition of size $(N,2)$
 if $j\equiv m-1\ (mod \ K)$.  In all other
cases $rank(z_m^{\phi _j}) < j$.
 \end{enumerate}

\item [(2)] If $j=r+L$, $0<r\leq K,\ (w_1\ldots w_k)\in
Sub_k(X^{\pm\gamma_K \ldots\gamma_{r+1}})$ then either $(w_1\ldots
w_k)^{\phi _j}=(w_1\ldots w_k)^{\phi _{j-1}},$ or $ (w_1\ldots
w_k)^{\phi _j}$ has a canonical $N$-large $A^*_j$-decomposition.
In any case, $ (w_1\ldots w_k)^{\phi _j}$   has a canonical
$N$-large $A^*_s$-decomposition in some rank $s,$ $j-K+1\leq s\leq
j.$
\end{enumerate}
\end{lm}

\begin{proof} (1) Consider  $y_i^{\phi _{L+m+4i}}:$

$$y_i^{\phi _{L+m+4i}}=(   x_{i+1}^{\phi _{L}} y_i^{-\phi _{L+
m+4i-1}})^{q_4-1}x_{i+1}^{\phi _{L}} (y_i^{\phi _{ L+m+4i-1}}
x_{i+1}^{-\phi _{L}})^{q_4},$$

In this case $A^*(\phi _{L+m+4i})=x_{i+1}^{\phi
_{L+m+4i-1}}y_i^{-\phi _{L+m+4i-1}}$.

To write a formula for  $x_i^{\phi _{L+m+4i}}$, denote $\tilde
y_{i-1}=y_{i-1}^{\phi _{L+m+4i-5}},\ \bar x_i=x_i^{\phi _{L}}, \
\bar y_i=y_i^{\phi _{L}}$. Then

\medskip
\begin{multline*} x_i^{\phi _ {L+m+4i}}=( \bar x_{i+1} y_i^{-\phi _
{L+m+4i-1}})^{q_4-1} \bar x_{i+1}\\ ( ((
 \bar x_i\tilde  y_{i-1}^{-1})^{q_0}    \bar x_i^{q_1}    \bar y_i)^{q_2-1} (
   \bar x_i\tilde
y_{i-1}^{-1})^{q_0}    \bar x_i^{q_1+1}    \bar y_i)^{-q_3+1} \bar
y_i^{-1} \bar x_i^{-q_1}(\tilde y_{i-1} \bar
x_i^{-1})^{q_0}.\end{multline*}

 Similarly we consider
 $ z_i^{\phi _{L+i}}$.

(2)  If in a word $(w_1\cdots w_k)^{\phi _j}$ all the powers of
$A_j^{p_j}$ are cancelled (by Lemma
\ref{le:properties-gamma-words} they can only cancel completely
and the process of cancellations does not depend on $p$) then if
we consider an $A_j^*$-decomposition of $(w_1\cdots w_k)^{\phi
_j}$, all the powers of $A_j^*$ are also completely cancelled. By
construction of the automorphisms $\gamma _j$, this implies that
$(w_1\cdots w_k)^{\gamma _j\phi _{j-1}}= (w_1\cdots w_k)^{\phi
_{j-1}}.$
\end{proof}

\section{Generic  solutions of orientable quadratic equations}
\label{se:7.2}

Let $G$ be a finitely generated fully residually free group and $S
= 1$ a standard quadratic orientable equation over $G$ which has a
solution in $G$. In this section we effectively construct
discriminating  sets of solutions of $S = 1$ in  $G$. The main
tool in this construction is an embedding
 $$\lambda: G_{R(S)} \rightarrow G(U,T)$$
of the coordinate group $G_{R(S)}$ into a group $G(U,T)$ which is
obtained from  $G$ by finitely many extensions of centralizers.
There is a nice  set $\Xi_P$ (see Section 1.4 in \cite{JSJ}) of
discriminating $G$-homomorphisms from $G(U,T)$ onto $G$. The
restrictions of homomorphisms from $\Xi_P$
 onto  the image $G_{R(S)}^\lambda$ of $G_{R(S)}$ in $G(U,T)$
 give a discriminating  set of
 $G$-homomorphisms from $G_{R(S)}^\lambda$ into $G$, i.e.,
solutions of $S = 1$ in $G$. This idea was introduced in
\cite{KMNull}  to describe the radicals of quadratic equations.

It has been shown in \cite{KMNull} that the coordinate groups of
non-regular standard quadratic equations $S = 1$ over $G$ are
already extensions of centralizers of $G$, so in this case we can
immediately put $G(U,T) = G_{R(S)}$ and the result follows. Hence
we can assume from the beginning that $S= 1$ is regular.

Notice, that all regular quadratic  equations have solutions in
general position, except for the equation $[x_1,y_1][x_2,y_2] = 1$
(see \cite{KMIrc}, Section 2).

 For the equation $[x_1,y_1][x_2,y_2]
= 1$ we do the following trick. In this case we view the
coordinate group $G_{R(S)}$  as the coordinate group of the
equation $[x_1,y_1]= [y_2,x_2]$ over the group of constants $G
\ast F(x_2,y_2)$. So the commutator $[y_2,x_2] = d$ is a
non-trivial constant  and the new equation is of the form $[x,y] =
d$, where all solutions are in general position. Therefore, we can
assume that  $S = 1$ is  one of the following types (below $d,c_i$
are nontrivial elements from $G$):
\begin{equation}\label{eq:st7}
\noindent \prod_{i=1}^{n}[x_i,y_i] = 1, \ \ \ n \geqslant 3;
\end{equation}
\begin{equation}\label{eq:st8}
\prod_{i=1}^{n}[x_i,y_i] \prod_{i=1}^{m}z_i^{-1}c_iz_i d = 1,\ \ \
n \geqslant 1, m \geqslant 0;
\end{equation}
\begin{equation}\label{eq:st9}
 \prod_{i=1}^{m}z_i^{-1}c_iz_i d = 1,\ \ \  m \geqslant 2,
\end{equation}
and it has a solution in $G$ in general position.

Observe, that since $S= 1$ is regular then Nullstellenzats  holds
for $S = 1$, so $R(S) = ncl(S)$ and $G_{R(S)} = G[X]/ncl(S) =
G_S$.

For a group $H$ and an element $u \in H$ by  $H(u,t)$ we denote
the  extension of the centralizer $C_H(u)$ of $u$:
$$H(u,t) = \langle H,t \ | \ t^{-1}xt = x  \ \ (x \in C_H(u)) \rangle.$$
If
 $$ G = G_1  \leqslant G_1(u_1,t_1) = G_2 \leqslant  \ldots \leqslant G_n(u_n,t_n) = G_{n+1}
$$
   is a chain of extensions of centralizers of elements $u_i \in G_i$,
then we denote the resulting group $G_{n+1}$ by $G(U,T)$, where $U
= \{u_1, \ldots, u_n\}$ and $T = \{t_1, \ldots, t_n\}$.

  Let   $\beta: G_{R(S)} \rightarrow G$ be  a solution of the equation
  $S(X) = 1$ in the group $G$ such that
 $$x_i^\beta = a_i, y_i^\beta = b_i, z_i^\beta = e_i. $$
Then
 $$d = \prod_{i=1}^{m}e_i^{-1}c_ie_i \prod_{i=1}^{n}[a_i,b_i].$$
Hence we can rewrite the equation $S = 1$  in the following form
(for appropriate $m$ and $n$):
\begin{equation}\label{6-general}
\prod_{i=1}^{m}z_i^{-1}c_iz_i\prod_{i=1}^{n}[x_i,y_i]  =
\prod_{i=1}^{m}e_i^{-1}c_ie_i \prod_{i=1}^{n}[a_i,b_i].
\end{equation}

\begin{prop}
\label{prop:hom-lambda} Let $S = 1$ be a regular quadratic
equation {\rm (\ref{6-general})}  and $\beta:G_{R(S)} \rightarrow G$ a
solution of $S= 1$ in $G$ in a general position. Then one can
effectively construct a sequence of extensions of centralizers
$$ G = G_1  \leqslant G_1(u_1,t_1) = G_2 \leqslant  \ldots \leqslant G_n(u_n,t_n) = G(U,T)$$
 and a $G$-homomorphism $\lambda_\beta : G_{R(S)} \rightarrow G(U,T)$.
 \end{prop}
\begin{proof}  By induction we define a sequence of  extensions of
centralizers and a sequence of group homomorphisms in the
following way.

{\it Case: $m \neq 0, n = 0$.} In this event for each $i = 1,
\ldots, m-1$ we define by induction  a pair $(\theta_i, H_i)$,
consisting of a group $H_i$ and a $G$-homomorphism $\theta_i:G[X]
\rightarrow H_i$.

Before we will go into formalities let us explain the idea that
lies behind this. If $z_1 \rightarrow e_1, \ldots, z_m \rightarrow
e_m$ is a solution of an equation
\begin{equation} \label{eq:zz}
z_1^{-1}c_1z_1 \ldots z_m^{-1}c_mz_m = d,
\end{equation}
then  transformations
\begin{equation}\label{eq:ztr2}
e_i \rightarrow e_i\left(c_i^{e_i}c_{i+1}^{e_{i+1}}\right)^q, \ \ e_{i +1}
\rightarrow e_{i+1}\left(c_i^{e_i}c_{i+1}^{e_{i+1}}\right)^q,  \ \ e_j
\rightarrow e_j  \ \ \ (j \neq i, i+1),
\end{equation}
produce a new solution of the equation (\ref{eq:zz}) for an
arbitrary integer $q$. This solution is composition of the
automorphism $\gamma _{i}^q$ and the solution $e$. To avoid
collapses under cancellation of the periods
$(c_i^{e_i}c_{i+1}^{e_{i+1}})^q$ (which is an important part of
the construction of the discriminating set of homomorphisms
$\Xi_P$ in Section 1.4 in \cite{JSJ}) one might want to have
number $q$ as big as possible, the best way would be to have $q =
\infty$. Since there are no infinite powers in $G$, to realize
this idea one should go outside the group $G$ into a bigger group,
 for example, into an ultrapower $G^\prime$ of $G$, in which
 a non-standard power, say $t$,  of the element $c_i^{e_i}c_{i+1}^{e_{i+1}}$
 exists. It is not hard to see that the subgroup $\langle G,t\rangle \leqslant G^\prime$
 is an extension of the centralizer $C_G(c_i^{e_i}c_{i+1}^{e_{i+1}})$ of
 the element $c_i^{e_i}c_{i+1}^{e_{i+1}}$ in  $G$. Moreover, in the group
 $\langle G,t\rangle $
 the transformation  (\ref{eq:ztr2}) can be described as
\begin{equation}\label{eq:ztr3}
e_i \rightarrow e_it, \ \ e_{i +1} \rightarrow e_{i+1}t, \ \ e_j
\rightarrow e_j  \ \ \ (j \neq i, i+1),
\end{equation}
Now, we are going to construct formally the subgroup  $\langle
G,t\rangle $ and the corresponding homomorphism using
(\ref{eq:ztr3}).

 Let $H$ be an arbitrary group and $\beta: G_{S} \rightarrow H$
a homomorphism. Composition of the canonical projection $G[X]
\rightarrow G_S$ and $\beta$ gives a homomorphism $\beta_0: G[X]
\rightarrow H$.   For $i =0$ put
$$H_0 = H, \ \ \ \theta _0 = \beta_0$$
 Suppose now, that a group $H_{i}$ and
a homomorphism $\theta _{i}: G[X] \rightarrow H_i$ are already
defined. In this event we define $H_{i+1}$ and $\theta_{i+1}$ as
follows
$$
H_{i+1}=\left< H_{i},r_{i+1} \mid
\left[C_{H_{i}}(c_{i+1}^{z_{i+1}^{\theta
_{i}}}c_{i+2}^{z_{i+2}^{\theta_{i}}}),r_{i+1}\right]=1\right>, $$
$$z_{i+1}^{\theta_{i+1}}=z_{i+1}^{\theta_{i}}r_{i+1}, \ \
z_{i+2}^{\theta_{i+1}}=z_{i+2}^{\theta_{i}}r_{i+1}, \ \
z_j^{\theta_{i+1}}= z_j^{\theta_{i}}, \ \ \ (j \neq i+1, i+2).
$$
By induction we constructed a series of extensions of centralizers
$$
G = H_0 \leqslant H_1 \leqslant \ldots \leqslant H_{m-1} =
H_{m-1}(G)
$$
and a homomorphism
$$\theta_{m-1,\beta} = \theta _{m-1}: G[X] \rightarrow H_{m-1}(G).$$
 Observe, that,
  $$c_{i+1}^{z_{i+1}^{\theta
_{i}}}c_{i+2}^{z_{i+2}^{\theta_{i}}} =
c_{i+1}^{e_{i+1}r_{i}}c_{i+2}^{e_i}$$
 so  the element $ r_{i+1}$ extends the
centralizer of the element $c_{i+1}^{e_{i+1}r_{i}}c_{i+2}^{e_i}$.
In particular, the following equality holds in the group
$H_{m-1}(G)$ for each $i = 0, \ldots, m-1$:
\begin{equation}\label{eq:r}
[r_{i+1},c_{i+1}^{e_{i+1}r_{i}}c_{i+2}^{e_i}] = 1.
\end{equation}
(where $r_0 = 1$). Observe also, that
\begin{equation}
\label{eq:z} z_1^{\theta_{m-1}} =  e_1r_1, \ \ z_i^{\theta_{m-1}}
=  e_ir_{i-1}r_{i}, \ \ z_m^{\theta_{m-1}} =  e_mr_{m-1} \ \ \
(0<i<m).
\end{equation}

From (\ref{eq:r}) and (\ref{eq:z}) it readily follows that
\begin{equation}
\label{eq:homz} \left(\prod_{i=1}^{m}z_i^{-1}c_iz_i \right) ^{\theta_{m- 1}}
= \prod_{i=1}^{m}e_i^{-1}c_ie_i,
\end{equation}
so   $\theta_{m-1}$ gives rise to a homomorphism (which we again
denote by $\theta_{m-1}$ or $\theta_\beta$)
 $$\theta_{m-1} : G_S \longrightarrow H_{m-1}(G). $$
  Now we iterate the construction one more time replacing $H$ by $H_{m-1}(G)$
  and $\beta$ by $\theta_{m-1}$ and put:
  $$H_\beta(G) = H_{m-1}(H_{m-1}(G)), \ \ \
  \lambda_\beta = \theta_{\theta_{m-1}}: G_S \longrightarrow H_\beta(G).$$
The group $H_\beta(G)$ is union of a chain of extensions of
centralizers which starts at the group $H$.

If $H = G$ then  all the homomorphisms
 above are $G$-homomorphisms. Now  we can write
  $$H_\beta(G) = G(U,T)$$
  where $U = \{u_1, \ldots, u_{m-1},
\bar u_1, \ldots, \bar u_{m-1}\}$, $T = \{r_1, \ldots, r_{m-1},
\bar r_1, \ldots, \bar r_{m-1}\}$ and   $ \bar u_i, \bar r_i$ are
the corresponding elements when
 we iterate the construction:
 $$
 u_{i+1} =
c_{i+1}^{e_{i+1}r_{i}}c_{i+2}^{e_{i+2}}, \ \  \bar u_{i+1} =
c_{i+1}^{e_{i+1}r_{i}r_{i+1}\bar
r_i}c_{i+2}^{e_{i+2}r_{i+1}r_{i+2}}.$$

{\it Case:  $m = 0, n > 0$. } In this case $S = [x_1,y_1] \cdots
[x_n,y_n]d^{-1}$. Similar to the case above we start with the
principal automorphisms. They consist of two Dehn's twists:
\begin{equation}\label{t1}
x \rightarrow y^px, \ \ y \rightarrow y;
\end{equation}
\begin{equation}\label{t2}
x  \rightarrow x, \ \ y \rightarrow x^py;
\end{equation}
which fix the commutator $[x,y]$, and the third transformation
which ties two consequent commutators
$[x_i,y_i][x_{i+1},y_{i+1}]$:
\begin{equation}\label{t3}
x_i \rightarrow (y_ix_{i+1}^{-1})^{-q}x_i, \ \ y_i \rightarrow
(y_ix_{i+1}^{- 1})^{-q}y_i(y_ix_{i+1}^{-1})^q,
\end{equation}
$$x_{i+1} \rightarrow (y_ix_{i+1}^{-1})^{-q}x_{i+1}(y_ix_{i+1}^{- 1})^q, \ \
y_{i+1} \rightarrow (y_ix_{i+1}^{-1})^{-q}y_{i+1}.$$

 Now we define
by induction on $i$,  for $i = 0, \ldots, 4n-1$,   pairs
$(G_i,\alpha_i)$ of groups $G_i$ and $G$-homomorphisms $\alpha_i:
G[X] \rightarrow G_i$. Put
 $$G_0 = G, \ \ \alpha_0  = \beta.$$
For each commutator $[x_i,y_i] $ in $S = 1$ we  perform
consequently three Dehn's twists (\ref{t2}), (\ref{t1}),
(\ref{t2}) (more precisely, their analogs for an extension of a
centralizer) and an analog of the connecting transformation
(\ref{t3}) provided the next commutator exists. Namely, suppose
$G_{4i}$ and $\alpha _{4i}$ have been already defined. Then

\bea G_{4i+1} &=& \left<G_{4i},t_{4i+1}\mid [C_{G_{4i}}(x_{i+1}^{\alpha
_{4i}}),t_{4i+1}]=1\right>;\\
y_{i+1}^{\alpha
_{4i+1}} &=& t_{4i+1}y_{i+1}^{\alpha _{4i}}, \ \ s^{\alpha
_{4i+1}}=s^{\alpha _{4i}} \ \ (s \neq y_{i+1}).\\
G_{4i+2} &=& \left<G_{4i+1},t_{4i+2} \mid [C_{G_{4i+1}}(y_{i+1}^{\alpha
_{4i+1}}),t_{4i+2}]=1\right>;\\
x_{i+1}^{\alpha _{4i+2}} &=& t_{4i+2}x_{i+1}^{\alpha _{4i+1}}, \ \
s^{\alpha _{4i+2}}=s^{\alpha _{4i+1}} \ \ \ (s \neq x_{i+1});\\
 G_{4i+3} &=& \left< G_{4i+2},t_{4i+3} \mid \left[C_{G_{4i+2}}(x_{i+1}^{\alpha
_{4i+2}}),t_{4i+3}\right]=1\right>;\\
 y_{i+1}^{\alpha
_{4i+3}} &=& t_{4i+3}y_{i+1}^{\alpha _{4i+2}}, \ \ s^{\alpha
_{4i+3}}=s^{\alpha _{4i+2}} \ \ \ (s \neq y_{i+1});\\
G_{4i+4} &=& \left< G_{4i+3},t_{4i+4}\mid \left[C_{G_{4i+3}}\left(y_{i+1}^{\alpha
_{4i+3}}x_{i+2}^{-\alpha _{4i+3}}\right),t_{4i+4}\right]=1\right>;\\
x_{i+1}^{\alpha _{4i+4}} &=& t_{4i+4}^{-1}x_{i+1}^{\alpha _{4i+3}},
y_{i+1}^{\alpha _{4i+4}}=y_{i+1}^{\alpha _{4i+3}t_{4i+4}},
x_{i+2}^{\alpha _{4i+4}}=x_{i+2}^{\alpha _{4i+3}t_{4i+4}},\\
y_{i+2}^{\alpha _{4i+4}} &=& t_{4i+4}^{-1}y_{i+2}^{\alpha _{4i+3}};\\
s^{\alpha _{4i+4}} &=& s^{\alpha _{4i+3}} \ \  (s \neq
x_{i+1},y_{i+1},x_{i+2},y_{i+2}).
\eea
 Thus we have defined groups
$G_i$ and  mappings $\alpha_i$ for all $i = 0, \dots, 4n-1$.    As
above, the straightforward verification shows that the mapping
$\alpha_{4n-1}$ gives rise to a $G$-homomorphism $\alpha_{4n-1}:
G_S \longrightarrow G_{4n-1}.$ We repeat now the above
construction once more  time with $G_{4n-1}$ in the place of
$G_0$, $\alpha _{4n-1}$ in the place of $\beta$, and  $\bar t_j$
in the place of $ t_j$. We denote the corresponding groups and
homomorphisms by ${\bar G_i}$ and ${\bar \alpha_i}: G_S
\rightarrow {\bar G_i}$.

 Put
 $$G(U,T) = \bar G_{4n-1}, \ \ \ \lambda_\beta = \bar \alpha_{4n-1},$$
By induction we have constructed  a $G$-homomorphism
$$\lambda_\beta : G_S \longrightarrow G(U,T).$$

{\it Case:}  $m > 0, n > 0$.  In this case we combine the two
previous cases together. To this end we take the group $H_{m-1}$
 and the homomorphism $\theta_{m-1}: G[X] \rightarrow H_{m-1}$
 constructed in the first case and put them as the input for
 the construction in the second case. Namely, put
$$G_0=\left<H_{m-1},r_{m}| [C_{H_{m-1}}(c_m^{z_m^{\theta_{m-1}}}
x_1^{-\theta_{m- 1}}),r_{m}]=1\right>,$$
 and define the homomorphism $\alpha_0$ as follows
 $$z_m^{\alpha _0}=z_m^{\theta
_{m-1}}r_{m}, \ \ x_1^{\alpha _0}=a_1^{r_{m}}, \ \ y_1^{\alpha
_0}=r_{m}^{-1}b_1, \ \  s^{\alpha _0}=s^{\theta _{m- 1}} \ \ (s\in
X, s \neq z_m, x_1,y_1).$$
  Now we  apply the construction from the second case.
 Thus we have defined groups
$G_i$ and  mappings $\alpha_i: G[X] \rightarrow G_i$ for all $i =
0, \dots, 4n-1$. As above, the straightforward verification shows
that the mapping $\alpha_{4n-1}$ gives rise to a $G$-homomorphism
$\alpha_{4n-1}: G_S \longrightarrow G_{4n-1}.$

We repeat now the above construction once more time with
$G_{4n-1}$ in place of $G_0$ and  $\alpha _{4n-1}$ in place of
$\beta$. This results in a group $\bar G_{4n-1}$ and a
homomorphism $\bar \alpha_{4n-1}: G_S \rightarrow \bar G_{4n-1}$.

Put
 $$G(U,T) = \bar G_{4n-1}, \ \ \ \lambda_\beta = \bar \alpha_{4n-1}.$$
We have constructed  a $G$-homomorphism
$$\lambda_\beta : G_S \longrightarrow G(U,T).$$

 We proved the proposition for all three types of equations (\ref{eq:st7}),
 (\ref{eq:st8}), (\ref{eq:st9}), as required.
 \end{proof}

\begin{prop}
\label{prop:hom-lambda-monic} Let $S = 1$ be a regular quadratic
equation {\rm (\ref{eq:st2})}
 and $$\beta:G_{R(S)} \rightarrow G$$ a
solution of $S= 1$ in $G$ in a general position. Then  the
homomorphism $$\lambda_\beta : G_{R(S)} \rightarrow G(U,T)$$ is a
monomorphism.
\end{prop}
\begin{proof} In the proof of this proposition we use induction on
the atomic rank of the equation in the same way as in the proof of
Theorem 1 in \cite{KMNull}.

Since all the intermediate groups are also fully residually free
 by induction it suffices to prove the following:

1. $n=1$, $m=0$; prove that $\psi =\alpha _3$ is an embedding of
$G_S$ into $G_3$.

2. $n=2$, $m=0$;  prove that $\psi =\alpha _4$ is a monomorphism on
$H=\left<G,x_1,y_1\right>.$

3. $n=1$, $m=1$; prove that $\psi=\alpha _{3}\bar \alpha _0$ is a
monomorphism on $H=\left<G,z_1\right>.$

4. $n=0$, $m>2$; prove that $\theta _{2}\bar \theta _2$ is an
embedding of $G_S$ into $\bar H_2$.

Now we consider all these cases one by one.

 Case 1. Choose an
arbitrary nontrivial element $h\in G_S$. It can be written in the
form $$ h = g_1\ v_1(x_1,y_1)\ g_2\ v_2(x_1,y_1) \ g_3 \ldots
v_n(x_1,y_1)\ g_{n+1}, $$ where $1 \neq v_i(x_1,y_1) \in
F(x_1,y_1)$ are words in $x_1, y_1$, not belonging to the subgroup
$\langle [x_1,y_1]\rangle ,$ and $1 \neq g_i \in G, g_i\not\in
\langle [a,b]\rangle$ (with the exception of $g_1$ and $g_{n+1}$,
they could be trivial). Then \beq \label{66} h^\psi = g_1\
v_1(t_3t_1a,t_2b)\ g_2\ v_2(t_3t_1a,t_2b) \ g_3 \cdots
v_n(t_3t_1a,t_2b)\ g_{n+1}. \eeq
 The group $G(U,T)$ is obtained from $G$ by three
HNN-extensions (extensions of centralizers), so every element in
$G(U,T)$ can be rewritten to its reduced form by making finitely
many pinches. It is easy to see that the leftmost occurrence of
either $t_3$ or $t_1$ in the product (\ref{66}) occurs in the
reduced form of $h^\psi$ uncancelled.

Case 2. $x_1\rightarrow t_4^{-1}t_2a_1,\ y_1\rightarrow
t_4^{-1}t_3t_1b_1t_4,\ x_2\rightarrow t_4^{-1}a_2t_4, \
y_2\rightarrow t_4^{-1}b_2.$ Choose an arbitrary nontrivial
element $h\in H=G*F(x_1,y_1)$. It can be written in the form $$ h
= g_1\ v_1(x_1,y_1)\ g_2\ v_2(x_1,y_1) \ g_3 \ldots v_n(x_1,y_1)\
g_{n+1}, $$ where $1 \neq v_i(x_1,y_1) \in F(x_1,y_1)$ are words
in $x_1, y_1$, and $1 \neq g_i \in G$ (with the exception of
$g_1$ and $g_{n+1}$, they could be trivial). Then \beq \label{7}
h^\psi = g_1\ v_1(t_4^{-1}t_2a,(t_3t_1b)^{t_4})\ g_2\
v_2(t_4^{-1}t_2a,(t_3t_1b)^{t_4}) \ g_3 \cdots
v_n(t_4^{-1}t_2a,(t_3t_1b)^{t_4})\ g_{n+1}. \eeq
 The group $G(U,T)$ is obtained from $G$ by four
HNN-extensions (extensions of centralizers), so every element in
$G(U,T)$ can be rewritten to its reduced form by making finitely
many pinches. It is easy to see that the leftmost occurrence of
either $t_4$ or $t_1$ in the product (\ref{7}) occurs in the
reduced form of $h^\psi$ uncancelled.

Case 3. We have an equation $c^z[x,y]=c[a,b]$, $z\rightarrow
zr_1\bar r_1,\ x\rightarrow (t_2a^{r_1})^{\bar r _1},\
y\rightarrow \bar r_1^{-1}t_3t_1r_1^{-1}b,$ and $[r_1,ca^{-1}]=1,\
[\bar r_1, (c^{r_1}a^{-r_1}t_2^{-1})]=1.$ Here we can always
suppose, that $[c,a]\not =1$, by changing a solution, hence
$[r_1,\bar r_1]\not =1.$ The proof for this case is a repetition
of the proof of Proposition 11 in \cite{KMNull}.

Case 4. We will consider the case when $m=3$; the general case can
be considered similarly. We have an equation
$c_1^{z_1}c_2^{z_2}c_3^{z_3}=c_1c_2c_3$, and can suppose
$[c_i,c_{i+1}]\neq 1.$

We will prove that $\psi=\theta _2\bar\theta _1$ is an embedding.
The images of $z_1,z_2,z_3$ under $\theta _2\bar\theta _1$ are the
following:
$$z_1\rightarrow c_1r_1\bar r_1,\ z_2\rightarrow c_2r_1r_2\bar
r_1,\ z_3\rightarrow c_3r_2,$$ where
$$
[r_1,c_1c_2]=1,\ [r_2, c_2^{r_1}c_3]=1,\ [\bar r_1,
c_1^{r_1}c_2^{r_1r_2}]=1.$$

Let $w$ be a reduced word in  $G*F(z_i, i=1,2,3),$ which does not
have subwords $c_1^{z_1}$. We will prove that if $w^{\psi}=1$ in
$\bar H_1$, then $w\in N,$ where  $N$ is  the normal closure of
the element $c_1^{z_1}c_2^{z_2}c_3^{z_3}c_3^{-1}c_2^{-1}c_1^{-1}.$
We use induction on the number of occurrences of $z_1^{\pm 1}$ in
$w$. The induction basis is obvious, because homomorphism $\psi$
is injective on the subgroup $<F,z_2,z_3>.$

Notice, that the homomorphism $\psi$ is also injective on the
subgroup $K=<z_1z_2^{-1}, z_3, F>.$

Consider $\bar H_1$ as an HNN-extension by letter $\bar r_1$.
Suppose $w^{\psi}=1$ in $\bar H_1$. Letter $\bar r_1$ can
disappear in two cases: 1) $w\in KN,$  2) there is a pinch between
$\bar r_1^{-1}$ and $\bar r_1$ (or between $\bar r_1$ and $\bar
r_1^{-1}$) in $w^{\psi}.$ This pinch corresponds to some element
$z_{1,2}^{-1}uz^{\prime}_{1,2}$ (or $z_{1,2}u(z^{\prime}
_{1,2})^{-1}$), where $z_{1,2}, z^{\prime}_{1,2}\in\{z_1,z_2\}.$

In the first case $w^{\psi}\neq 1$, because $w\in K$ and $w\not\in
N$.

In the second case, if the pinch happens in
$(z_{1,2}u(z^{\prime}_{1,2})^{-1})^{\psi}$, then
$z_{1,2}u(z^{\prime}_{1,2})^{-1}\in KN,$ therefore it has to be at
least one pinch that corresponds to
$(z_{1,2}^{-1}uz^{\prime}_{1,2})^{\psi}$. We can suppose, up to a
cyclic shift of $w$, that $z_{1,2}^{-1}$ is the first letter, $w$
does not end with some $z_{1,2}^{\prime\prime}$, and $w$ cannot be
represented as
$z_{1,2}^{-1}uz^{\prime}_{1,2}v_1z^{\prime\prime}_{1,2}v_2,$ such
that $z_{1,2}^{\prime}v_1\in KN.$ A pinch can only happen if
$z_{1,2}^{-1}uz^{\prime}_{1,2}\in <c_1^{z_1}c_2^{z_2}>$.
Therefore, either $z_{1,2}=z_1,$ or $z_{1,2}^{\prime}=z_1$, and
one can replace $c_1^{z_1}$ by  $c_1c_2c_3c_3^{-z_3}c_2^{-z_2}$,
therefore replace $w$ by $w_1$ such that $w=uw_1$, where $u$ is in
the normal closure of the element
$c_1^{z_1}c_2^{z_2}c_3^{z_3}c_3^{-1}c_2^{-1}c_1^{-1},$ and apply
induction.
  \end{proof}

The embedding $\lambda_\beta : G_S \rightarrow G(U,T)$ allows one
to construct effectively discriminating sets of solutions in $G$
of the equation $S = 1$. Indeed, by the construction above  the
group $G(U,T)$ is union of the following  chain of length $2K =
2K(m,n)$ of extension of centralizers:
 $$
  G = H_0 \leqslant H_1 \ldots \leqslant H_{m-1} \leqslant G_0 \leqslant G_1 \leqslant \ldots
  \leqslant G_{4n-1} = $$
   $$ = \bar H_0 \leqslant \bar H_1 \leqslant \ldots \leqslant  \bar H_{m-1 }=  \bar G_0 \leqslant \ldots \leqslant \bar G_{4n-1}
 = G(U,T).$$
 Now, every $2K$-tuple $p \in {\mathbb N}^{2K}$ determines  a
$G$-homomorphism
 $$\xi_p: G(U,T) \rightarrow G.$$
  Namely, if $Z_i$ is the $i$-th term of the chain above then $Z_i$ is an extension of
  the centralizer of some element $g_i \in Z_{i-1}$ by a stable letter $t_i$. The $G$-homomorphism
$\xi_p$  is defined as composition
 $$\xi_p= \psi_1 \circ \ldots \circ \psi_K$$
 of homomorphisms
$\psi_i:Z_i \rightarrow Z_{i-1}$ which are identical on $Z_{i-1}$
 and such that $t_i^{\psi_i} = g_i^{p_i}$, where $p_i$ is the $i$-th component of
$p$.

 It follows
 (see   \cite[Section~1.4]{JSJ})
  that for every unbounded set of tuples
$P \subset {\mathbb N}^{2K}$  the set of homomorphisms
 $$\Xi_P =  \{\xi_p \mid p \in P\}$$
 $G$-discriminates $G(U,T)$ into $G$. Therefore, (since $\lambda_\beta$ is monic),
 the family of $G$-homomorphisms
 $$ \Xi_{P,\beta} = \{\lambda_\beta  \xi_p \mid \xi_p \in \Xi_P
 \}$$
$G$-discriminates $G_S$ into $G$.

 One can give another description of the set $\Xi_{P,\beta}$ in terms of the basic
 automorphisms from the  basic sequence $\Gamma$. Observe first that
 $$ \lambda_\beta  \xi_p = \phi_{2K,p} \beta ,$$
 therefore
 $$\Xi_{P,\beta} = \{\phi_{2K,p} \beta \mid \ p \in P \}.$$

  We summarize the discussion above as follows.

\begin{theorem} \label{cy2} Let $G$ be a finitely generated
fully residually free group, $S = 1$  a regular standard quadratic
orientable equation,  and   $\Gamma$   its basic sequence of
automorphisms. Then for  any  solution $\beta:G_S \rightarrow G $
in general position, any positive integer $J \geq 2$,  and  any
unbounded set $P \subset {\mathbb N}^{JK}$  the set of
$G$-homomorphisms $\Xi_{P,\beta}$ $G$-discriminates  $G_{R(S)}$
into $G$. Moreover, for any fixed tuple $p'\in {\mathbb N}^{tK}$
the family
 $$\Xi_{P,\beta,p'} = \{\phi _{tK,p'}\theta \mid \theta\in\Xi_{P,\beta}\}$$
$G$-discriminates  $G_{R(S)}$ into $G$.
\end{theorem}

For tuples $f = (f_1, \ldots,f_k)$ and $g = (g_1, \ldots,g_m)$
 denote the tuple
  $$fg = (f_1, \ldots,f_k,g_1, \ldots,g_m).$$
Similarly, for a set of tuples $P$ put
 $$f Pg =\{fpg \mid p \in P\}.$$

\begin{cy} \label{cy:cy2}
Let $G$ be a finitely generated fully residually free group, $S =
1$ a regular standard quadratic orientable equation,    $\Gamma$
the basic sequence of automorphisms of $S$, and $\beta:G_S
\rightarrow G $ a solution of $S=1$ in general position. Suppose
$P \subseteq \mathbb{N}^{2K}$ is unbounded set, and $f \in
\mathbb{N}^{Ks}$, $g \in \mathbb{N}^{Kr}$ for some $r,s \in
\mathbb{N}$. Then there exists a number $N$ such that if $f$ is
$N$-large and $s\geq 2$ then the family
 $$\Phi_{P,\beta,f,g} = \{\phi _{K(r+s+2),q} \beta  \mid q \in fPg\}$$
$G$-discriminates  $G_{R(S)}$ into $G$.
\end{cy}
 \begin{proof} By Theorem \ref{cy2} it suffices to show that if
 $f$ is $N$-large for some $N$ then
 $\beta_{f} = \phi_{2K,f}\beta$ is a solution of $S = 1$ in general
 position, i.e., the images of some particular finitely many  non-commuting
 elements from $G_{R(S)}$ do not commute in $G$. It has been shown above  that the set
 of solutions $\{\phi_{2K,h}\beta \mid h \in \mathbb{N}^{2K}\}$ is
 a discriminating set for $G_{R(S)}$. Moreover, for any finite set
 $M$ of non-trivial elements from $G_{R(S)}$ there exists a number $N$
 such that for any $N$-large tuple $h \in \mathbb{N}^{2K}$ the
 solution $\phi_{2K,h}\beta$ discriminates all elements from $M$ into
 $G$. Hence the result.
\end{proof}

\section{Small cancellation solutions of standard orientable
equations} \label{subsec:7.2-b}

Let $S(X) = 1$ be  a standard regular orientable quadratic
equation over $F$ written in the form (\ref{6-general}):
 $$\prod_{i=1}^{m}z_i^{-1}c_iz_i\prod_{i=1}^{n}[x_i,y_i]  =
\prod_{i=1}^{m}e_i^{-1}c_ie_i \prod_{i=1}^{n}[a_i,b_i].$$ In this
section we construct solutions in $F$ of $S(X) = 1$ which satisfy
some small cancellation conditions.

\begin{df}
\label{de:smallcan} Let $S = 1$ be a standard regular orientable
quadratic equation written in the form (\ref{6-general}). We say
that a solution $\beta: F_S \rightarrow F$ of $S= 1$  satisfies
the small cancellation condition $(1/\lambda)$ with respect to the
set $\bar {\mathcal W}_{\Gamma ,L}$ if the following conditions
are satisfied:
 \begin{enumerate}
  \item [1)] $\beta$ is in general position;
  \item [2)] for any 2-letter word $uv \in {\mathcal
W}_{\Gamma ,L}$ (in the alphabet $Y$) cancellation in the word
$u^\beta v^\beta$ does not exceed $(1/\lambda) \min\{|u^\beta|,
|v^\beta|\}$ (we assume here and below that $u^\beta, v^\beta$ are
given by their reduced forms in $F$);
 \item [3)]   cancellation in a word $u^{\beta} v^{\beta}$
 does not exceed $(1/\lambda) \min\{|u^{\beta}|,|v^\beta|\}$
 provided $u, v$ satisfy one of the conditions below:
  \begin{itemize}
   \item [a)] $u=z_i, v=(z_{i-1}^{-1}c_{i-1}^{-1}z_{i-1})$,
    \item [b)] $u=c_i, v=z_i$,
    \item [c)]  $u=v=c_i$.
\end{itemize}
 \end{enumerate}
\end{df}

\begin{notation} For a homomorphism  $\beta : F[X] \rightarrow F$ by
$C_\beta$ we denote the set of all  elements that cancel in $u^\beta
v^\beta$  where $u, v$ are as in 2), 3) from Definition
\ref{de:smallcan}.
\end{notation}

\begin{lm}
\label{lm:32} Let $u,v$ be cyclically reduced elements of $G\ast
H$ such that $|u|,|v|\geqslant 2$. If for some $m, n > 1$ elements
$u^{m}$ and $v^{n}$ have a common initial segment of length
$|u|+|v|$, then $u$ and $v$ are both powers of the same element
$w\in G \ast H$. In particular, if both $u$ and $v$ are not proper
powers then $u = v$.
\end{lm}
 \begin{proof}  The same argument as in the case of free groups.

\begin{cy} If $u,v\in F,\ [u,v]\neq 1,$ then for any
$\lambda \geqslant 0$ there exist $m_0, n_0$ such that
for any $m \geqslant m_0, n \geqslant n_0$ cancellation
between $u^m$ and $v^n$ is less than
$\frac{1}{\lambda} \max\{|u^m|,|v^n|\}.$\end{cy}

\begin{lm}
\label{le:7.1.beta} Let $S(X) = 1$ be  a standard regular
orientable quadratic equation written in the form {\rm
(\ref{6-general})}:

 $$\prod_{i=1}^{m}z_i^{-1}c_iz_i\prod_{i=1}^{n}[x_i,y_i]  =
\prod_{i=1}^{m}e_i^{-1}c_ie_i \prod_{i=1}^{n}[a_i,b_i], \ \
n\geqslant 1, $$ where  all  $c_i$ are cyclically reduced,
 and
 $$\beta_1: x_i \rightarrow a_i, y_i \rightarrow b_i, z_i \rightarrow e_i $$
  a solution of $S = 1$ in $F$ in
 general position. Then for any $\lambda \in \mathbb{N}$ there are positive
integers  $m_i, n_i, k_i, q_j$ and a tuple $p=(p_1,\ldots p_{m})$
such that  the map $\beta: F[X] \rightarrow F$
 defined by

$$x_1^{\beta}=(\tilde b_1^{n_1}\tilde a_1)^{[\tilde a_1,\tilde b_1]^{m_1}}, \ \ y_1^{\beta}=
((\tilde b_1^{n_1}\tilde a_1)^{k_1}\tilde b_1)^{[\tilde a_1,\tilde
b_1]^{m_1}}, \ \ where \ \tilde a_1=x_1^{\phi _{m}\beta_1},\
\tilde b_1=y_1^{\phi
 _{m}\beta_1}$$

  $$x_i^{\beta}=(b_i^{n_i}a_i)^{[a_i,b_i]^{m_i}}, \ \ y_i^{\beta}=
((b_i^{n_i}a_i)^{k_i}b_i)^{[a_i,b_i]^{m_i}}, \ \ i=2,\ldots n,$$

 $$z_i^{\beta}=c_i^{q_i}z_i^{\phi _{m}\beta _1}, \ \ i=1,\ldots m,$$
 is a  solution of $S= 1$ satisfying the small cancellation condition
 $(1/\lambda)$ with respect to $\bar {\mathcal W}_{\Gamma ,L}$.
Moreover, one can choose the solution $\beta_1$ such that
 if $u=c_i^{z_i}$ or $u=x_j^{-1}$ and $v=c_1^{z_1}$, then
the cancellation between $u^{\beta}$ and $v^{\beta}$ is less than
$(1/\lambda)\min \{|u|,|v|\}.$
 \end{lm}
\begin{proof}
 The solution $$x_i \rightarrow a_i, y_i \rightarrow b_i, z_i
\rightarrow e_i $$
  $i = 1,  \ldots, n, j = 1, \ldots, m$ is in general position,
therefore the neighboring items in the sequence
  $$c_1^{e_1}, \ldots, c_m^{e_m}, [a_1,b_1],\dots ,[a_n,b_n]$$
   do not commute.

There is a homomorphism $\theta _{\beta _1}:F_S\rightarrow \bar
F=F(\bar U,\bar T)$ into the group $\bar F$ obtained from $F$ by a
series of extensions of centralizers, such that $\beta=\theta
_{\beta _1}\psi _p$, where $\psi _p:\bar F\rightarrow F$. This
homomorphism $\theta _{\beta _1}$ is a  monomorphism on $F\ast
F(z_1,\dots, z_m)$ (this follows from the proof of Theorem 4 in
\cite{KMNull}, where the same sequence of extensions of centralizers
is constructed).

 The set  of
solutions $\psi_p$ for different tuples $p$ and numbers
$m_i,n_i,k_i, q_j$ is a discriminating family for $\bar F.$ We
just have to show that the small cancellation condition for
$\beta$ is equivalent to a finite number of inequalities in the
group $\bar F$.

 We have  $z_i^\beta =c_i^{q_i}z_i^{\phi _{m}\beta
_1}$ such that $\beta _1(z_i)=e_i$, and $p=(p_1,\dots ,p_{m})$ is
a large tuple.  Denote $\bar A_j= A_j^{\beta _1},\ j=1,\dots ,m.$
Then it follows from Lemma \ref{le:7.1.zforms}
 that

\medskip
$z_i^{\beta}=c_i^{q_i+1}e_i\bar
A_{i-1}^{p_{i-1}}c_{i+1}^{e_{i+1}}\bar A_{i}^{p_{i}-1}$, where
$i=2,\dots ,m-1$

\medskip
$z_m^{\beta}=c_m^{q_m+1}e_m\bar A_{m-1}^{p_{m-1}}a_1^{-1}\bar
A_{m}^{p_{m}-1}$,

where
\bea \bar A_1 &=&c_1^{e_1}c_2^{e_2},\\
\bar A_2 &=& \bar A_1(p_1)= \bar
A_1^{-p_1}c_2^{e_2}\bar A_1^{p_1}c_3^{e_3},\\ && \vdots\\
\bar A_{i} &=& \bar A_{i-1}^{-p_{i-1}}c_i^{e_i}\bar
A_{i-1}^{p_{i-1}}c_{i+1}^{e_{i+1}},\ i=2,\dots ,m-1,\\
\bar A_{m} &=& \bar A_{m-1}^{-p_{m-1}}c^{e_m}\bar
A_{m-1}^{p_{m-1}}a_1^{-1}.
\eea

 One can choose $p$ such that
$[\bar A_{i},\bar A_{i+1}]\neq 1, [\bar
A_{i-1},c_{i+1}^{e_{i+1}}]\neq 1, [\bar A_{i-1},c_{i}^{e_{i}}]\neq
1$ and $[\bar A_{m},[a_1,b_1]]\neq 1$, because their pre-images do
not commute in $\bar F$. We need the second and third inequality
here to make sure that $\bar A_{i}$ does not end with a power of
$\bar A_{i-1}$. Alternatively, one can prove by induction on $i$
that $p$ can be chosen to satisfy these inequalities.
 Then $c_i^{z_i^{\beta}}$
and $c_{i+1}^{z_{i+1}^{\beta}}$ have small cancellation, and
$c_m^{z_m^{\beta}}$ has small cancellation with $x_1^{\pm\beta},
y_1^{\pm\beta}$.

    Let
   $$ x_i^{\beta}=(b_i^{n_i}a_i)^{[a_i,b_i]^{m_i}}, \ \
   y_i^{\beta}=((b_i^{n_i}a_i)^{k_i}b_i)^{[a_i,b_i]^{m_i}}, \ \ i=2,\dots ,n
$$
 for some positive integers  $m_i, n_i, k_i, s_j$ which values we
will specify in a due course.  Let $uv \in\bar{\mathcal
W}_{\Gamma}$. There are several cases to consider.

1) $uv = x_ix_i$. Then
 $$u^\beta v^\beta = (b_i^{n_i}a_i)^{[a_i,b_i]^{m_i}}
 (b_i^{n_i}a_i)^{[a_i,b_i]^{m_i}}.$$
Observe that the cancellation  between  $(b_i^{n_i}a_i)$ and
$(b_i^{n_i}a_i)$ is not more then $|a_i|$. Hence the cancellation
in $u^\beta v^\beta$ is not more then $|[a_i,b_i]^{m_i}| + |a_i|$.
We chose $n_i \gg m_i$ such that

$$|[a_i,b_i]^{m_i}| + |a_i| < (1/\lambda) |(b_i^{n_i}a_i)^{[a_i,b_i]^{m_i}}|$$
which is obviously possible. Similar arguments prove the cases $uv
= x_iy_i$ and $uv = y_ix_i.$

2) In all other cases the cancellation in
 $u^\beta v^\beta$ does not exceed the cancellation between $[a_i,b_i]^{m_i}$ and
 $[a_{i+1},b_{i+1}]^{m_{i+1}}$, hence by Lemma  \ref{lm:32} it is not greater than
 $|[a_i,b_i]| +  |[a_{i+1},b_{i+1}]|.$

Let $u=z_i^{\beta}, v=c_{i-1}^{-z_{i-1}^{\beta}}.$ The
cancellation is the same as between $\bar A_{2i}^{p_{2i}}$ and
$\bar A_{i-1}^{-p_{i-1}}$ and, therefore, small.

Since $c_i$ is cyclically reduced, there is no cancellation
between $c_i$ and $z_{i}^{\beta}$.

The first statement of the lemma is proved.

We now will prove the second statement of the lemma.
 We can  choose the  initial solution $e_1,\dots, e_m, a_1,b_1,\dots ,a_n,b_n$ so
   that $[c_1^{e_1}c_2^{e_2}, c_3^{e_3}\ldots c_i^{e_i}]\neq 1$ ( $i\geqslant    3$),
   $[c_1^{e_1}c_2^{e_2},[a_i,b_i]]\neq 1, (i=2,\dots ,n)$ and
   $[c_1^{e_1}c_2^{e_2}, b_1^{-1}a_1^{-1}b_1]\neq 1.$
Indeed, the equations $[c_1^{z_1}c_2^{z_2}, c_3^{z_3}\ldots
   c_i^{z_i}]=1$, $[c_1^{z_1}c_2^{z_2},[x_i,y_i]]=1, (i=2,\dots ,n)$ and
   $[c_1^{z_1}c_2^{z_2}, y_1^{-1}x_1^{-1}y_1]=1$ are not  consequences of the equation $S=1$, and,
therefore, there is a solution of $S(X)=1$ which does not satisfy
 any of  these equations.

To show that $u=c_i^{z_i^{\beta}}$  and $v=c_1^{z_1^{\beta}}$,
have small cancellation, we have to show  that $p$ can be chosen
so that $[\bar A_1, \bar A_{i}]\neq 1$ (which is obvious, because
the pre-images in $\bar G$ do not commute), and that $\bar
A_{i}^{-1}$ does not begin with a power of $\bar A_1$. The period
$\bar A_{i}^{-1}$ has form ($c_{i+1}^{-z_{i+1}}\ldots
c_3^{-z_3}\bar A_1^{-p_2}\ldots ).$ It begins with a power of
$\bar A_1$ if and only if $[\bar A_1, c_3^{e_3}\ldots
c_i^{e_i}]=1$, but this equality does not hold.

Similarly one can show, that the cancellation between
 $u=x_j^{-\beta}$ and $v=c_1^{z_1^{\beta}}$ is small. \end{proof}

\begin{lm} \label{sc2}  Let $S(X) = 1$ be  a standard regular
orientable quadratic equation of the type {\rm  (\ref{eq:st9})}
$$\prod _{i=1}^{m}z_i^{-1}c_iz_i=c_1^{e_1}\ldots c_m^{e_m}=d,$$
 where  all  $c_i$ are cyclically reduced,
 and
 $$\beta_1:  z_i \rightarrow e_i $$
  a solution of $S = 1$ in $F$ in
 general position. Then for any $\lambda \in \mathbb{N}$ there is a  positive
integer $s$   and a tuple $p=(p_1,\ldots p_{K})$ such that the map
$\beta: F[X] \rightarrow F$
 defined by
$$z_i^{\beta}=c_i^{q_i}z_i^{\phi _K\beta _1}d^s, $$
  is a  solution of $S= 1$ satisfying the small cancellation condition
 $(1/\lambda)$ with respect to $\bar {\mathcal W}_{\Gamma ,L}$ with
 one exception when $u = d$ and $v =c_{m-1}^{-z_{m-1}}$ (in this
 case $d$ cancels out in $v^\beta$). Notice, however, that such
 word $uv$ occurs only in the product $wuv$ with  $w = c_2^{z_2}$,
 in which case cancellation
between $w^{\beta}$ and $dv^{\beta}$ is less than
$\min\{|w^{\beta}|,|dv^{\beta}|\}.$
\end{lm}
\begin{proof}  Solution $\beta$ is chosen   the same way as in the
previous lemma (except for the multiplication by  $d^s$) on the
elements $z_i,\ i\neq m$. We do not take $s$ very large, we just
need it to avoid cancellation between $z_2^{\beta}$ and $d$.
Therefore the cancellation between $c_i^{z_i^{\beta}}$ and
$c_{i+1}^{\pm z_{i+1}^{\beta}}$ is small for $i<m-1$. Similarly,
for $u=c_2^{z_2},\  v=d,\ w=c_{m-1}^{-z_{m-1}},$ we can make the
cancellation between $u^{\beta}$ and $dw^{\beta}$  less than
$\min\{|u^{\beta}|,|dw^{\beta}|\}.$ \end{proof}

\begin{lm}
\label{le:small-cancellation-beta-double} Let $U, V \in
\bar{\mathcal W}_{\Gamma, L}$ such that $UV = U \circ V$ and $UV
\in \bar{\mathcal W}_{\Gamma, L}$.

1. Let $n\neq 0.$ If $u$ is the last letter of $U$ and $v$ is the
first letter of $V$ then cancellation between  $U^\beta$ and
$V^\beta$ is equal to the cancellation between  $u^\beta$ and
$v^\beta$.

2. Let $n=0.$ If $u_1u_2$ are the last two letters of $U$ and
$v_1,v_2$ are the first two letters of $V$ then cancellation
between $U^\beta$ and $V^\beta$ is equal to the cancellation
between $(u_1u_2)^\beta$ and $(v_1v_2)^\beta$.

  \end{lm}
\begin{proof} Since  $\beta$ has the small
  cancellation property with respect to
  $\bar {\mathcal W}_{\Gamma, L},$
  this implies that the cancellation in $U^\beta V^{\beta}$ is equal to the cancellation in
  $u^\beta v^{\beta}$, which is equal to some element in
  $C_\beta$.
This proves the lemma.
 \end{proof}

Let $w\in \bar{\mathcal W}_{\Gamma ,L}, W=w^{\phi _j}.$ We start
with the canonical $N$-large $A$-representation of $W$:
  \begin{equation}
\label{eq:Adecomp3} W = B_1 \circ A^{q_1} \circ \cdots \circ B_k
\circ A^{q_k} \circ B_{k+1}
\end{equation}
where $|q_i| \geqslant N$ and $\max_j(B_i) \leqslant r$.

Since the occurrences $A^{q_i}$ above are stable we have $$B_1 =
{\bar B}_1 \circ A^{sgn(q_1)}, \ \ B_i = A^{sgn(q_{i-1})} \circ
{\bar B}_i \circ A^{sgn(q_i)} \ (2 \leqslant i \leqslant k), \ \
B_{k+1} = A^{sgn(q_k)} \circ {\bar B}_{k+1}.$$ Denote
$A^{\beta}=c^{-1}A^{\prime} c,$ where $A^{\prime}$ is cyclically
reduced, and $c\in C_{\beta}.$ Then
$$B_1^\beta = {\bar B}_1^\beta c^{-1} (A^\prime)^{sgn(q_1)}c, \ \
B_i^\beta = c^{-1}(A^\prime)^{sgn(q_{i-1})} c {\bar B}_i^\beta
c^{-1} (A^\prime)^{sgn(q_i)}c, $$

$$ \ B_{k+1}^\beta  =c^{-1} (A^\prime)^{sgn(q_k)} c{\bar
B}_{k+1}^\beta.$$ By Lemma \ref{le:small-cancellation-beta-double}
we can assume that the cancellation in the words above is small,
i.e., it does not exceed a fixed number $\sigma$ which is the
maximum length of words from $C_\beta$. To get an $N$-large
canonical $A^\prime$-decomposition of $W^\beta$ one has to take
into account stable occurrences of $A^\prime$. To this end, put
$\varepsilon_i = 0$  if ${A^\prime}^{sgn(q_i)}$ occurs in the
reduced form of ${\bar B}_i^\beta c^{-1}(A^\prime)^{sgn(q_i)}$ as
written (the cancellation does not touch it), and put
$\varepsilon_i = sgn(q_i)$
 otherwise. Similarly, put $\delta_i = 0$ if
${A^\prime}^{sgn(q_i)}$ occurs in the reduced form of
$(A^\prime)^{sgn(q_{i})} c {\bar B}_{i+1}^\beta$ as written, and
put $\delta_i = sgn(q_{i})$ otherwise.

Now one can rewrite $W^\beta$ in the following form

\begin{equation}
\label{eq:7.3.can}
 W^\beta  = E_1 \circ
   (A^\prime)^{q_1 - \varepsilon_1 -\delta_1}
   \circ
  E_2
  \circ (A^\prime)^{q_2 - \varepsilon_2 -\delta_2} \circ \cdots
 \circ (A^{\prime})^{q_k -\varepsilon_k -\delta_k}  \circ E_{k+1},
 \end{equation}
 where $E_1=(B_1^\beta c^{-1}(A^\prime)^{\varepsilon_1}),\
 E_2=((A^\prime)^{\delta_1}cB_2^\beta
 c^{-1}(A^\prime)^{\varepsilon_2}),\ E_{k+1}=((A^\prime)^{\delta_{k}}
   cB_{k+1}^\beta).$

   Observe, that $d_i$ and $\varepsilon_i, \delta_i$ can be effectively computed from
 $W$ and $\beta$. It follows that one can effectively rewrite  $W^\beta$
 in the form  (\ref{eq:7.3.can}) and the form is unique.

The decomposition  (\ref{eq:7.3.can}) of  $W^\beta$ induces a
corresponding  $A^*$-decomposition of  $W$. Namely, if the
canonical $N$-large $A^*$-decomposition of $W$ has the form:
$$
D_1 (A^*)^{q_1} D_2 \cdots D_k(A^*)^{q_k} D_{k+1}$$ then the
induced one has the form: $W=$
\begin{equation}
\label{eq:induced} (D_1A^{* \varepsilon_1}) A^{* q_1 -
\varepsilon_1 - \delta_1} (A^{*  \delta_1} D_2 A^{*\varepsilon_2})
\cdots (A^{*  \delta_{k-1}} D_k A^{*\varepsilon_k}) A^{* q_k -
\varepsilon_k - \delta_k} (A^{*\delta_k}D_{k+1}).\end{equation}

We call this decomposition the {\it induced} $A^*$-decomposition
of $W$ with respect to $\beta$ and write it in the form:
\begin{equation}
\label{eq:induced*} W = D_1^* {(A^*)}^{q_1^*} D_2^*
 \cdots D_k^* {(A^*)}^{q_k^*}D_{k+1}^*,\end{equation}
where $D_i^* = {(A^*)}^{\delta_{i-1}} D_i
{(A^*)}^{\varepsilon_i}$, $q_i^* = q_i - \varepsilon_i -
\delta_i$, and, for uniformity, $\delta_1 = 0$ and
$\varepsilon_{k+1} = 0$.

\begin{lm}
\label{le:D-star-beta} For given positive integers $j$,  $M$, $N$
there is a constant  $C = C(j,M,N) >0$  such that  if $p_{t+1}-p_t
> C$ for every $t = 1, \ldots, j-1,$ and a word $W\in\bar{\mathcal
W}_{\Gamma ,L}$ has a canonical
 $N$-large $A^*$-decomposition (\ref{eq:induced*}), then this decomposition satisfies the following conditions:
\begin{equation} \label{eq:E-D} (D_1^*)^\beta = E_1
\circ_{\theta} (c R^\beta), \ \ (D_i^*)^\beta = (R^{-\beta}c^{-1})
\circ_{\theta} E_i \circ_{\theta} (cR^\beta), \ \
(D_{k+1}^*)^\beta = (R^{-\beta}c^{-1}) \circ_{\theta} E_{k+1},
\end{equation}
where $\theta  < |A| - M$. Moreover, this constant $C$ can be
found effectively.
\end{lm}
\begin{proof} Applying homomorphism $\beta$ to the reduced
 $A^*$-decomposition of $W_\sigma$ (\ref{eq:induced*}) we can see
 that

 $$
 W_{\sigma}^\beta = \left ((D_1^*)^\beta R^\beta c \right )
 {(A^\prime)}^{q_1^*} \left (c R^{\beta}(D_2^*)^\beta
 R^{-\beta}c^{-1}\right ){(A^\prime)}^{q_2^*} \ldots $$ $$ \left (cR^\beta
 (D_k^*)^\beta R^{-\beta}c^{-1}\right ){(A^\prime)}^{q_k^*} \left
 (cR^\beta (D_{k+1}^*)^\beta\right ). $$
 Observe that this decomposition has the same powers of $A^\prime$
 as the canonical $N$-large $A^\prime$-decomposition
 (\ref{eq:7.3.can}). From the uniqueness of such decompositions we
 deduce that
$$E_1 = (D_1^*)^\beta R^\beta c, \ \ E_i = c R^{\beta}(D_i^*)^\beta
 R^{-\beta}c^{-1}, \ \ E_{k+1} = cR^\beta (D_{k+1}^*)^\beta$$
 Rewriting these equalities one can get
 $$
(D_1^*)^\beta = E_1 \circ_{\theta} (c R^\beta), \ \ (D_i^*)^\beta
= (R^{-\beta}c^{-1}) \circ_{\theta} E_i \circ_{\theta} (cR^\beta),
\ \ (D_{k+1}^*)^\beta = (R^{-\beta}c^{-1}) \circ_{\theta}
E_{k+1}$$ and $\theta \ll |A|$.  Indeed,  in the decomposition
(\ref{eq:7.3.can}) every occurrence $(A^\prime)^{q_i -
\varepsilon_i -\delta_i}$ is stable hence $E_i$ starts (ends) on
$A^\prime$. The $N$-large rank of $R$ is at most $rank _N(A)$, and
$\beta$ has small cancellation. Taking $p_{j+1}\gg p_j$ we may
assume that $|A^\prime| \gg |c|, |R^\beta|$.
\end{proof}

 Notice, that one can
effectively write down the induced $A^*$-decomposition of $W$ with
respect to $\beta$.

We summarize the discussion above in the following statement.

\begin{lm}\label{Claim1.b} For given positive integers $j$,  $N$
there is a constant  $C =  C(j,N)$  such that if $p_{t+1}-p_t > C$,
for  every $t = 1, \ldots, j-1,$ then for any   $W\in \bar{\mathcal
W}_{\Gamma , L}$    the following conditions are  equivalent:

\begin{enumerate}
\item  Decomposition {\rm (\ref{eq:Adecomp3})} is the canonical
(the canonical $N$-large) $A$-decomposition of $W$,

\item Decomposition {\rm (\ref{eq:7.3.can})} is the canonical (the
canonical $N$-large) $A^\prime$-decomposition of $W^\beta$,

\item  Decomposition {\rm (\ref{eq:induced})} is the canonical
(the canonical $N$-large) $A^*$-decomposition of $W.$

\end{enumerate}
\end{lm}

\section{Implicit function theorem for quadratic equations}\label{se:7.3}
In this section we prove Theorem A for orientable quadratic
equations over a free group $F = F(A)$. Namely, we prove the
following statement.

\medskip
 {\it Let $S(X,A)=1$ be a regular standard orientable
quadratic equation over $F$. Then every equation $T(X,Y,A) = 1$
compatible with $S(X,A) = 1$ admits an effective complete
$S$-lift.}

 \medskip
{\bf A special discriminating set of solutions
 ${\mathcal L}$ and the corresponding   cut equation $\Pi$.}

 \medskip
Below we continue to use notations from the previous sections. Fix
a  solution $\beta$ of $S(X,A) = 1$ which satisfies the
cancellation condition $(1/\lambda)$ (with $\lambda
> 10$) with respect to $\bar {\mathcal W}_{\Gamma}$.

 Put
$$x_i^{\beta}=\tilde a_i,y_i^{\beta}=\tilde b_i,z_i^{\beta}=\tilde c_i. $$
 Recall that $$\phi_{j,p} = \gamma_j^{p_j} \cdots \gamma_1^{p_1} = \stackrel{\leftarrow}{\Gamma}_j^p$$
 where $j \in {\mathbb N}$,  $\Gamma_j = (\gamma_1,
 \ldots, \gamma_j)$ is the initial subsequence of length $j$ of the
sequence
 $\Gamma^{(\infty)}$, and $p = (p_1, \ldots,p_j) \in {\mathbb N}^j$.
  Denote by $\psi_{j,p}$ the following solution of $S(X) = 1$:
  $$\psi_{j,p}=\phi_{j,p}\beta .$$ Sometimes we
 omit $p$ in $\phi_{j,p}, \psi_{j,p}$ and simply write
$\phi_j, \psi_j$.

Below we continue to use notation:
  $$A = A_j = A_j,  \ A^* =
A^*_j = A^*(\phi_j) = R_j^{-1} \circ A_j \circ R_j, \ d = d_j =
|R_j|.$$
 Recall that $R_j$ has rank $\leq j-K +2$ (Lemma
\ref{le:cyclicAA*}).  By $A^\prime$ we denote the cyclically
reduced form of $A^\beta$ (hence of $(A^*)^\beta$).
  Recall that $C_\beta$ is the finite  set of  all initial and terminal segments of
elements  in $(X^{\pm 1})^\beta$.

Let
 $$ \Phi = \{\phi_{j,p} \mid j
\in {\mathbb N}, p \in {\mathbb N}^j \}.$$
  For    an arbitrary  subset ${\mathcal L}$ of $\Phi$ denote
  $${\mathcal L}^\beta = \{\phi \beta \mid \phi \in {\mathcal L}\}.$$

Specifying step by step various subsets of $\Phi$ we will
eventually ensure a very particular choice of a set of solutions
of $S(X) = 1$ in $F$.

Let $K = K(m,n)$ and $J \in \mathbb{N}, J \geq 3,$  a sufficiently
large positive integer which will be specified precisely in  due
course.  Put $L = JK$ and define ${\mathcal P}_1 = \mathbb{N}^L$,
 $${\mathcal L_1} = \{\phi_{L,p} \mid  p \in {\mathcal P}_1 \}.$$
By Theorem \ref{cy2} the set ${\mathcal L_1}^\beta$  is a
discriminating set of solutions of $S(X) = 1$ in $F$. In fact, one
can replace the set ${\mathcal P}_1$ in the definition of
${\mathcal L_1}$ by any unbounded subset ${\mathcal P}_2 \subseteq
{\mathcal P}_1$, so that the new set is still discriminating. Now
we construct by induction a very particular unbounded subset
${\mathcal P}_2 \subseteq \mathbb{N}^L$.
 Let $a \in  \mathbb{N}$ be a natural number and $h: \mathbb{N} \times \mathbb{N} \rightarrow \mathbb{N}$ a  function. Define a tuple
  $$p^{(0)} = (p^{(0)}_1, \ldots, p^{(0)}_L)$$
where
 $$ p^{(0)}_1 = a, \ \ p^{(0)}_{j+1} = p^{(0)}_j + h(0,j).$$
Similarly, if a tuple $p^{(i)} = (p^{(i)}_1, \ldots, p^{(i)}_L)$
is defined then put $p^{(i+1)} = (p^{(i+1)}_1, \ldots,
p^{(i+1)}_L)$, where
  $$ p^{(i+1)}_1 = p^{(i)}_1 + h(i+1,0), \ \ p^{(i+1)}_{j+1} = p^{(i+1)}_j + h(i+1,j).$$
This defines by induction an infinite set
 $${\mathcal P}_{a,h} = \{p^{(i)} \mid i \in \mathbb{N}\} \subseteq
\mathbb{N}^L$$
 such that  any infinite subset of ${\mathcal P}_{f,h}$ is also unbounded.

From now on fix a recursive non-negative monotonically increasing
with respect to both variables function $h$ (which will be
specified in due course) and put
$${\mathcal P}_2 =  {\mathcal P}_{a,h}, \ \ \ {\mathcal L_2} = \{\phi_{L,p} \mid  p \in {\mathcal P}_2 \}.$$

\begin{prop} \label{pr:disc-p}
Let $r \geq 2$ and $K(r+2) \leq L$ then there exists a number
$a_0$ such that if $a \geq a_0$ and the function $h$ satisfies the
condition
 \begin{equation}
 \label{eq:h}
 h(i+1,j) > h(i,j) \ \ \text{for any} \  j = Kr+1,
 \ldots, K(r+2), i = 1,2, \ldots;
\end{equation}
then for any infinite subset ${\mathcal P} \subseteq {\mathcal
P}_2$ the set of solutions
 $${\mathcal L_{\mathcal P}}^\beta = \{\phi_{L,p}\beta \mid  p \in {\mathcal P} \}$$
is a discriminating set of solutions of $S(X,A) = 1$.
\end{prop}
\begin{proof}  The result  follows from Corollary \ref{cy:cy2}.
 \end{proof}

   Let  $\psi = \psi_p \in {\mathcal L}_2^\beta$. Denote by $U_\psi$ the solution
  $X^{\psi}$ of the equation $S(X) = 1$ in $F$.  Since  $T(X,Y) = 1$
   is compatible with $S(X) = 1$ in $F$ the equation
$T(U_\psi,Y) = 1$ (in variables $Y$)  has a solution in $F$, say
$Y = V_\psi$.  Set
  $$\Lambda = \{(U_\psi, V_\psi) \mid \psi \in {\mathcal
  L}_2^\beta\}.$$
  It follows that every pair
$(U_\psi,V_\psi) \in \Lambda $ gives a solution of  the system
    $$ R(X,Y) = (S(X) = 1 \  \wedge \ T(X,Y) = 1).$$
      By Theorem
 \ref{th:cut} there exists a finite set ${\mathcal CE}(R)$ of cut equations
  which describes all solutions of $R(X,Y) = 1$ in $F$, therefore there exists a cut equation
$ \Pi_{\mathcal L_3, \Lambda}  \in {\mathcal CE}(R)$ and an
infinite subset ${\mathcal L}_3 \subseteq {\mathcal L_2}$ such
that $\Pi_{\mathcal L_3, \Lambda}$ describes all solutions of the
type $(U_\psi, V_\psi)$, where  $\psi \in  {\mathcal L}_3$.
 We state the precise formulation of this result  in the following
 proposition which, as we have mentioned already,  follows from Theorem \ref{th:cut}.

\begin{prop}
\label{prop:cut-L3}
 Let ${\mathcal L_2}$ and $\Lambda$ be as above. Then
 there exists  an infinite subset ${\mathcal P}_3 \subseteq {\mathcal P}_2$ and the corresponding set
 ${\mathcal L}_3  = \{\phi_{L,p} \mid p \in {\mathcal P}_3\} \subseteq
{\mathcal L_2}$,   a cut equation $\Pi_{\mathcal L_3, \Lambda} =
({\mathcal E}, f_X, f_M) \in {\mathcal CE}(R)$, and a
  tuple of words $Q(M)$ such that the following  conditions hold:
 \begin{enumerate}
 \item [1)] $f_X({\mathcal E}) \subset X^{\pm 1}$;
 \item [2)]  for every $\psi \in {\mathcal L}_3^\beta$  there exists a
tuple of words $P_\psi = P_\psi(M)$ and a solution $\alpha_\psi: M
\rightarrow F$ of  $\Pi_{\mathcal L_3, \Lambda}$ with respect to
$\psi: F[X] \rightarrow F$ such that:
 \begin{itemize}
  \item the solution  $U_\psi = X^\psi$ of $S(X) =1$ can be
presented as $U_\psi = Q(M^{\alpha_\psi})$ and the word
$Q(M^{\alpha_\psi})$ is reduced as written,
 \item $V_\psi = P_\psi(M^{\alpha_\psi})$.
 \end{itemize}
 \item [3)] there exists a tuple of words
$P$ such that for any solution (any group solution) $(\beta,
\alpha)$ of $\Pi_{\mathcal L_3, \Lambda}$  the pair $(U,V),$ where
$U = Q(M^\alpha)$ and $V = P(M^\alpha),$ is a solution of $R(X,Y)
= 1$ in $F$.
 \end{enumerate}
\end{prop}

Put
 $${\mathcal P} = {\mathcal P}_3, \ \ \ {\mathcal L} = {\mathcal L}_3, \ \ \  \Pi_{\mathcal L} = \Pi_{\mathcal L_3, \Lambda}.$$
   By Proposition
\ref{pr:disc-p}
 the set ${\mathcal L}^\beta$  is a discriminating
set of solutions of $S(X) = 1$ in $F$.

 \medskip
 {\bf The initial cut equation $\Pi_{\phi}$.}

 \medskip
 Now fix a tuple $p \in {\mathcal P}$ and the  automorphism $\phi= \phi_{L,p}  \in {\mathcal L}$.
Recall,  that for every $j \leq L$ the automorphism $\phi_j$ is
defined by $\phi_j = \stackrel{\leftarrow}{\Gamma}_j^{p_j}$, where
$p_j$ is the initial subsequence of $p$ of length $j$. Sometimes
we use  notation $\psi = \phi \beta, \psi_j = \phi_j  \beta$.

Starting with the cut equation $\Pi_{\mathcal L}$ we construct a
cut equation $\Pi_{\phi} = ( {\mathcal E}, f_{\phi,X}, f_M) $
which  is obtained from $\Pi_{\mathcal L}$  by replacing the
function $f_X: {\mathcal E} \rightarrow F[X] $ by a new function
$f_{\phi,X}: {\mathcal E} \rightarrow F[X]$, where $f_{\phi,X}$ is
the composition of $f_X$ and the automorphism $\phi$. In other
words, if an interval $e \in {\mathcal E}$ in $\Pi_{\mathcal L}$
has a label $x \in X^{\pm 1}$ then its label in $\Pi_{\phi}$ is
$x^\phi$.

Notice, that $\Pi_{\mathcal L}$ and $\Pi_{\phi}$ satisfy  the
following conditions:
 \begin{enumerate}
  \item [a)]  $\sigma ^{f_X \phi \beta} = \sigma ^{f_{\phi,X}
  \beta}$ for every $ \sigma \in {\mathcal E}$;
  \item [b)] the solution  of $\Pi_{\mathcal L}$ with respect to
  $\phi  \beta$ is also a solution of $\Pi_{\phi}$ with respect to
  $\beta$;
  \item [c)] any solution (any group solution) of $\Pi_\phi$ with
  respect to $\beta$ is a solution (a group solution) of $\Pi_{\mathcal L}$
  with respect to $\phi \beta$.
  \end{enumerate}

 The cut equation $\Pi_\phi$ has a very particular type. To deal with such cut
 equations we need the following definitions.

\begin{df}
Let $\Pi = ({\mathcal E}, f_X, f_M)$ be a cut equation. Then the
number $$length(\Pi) = \max \{|f_M(\sigma)| \mid \sigma \in
{\mathcal E}\}$$ is called  the length of $\Pi$. We denote it by
$length(\Pi)$  or simply by $N_{\Pi}$.
\end{df}

Notice,  by construction,  $length(\Pi_\phi) =
length(\Pi_{\phi^\prime})$ for every $\phi, \phi^\prime \in
{\mathcal L}$. Denote
 $$N_{\mathcal L} = length(\Pi_\phi).$$

 \begin{df}
\label{df:7.3.cutgamma} A  cut equation $\Pi = ( {\mathcal E},
f_{X}, f_M)$ is called a $\Gamma$-cut equation in {\em  rank} $j$
($rank(\Pi) = j$) and size $l$ if it satisfies the following
conditions:
\begin{enumerate}
\item [1)] let $W_\sigma = f_X(\sigma)$ for $\sigma \in {\mathcal
E}$ and  $N = (l+2)N_{\Pi}$.
  Then for every $\sigma \in {\mathcal E}$ $W_\sigma \in \bar {\mathcal W}_{\Gamma,L}$
  and one of the following conditions
  holds:
\begin{enumerate}
  \item [1.1)] $W_\sigma$ has $N$-large rank $j$ and its canonical
 $N$-large  $A_j$-decomposition has  size
$(N,2)$ i.e., $W_\sigma$ has the canonical $N$-large
$A_j$-decomposition
\begin{equation}
 \label{eq:Adecomp-1}
 W_{\sigma} = B_1 \circ A^{q_1}_j \circ \ldots B_k \circ A^{q_k}_j \circ B_{k+1},
 \end{equation}
with   $max_j(B_i) \leq 2$ and $q_i \geq N$;
 \item  [1.2)] $W_\sigma$ has rank $j$ and $\max_j(W_\sigma) \leq 2$;
\item [1.3)]  $W_\sigma$ has rank $< j$.
\end{enumerate}
Moreover,  there exists  at least one interval $\sigma \in
{\mathcal E}$ satisfying the condition 1.1).
 \item [2)] there exists a solution $\alpha : F[M]
\rightarrow F$ of the cut equation $\Pi$ with respect to the
homomorphism   $\beta:F[X] \rightarrow F$.
\end{enumerate}
\end{df}

\begin{lm} Let $l\geq 3$.
The cut equation $\Pi_\phi$ is a $\Gamma$-cut equation in rank $L$
and size $l$, provided  $$ p_L \geq (l+2)N_{\Pi_\phi}+3.$$
\end{lm}
\begin{proof} By construction the labels of intervals from $\Pi_\phi$ are
precisely the words of the type $x^{\phi_L}$ and every such word
appears as a label. Observe, that $rank(x_i^{\phi_L}) < L$ for
every $i, 1 \leq i \leq n$ (Lemma \ref{le:7.1.gammawords}, 1a).
Similarly, $rank(x_i^{\phi_L}) < L$ for every $i < n$ and
$rank(y_n^{\phi_L}) = L$ (Lemma \ref{le:7.1.gammawords} 1b). Also,
$rank(z_i^{\phi_L}) < L$ unless $n = 0$ and  $i = m$, in the
latter case $z_m^{\phi_L}) = L$ (Lemma \ref{le:7.1.gammawords} 1c
and 1d).
  Now consider the labels $y_n^{\phi_L}$ and $z_m^{\phi_L})$ (in the case $n =
  0$) of  rank $L$.  Again, it has been shown in Lemma
  \ref{le:7.1.gammawords} 1) that these labels have
 $N$-large  $A_L$-decompositions of size $(N,2)$, as required in 1.1)
 of the definition of a $\Gamma$-cut equation of rank $L$ and size
 $l$.

\end{proof}

\medskip

 {\bf Agreement 1 on ${\mathcal P}.$} Fix an arbitrary integer $l$, $l \geq 5$.
 We may assume, choosing the
 constant $a$ to satisfy the condition
   $$a \geq (l+2)N_{\Pi_\phi}+3,$$
   that all tuples in the set ${\mathcal P}$ are
$[(l+2)N_{\Pi_\phi}+3]$-large. Denote $N = (l+2)N_{\Pi_\phi}.$

Now we introduce one  technical restriction on the set ${\mathcal
P}$, its real meaning will be clarified later.

\medskip
 {\bf Agreement 2 on ${\mathcal P}.$}
Let $r$ be an arbitrary fixed positive integer with $Kr \leq L$
and $q$ be a fixed tuple of length $Kr$ which is an initial
segment of some tuple from
 ${\mathcal P}$. The choice of $r$ and $q$ will be clarified later.
  We may assume (suitably choosing the  function
  $h$) that all tuples from ${\mathcal P}$
 have $q$ as their initial segment. Indeed, it suffices to define
 $h(i,0) = 0$ and $h(i,j) = h(i+1,j)$ for all $i \in \mathbb{N}$ and
 $j = 1, \ldots, Kr$.

\medskip
 {\bf Agreement 3 on ${\mathcal P}.$} Let $r$ be the number from  Agreement
 2. By Proposition \ref{pr:disc-p} there exists a number $a_0$ such
 that for every infinite subset of ${\mathcal P}$ the corresponding set
 of solutions is a discriminating set. We may assume that $a > a_0$.

\medskip
{\bf Transformation $T^*$ of $\Gamma$-cut equations.}

\medskip
 Now we describe a  transformation $T^*$ defined on
  $\Gamma$-cut
 equations and their solutions,  namely, given a $\Gamma$-cut equation $\Pi$ and its solution
 $\alpha$ (relative to the  fixed map $\beta:F[X] \rightarrow F$ defined above) $T^*$
 transforms $\Pi$ into a new $\Gamma$-cut equation $\Pi^* = T^*(\Pi)$ and $\alpha$ into a solution
 $\alpha^* = T^*(\alpha)$ of $T^*(\Pi)$ relative to $\beta$.

 Let $\Pi = ( {\mathcal E}, f_{X}, f_M)$ be a $\Gamma$-cut equation in rank $j$ and size
 $l$. The cut equation
 $$T^*(\Pi) = ( {\mathcal E^*}, f^*_{X^*}, f^*_{M^*})$$
 is defined as follows.

\medskip
{\bf Definition of the set ${\mathcal E}^*$.}

\medskip
 For $\sigma \in {\mathcal E}$ we denote $W_{\sigma} = f_X(\sigma)$. Put
 $${\mathcal E}_{j,N} = \{\sigma \in {\mathcal E} \mid W_\sigma \ \mbox{satisfies} \ 1.1) \}.$$
 Then ${\mathcal E} = {\mathcal E}_{j,N} \cup {\mathcal E}_{< j,N}$ where ${\mathcal E}_{< j,N}$ is
 the complement
 of $ {\mathcal E}_{j,N}$ in ${\mathcal E}$.

 Now let $\sigma \in  {\mathcal E}_{j,N}$.
Write the  word $W_\sigma^\beta$ in its canonical $A^\prime$
decomposition:
\begin{equation}
\label{eq:can-W-sigma-beta}
 W_\sigma^\beta = E_1 \circ {A^\prime}^{q_1} \circ E_2 \circ \cdots \circ E_k \circ
 {A^\prime}^{q_k} \circ E_{k+1}
 \end{equation}
 where $|q_i| \geqslant 1$, $E_i \neq 1$ for $2 \leqslant i \leqslant k$.

 Consider the partition $$f_M(\sigma)  = \mu_1 \ldots \mu_n$$ of $\sigma$.
  By the condition 2) of the definition of $\Gamma$-cut equations
   for the solution  $\beta:F[X] \rightarrow F$ there exists a solution
  $\alpha : F[M] \rightarrow F$ of the cut equation $\Pi$ relative to $\beta$. Hence
 $W_\sigma^\beta = f_M(M^\alpha)$ and the element
 $$f_M(M^\alpha) = \mu_1^\alpha \ldots \mu_n^\alpha$$ is reduced as written.
   It  follows that
   \begin{equation}
   \label{eq:7.3.13}  W_\sigma^\beta = E_1 \circ {A^\prime}^{q_1} \circ E_2 \circ \cdots  \circ E_k \circ
 {A^\prime}^{q_k} \circ E_{k+1}  = \mu_1^{\alpha} \circ  \cdots \circ \mu_n^{\alpha}.
\end{equation}

   We say that a variable $\mu_i$ is {\em long} if ${A^\prime}^{\pm (l+2)}$  occurs in
   $\mu_i^\alpha$ (i.e., $\mu_i^\alpha$ contains a stable occurrence of
   ${A^\prime}^{l}$),   otherwise it is called {\em short}.
   Observe, that the definition of long  (short) variables $\mu
   \in M$ does not depend on a choice of $\sigma$,  it depends only
   on the given homomorphism $\alpha$.
The graphical  equalities  (\ref{eq:7.3.13}) (when $\sigma$ runs
over ${\mathcal E}_{j,N}$)  allow one to effectively recognize
long and short variables in $M$. Moreover, since for every $\sigma
\in {\mathcal E}$ the length of the word $f_M(\sigma)$ is bounded
by $length(\Pi) = N_{\Pi}$ and $N = (l+2)N_{\Pi}$, every word
$f_M(\sigma)$ ($\sigma \in {\mathcal E}_j$) contains long
variables. Denote by $M_{\rm short}$, $M_{\rm long}$ the sets of short
and long variables in $M$. Thus, $M = M_{\rm short} \cup M_{\rm long}$
is a non-trivial partition of $M$.

Now we define the following property $P = P_{long,l}$ of
occurrences of powers of $A^\prime$ in $W_\sigma^\beta$: a given
stable occurrence ${A^\prime}^q$ satisfies $P$ if it occurs in
$\mu^\alpha$ for some long variable $\mu \in M_{long}$ and $q
\geqslant l$. It is easy to see that $P$ preserves correct
overlappings. Consider the set of stable occurrences ${\mathcal
O}_{P}$ which are maximal with respect to $P$. As we have
mentioned already in Section \ref{se:7.2.5}, occurrences from
${\mathcal O}_{P}$ are pair-wise disjoint and this set is uniquely
defined. Moreover, $W_\sigma^\beta$ admits the  unique
$A^\prime$-decomposition relative  to the set ${\mathcal O}_{P}$:
 \begin{equation}
  \label{eq:7.3.can*}
   W_\sigma ^\beta = D_1 \circ (A^\prime)^{q_1}
   \circ D_2 \circ \cdots \circ D_k
 \circ (A^{\prime})^{q_k}  \circ
   D_{k+1},
  \end{equation}
where $D_i \neq 1$ for $i = 2, \ldots, k$. See Figure 1.
\begin{figure}[here]
\centering{\mbox{\psfig{figure=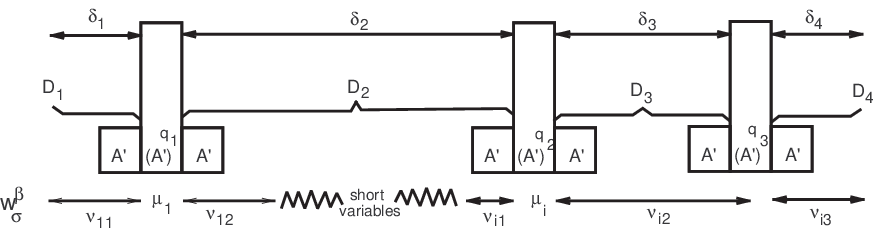,height=1.3in}}}
\caption{ Decomposition
(\ref{eq:7.3.can*})}\label{picture-decomposition-997}
\end{figure}

 Denote by $k(\sigma)$
the number of nontrivial elements among $D_1, \ldots, D_{k+1}$.

According to Lemma \ref{Claim1.b}  the $A^\prime$-decomposition
\ref{eq:7.3.can*} gives rise to the  unique associated
$A$-decomposition of $W_\sigma$ and hence the  unique associated
$A^*$-decomposition of $W_\sigma$.

Now with a given  $\sigma \in {\mathcal E}_j$ we
  associate a finite set of new intervals $E_\sigma$ (of the equation
  $T^*(\Pi)$):
  $$E_\sigma = \{\delta_1, \ldots , \delta_{k(\sigma)}\}$$
   and put
 $${\mathcal E}^* = {\mathcal E}_{<j} \cup \bigcup_{\sigma \in {\mathcal E}_j} E_\sigma.$$

\medskip
{\bf Definition of the set $M^*$}

 \medskip
 Let $\mu \in M_{long}$ and
 \begin{equation}
 \label{eq:7.3.14}
 \mu^\alpha = u_1\circ (A^\prime)^{s_1} \circ u_2 \circ \cdots \circ u_t\circ
 (A^{\prime})^{s_t}
 \circ u_{t+1}
 \end{equation}
  be the canonical $l$-large $A^\prime$-decomposition of $\mu^\alpha$.  Notice that
if $\mu$ occurs in  $f_M(\sigma)$ (hence $\mu^\alpha$ occurs in
$W_\sigma^\beta$) then this decomposition (\ref{eq:7.3.14})
   is precisely the $A^\prime$-decomposition of $\mu^\alpha$
  induced on $\mu^\alpha$ (as a subword of $W_\sigma^\beta$)
  from the $A^\prime$-decomposition (\ref{eq:7.3.can*}) of $W_\sigma^\beta$  relative to ${\mathcal O}_P$.

   Denote by  $t(\mu)$ the number of non-trivial elements among  $u_1, \ldots,  u_{t+1}$
 (clearly, $u_i \neq 1$ for $2 \leqslant i \leqslant t$).

We associate with each long variable $\mu$ a sequence of new
variables (in the equation $T^*(\Pi)$)
 $S_\mu = \{\nu_1, \ldots, \nu_{t(\mu)} \}$. Observe, since the decomposition
(\ref{eq:7.3.14}) of $\mu^\alpha$ is unique, the set $S_\mu$ is
well-defined (in particular, it does not depend on   intervals
$\sigma$).

 It is convenient to define here two
functions $\nu_{\rm left}$ and $\nu_{\rm right}$  on the set
$M_{long}$: if $\mu \in M_{long}$ then
 $$\nu_{\rm left}(\mu) =\nu_1, \ \ \ \nu_{\rm right}(\mu) =\nu_{t(\mu)}.$$

Now we define a new set of variable $M^*$ as follows:
$$ M^* = M_{\rm short}\cup \bigcup_{\mu \in M_{long}} S_\mu .$$

\medskip
{\bf Definition of the labelling function $f^*_{X^*}$}

\medskip
Put $X^* = X$. We define the labelling function $f^*_{X^*} :
{\mathcal E}^* \rightarrow F[X]$ as follows.

Let $\delta \in {\mathcal E}^*$. If $\delta \in {\mathcal
E}_{<j}$, then put $$f^*_{X^*}(\delta) = f_X(\delta).$$

 Let now  $\delta = \delta_i \in
E_\sigma$ for some $\sigma \in M_{\rm long}$. Then there are three
cases to consider.

a) $\delta$ corresponds to the consecutive occurrences of powers
${A^\prime}^{q_{j-1}}$ and ${A^\prime}^{q_{j}}$ in the
$A^\prime$-decomposition (\ref{eq:7.3.can*}) of $W_\sigma^\beta$
relative to ${\mathcal O}_P$. Here $j = i$ or $j = i-1$ with
respect to whether $D_1 = 1$ or $D_1 \neq 1$.

As we have mentioned before, according to Lemma \ref{Claim1.b} the
$A^\prime$-decomposition (\ref{eq:7.3.can*}) gives rise to the
unique associated $A^*$-decomposition of $W_\sigma$:
$$W_\sigma = D_1^* \circ_d (A^*)^{q^*_1} \circ_d D^*_2 \circ \cdots \circ_d D^*_k \circ_d (A^*)^{q^*_{k}} \circ_d
D^*_{k+1}.$$

Now put
 $$f^*_X(\delta_i) = D_j^* \in F[X]$$
where $j = i$ if $D_1 = 1$ and  $j = i-1$ if  $D_1 \neq 1$. See
Figure 2.
\begin{figure}[here]
\centering{\mbox{\psfig{figure=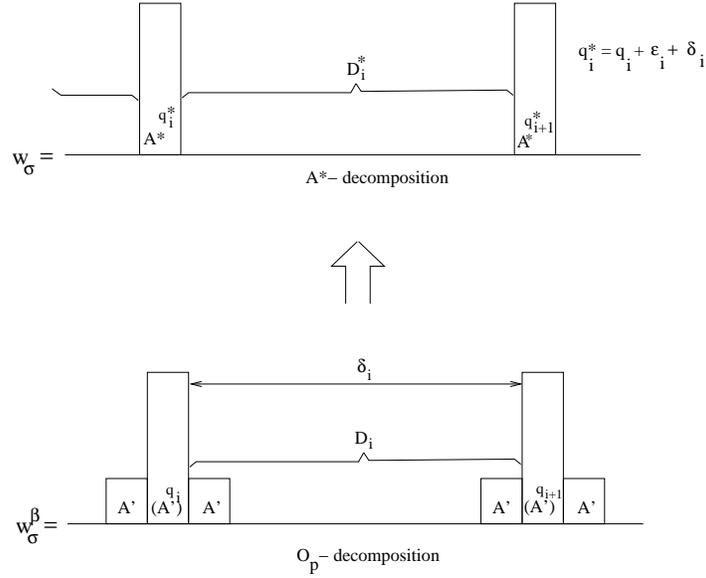,height=3in}}}
\caption{Defining $f^*_{X^*}$.} \label{fx-star}
\end{figure}

The other two cases are treated similarly to  case a).

 b) $\delta$ corresponds to the interval from the beginning of
$\sigma$ to the first $A^\prime$ power ${A^\prime}^{q_1}$ in the
decomposition (\ref{eq:7.3.can*}) of $W_\sigma^\beta$. Put
$$f^*_X(\delta )=D_1^*.$$

c) $\delta$ corresponds to the interval from the last occurrence
of  a power ${A^\prime}^{q_k}$ of $A^\prime$   in the
decomposition (\ref{eq:7.3.can*}) of $W_\sigma^\beta$ to the end
of the interval. Put
$$f^*_X(\delta )=D_{k+1}^*.$$

\medskip
{\bf   Definition of the function $f^*_{M^*}$.}

\medskip
Now we define the function $f^*: {\mathcal E}^* \rightarrow
F[M^*]$.

Let $\delta \in {\mathcal E}^*$. If $\delta \in {\mathcal
E}_{<j}$, then put
 $$f^*_{M^*}(\delta) = f_M(\delta)$$
 (observe that all variables in
$f_M(\delta)$ are short, hence they belong to $M^*$).

Let $\delta = \delta_i \in E_\sigma$ for some $\sigma \in
M_{long}$. Again, there are three cases to consider.

a) $\delta$ corresponds to the consecutive occurrences of powers
${A^\prime}^{q_s}$ and ${A^\prime}^{q_{s+1}}$ in the
$A^\prime$-decomposition (\ref{eq:7.3.can*}) of $W_\sigma^\beta$
relative to ${\mathcal O}_P$.  Let the stable occurrence
${A^\prime}^{q_s}$ occur in $\mu_i^\alpha$ for a long variable
$\mu_i$, and the stable occurrence ${A^\prime}^{q_{s+1}}$ occur in
$\mu_j^\alpha$ for a long variable $\mu_j$.

Observe that
$$D_s = right(\mu_i)\circ  \mu_{i+1}^\alpha \circ \cdots \circ
\mu_{j-1}^\alpha  \circ left(\mu_j),$$ for some elements
$right(\mu_i), left(\mu_j) \in F$.

 Now put
$$f^*_{M^*}(\delta) = \nu_{i,\rm right} \mu_{i+1} \ldots
\mu_{j-1} \nu_{j,left},$$
  See Figure 3.
 \begin{figure}[here]
\centering{\mbox{\psfig{figure=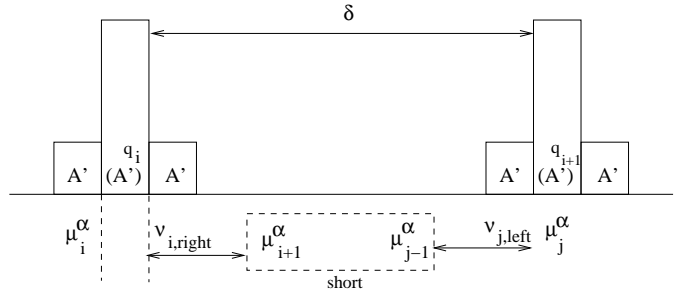,height=1.5in}}}
\caption{Defining $f^*_{M^*},$ case a)}\label{case-a}
\end{figure}

The other two cases are treated similarly to case a).

 b) $\delta$ corresponds to the interval from the beginning of
$\sigma$ to the first $A^\prime$ power ${A^\prime}^{q_1}$ in the
decomposition (\ref{eq:7.3.can*}) of $W_\sigma^\beta$. Put
 $$f^*_{M^*}(\delta) =  \mu_{1} \ldots
\mu_{j-1}\nu _{j,left}.$$

c) $\delta$ corresponds to the interval from the last occurrence
of  a power ${A^\prime}^{q_k}$ of $A^\prime$   in the
decomposition (\ref{eq:7.3.can*}) of $W_\sigma^\beta$ to the end
of the interval.

The cut equation $T^*(\Pi) = ({\mathcal E}^*, f^*_X, f^*_{M^*}) $
has been  defined.

We define now  a sequence
\begin{equation}
\label{eq:7.3.15^*}
 \Pi_L \stackrel{T^*}{ \rightarrow} \Pi_{L-1}  \stackrel{T^*}{
\rightarrow} \ldots \stackrel{T^*}{ \rightarrow} \Pi_1
\end{equation} of $N$-large $\Gamma$-cut equations,
 where  $\Pi_L = \Pi_\phi$, and $\Pi_{i-1} = T^*(\Pi_i)$.
 In Claims  \ref{2b} and  \ref{2} below we show that in this case if
 $\Pi$ is a $\Gamma$-cut  equation then
 $T^*(\Pi)$ is also a $\Gamma$-cut  equation of the corresponding
 rank and size, so the sequence is well-defined.
 However, it is convenient to assume this as a fact now and
  introduce some notation and agreements before proving the claims.

\begin{Claim} \label{2.0.1} Let $\Pi_j$ be a cut equation from the
sequence (\ref{eq:7.3.15^*}). Then there exists an infinite subset
${\mathcal P}' \subseteq {\mathcal P}$ such that the cut equation
$\Pi_{j-1} = T^*(\Pi_j)$ satisfies the following conditions:

 \begin{enumerate}
  \item  the  words $f_{X^*}(\sigma)
\in F[X]$, as parametric words in the parameters from $p$, are the
same for every $p \in {\mathcal P}'$, i.e., they differ only in
exponents corresponding to components  of the tuples  $p$.
 \item the words $f_{M^*}(\sigma)$ are the same for every $p \in
 {\mathcal P}'$.
 \end{enumerate}
\end{Claim}
\begin{proof}  The claim follows from the construction of $T^*(\Pi )$.
 \end{proof}

\medskip
{\bf Agreement 4 on the set ${\mathcal P}$:} we assume (replacing
$P$ with a suitable infinite subset) that every tuple $p \in
{\mathcal P}$ satisfies the conditions of  Claim \ref{2.0.1}.
Thus, every $\Pi = \Pi_i$ from the sequence (\ref{eq:7.3.15^*})
satisfies the conclusion of Claim \ref{2.0.1} for ${\mathcal P}' =
{\mathcal P}$.

\begin{Claim} \label{2.0.-1.}
 The homomorphism $\alpha^*: F[M^*]
\rightarrow F$  defined as (in the notations above):
$$\alpha^*(\mu) = \alpha(\mu) \ \ (\mu \in M_{\rm short}), $$
$$\alpha^*(\nu_{i,right}) = R^{-\beta}c^{-1} right(\mu_i) \ \ (\nu_i \in S_{\mu} \ for \ \mu \in M_{long})$$
$$\alpha^*(\nu_{i,left}) =  left(\mu_i)c R^{\beta}$$
is a solution of the cut equation $T^*(\Pi)$ with respect to
$\beta: F[X] \rightarrow F$.
\end{Claim}
\begin{proof}   Indeed,  by Lemma \ref{le:D-star-beta}
$$(D_s^*)^\beta = (R^{-\beta}c^{-1}) \circ_{\theta} D_s \circ_{\theta}
(cR^\beta)$$
 where $\theta << |A^\prime|$. Therefore, $\mu_{i+1}^\alpha \ldots
 \mu_{j-1}^{\alpha}$ occurs in $D_s$ without cancellation.
 Therefore $\alpha^*$ is a required solution. \end{proof}

\medskip

\medskip
{\bf Agreement 5 on the set ${\mathcal P}$:} we assume (by
choosing the function $h$ properly, i.e., $h(i,j) > C(L,N+3)$, see
Lemma )
 that every tuple $p
\in {\mathcal P}$ satisfies the conditions of Lemma
\ref{le:D-star-beta}, so Claim \ref{2.0.-1.} holds for every $p
\in {\mathcal P}$. Thus, for every $\Pi = \Pi_i$ from the sequence
(\ref{eq:7.3.15^*}) with a solution $\alpha$ (relative to $\beta$)
the solution $\alpha^*$ of the equation $T^*(\Pi)$ is defined as
in Claim \ref{2.0.-1.}.

\begin{Claim} \label{2.0.0}  Let $\Pi =({\mathcal E},f_X,f_M)$ be a
$\Gamma$-cut equation in rank $j \geq 1$ from the sequence
(\ref{eq:7.3.15^*}). Then for every variable $\mu \in M$ there
exists a word ${\mathcal M}_\mu (M_{T(\Pi)}, X^{\phi_{j-1}},F)$
such that the following equality holds in the group $F$
$$\mu^\alpha = {\mathcal M}_\mu (M_{T(\Pi)}^{\alpha^*}, X^{\phi_{j-1}})^\beta.$$
Moreover, there exists an infinite subset $P' \subseteq P$ such
that
 the words ${\mathcal M}_\mu
(M_{T(\Pi)}, X)$ depend only on exponents $s_1, \ldots, s_t$ of
the canonical $l$-large decomposition (\ref{eq:7.3.14}) of the
words $\mu^\alpha$.
\end{Claim}
\begin{proof}   The claim follows from the construction. Indeed, in
constructing
 $T(\Pi)$ we cut out leading periods of the type $(A_j^\prime)^s$
 from $\mu^\alpha$ (see (\ref{eq:7.3.14})). It follows that to get
 $\mu^\alpha$ back from $M_{T(\Pi)}^{\alpha^*}$ one needs to put
 the exponents $(A_j^\prime)^s$ back. Notice, that
 $$ A_j = A(\gamma_j)^{\phi_{j-1}}$$
  Therefore,  $$ (A_j)^s = A(\gamma_j)^{\phi_{j-1}\beta}$$
 Recall that $A_j^\prime$ is the cyclic reduced form of
 $A_j^\beta$, so
 $$ (A_j^\prime)^s = uA(\gamma_j)^{\phi_{j-1}\beta}v$$
  for some  constants $u, v \in C_{\beta}\subseteq F$. To see
  existence of the subset $P' \subseteq P$ observe that the length
  of the words $f_M(\sigma)$ does not depend on $p$, so there are
  only finitely many ways to cut out the leading periods  $(A_j^\prime)^s$
 from $\mu^\alpha$.  This proves the claim.
\end{proof}

\medskip
{\bf Agreement 6 on the set ${\mathcal P}$:} we assume (replacing
$P$ with a suitable infinite subset) that every tuple $p \in
{\mathcal P}$ satisfies the conditions of  Claim \ref{2.0.0}.
Thus, for every $\Pi = \Pi_i$ from the sequence
(\ref{eq:7.3.15^*}) with a solution $\alpha$ (relative to $\beta$)
the solution $\alpha^*$ satisfies the conclusion of Claim
\ref{2.0.0}.

\begin{df}
We define a new transformation $T$ which is a modified version of
$T^\ast$. Namely, $T$ transforms cut equations and their solutions
$\alpha$ precisely as the transformation $T^\ast$, but it also
transforms the set of tuples ${\mathcal P}$ producing an infinite
subset ${\mathcal P}^\ast \subseteq {\mathcal P}$ which satisfies
the Agreements 1-6.
\end{df}

Now we define a sequence
\begin{equation}
\label{eq:7.3.15}
 \Pi_L \stackrel{T}{ \rightarrow} \Pi_{L-1}  \stackrel{T}{
\rightarrow} \ldots \stackrel{T}{ \rightarrow} \Pi_1
\end{equation} of $N$-large $\Gamma$-cut equations,
 where  $\Pi_L = \Pi_\phi$, and $\Pi_{i-1} = T(\Pi_i)$.
 From now on we fix the sequence (\ref{eq:7.3.15}) and refer to it as the
 {\em $T$-sequence.}

\begin{Claim}\label{2b} The following statements are true:

1)  for every $i = 1, \ldots, L/K$ and every interval
 $\sigma$ of the cut equation $\Pi _{L-iK}$  from the $T$-sequence (\ref{eq:7.3.15})
there exists a word $w = w_\sigma \in \bar{\mathcal W}_{\Gamma, L}$
without
  $N$-large powers of elementary periods
 such that  $ f_X(\sigma) = w^{\phi _{L-iK}};$

 2)  for every $j = 1, \ldots, L$ and every interval
 $\sigma$ of the cut equation $\Pi_{L-j}$ from the $T$-sequence (\ref{eq:7.3.15})
the label $f_X(\sigma)$ of $\sigma$  belongs to $\bar{\mathcal
W}_{\Gamma, L}$.
 \end{Claim}
 \begin{proof} We  prove the claim by induction on $i$.

 Let  $i=1$.
For every $x \in X^{\pm 1}$ one can represent
 the  element $x^{\phi _L}$ as a product of elements of the type
 $y^{\phi_{L-K}}, y \in X^{\pm 1}$ (in this event we say that the
 element $x^{\phi _L}$ is a word in the alphabet $X^{\phi _{L-K}}$).
 Indeed,
  $$x^{\phi _L}= (x^{\phi_K})^{\phi _{L-K}} = w^{\phi _{L-K}},$$
   where $w = x^{\phi_K}$ is a word in $X$.

 Now consider the first $K$ terms in the $T$-sequence:
$$\Pi_{L}\rightarrow \ldots\rightarrow \Pi_{L-K}.$$
 We use induction on $m$ to prove that for every interval
 $\sigma \in \Pi_{L-m} = \newline({\mathcal E}^{(L-m)}, f_{X}^{(L-m)}, f_M^{(L-m)})$
  the label $f_X^{(L-m)}(\sigma)$ is of the form $u^{\phi_{L-K}}$
  for some  $u \in Sub(X^{\phi_K})$.

For $m=1$  by Lemma \ref{le:A-A^*} for $j=L, r=K,$ there is a
precise correspondence between stable $A_L^*$-decompositions of
 $$x^{\phi _L}= w^{\phi
_{L-K}}= D_1^{\phi _{L-K}}\circ _d A_L^{*q_1}\circ _dD_2^{\phi
_{L-K}}\ldots D_k^{\phi _{L-K}}\circ_d A_L^{*q_k}\circ D_{k+1}^{\phi
_{L-K}}$$ and stable $A_K$-decompositions of $w$
 $$w=D_1\circ
{A_K}^{q_1}\circ D_2\ldots D_k\circ A_K^{q_k}\circ D_{k+1}.$$  By
construction, application of the transformation $T$ to $\Pi _L$
removes  powers  $A_L^{*q_s}=A_K^{q_s\phi _{L-K}}$ which are
subwords of the word $w^{\phi _{L-K}}$ written in the alphabet
$X^{\phi _{L-K}}$. By construction the words $ D_s^{\phi _{L-K}}$
are the labels of the new intervals of the equation $\Pi _{L-1}$.
Suppose by induction  that for an interval $\sigma $ of the cut
equation $\Pi _j$ (for $m = L-j$) $f_X^{(j)}(\sigma)=u^{\phi
_{L-K}}$ for some  $u\in Sub(X^{\pm\phi _K}).$ Then either $\sigma$
does not change under $T$  or $f_X^{(j)}(\sigma )$ has
 a stable  $(l+2)$-large ${A_j}^*$-decomposition  in  rank
$j=r+(L-K)$ associated with long variables in $f_M^{(j)}(\sigma)$:
$$u^{\phi _{L-K}}=\bar D_1^{\phi _{L-K}}\circ _d A_j^{*q_1}\circ
_d\bar D_2^{\phi _{L-K}}\ldots \bar D_k^{\phi _{L-K}}\circ_d
A_j^{*q_k}\circ \bar D_{k+1}^{\phi _{L-K}},$$ and $\sigma$ is an
interval in $\Pi_j$. By Lemma \ref{le:A-A^*}, in this case there is
a stable $A_r$-decomposition of $u$:
$$u=\bar D_1\circ
A_r^{q_1}\circ \bar D_2\ldots \bar D_k\circ A_r^{q_k}\circ \bar
D_{k+1}.$$
 The application of the transformation
$T$ to $\Pi _j$ removes powers $A_j^{*q_s} = A_r^{q_s\phi _{L-K}}$
 (since ${A_j}^*=A_r^{\phi _{L-K}}$)  which are subwords of the word $u^{\phi
_{L-K}}$ written in the alphabet $X^{\phi _{L-K}}$. By construction
the words $\bar D_s^{\phi _{L-K}}$ are the labels of the new
intervals of the equation $\Pi _{j-1}$, so they have  the required
form. By induction the statement holds for $m = K$, so the label
$f_X^{(L-K)}(\sigma)$ of an interval $\sigma$ in $\Pi _{L-K}$ is of
the form $u^{\phi _{L-K}},$ for some $u \in Sub(X^{\pm\phi _K})$.
 Notice that  $ Sub(X^{\pm\phi _K})\subseteq {\mathcal
W}_{\Gamma,L}$ which proves  statement 1) of the Claim for $i=1$
and proves the  statement 2) for all $j=1,\ldots ,K.$

Suppose, by induction, that labels of intervals in the cut equation
$\Pi _{L-Ki}$ have form $w^{\phi _{L-Ki}},$ $w\in\bar {\mathcal
W}_{\Gamma ,L}.$ We can rewrite each label in the form $v^{\phi
_{L-K(i+1)}},$ where $v=w^{\phi _K}\in \bar{\mathcal W}_{\Gamma
,L}$. In the $T$-sequence
$$\Pi_{L-Ki}\rightarrow \ldots\rightarrow \Pi_{L-K(i+1)}$$
 each application of the transformation $T$ removes subwords in the alphabet \newline $
X^{\phi _{L-K(i+1)}}$. The argument above shows that the  labels
of the  new intervals in all cut equations  $\Pi _{L-Ki-1)},\ldots
,\Pi _{L-K(i+1)}$ are of the form $v^{\phi _{L-K(i+1)}},$ where
$v\in\bar{\mathcal W}_{\Gamma ,L}.$ Following the proof it is easy
to see
 that the word $v$ does not contain $N$-large powers of $e^{\phi
_{L-K(i+1)}}$ for an elementary period $e$.

\end{proof}

\begin{Claim} \label{2} Let $l\geq 3$, $p_{j-1}\geq (l+2)N_{\Pi}+3$.
The cut equation $T(\Pi)$ is a $\Gamma$-cut equation in  rank
$\leq j-1$ of size $l$.
\end{Claim}
\begin{proof}  The claim follows from the construction of $T(\Pi)$. More
precisely,  we show first that $T(\Pi)$ has a solution relative to
$\beta$. It has been shown in Claim 1 that $T^*(\Pi)$ has a
solution $\alpha^*$ relative to $\beta$.
 This
 proves  condition 2) in the definition of the $\Gamma$-cut
 equation.

 Observe also, that to show
1) it suffices to show that 1.1) in rank $j$ does not hold for
$T^*(\Pi)$. It is not hard to see that it suffices to prove the
required inequalities for $A^\prime$-decompositions (see Lemma
\ref{Claim1.b}).

Let $\delta \in {\mathcal E}^*$.  By the construction
$(A^\prime)^{l+2}$ does not occur in $\mu^\alpha$ for any $\mu \in
M^*$. Therefore  the maximal power of $A^\prime$ that can occur in
$f_{M^*}^*(\delta)^\alpha$ is bounded from above by
$(l+1)|f_{M^*}^*(\delta)|$ which is less then
$(l+1)length(T^*(\Pi))$, as required. Let $t$ be the rank of
$T(\Pi ),\ t\leq j-1.$ It follows from the construction that if
conditions 1.1) and 1.3) for rank $t$ are not satisfied for an
interval in $T(\Pi ),$ then condition 1.2) is satisfied.
\end{proof}

\begin{df} Let $\Pi = ({\mathcal E}, f_X,f_M)$ be a cut equation.
For a positive  integer $n$ by  $k_n(\Pi)$ we denote the number of
intervals $\sigma \in {\mathcal E}$ such that $|f_M(\sigma)| = n$.
The following finite sequence of integers
 $$Comp(\Pi) = (k_2(\Pi), k_3(\Pi), \ldots, k_{length(\Pi)}(\Pi))$$
is called the {\em complexity} of $\Pi$.
 \end{df}

 We well-order complexities of cut equations in the (right) shortlex order:
 if $\Pi$ and $\Pi^\prime$ are two cut equations then
  $Comp(\Pi) <  Comp(\Pi^\prime)$ if and only if $length(\Pi) <
 length(\Pi^\prime)$ or $length(\Pi) =  length(\Pi^\prime)$ and there exists $1 \leqslant i\leqslant length(\Pi)$ such
 that
  $k_j(\Pi) = k_j(\Pi^\prime)$
 for all $j > i$ but $k_i(\Pi) < k_i(\Pi^\prime)$.

Observe that intervals $\sigma \in {\mathcal E}$ with
$|f_M(\sigma)| = 1$ have no input into the complexity of a  cut
equation $\Pi = ({\mathcal E}, f_X,f_M)$. In particular, equations
with $|f_M(\sigma)| = 1$ for every $\sigma \in {\mathcal E}$ have
the minimal possible complexity among equations of a given length.
We will write $Comp(\Pi) = {\bf 0}$ in the case when $k_i(\Pi) =
0$ for every $i = 2, \ldots, length(\Pi)$.

\begin{Claim}\label{3.} Let $\Pi = ({\mathcal E}, f_X,f_M)$. Then the following holds:
 \begin{enumerate}
 \item $length(T(\Pi)) \leqslant length(\Pi)$;
  \item $Comp(T(\Pi)) \leqslant Comp(\Pi)$.

  \end{enumerate}
\end{Claim}

\begin{proof}  By straightforward verification. Indeed, if $\sigma
\in {\mathcal E}_{<j}$
  then  $f_M(\sigma) = f_{M^*}^*(\sigma)$. If $\sigma \in {\mathcal E}_{j}$ and
   $\delta_i \in E_\sigma$ then   $$ f_{M^*}^*(\delta_i) = \mu_{i_1}^* \mu_{i_1+1} \ldots
\mu_{i_1+r(i)}^*,$$ where $\mu_{i_1} \mu_{i_1+1} \ldots
\mu_{i_1+r(i)}$ is a subword  of $ \mu_1 \ldots \mu_n$ and hence
$|f_{M^*}^*(\delta_i)| \leqslant |f_M(\sigma)|$, as required.
\end{proof}

  We need a few definitions related to the sequence (\ref{eq:7.3.15}). Denote by $M_j$ the set
of variables in the equation $\Pi_j$. Variables from $\Pi_L$ are
called {\em initial } variables. A variable $\mu$ from $M_j$ is
called {\em essential} if it occurs in some $f_{M_j}(\sigma)$ with
$|f_{M_j}(\sigma)| \geqslant 2$, such occurrence of $\mu$ is
called {\em essential}. By $n_{\mu,j}$ we denote the total number
of all essential occurrences of $\mu$ in $\Pi_j$. Then
 $$S({\Pi_j})=\sum_{i=2}^{N_{\Pi_j}} ik_i(\Pi_j) = \sum_{\mu \in M_j} n_{\mu,j} $$
is the total number of all essential occurrences of variables from
$M_j$ in $\Pi_j$.

\begin{Claim} \label{claim:S-Pi}
If  $1 \leqslant j \leqslant L$ then $S(\Pi_j) \leqslant
2S(\Pi_L)$.
\end{Claim}
\begin{proof}
 Recall, that every variable
$\mu$ in $M_j$ either belongs to $M_{j+1}$ or it is replaced in
$M_{j+1}$ by the  set $S_\mu$ of new variables (see definition of
the function $f^*_{M^*}$ above). We refer to  variables from
$S_\mu$ as to {\em children} of $\mu$. A given occurrence of $\mu$
in some $f_{M_{j+1}}(\sigma)$, $\sigma \in {\mathcal E}_{j+1}$, is
called a {\em side occurrence} if it is either the first variable
or the last variable (or both) in  $f_{M_{j+1}}(\sigma)$. Now we
formulate several properties of variables from the sequence
(\ref{eq:7.3.15})  which come directly from the construction.  Let
$\mu \in M_j$. Then the following conditions hold:
 \begin{enumerate}
  \item every child of $\mu$ occurs only as a side variable in
  $\Pi_{j+1}$;
   \item every side variable $\mu$ has at most one essential
   child, say $\mu^*$. Moreover, in this event $n_{\mu^*,j+1}
   \leqslant n_{\mu,j}$;
   \item every initial variable $\mu$ has at most two essential
   children, say $\mu_{\rm left}$ and $\mu_{\rm right}$. Moreover, in this
   case $n_{\mu_{\rm left},j+1} + n_{\mu_{\rm right},j+1} \leqslant 2n_\mu$.
\end{enumerate}
Now the claim follows from the properties listed above. Indeed,
every initial variable from $\Pi_j$ doubles, at most, the number
of essential occurrences of its children in the next equation
$\Pi_{j+1}$, but all other variables (not the initial ones) do not
increase this number. \end{proof}

 Denote by $width(\Pi)$ the {\em width} of $\Pi$ which is
   defined as
    $$width(\Pi) = \max_i {k_i(\Pi)}.$$

\begin{Claim} \label{claim:width}
For every $1 \leqslant j \leqslant L$ $width(\Pi_j) \leqslant
2S(\Pi_L)$
\end{Claim}
 \begin{proof}  It follows directly from Claim \ref{claim:S-Pi}.
 \end{proof}

 Denote by $\kappa(\Pi)$ the number of all $(length(\Pi)-1)$-tuples
of non-negative integers which are bounded by $2S(\Pi_L)$.

\begin{Claim}\label {4.} $Comp(\Pi_L) = Comp(\Pi_{\mathcal L})$.
\end{Claim}

\begin{proof}   The complexity $Comp(\Pi_L)$ depends only on the
function $f_M$ in $\Pi_L$. Recall that $\Pi_L = \Pi_{\phi}$ is
obtained from the cut equation $\Pi_{\mathcal L}$ by changing only
the labelling function $f_X$, so $\Pi_{\mathcal L}$ and $\Pi_L$
have the same functions $f_M$, hence the same complexities.
\end{proof}

We say that a $T$-sequence
 has {\em $3K$-stabilization} at $K(r+2)$ , where  $2 \leqslant r
\leqslant L/K$,
  if
$$Comp(\Pi_{K(r+2)})= \ldots = Comp(\Pi_{K(r-1)}).$$
 In this event we denote
 $$K_0 = K(r+2), \ \ \ K_1 = K(r+1), \ \ \ K_2 = Kr, \ \ \ K_3 = K(r-1).$$
For the cut equation $\Pi_{K_1}$ by $M_{\rm veryshort}$ we denote
the subset of variables from $M(\Pi_{K_1})$  which occur unchanged
in $\Pi_{K_2}$ and are short in $\Pi_{K_2}$.

\begin{Claim}\label{5.} For a given $\Gamma$-cut equation $\Pi$ and a
positive integer $r_0 \geqslant 2$  if $L \geqslant Kr_0 +
\kappa(\Pi)4K$ then for some $r \geqslant r_0$ either the sequence
{\rm (\ref{eq:7.3.15})} has $3K$-stabilization at $K(r+2)$ or
$Comp(\Pi_{K(r+1)}) = {0}$.
\end{Claim}
 \begin{proof}  Indeed, the claim follows by the ``pigeon  hole" principle  from
Claims \ref{3.} and \ref{claim:width}  and the fact that there are
not more than $\kappa(\Pi)$ distinct complexities which are less
or equal to $Comp(\Pi)$.
  \end{proof}

  Now we define a special set of solutions of the  equation $S(X) =
1$. Let $L = 4K + \kappa(\Pi)4K$,  $p$ be  a fixed $N$-large tuple
from ${\mathbb N}^{L-4K}$, $q$ be an arbitrary fixed $N$-large
tuple from ${\mathbb N}^{2K}$, and $p^*$ be an arbitrary $N$-large
tuple from ${\mathbb N}^{2K}$. In fact, we need $N$-largeness of
$p^*$ and $q$ only to formally satisfy the conditions of the
claims above. Put
$${\mathcal B}_{p,q,\beta} = \left\{\phi_{L-4K,p}  \phi_{2K,p^*}  \phi_{2K,q}  \beta
\mid p^* \in  {\mathbb N}^{2K},pp^*q\in{\mathcal P} \right\}.$$

It follows from Theorem \ref{cy2} that ${\mathcal B}_{p,q,\beta}$
is a discriminating family of solutions of $S(X) = 1$.

 Denote $\beta_q =\phi_{2K,q} \circ
  \beta$. Then $\beta_q$ is a solution of $S(X) = 1$ in  general
  position and
   $${\mathcal B}_{q,\beta} = \{ \phi_{2K,p^*}  \beta_q
\mid p^* \in  {\mathbb N}^{2K} \}$$ is also a discriminating
family
 by Theorem \ref{cy2}.

Let $${\mathcal B} = \{\psi _{K_1}=\phi_{K(r-2),p'}  \phi_{2K,p^*}
 \phi_{2K,q}  \beta \mid p^* \in  {\mathbb N}^{2K} \},$$ where $p'$ is a
beginning of $p$.
\begin{prop}
\label{Or}

Let $L = 2K + \kappa(\Pi)4K$ and $\phi_L \in {\mathcal
B}_{p,q,\beta}$. Suppose  the $T$-~sequence
  of cut equations {\rm (\ref{eq:7.3.15})} has
$3K$-stabilization at $K(r+2), r\geqslant 2$.
 Then the set of  variables $M$ of the cut equation $\Pi_{K(r+1)}$ can be partitioned
 into three disjoint subsets
 $$M = M_{\rm veryshort} \cup M_{\rm free} \cup M_{\rm useless}$$
 for which the following holds:
\begin{enumerate}
  \item   there exists a finite  system of  equations
$\Delta(M_{\rm veryshort}) = 1$ over $F$ which has a solution in
$F$;
 \item   for every $\mu \in M_{\rm useless}$ there exists a word
$V_\mu \in F[X \cup M_{\rm free}\cup M_{\rm veryshort}]$ which
does not depend on tuples $p^*$ and $q$;
 \item  for every solution $\delta \in {\mathcal B}$, for every map
$\alpha_{\rm free} : M_{\rm free}
 \rightarrow F$, and every solution $\alpha_s: F[M_{\rm veryshort}] \rightarrow F$
  of the system $\Delta(M_{\rm veryshort}) = 1$  the map $\alpha: F[M] \rightarrow F$
  defined by
   \[ \mu^\alpha = \left\{\begin{array}{ll}
   \mu^{\alpha_{\rm free}} &\mbox{ if  $\mu \in M_{\rm free}$;}\\
   \mu^{\alpha_{s}} & \mbox{ if $\mu \in M_{\rm veryshort}$;}\\
     V_\mu(X^\delta, M_{\rm free}^{\alpha_{\rm free}},
   M_{\rm veryshort}^{\alpha_s}) & \mbox{ if $\mu \in M_{\rm useless}$.}
   \end{array}
 \right. \]
  is a group solution of $\Pi_{K(r+1)}$ with respect to $\beta$.
   \end{enumerate}
 \end{prop}
\begin{proof}     Below we describe (in a series of claims
\ref{6.}-\ref{21}) some properties of partitions of intervals of
cut equations from the sequence (\ref{eq:7.3.15}):
  $$ \Pi_{K_1} \stackrel{T}{ \rightarrow} \Pi_{K_1-1}
\stackrel{T}{ \rightarrow} \ldots \stackrel{T}{
\rightarrow}\Pi_{K_2}.$$

Fix an arbitrary integer $s$  such that $K_1 \geqslant s \geqslant
K_2$.

\begin{Claim}\label{6.}  Let $f_M(\sigma) = \mu_1 \cdots \mu_k$ be  a
partition of an interval $\sigma$ of rank $s$ in $\Pi_s$. Then:
 \begin{enumerate}
 \item
the variables  $\mu_2, \ldots, \mu_{k-1}$ are very short;
 \item  either $\mu_1$ or $\mu_k$, or both, are long variables.
 \end{enumerate}
\end{Claim}
 \begin{proof}  Indeed,  if any of the variables $\mu_2, \ldots, \mu_{k-1}$
 is long then the interval $\sigma$ of $\Pi_s$ is replaced in $T(\Pi_s)$ by a set of intervals
 $E_\sigma$ such that $|f_M(\delta)| < |f_M(\sigma)|$ for every $\delta \in
 E_\sigma$.  This implies that complexity of $T(\Pi_s)$ is
smaller than of $\Pi_s$ - contradiction. On the other hand, since
$\sigma$ is a partition of rank $s$ some variables must be long -
hence the result. \end{proof}

Let $f_M(\sigma) = \mu_1 \ldots \mu_k$ be  a partition of an
interval $\sigma$ of rank $s$ in $\Pi_s$. Then the variables
$\mu_1$ and $\mu_k$ are called {\em side variables}.

\begin{Claim} \label{7.} Let $f_M(\sigma) = \mu_1 \ldots \mu_k$ be  a
partition of an interval $\sigma$ of rank $s$ in $\Pi_s$.  Then
this partition will induce a partition of the form
$\mu_1'\mu_2\ldots \mu_{k-1}\mu_k'$ of some interval in rank $s-1$
in $\Pi_{s-1}$ such that  if $\mu_1$ is short in rank $s$ then
$\mu_1' = \mu_1$, if $\mu_1$ is long in $\Pi_s$ then  $\mu_1'$ is
a new variable which does not appear in the previous ranks.
Similar conditions hold for $\mu_k$.
\end{Claim}
 \begin{proof}  Indeed, this follows from the
construction of the transformation $T$.\end{proof}
\begin{Claim}\label{8.} Let $\sigma_1$ and
$\sigma_2$  be two intervals of ranks $s$ in $\Pi_s$ such that
$f_X(\sigma_1) = f_X(\sigma_2)$
 and
 $$f_M(\sigma_1) = \mu_1\nu_2\ldots \nu_k,\ \  f_M(\sigma_2) = \mu_1\lambda_2\ldots \lambda_l.$$
  Then   for any solution $\alpha$ of
$\Pi_s$ one has
$$\nu_k^\alpha = \nu_{k-1}^{-\alpha} \ldots \nu_2^{-\alpha}\lambda_2^{-\alpha}
\ldots \lambda_{l-1}^{-\alpha}\lambda_l^{-\alpha}$$ i.e,
$\nu_k^\alpha$ can be expressed via $\lambda_l^{\alpha}$ and a
product of images of short variables.\end{Claim}

\begin{Claim}\label{9.} Let $f_M(\sigma) = \mu_1 \ldots \mu_k$ be  a
partition of an interval $\sigma$ of rank $s$ in $\Pi_s$. Then for
any $u \in X \cup E(m,n)$ the word  $\mu_2^\alpha \ldots
\mu_{k-1}^\alpha $ does not contain a subword of the type
$c_1(M_u^{\phi_{K_1}})^{\beta}c_2,$ where $c_1,c_2\in C_\beta$,
and $M_u^{\phi_{K_1}}$ is the middle of $u$ with respect to
$\phi_{K_1}$.
\end{Claim}
 \begin{proof}  By Corollary \ref{cy:middles} every word $M_u^{\phi_{K_1}}$
contains a big power (greater than $(l+2)N_{\Pi_s}$) of a period
in rank strictly  greater than $K_2$. Therefore, if
$(M_u^{\phi_{K_1}})^\beta$ occurs in the word $\mu_2^\alpha \ldots
\mu_{k-1}^\alpha $ then some of the variables $\mu_2, \ldots,
\mu_{k-1}$ are not short in some rank greater than $K_2$ -
contradiction. \end{proof}

\begin{Claim}\label{15.} Let $\sigma$ be an interval   in $\Pi_{K_1}$ and  $\phi
_{K_1}={\phi_{K_1,p}}$.  Then $f_X(\sigma) = W_{\sigma}$
 written in the form
  $$W_\sigma = w^{\phi_{K_1}},$$
 and the following
holds:
 \begin{enumerate}
\item [(1)] the word $w$  can be uniquely written as $w=v_1\ldots
v_e,$ where  $v_1,\ldots v_e\in X^{\pm 1}\cup E(m,n)^{\pm 1}$, and
$v_iv_{i+1}\not \in E(m,n)^{\pm 1}$.

 \item [(2)]  $w$ is either a subword of a word from the list in Lemma
\ref{le:xyu} or  there exists $i$ such that
 $v_1\cdots v_i$, $v_{i+1}\cdots v_e$ are
 subwords of words from the list in Lemma \ref{main11}. In addition,
$(v_iv_{i+1})^{\phi _K}=v_i^{\phi _K}\circ v_{i+1}^{\phi _K}.$

 \item  [(3)] if $w$ is  a subword of a word from the list in Lemma
\ref{le:xyu}, then  at most for two indices $i,j$ elements $v_i,
v_j$ belong to $ E(m,n)^{\pm 1},$ and, in this case $j=i+1.$
 \end{enumerate}
\end{Claim}

\begin{proof}  The fact that $W_{\sigma}$ can be written in such a form
follows from Claim \ref{2b} for $r=K.$  Indeed, by Claim \ref{2b},
$W_{\sigma}=w^{\phi _{K_1}},$ where $w\in {\mathcal W}_{\Gamma
,L},$
 therefore it is either a subword of a word
from the list in Lemma \ref{le:xyu} or contains a subword from the
set $Exc$ from  Lemma \ref{main11}. It can contain only one such
subword, because two such subwords of a word from $X^{\pm \phi
_L}$ are separated by big (unbounded) powers of elementary
periods.
 The uniqueness
of $w$ in the first statement follows from the fact that
$\phi_{K_1,p^\prime}$ is an automorphism. Obviously, $w$ does not
depend on $p$.
 Property (3) follows from the comparison of
the set $E(m,n)$
 with the list from Lemma \ref{le:xyu}.\end{proof}

 We say that the decomposition $w=v_1\cdots v_e,$
above  is the  {\em canonical decomposition} of $w$ and
$(v_1\ldots v_s)^{\phi _{K_1}}$ is a canonical decomposition of
$w^{\phi _{K_1}}. $

\begin{Claim}
\label{claim:16} Let $ \Pi_{K_1}=({\mathcal E},f_X, f_M)$ and $\mu
\in M$ be  a long variable (in rank $K_1$) such that $f_M(\delta)
\neq \mu$  for any $\delta \in  {\mathcal E}$. If $\mu$ occurs as
the left variable in $f_M(\sigma)$ for some  $\delta \in {\mathcal
E}$ then it does not occur as the right variable in $f_M(\delta)$
for any $\delta \in {\mathcal E}$ (however,  $\mu^{-1}$ can occur
as the right variable). Similarly, If $\mu$ occurs as the right
variable in $f_M(\sigma)$  then it does not occur as the right
variable in any $f_M(\delta)$.

\end{Claim}
\begin{proof}   Notice, that in this case if $\mu _1$ is not a single
variable, it cannot be a right side variable of $f_M(\bar\sigma )$
for some interval $\bar\sigma$. Indeed, suppose $W_{\bar\sigma}$
ends with $\mu _1$.  If $v_{left}\neq z_i, y_n^{-1}$, $W_{\sigma}$
begins with a big power of some period $A_j^{*\beta}, \ j>K_2$,
therefore $\mu _1$ begins with this big power, and the complexity
of $\bar\sigma$ would decrease when we apply $T$ to the cut
equation in rank $j$. If $v_{left}=z_i$, $\mu _1$ cannot be the
right side variable, because $c_i^N$ can occur only in the
beginning of labels of intervals. If $v_{left}=y_n^{-1}$, then
$W_{\bar\sigma}=\cdots x_n^{-1}\circ y_n^{-1}$, and the complexity
would also decrease when $T$ is applied in rank
 $K_2+m+4n-4$.
\end{proof}

Our next goal is to  transform further the cut equation $\Pi
_{K_1}$ to the form where all  intervals are labelled by elements
$x^{\phi _{K_1}},\ x\in (X\cup E(m,n))^{\pm 1}$. To this end we
introduce several  new  transformations of $\Gamma$-cut equations.

Let $\Pi = ({\mathcal E},f_X, f_M)$ be a $\Gamma$-cut equation in
rank $K_1$ and size $l$ with a solution $\alpha:F[M] \rightarrow
F$ relative to $\beta : F[X] \rightarrow F$. Let $\sigma \in
{\mathcal E}$ and
   $$W_\sigma=(v_1\cdots v_e)^{\phi
_{K_1}}, \ \ \ e \geq 2,$$  be the canonical decomposition of
$W_\sigma $. For $i, 1 \leq i < e,$ put
 $$v_{\sigma,i,left }=v_1\cdots v_i,\ v_{\sigma,i,right}  = v_{i+1}\cdots v_e.$$
 Let, as usual,
 $$f_M(\sigma) = \mu_1 \cdots \mu_k.$$

We start with a {\bf transformation $T_{1,left}$}.  For $\sigma
\in {\mathcal E}$ and $1 \leq i < e$ denote by $\theta$ the
boundary between $v_{\sigma,i,left}^{\phi _{K_1}\beta}$ and
$v_{\sigma,i,right}^{\phi _{K_1}\beta}$  in the reduced form of
the product $v_{\sigma,i,left}^{\phi _{K_1}\beta}
v_{\sigma,i,right}^{\phi _{K_1}\beta}$. Suppose now that there
exist $\sigma$ and $i$ such that the following two conditions
hold:

\begin{enumerate}
 \item [C1)]
$\mu _1^\alpha$ almost   contains the beginning of the word
$v_{\sigma,i,left}^{\phi _{K_1}\beta}$ till the boundary $\theta$
(up to a very short end of it), i.e.,  there are elements $u_1,
u_2, u_3, u_4 \in F$ such that $v_{\sigma,i,left}^{\phi
_{K_1}\beta} = u_1\circ u_2 \circ u_3$, $ v_{i+1}^{\phi
_{K_1}\beta} = u_3^{-1}\circ u_4$, $u_1u_2u_4 = u_1\circ u_2 \circ
u_4$, and $\mu_1^\alpha $ begins with $u_1$, and $u_2$ is very
short (does not contain $A^{\pm l} _{K_2}$) or trivial.

\item [ C2)] the boundary $\theta$  does not lie inside
$\mu_1^{\alpha}$.
\end{enumerate}

In this event the transformation $T_{1,left}$ is applicable to
$\Pi$ as  described below. We consider three cases with respect to
the location of $\theta$ on $f_M(\sigma)$.

\begin{figure}[here]

\vspace{3ex} \centering{\mbox{\psfig{figure=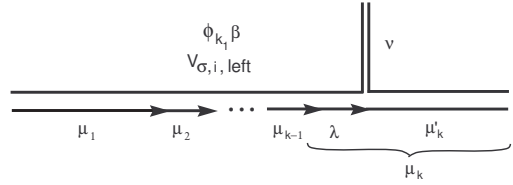,height=3in}}}

\vspace{-25ex} \caption{ T2, Case 1)} \label{T2}

\end{figure}
\begin{enumerate}

\item [Case 1)] $\theta$ is inside $\mu_k^\alpha$ (see Fig.
\ref{T2}). In this case we perform the following:

 a) Replace  the interval
$\sigma$ by two new intervals $\sigma_1, \sigma _2$ with the
labels $v_{\sigma,i,left}^{\phi _{K_1}},\ v_{\sigma,i,right}^{\phi
_{K_1}}$;

b) Put $f_M(\sigma _1)=\mu _1\ldots \mu _{k-1}\lambda\nu,$
$f_M(\sigma _2)=\nu ^{-1}\mu' _k,$  where $\lambda$ is a new very
short variable, $\nu$ is a new variable.

c) Replace everywhere $\mu _k$ by $\lambda\mu' _{k}$. This
finishes the description of the cut equation  $T_{1,left}(\Pi)$.

d) Define a solution $\alpha^\ast$ (with respect to $\beta$)
  of  $T_{1,left}(\Pi)$ in the natural way.
  Namely, $\alpha^\ast(\mu) = \alpha(\mu)$ for all variables $\mu$ which
  came unchanged from $\Pi$. The values
  $\lambda^{\alpha^\ast},  {\mu'}_{k}^{\alpha^\ast}$, $\nu ^{\alpha^\ast}$ are defined
  in the natural way, that is
${\mu'}_{k}^{\alpha^\ast}$ is the whole end part of $\mu_k^\alpha$
after the boundary $\theta$,
 $(\nu ^{-1}\mu'
_k)^{\alpha^\ast}=v_{\sigma,i,{\rm right}}^{\phi_{K_1}\beta}$,
$\lambda^{\alpha^\ast} = \mu_k^\alpha ({\mu'}_{k}^\alpha)^{-1}$.

\item [Case 2)] $\theta$ is on the  boundary between  $\mu
_j^{\alpha}$ and $\mu _{j+1}^\alpha$ for some $j$. In this case we
perform the following:

a) We split the interval $\sigma$ into two new intervals $\sigma
_1$ and $\sigma _2$ with labels $v_{\sigma,i, left}^{\phi_{K_1}}$
and $v_{\sigma,i,right}^{\phi _{K_1}}$.

b) We introduce a new variable $\lambda$ and put $f_M(\sigma_1) =
\mu_1\ldots \mu_j\lambda$, $f_M(\sigma _2) =
\lambda^{-1}\mu_{j+1}\ldots \mu_k$.

 c) Define $\lambda^{\alpha^\ast}$ naturally.

\item [Case 3)] The boundary $\theta$  is contained inside
$\mu_i^\alpha$ for some $i (2 \leq i \leq r-1)$. In this case we
do the following:

 a) We split the interval $\sigma$ into two intervals $\sigma _1$ and
$\sigma _2$ with labels $v_{left}^{\phi_{K_1}}$ and
$v_{\sigma,i,right}^{\phi _{K_1}}$, respectively.

 b) Then we introduce three new variables $\mu _j', \mu _j'',
\lambda$, where  $\mu _j', \mu _j''$ are  ``very short", and add
equation $\mu _j = \mu _j'\mu _j''$ to the system $\Delta_{\rm
veryshort}$.

 c) We define
$f_M(\sigma _1) = \mu _1\cdots \mu_j'\lambda $, $f_M(\sigma _2) =
\lambda ^{-1}\mu _j''\mu _{i+1}\cdots \mu _k$.

 d) Define values of
$\alpha^\ast$ on the new variables naturally. Namely, put
$\lambda^{\alpha^\ast} $ to be equal to the terminal segment of
$v_{left}^{\phi_{K_1}\beta}$ that cancels in the product
$v_{left}^{\phi_{K_1}\beta}v_{\sigma,i,{\rm right}}^{\phi
_{K_1}\beta}$. Now the values $\mu _j'^{\alpha^\ast}$ and $\mu
_j''^{\alpha^\ast}$ are defined to satisfy the equalities
 $$f_X(\sigma_1)^\beta = f_M(\sigma_1)^{\alpha^\ast}, f_X(\sigma_2)^\beta = f_M(\sigma_2)^{\alpha^\ast}.$$

\end{enumerate}

We described the transformation $T_{1,left}$. The transformation
$T_{1, right}$ is defined similarly. We denote both of them by
$T_1$.

Now we describe a {\bf transformation $T_{2,left}$}.

Suppose again that a cut equation $\Pi$ satisfies C1). Assume in
addition that for  these $\sigma$ and $i$ the following condition
holds:
\begin{enumerate}
  \item [C3)] the boundary $\theta$ lies inside $\mu_1^\alpha$.
  \end{enumerate}
Assume also that  one of the following three conditions holds:
 \begin{enumerate}
  \item [C4)] there are no intervals $\delta \neq \sigma$ in $\Pi$ such that
  $f_M(\delta)$ begins with $\mu_1$ or ends on $\mu_1^{-1}$;
  \item  [C5)] $v_{\sigma,i,left} \neq x_n$ (i.e., either  $i > 1$ or $i = 1$ but $v_1 \neq
 x_n$) and for every $\delta \in \mathcal{E}$ in $\Pi$ if $f_M(\delta)$ begins
with $\mu_1$ (or ends on $\mu_1^{-1}$) then  the canonical
decomposition of $f_X(\delta)$ begins with
$v_{\sigma,i,left}^{\phi_{K_1}}$ (ends with $v_{\sigma,i,left}^{-
\phi_{K_1}}$);
 \item  [C6)] $v_{\sigma,i,left} =  x_n$ ($i = 1$ and $v_1 = x_n$)
 and for every $\delta \in \mathcal{E}$ if $f_M(\delta)$ begins
with $\mu_1$ (ends with $\mu_i^{-1}$) then  the canonical
decomposition of $f_X(\delta)$ begins with $x_n^{\phi_{K_1}}$ or
with $y_n^{\phi_{K_1}}$ (ends with $x_n^{-\phi_{K_1}}$ or
$y_n^{-\phi_{K_1}}$).
 \end{enumerate}

In this event the transformation $T_{2,left}$ is applicable to
$\Pi$ as described below.

\begin{enumerate}
 \item [Case C4)] Suppose the condition C4) holds. In this case we
 do the following.

a)  Replace $\sigma$ by two new intervals $\sigma_1, \sigma _2$
with the labels $v_{\sigma,i,left}^{\phi _{K_1}},\
v_{\sigma,i,right}^{\phi _{K_1}}$;

b) Replace $\mu_1$ with  two new variables  $\mu' _1, \mu_1''$ and
put $f_M(\sigma_1) = \mu_1'$, $f_M(\sigma_2) = \mu_1''\mu_2 \ldots
\mu_k.$

c) Define $(\mu_1')^{\alpha^\ast}$ and $(\mu_1'')^{\alpha^\ast}$
such that $f_M(\sigma_1)^{\alpha^\ast} = v_{\sigma,i,left}^{\phi
_{K_1}\beta}$ and $f_M(\sigma_2)^{\alpha^\ast}
\newline =v_{\sigma,i,right}^{\phi _{K_1}\beta}$.

 \item [Case C5)]  Suppose $v_{left}\neq x_n$. Then do the
 following.

 a) Transform $\sigma$ as described in C4).

 b) If for some interval $\delta \neq \sigma$ the word
$f_M(\delta)$ begins with $\mu_1$ then replace $\mu _1$ in
$f_M(\delta)$ by the variable $\mu''_1$  and replace $f_X(\delta
)$ by $v_{\sigma,i,left}^{-\phi_{K_1}}f_X(\delta).$ Similarly
transform intervals $\delta$ that end with $\mu_1^{-1}$.

\item [Case C6)] Suppose $v_{left}=x_n$. Then do the following.

 a)  Transform $\sigma$  as described in C4).

 b) If for some $\delta$ the word $f_M(\delta)$ begins with $\mu_1$  and $f_X(\delta)$ does not begin
with $y_n$ then  transform $\delta$ as described in Case C5).

 c)  Leave all other intervals  unchanged.

\end{enumerate}
We described the transformation $T_{2,left}$. The transformation
$T_{2, right}$ is defined similarly. We denote both of them by
$T_2$.

Suppose now that $\Pi = \Pi_{K_1}$.
 Observe that  the transformations
$T_1$ and $T_2$ preserve the properties described in Claims
\ref{3.}--\ref{4.} above.
 Moreover, for the
homomorphism $\beta : F[X] \rightarrow F$ we have constructed  a
solution $\alpha^* :F[M] \rightarrow F$ of $T_n(\Pi _{K_1})$  ($n
= 2,3$) such that the initial solution $\alpha$ can be
reconstructed from $\alpha^\ast$ and the equations $\Pi$ and
$T_n(\Pi)$.
 Notice also that the
length of the elements $W_{\sigma'}$ corresponding to new
intervals $\sigma$ are shorter than the length of the words
$W_\sigma$ of the original intervals $\sigma$ from which $\sigma'$
were obtained. Notice also that the transformations $T_1, T_2$
preserves the property of intervals formulated in the Claim
\ref{6.}.

\begin{Claim} \label{cl:T2-T3}

Let $\Pi$ be a cut equation which satisfies the conclusion of the
Claim~\ref{6.}. Suppose  $\sigma$ is an interval in $\Pi$ such
that $W_{\sigma}$ satisfies the conclusion of Claim~\ref{15.}. If
for some $i$
 $$(v_1 \ldots v_e)^{\phi_{K}} = (v_1 \ldots v_i)^{\phi_{K}} \circ
(v_{i+1} \ldots
 v_e)^{\phi_{K}}$$
  then either $T_1$ or $T_2$ is applicable to given
  $\sigma$ and $i$.
\end{Claim}

\begin{proof} By Corollary \ref{R} the automorphism $\phi_{K_1}$ satisfies
the
 Nielsen property with respect to $\bar {\mathcal W}_{\Gamma}$ with
 exceptions $E(m,n)$. By Corollary 12, equality
 $$(v_1 \ldots v_e)^{\phi_{K}} = (v_1 \ldots v_i)^{\phi_{K}} \circ
(v_{i+1} \ldots
 v_e)^{\phi_{K}}$$
implies that the element that is cancelled between $(v_1 \ldots
v_i)^{\phi_{K}\beta}$ and   $(v_{i+1} \ldots
 v_e)^{\phi_{K}\beta}$ is short in rank $K_2$. Therefore either $\mu
_1^{\alpha}$ almost
 contains
$(v_1 \ldots v_i)^{\phi_{K}\beta}$
 or $\mu _k^{\alpha}$ almost contains $(v_{i+1} \ldots
 v_e)^{\phi_{K}\beta}$. Suppose $\mu _1^{\alpha}$ almost contains
$(v_1 \ldots v_i)^{\phi_{K}\beta}$. Either we can apply
$T_{1,left}$, or the boundary $\theta$ belongs to $\mu
_1^{\alpha}$. One can verify using formulas from Lemmas
\ref{le:7.1.zforms}-\ref{le:7.1.xiforms} and \ref{main} directly
that in this case one of the conditions $C4)-C6)$ is satisfied,
and, therefore $T_{2,left}$ can be applied.
\end{proof}

\begin{lm} \label{n12} Given a cut equation $\Pi_{K_1}$ one can effectively
find a finite sequence of transformations $Q_1, \ldots, Q_s$ where
$Q_i \in \{T_1, T_2\}$ such that for every interval $\sigma$  of
the cut equation $\Pi_{K_1}^\prime = Q_s \ldots Q_1 (\Pi_{K_1})$
the label $f_X(\sigma)$ is of
 the form $u^{\phi _{K_1}}$, where $u \in X^{\pm 1}\cup E(m,n)$.

 Moreover, there exists an infinite subset $P'$ of the solution set $P$ of
$\Pi _{K_1}$ such that this sequence
 is the same for any solution in $P'$.
\end{lm}
\begin{proof} Let $\sigma$ be an interval of the equation $\Pi_{K_1}$. By
Claim \ref{15.} the word $W_\sigma$ can be uniquely written in the
canonical decomposition form
 $$W_\sigma = w^{\phi_{K_1}}=(v_1 \ldots v_e)^{\phi_{K_1}},$$
 so that the conditions 1), 2), 3) of Claim \ref{15.} are satisfied.

It follows from the construction of $\Pi _{K_1}$ that either $w$
is a subword of a word between two elementary squares $x\neq c_i$
or begins and (or) ends with some power $\geq 2$ of an elementary
period. If $u$ is an elementary period,
 $u^{2\phi _K}=u^{\phi _K}\circ u^{\phi _K}$, except $u=x_n$,
 when the middle is exhibited in the proof of Lemma \ref{main}.
 Therefore, by Claim \ref{cl:T2-T3}, we can apply $T_1$ and $T_2$ and cut
$\sigma$ into
 subintervals $\sigma _i$
such that for any $i$ $f_X(\sigma _i)$
 does not contain powers $\geq 2$ of elementary periods.
All possible values of $u^{\phi _K}$ for
 $u\in E(m,n)^{\pm 1}$ are shown in the proof of Lemma
 \ref{main}. Applying $T_1$ and $T_2$ as in Claim \ref{cl:T2-T3}
we can split intervals (and their
 labels) into parts with labels of the form $x^{\phi _{K_1}},\
 x\in (X\cup E(m,n)),$
except for the following cases:

1. $w=uv$, where $u$ is $x_i^2, i<n,\ v\in E_{m,n},$ and $v$ has
at
 least three letters,

2.
$w=x_{n-2}^2y_{n-2}x_{n-1}^{-1}x_nx_{n-1}y_{n-2}^{-1}x_{n-2}^2,$

3. $w=x_{n-1}^2y_{n-1}x_n^{-1}x_{n-1}y_{n-2}^{-1}x_{n-2}^{-2},$

4. $y_{r-1}x_r^{-1}y_r^{-1},\ r<n,$

5. $w=uv$, where $u=(c_1^{z_1}c_2^{z_2})^2,\ v\in E(m,n),$ and $v$
is one of the following: $v=\prod _{t=1}^mc_t^{z_t}x_1^{\pm 1},$
$v=\prod _{t=1}^mc_t^{z_t}x_1^{\pm 1}\prod _{t=m}^1c_t^{-z_t}$, $v
= \prod_{t=1}^mc_t^{z_t}x_1\prod
_{t=m}^1c_t^{-z_t}(c_1^{z_1}c_2^{z_2})^{-2},$

6. $w=uv$, where $u=(c_1^{z_1}c_2^{z_2})^2,\ v\in E(m,n),$ and $v$
is one of the following: $v=\prod _{t=1}^mc_t^{z_t}x_1^{-1
}x_2^{-1}$ or $v=\prod _{t=1}^mc_t^{z_t}x_1^{-1}y_1^{-1}.$

7. $w=z_iv.$

Consider the first case. If $f_M(\sigma)=\mu _1\cdots\mu _k,$ and
$\mu _1^{\alpha}$ almost contains
 $$x_i^{\phi
_{K_1}}(A_{m+4i+K_2}^*)^{-p_{m+4i+K_2}+1}x_{i+1}^{\phi
_{K_2}\beta}$$ (which is a non-cancelled initial peace of
$x_i^{2\phi _{K_1}\beta}$ up to a very short part of it), then
either $T_{1,\rm left}$ or $T_{2,\rm left}$ is applicable and we
split $\sigma$ into two intervals $\sigma _1$ and $\sigma _2$ with
labels $x_i^{2\phi _{K_1}}$ and $v^{\phi _{K_1}}$.

Suppose $\mu _1^{\alpha}$ does not contain $x_i^{\phi
_{K_1}}(A_{m+4i+K_2}^*)^{-p_{m+4i+K_2}+1}x_{i+1}^{\phi
_{K_2}\beta}$ up to a very short part. Then $\mu _k^{\alpha}$
contains the non-cancelled left end $E$ of $v^{\phi _{K+1}\beta},$
and $\mu _k^{\alpha}E^{-1}$ is not very short. In this case
$T_{2,\rm right}$ is applicable.

We can similarly consider all Cases 2-6.

Case 7. Letter $z_i$  can  appear only
 in the beginning of $w$ (if  $z_i^{-1}$ appears at the
 end of $w$, we can replace $w$ by $w^{-1}$)
 If $w=z_it_1\cdots t_s$ is the canonical decomposition, then
 $t_k=c_j^{\pm z_j}$ for each $k$. If $\mu _1^{\alpha}$ is longer
 than
the non-cancelled part of  $(c_i^{p}z_i)^{\beta}$, or
 the difference between $\mu _1^{\alpha}$
and $(c_i^{p}z_i)^{\beta}$ is very short, we can split $\sigma$
into two parts, $\sigma _1$ with label $f_X(\sigma _1)=z^{\phi
_{K_1}}$ and $\sigma _2$ with label $f_X(\sigma _2)=(t_1\ldots
t_s)^{\phi _{K_1}}.$

If the difference between $\mu _1^{\alpha}$ and
$(c_i^{p}z_i)^{\beta}$ is not very short, and $\mu _1^{\alpha}$ is
shorter than the non-cancelled part of  $(c_i^{p}z_i)^{\beta}$,
then there is no interval $\delta$ with $f(\delta )\neq f(\sigma
)$ such that $f_M(\delta )$ and $f_M(\sigma )$ end with $\mu _k,$
and we can split $\sigma$ into two parts using $T_1$, $T_2$  and
splitting $\mu _k.$

We have considered all possible cases. \end{proof}

Denote the resulting cut equation by $\Pi '_{K_1}.$

\begin{cy} The intervals of $\Pi '_{K_1}$ are labelled by elements
$u^{\phi _{K_1}},$ where

 for $n=1$
$$u\in\{z_i,\ x_i,\ y_i,\  \prod c_s^{z_s},\  x_1\prod
_{t=m}^{1}c_t^{-z_t},\}
$$

for $n=2$ \begin{multline*}u\in\{z_i,\ x_i,\ y_i,\  \prod
c_s^{z_s},\ y_1x_1\prod _{t=m}^{1}c_t^{-z_t},\  y_1x_1,\ \prod
_{t=1}^{m}c_t^{z_t}x_1\prod _{t=m}^{1}c_t^{-z_t}\ , \prod
_{t=1}^{m}c_t^{z_t}x_1^{-1}x_2^{\pm 1},\\
\prod _{t=1}^{m}c_t^{z_t}x_1^{-1}x_2x_1,\  \prod
_{t=1}^{m}c_t^{z_t}x_1^{-1}x_2x_1\prod _{t=m}^{1}c_t^{-z_t},\
x_1^{-1}x_2x_1\prod _{t=m}^{1}c_t^{-z_t},\  x_2x_1\prod
_{t=m}^{1}c_t^{-z_t},\\ x_1^{-1}x_2,\ x_2x_1\},\end{multline*} and
for $n\geq 3$, \begin{multline*}u\in \{z_i,\ x_i,\ y_i,\
c_s^{z_s},\ y_1x_1\prod _{t=m}^{3}c_t^{-z_t},\ \prod
_{t=1}^{m}c_t^{z_t}x_1^{-1}x_2^{-1},\ y_rx_r, \ x_1\prod
_{t=m}^{1}c_t^{-z_t},\\ y_{r-2}x_{r-1}^{-1}x_r^{-1},\
y_{r-2}x_{r-1}^{-1},\ x_{r-1}^{-1}x_r^{-1},\ y_{r-1}x_r^{-1},\
r<n,\ x_{n-1}^{-1}x_nx_{n-1},\\
y_{n-2}x_{n-1}^{-1}x_nx_{n-1}y_{n-2}^{-1},\
y_{n-2}x_{n-1}^{-1}x_n^{\pm 1},\  x_{n-1}^{-1}x_n,\ x_nx_{n-1},\\
y_{n-1}x_n^{-1}x_{n-1}y_{n-2}^{-1},\ y_{n-1}x_n^{-1},
y_{r-1}x_r^{-1}y_r^{-1}\}.\end{multline*}
\end{cy}

\begin{proof} Direct inspection from Lemma \ref{n12}.
 \end{proof}

 Below we
suppose $n>0$. We still want to reduce the variety of possible
labels of intervals in $\Pi '_{K_1}$. We cannot apply $T_1$, $T_2$
to some of the intervals labelled by $x^{\phi _{K_1}}$, $x\in
X\cup E(m,n)$, because there are some cases when $x^{\phi _{K_1}}$
is completely cancelled in  $y^{\phi _{K_1}}$, $x,y\in (X\cup
E(m,n))^{\pm 1}.$

We will change the basis of $F(X\cup C_S)$, and then apply
transformations $T_1$, $T_2$ to the labels written in the new
basis. Replace, first, the basis $(X\cup C_S)$ by a new basis
$\bar X\cup C_S$ obtained by replacing each variable $x_s$ by
$u_s=x_sy_{s-1}^{-1}$ for $s>1$, and replacing $x_1$ by
$u_1=x_1c_m^{-z_m}$.

Consider  the case $n\geqslant 3.$ Then the labels of the
intervals will be rewritten as $u^{\phi _{K_1}}$, where
\begin{multline*}
u\in\{z_i,\ u_iy_{i-1},\ y_i,\ \prod _sc_s^{z_s},\
 y_1u_1\prod _{j=n-1}^1c_j^{-z_j}, \  u_1^{-1}y_1^{-1}u_2^{-1},\\ y_ru_ry_{r-1}, \ u_r,\
 u_{r-1}^{-1}y_{r-1}^{-1}u_r^{-1},\  u_ry_{r-1}u_{r-1}y_{r-2},\
 u_2y_1u_1\prod _{j=n-1}^1c_j^{-z_j},\ r<n;\\
y_{n-2}^{-1}u_{n-1}^{-1}u_ny_{n-1}u_{n-1}y_{n-2}, \
 u_{n-1}^{-1}u_ny_{n-1}u_{n-1},\ u_{n-1}^{-1}u_ny_{n-1}, \\
 u_{n-1}^{-1}y_{n-1}^{-1}u_n^{-1},
\left. y_{n-2}^{-1}u_{n-1}^{-1}u_ny_{n-1},  \ u_ny_{n-1}u_{n-1}y_{n-2},\
 u_n^{-1}u_{n-1},\ u_n\right\}.
 \end{multline*}

In the cases $n=1,2$ some of the labels above do not appear, some
coincide. Notice, that $x_n^{\phi _K}=u_n^{\phi _K}\circ
y_{n-1}^{\phi _K},$ and that the first letter of $y_{n-1}^{\phi
_K}$ is not cancelled in the products
$(y_{n-1}x_{n-1}y_{n-2}^{-1})^{\phi _K},$ $(y_{n-1}x_{n-1})^{\phi
_K}$ (see Lemma \ref{le:7.1.x1formsmneq0}). Therefore,  applying
transformations similar to $T_1$ and $T_2$ to the cut equation
$\Pi'_{K_1}$ with labels written in the basis $\bar X$, we can
split all the intervals with labels containing $(u_ny_{n-1})^{\phi
_{K_1}}$ into two parts and obtain a cut equation with the same
properties and intervals labelled by $u^{\phi _{K_1}},$ where
\begin{multline*}
u\in\{z_i,\ u_iy_{i-1},\ y_i,\ \prod _sc_s^{z_s},\
 y_1u_1\prod _{j=n-1}^1c_j^{-z_j}, \  u_1^{-1}y_1^{-1}u_2^{-1},\\ y_ru_ry_{r-1}, u_r,\
 u_{r-1}^{-1}y_{r-1}^{-1}u_r^{-1}, u_ry_{r-1}u_{r-1}y_{r-2},\
 u_2y_1u_1\prod _{j=n-1}^1c_j^{-z_j},\ r<n;\\
 y_{n-2}^{-1}u_{n-1}^{-1}u_n,\ y_{n-1}u_{n-1}y_{n-2},\
 u_{n-1}^{-1}u_n,\ y_{n-1}u_{n-1},
\ u_n\}.
\end{multline*}
Consider for $i<n$ the expression for

\begin{multline*}(y_iu_i)^{\phi
_K} =A_{m+4i}^{-p_{m+4i}+1} \circ x_{i+1}\circ A
_{m+4i-4}^{-p_{m+4i-4}}\\ \circ x^{p_{m+4i-3}}\circ y_i\circ A
_{m+4i-2}^{p_{m+4i-2}-1}\circ x_i \circ \tilde y_{i-1}^{-1}.
\end{multline*}

Formula 3.a) from Lemma \ref{main} shows that $u_i^{\phi _K}$ is
completely cancelled in  the product $y_i^{\phi _K}u_i^{\phi _K}$.
This implies that $y_i^{\phi _{K}}=v_i^{\phi _{K}}\circ u_i^{-\phi
_{K}}$.

Consider also the product
\begin{multline*}
y_{i-1}^{-\phi _K}u_i^{-\phi _K}\\=
 \left({\bf A_{m+4i-4}^{-p_{m+4i-4}+1}\circ x_i\circ
\tilde y_{i-1}} \circ x_i^{-1}A_{m+4i-4}^{p_{m+4i-4}-1}\right)\\
 \left ( A_{m+4i-4}^{-p_{m+4i-4}+1}x_i\circ {\bf
(x_{i}^{p_{m+4i-3}}y_{i-1}\ldots
*)^{p_{m+4i-1}-1}x_i^{p_{m+4i-3}}y_ix_{i+1}^{-1}A_{m+4i}^{p_{m+4i}-1}}\right
),
\end{multline*}

 where the non-cancelled part is made bold.

Notice that $(y_{r-1}u_{r-1})^{\phi _{K}}y_{r-2}^{\phi
_{K}}=(y_{r-1}u_{r-1})^{\phi _{K}}\circ y_{r-2}^{\phi _{K}},$
because $u_{r-1}^{\phi _{K}}$ is completely cancelled in the
product $y_i^{\phi _{K}}u_i^{\phi _{K}}$.

Therefore, we can again apply the transformations similar to $T_1$
and $T_2$ and split the intervals into the ones with labels
$u^{\phi _{K_1}}$, where \begin{multline*}u\in\{z_s,\  y_i,\ u_i,\
\prod _s c_s^{z_s},\ y_ru_r,\  y_1u_1\prod _{j=m-1}^1c_j^{-z_j}, \
u_{n-1}^{-1}u_n=\bar u_n,\\ 1\leqslant i\leqslant n, \ 1\leqslant
j\leqslant m,\ 1\leqslant r <n\}.\end{multline*}

We change the basis again replacing $y_r, 1<r<n$ by a new variable
$v_r=y_ru_r$, and replacing  $y_1u_1\prod _{j=m-1}^1c_j^{-z_j}$ by
$v_1$. Then $y_r^{\phi _{K}}=v_r^{\phi _{K}}\circ u_r^{-\phi
_{K}},$ and $y_1^{\phi _K}=v_1^{\phi _K}\circ c_1^{z_1^{\phi
_K}}\circ c_{m-1}^{z_{m-1}^{\phi _K}}\circ u_1^{-\phi _K}$ (if
$n\neq 1$). Formula 2.c) shows that $u_n^{\phi _{K}}=u_{n-1}^{\phi
_{K}}\circ (u_{n-1}^{-1}u_n)^{\phi _{K}}.$

Apply transformations similar to $T_1$ and $T_2$ to the intervals
with labels written in the new basis $$\hat X=\{z_j,\ u_i,\ v_i,\
y_n,\  \bar u_n=u_{n-1}u_n,\ 1\leqslant j\leqslant m,\ 1\leqslant
i<n,\ j\leqslant m \},$$ and obtain intervals with labels $u^{\phi
_{K_1}},$ where
$$u\in \hat X\cup \{c_m^{z_m}\}.$$
Denote the resulting cut equation by ${\bar \Pi} _{K_1} = (\bar{
\mathcal E}, f_{\bar X},f_{\bar M})$. Let  $\alpha$ be the
corresponding  solution of ${\bar \Pi} _{K_1}$ with respect to
$\beta .$

Denote by $\bar M_{\rm side}$ the set of long variables in $\bar \Pi
_{K_1}$, then $ \bar M =  \bar M_{\rm veryshort}  \cup {\bar
M}_{\rm side}$.

Define a binary relation $\sim_{\rm left}$ on $\bar M_{\rm side}^{\pm 1}$
as follows.  For $\mu_1, \mu'_1 \in \bar M_{\rm side}^{\pm 1}$ put
$\mu_1 \sim_{\rm left} \mu_1^\prime$ if and only if there exist two
intervals $\sigma, \sigma'  \in \bar{E}$ with $f_{\bar X}(\sigma
)=f_{\bar X}(\sigma')$ such that
 $$f_{\bar M}(\sigma)=\mu _1\mu_2\cdots \mu_{r}, \ \ \ f_{\bar M}(\sigma ')=\mu _1'\mu _2'\cdots
 \mu_{r'}'$$
and either $\mu_r = \mu'_{r'}$ or $\mu_r,  \mu'_{r'} \in  M_{\rm
veryshort}.$ Observe that if $\mu_1 \sim_{\rm left} \mu_1^\prime$ then
 $$\mu_1 = \mu_1^\prime \lambda_1 \cdots \lambda_t$$
for some $\lambda_1, \ldots, \lambda_t \in M_{\rm veryshort}^{\pm
1}.$ Notice, that $\mu \sim_{\rm left} \mu$.

Similarly, we define a binary relation $\sim_{\rm right}$ on $\bar
M_{\rm side}^{\pm 1}$.  For $\mu_r, \mu'_{r'} \in \bar M_{\rm side}^{\pm
1}$ put $\mu_r \sim_{\rm right} \mu_{r'}^\prime$ if and only if there
exist two intervals $\sigma, \sigma'  \in \bar{E}$ with $f_{\bar
X}(\sigma )=f_{\bar X}(\sigma')$ such that
 $$f_{\bar M}(\sigma)=\mu _1\mu_2\cdots \mu_{r}, \ \ \ f_{\bar M}(\sigma ')=\mu _1'\mu _2'\cdots
 \mu_{r'}'$$
and either $\mu_1 = \mu'_1$ or $\mu_1,  \mu'_1 \in  M_{\rm
veryshort}.$  Again, if  $\mu_r \sim_{\rm right} \mu_{r'}^\prime$ then
 $$\mu_r =  \lambda_1 \ldots \lambda_t\mu_{r'}^\prime$$
for some $\lambda_1, \ldots, \lambda_t \in M_{\rm veryshort}^{\pm
1}.$

  Denote by $\sim$ the transitive closure of
   $$\{(\mu,\mu') \mid \mu \sim_{\rm left} \mu'\} \cup \{(\mu,\mu') \mid \mu
   \sim_{\rm right} \mu'\} \cup \{(\mu,\mu^{-1}) \mid \mu \in \bar M_{\rm side}^{\pm 1} \}.$$
   Clearly, $\sim$ is an equivalence relation on $\bar
  M_{\rm side}^{\pm 1}$.
  Moreover, $\mu \sim \mu'$ if and only if there exists a sequence
  of variables
   \begin{equation}
   \label{eq:seq-mu}
   \mu = \mu_0, \mu_1, \ldots, \mu_k = \mu'
    \end{equation}
 from $\bar M_{\rm side}^{\pm 1}$  such that either $\mu_{i-1} = \mu_i$, or $\mu_{i-1} = \mu_i^{-1}$,
   or  $\mu_{i-1} \sim_{\rm left} \mu_i$, or $\mu_{i-1} \sim_{\rm right} \mu_i$
    for $i= 1, \ldots, k.$ Observe that if $\mu_{i-1}$ and $\mu_i$ from (\ref{eq:seq-mu})
    are side variables of ``different sides" (one is on  the left, and the other is on the right)
     then $\mu_i = \mu_{i-1}^{-1}$. This implies that replacing in
     the sequence (\ref{eq:seq-mu}) some elements $\mu_i$ with
     their inverses one can get a new sequence
   \begin{equation}
   \label{eq:seq-mu-2}
   \mu = \nu_0, \nu_1, \ldots, \nu_k = (\mu')^\varepsilon
    \end{equation}
 for some $\varepsilon \in \{1,-1\}$ where  $\nu_{i-1} \sim \nu_i$ and all the variables $\nu_i$ are
 of the same side.   It follows  that if $\mu$ is a left-side variable and $\mu \sim \mu'$ then
     \begin{equation}
     \label{eq:sim-mu}
     (\mu')^\varepsilon =  \mu \lambda_1 \cdots \lambda_t
      \end{equation}
for some $\lambda_j \in M_{\rm veryshort}^{\pm 1}.$

It follows from (\ref{eq:sim-mu}) that for a variable $\nu \in
\bar M_{\rm side}^{\pm 1}$ all variables from the equivalence class
$[\nu]$ of $\nu$ can be expressed via $\nu$ and very short
variables from $M_{\rm veryshort}$. So if we fix a system of
representatives $R$ of $\bar M_{\rm side}^{\pm 1}$ relative to $\sim$
then all other variables from $\bar M_{\rm side}$ can be expressed as
in (\ref{eq:sim-mu}) via variables from $R$ and very short
variables.

This allows one to introduce a new transformation $T_3$ of cut
equations. Namely, if a set of representatives $R$ is fixed then
using (\ref{eq:sim-mu}) replace every variable $\nu$ in every word
 $f_M(\sigma)$ of a cut equation $\Pi$ by its expression via the
 corresponding representative variable from $R$ and a product of
 very short variables.

Now we repeatedly apply the transformation $T_3$  till the
equivalence relations $\sim_{\rm left}$ and $\sim_{\rm right}$
become trivial.  This process stops in finitely many steps since
the non-trivial relations decrease the number of side variables.

Denote the resulting equation again by ${\bar \Pi} _{K_1}$.

Now we introduce an equivalence relation on partitions of  ${\bar
\Pi} _{K_1}$. Two partitions $f_M(\sigma )$ and $f_M(\delta )$ are
equivalent ($f_M(\sigma ) \sim f_M(\delta )$) if $f_X(\sigma
)=f_X(\delta )$ and either the left side variables or the right
side variables of $f_M(\sigma )$ and $f_M(\delta )$ are
equivalent. Observe, that $f_X(\sigma )=f_X(\delta )$ implies
$f_M(\sigma )^\alpha=f_M(\delta )^\alpha$, so in this case the
partitions $f_M(\sigma )$ and $f_M(\delta )$ cannot begin with
$\mu$ and $\mu^{-1}$ correspondingly. It follows that if
$f_M(\sigma ) \sim f_M(\delta )$ then the left side variables and,
correspondingly,  the right side variables of $f_M(\sigma )$ and
$f_M(\delta )$ (if they exist) are equal. Therefore, the relation
$\sim$ is, indeed, an equivalence relation on the set of
partitions of ${\bar \Pi} _{K_1}$.

 If an equivalence class of partitions
contains two distinct elements $f_M(\sigma )$ and $f_M(\delta )$
then the equality
 $$ f_M(\sigma )^\alpha = f_M(\delta )^\alpha$$
implies the corresponding  equation on the variables $\bar M_{\rm
veryshort},$ which is obtained by deleting all side variables
(which are equal) from $f_M(\sigma )$ and $f_M(\delta )$ and
equalizing the resulting words in very short variables. Denote by
$\Delta (\bar M _{\rm veryshort})=1$ this system.

Now we describe a transformation $T_4$. Fix a set of
representatives $R_p$ of partitions of ${\bar \Pi} _{K_1}$ with
respect to the equivalence relation $\sim$. For a given class of
equivalent partitions we take as a representative an interval
$\sigma$ with $f_M(\sigma )=\mu _{\rm left}\ldots \mu _{\rm
right}.$
 Below we say that $\mu ^{\alpha}$ almost contains
$u^{\beta}$ if $\mu ^{\alpha}$ contains  a subword which is the
reduced form of $c_1u^{\beta}c_2$ for some $c_1,c_2\in C_{\beta}.$

{\bf Principal variables}  A long variable $\mu _{\rm left}$ or
$\mu _{\rm right}$ for the interval $\sigma$ which represents a
class of equivalent partitions is called {\em principal } {\bf in
$\sigma$} in the following cases.

1) Let $f_X(\sigma )=u_i\ (i\neq n)$, where $u_i=x_iy_{i-1}^{-1}$
for $i>1$ and $u_1=x_1c_m^{-z_m}$ for $m\neq 0$. Then (see Lemma
\ref{main})
\begin{multline*}u_i^{\phi
_{K_1}}=A_{K_2+m+4i}^{*-q_4+1}x_{i+1}^{\phi _{K_2}}y_i^{-\phi
_{K_2}}x_i^{-q_1\phi _{K_2}}\\ \left (x_i^{-\phi
_{K_2}}A_{K_2+m+4i-4}^{*q_0}A_{K_2+m+4i-2}^{*(-q_2+1)}y_i^{\phi
_{K_2}}x_i^{-q_1\phi _{K_2}}\right )^{q_3-1}
A_{K_2+m+4i-4}^{*q_0}.\end{multline*}

A right variable $\mu _{\rm right}$  is principal in $\sigma$ if
$\mu _{\rm right}^{\alpha}$  almost contains a cyclically reduced
part of
\begin{multline*}
\left (x_i^{-\psi
_{K_2}}A_{K_2+m+4i-4}^{*q_0\beta}A_{m+4i-2}^{*(-q_2+1)\beta
}y_i^{\psi _{K_2}}x_i^{-q_1\psi _{K_2}}\right )^{q}\\ =
(x_i^{q_1}y_i)^{\psi _{K_2}}
(A^{*\beta}_{K_2+m+4i-1})^{-q}(y_i^{-1}x_i^{-q_1})^{\psi
_{K_2}},
\end{multline*} for some $q>2$. If $\mu _{\rm right}$ is not principal,
thenwe define  $\mu _{\rm left}$ as principal.

2) Let $f_X(\sigma )=v_i$, where $v_i=y_iu_i \ (i\neq 1,n)$ and
$v_1=y_1u_1\prod _{j=m-1}^1c_j^{-z_j}$. Then (see formula 3.a)
from Lemma \ref{main})
$$v_i^{\phi _{K_1}}=A_{K_2+m+4i}^{*(-q_4+1)} x_{i+1}^{\phi _{K_2}} A_{K_2+m+4i-4}^{*(-q_0)}
x_i^{q_1\phi _{K_2}}y_i^{\phi _{K_2}} A_{K_2+m+4i-2}^{*(q_2-1)}
A_{K_2+m+4i-4}^{*-1},$$ if $i\neq 1$, and

$$v_1^{\phi _{K_1}}=A_{K_2+m+4}^{*(-q_4+1)} x_{2}^{\phi _{K_2}} A_{K_2+2m}^{*(-q_0)}
x_1^{q_1\phi _{K_2}}y_1^{\phi _{K_2}} A_{K_2+m+1}^{*(q_2-1)}x_1\Pi
_{j=n}^1c_j^{-z_j},$$ if $i=1$.

\noindent
 A side variable $\mu_{\rm right}$ or $\mu_{\rm left}$ is
principal if $\mu_{\rm right}^{\alpha}$ (correspondingly,
$\mu_{\rm left}^{\alpha}$) almost contains
$(A^{\beta}_{K_2+m+4i})^{-q},$ for some $q>2$.

 3) Let $f_X(\sigma )=\bar u_n$.
Formula 3.c) from Lemma \ref{main} gives $\bar u_n^{\phi
_{K_1}}=A_{K_2+m+4n-8}^{*}$
$$
A_{K_2+m+4n-6}^{-q_2+1}(y_{n-1}^{-1}x_n^{-q_1})^{\phi
_{K_1}}A_{K_2+m+4n-8}^{*q_0}(x_n^{q_5}y_n)^{\phi
_{K_1}}A_{K_2+m+4n-2}^{*q_6-1}A_{K_2+m+4n-4}^{*-1}.$$

 A side variable $\mu_{\rm right}$ or
$\mu_{\rm left}$ is principal if $\mu _{\rm right}^{\alpha}$
(correspondingly, $\mu _{\rm left} ^{\alpha}$) almost contains
$(A^{\beta}_{K_2+m+4n-2})^q,$ $q>2$.

\vspace{3mm}
 4) Let $f_X(\sigma )= y_n$. A side variable $\mu _{\rm right}$ or
$\mu_{\rm left}$ is principal if $\mu_{\rm right}^{\alpha}$
($\mu_{\rm left}^{\alpha}$) almost contains
$(A^{\beta}_{K_2+m+4n-1})^q,$ $2q>p_{K_1}-2$.

\vspace{3mm}

 5) Let $f_X(\sigma )=z_j$, $j=1,\dots ,m-1$.
Then (by Lemma \ref{le:7.1.zforms}) $$z_j^{\phi
_{K_1}}=c_jz_j^{\phi _{K_2}}A_{K_2+j-1}^{*\beta
p_{j-1}}c_{j+1}^{z_{j+1}^{\phi _{K_2}}}A_{K_2+j}^{*\beta p_j-1}.$$

 A variable $\mu _{\rm left}$
($\mu _{\rm right}$) is {\em principal} if $\mu _{\rm
right}^{\alpha}$ (correspondingly, $\mu _{\rm left} ^{\alpha}$)
almost contains $(A^{\beta}_{K_2+j})^q,$ for some  $|q|>2.$ Both
left and right side variables can be simultaneously principal.

 \vspace{3mm}
 6)
Let $f_X(\sigma )=z_m$. Then $z_m^{\phi_{K_1}} = c_m^{K_2}
z_m^{\phi_{K_2}} A_{K_2+m-1}^{\ast p_{m-1}} x_1^{-\phi_{K_2}}
A_{K_2 +m}^{\ast p_m-1}.$

 In this case  $\mu _{\rm left}$ is {\em principal} in
$\sigma$ if and only if $\mu_{\rm left}$ is long (i.e.,  it is not
very short), and we define $\mu_{\rm right}$ to be always
non-principal. Observe that if  $\mu_{\rm left}$ is very short
then $\mu _{\rm right}^{\alpha}= f z_m^{\phi _{K_1}\beta}$ for a
very short $f \in F$.

 Let $f_X(\sigma )=z_m^{-1}c_mz_m$. By Lemma \ref{le:7.1.zforms}
$f_X(\sigma )^{\phi _{K_1}}=A_{K_2+m}^{*-p_{m}+1}x_1^{\phi
_{K_2}}A_{K_2+m}^{*p_{m}}.$

The variable $\mu _{\rm left}$  is {\em principal} in $\sigma$
 if and only if the following two
conditions hold:  $\mu _{\rm left}^{\alpha}$ almost contains $(A^{
\beta}_{K_2+m})^q,$ for some $q$ with $ |q|>2$; $\mu_{\rm
left}^{-1}\neq fz_m^{\phi _{K_1}\beta}$ for a very short $f\in F$.

  Similarly, the variable $\mu _{\rm right}$  is
{\em principal} in $\sigma$  if and only if the following two
conditions hold:  $\mu _{\rm right}^{\alpha}$ almost contains
$(A^{\beta}_{K_2+m})^q,$ for some $q$ with $ |q|>2$;  $\mu _{\rm
right}^{\alpha}\neq fz_m^{\phi _{K_1}\beta}$ for a very short
$f\in F$.

 Observe, that in this case
the variables $\mu _{\rm left}$ and $\mu _{\rm right}$ can be
simultaneously principal in $\sigma$ and non-principal in
$\sigma$. The latter happens if and only if  $\mu _{\rm
right}^{\alpha}=f_1z_m^{\phi _{K_1}\beta}$ and
 $\mu _{\rm left}^{\alpha}=z_m^{-\phi _{K_1}\beta}f_2$ for
 some very short elements $f_1, f_2 \in F$. Therefore, if both
 $\mu _{\rm left}$ and $\mu _{\rm right}$ are non-principal then
 they can be expressed in terms of $z_m^{\phi _{K_1}}$ and very short variables.

\begin{Claim} \label{claim:one-p} Every partition has at least one principal variable,
unless this partition is of that particular  type from  Case 6).
\end{Claim}

\begin{Claim}
\label{claim:N-N} If both side variables of a partition of
$\bar\Pi _{K_1}$  are non-principal, then they are non-principal
in every partition of $\bar\Pi _{K_1}$.
\end{Claim}

\begin{Claim}\label{claim:p-n} Let $n\neq 0$. Then  a side
variable can be principal only in one class of equivalent
partitions.
\end{Claim}
\begin{proof}   Follows from the definition of principal variables.
\end{proof}

For the cut equation $\bar\Pi _{K_1}$ we construct a finite graph
$\Gamma = (V,E).$ Every vertex from $V$ is marked by variables
from $\bar M_{\rm side}^{\pm 1}$ and letters
from the alphabet $\{P,N\}$. Every edge from $E$ is colored either as red or blue.
The graph $\Gamma$ is constructed as follows.
Every partition $f_M(\sigma) = \mu_1 \cdots \mu_k$
of $\bar\Pi _{K_1}$ gives two
vertices $v_{\sigma,\rm left}$ and $v_{\sigma,\rm right}$ into $\Gamma$,
so
  $$V = \bigcup_{\sigma} \{v_{\sigma,\rm left},  v_{\sigma,\rm right}\}.$$
We mark $v_{\sigma,{\rm left}}$ by $\mu_1$ and $v_{\sigma,\rm
right}$ by $\mu_k$. Now we mark the  vertex $v_{\sigma,\rm left}$
 by a letter $P$ or letter $N$  if $\mu_1$
   is correspondingly principal or non-principal in
 $\sigma$. Similarly, we mark $v_{\sigma,\rm right}$  by $P$ or $N$ if
 $\mu_k$ is principal or non-principal in $\sigma$.

 For every $\sigma$ the
vertices $v_{\sigma,left}$ and $v_{\sigma,right}$ are connected by
a {\em red} edge. Also, we connect by a {\em blue} edge every pair
of
 vertices which are marked by variables $\mu, \nu$ provided $\mu =
 \nu$ or $\mu = \nu^{-1}$. This describes the graph $\Gamma$.

 Below we construct a new graph $\Delta$
  which is obtained from $\Gamma$ by deleting some  blue edges
  according to the following procedure.
Let $B$ be a maximal connected blue component of $\Gamma$, i.e., a
connected component of the graph obtained from $\Gamma$ by
deleting all red edges. Notice, that $B$ is a complete graph, so
every two vertices in $B$ are connected by a blue edge. Fix a
vertex $v$ in $B$ and consider the star-subgraph $Star_B$ of $B$
generated by all edges adjacent to $v$. If $B$ contains a vertex
marked by $P$ then we choose $v$ with label $P$, otherwise $v$ is
an arbitrary vertex of $B$. Now, replace $B$ in $\Gamma$ by the
graph $Star_B$, i.e., delete all edges in $B$ which are not
adjacent to $v$. Repeat this procedure for every maximal blue
component $B$ of $\Gamma$. If the blue component corresponds to
long bases of case 6) that are non-principal and equal to
$f_1z_m^{\phi _{K_1}}f_2$ for very short $f_1,f_2$, we remove all
the blue edges that produce  cycles if the red edge from $\Gamma$
connecting non-principal $\mu _{\rm left}$ and $\mu _{\rm right}$
is added to the component (if such a red edge exists). Denote the
resulting graph by $\Delta$.

In the next claim we describe connected components of the graph
$\Delta$.

\begin{Claim}\label{20.}
Let $C$ be a connected component of $\Delta$. Then one of the
following holds:
 \begin{enumerate}
   \item  [(1)]\label{item:4} there is a vertex  in $C$ marked by a variable
   which does not occur as a principal variable in any partition
    of\/ $\bar\Pi _{K_1}$. In particular, any component which
    satisfies one of the following conditions has such a vertex:
  \begin{itemize}
  \item  [a)] there is a vertex  in $C$ marked by a variable which is a short variable in some partition of
   $\bar\Pi _{K_1}$.
   \item [b)] there is a red  edge   in $C$ with both endpoints marked by $N$ {\rm (}it corresponds to a partition described
   in Case $6$ above{\rm )};
    \end{itemize}
   \item [(2)]\label{item:3} both endpoints of every red edge in $C$ are marked by $P$.
  In this case $C$ is an isolated vertex;
  \item [(3)] \label{item:5} there is a vertex  in $C$  marked by a  variable
  $\mu$ and $N$ and if $\mu$ occurs as a label of an endpoint  of some red edge in $C$
   then the other endpoint of this edge is marked by $P$.
  \end{enumerate}
\end{Claim}

\begin{proof}   Let $C$ be a connected component of $\Delta$.
Observe first, that if $\mu$ is a short variable  in
   $\bar\Pi _{K_1}$ then $\mu$ is not principle in $\sigma$ for
   any interval $\sigma$ from $\bar\Pi _{K_1}$, so there is no vertex in $C$
   marked by both $\mu$ and $P$.  Also, it follows from  Claim \ref{claim:N-N} that
if there is a red  edge  $e$ in $C$ with both endpoints marked by
$N$, then the variables assigned to endpoints of $e$ are
non-principle in any interval $\sigma$ of  $\bar\Pi _{K_1}$. This
proves the part ``in particular" of 1).

  Now assume that  the component $C$  does not satisfy any
of the conditions (1), (2). We need to show that $C$ has type (3).
It follows that every variable which occurs as a label of a vertex
in $C$ is long and it labels, at least, one vertex in $C$ with
label $P$. Moreover,  there are non-principle occurrences of
variables in $C$.

 We summarize some properties of $C$ below:
 \begin{itemize}
  \item  There are no blue edges in $\Delta$ between
vertices with labels $N$ and $N$ (by construction).
 \item  There are no blue edges between
 vertices labelled by $P$ and $P$ (Claim  \ref{claim:p-n}).
  \item   There are no  red edges in
  $C$ between vertices labelled by $N$ and $N$ (otherwise 1) would hold).
   \item  Any reduced path in $\Delta$ consists of
   edges of alternating color (by construction).
 \end{itemize}

We claim that $C$ is a tree. Let $p =  e_1 \ldots e_k$ be a simple
loop in $C$ (every vertex in $p$ has degree 2 and the terminal
vertex of $e_k$ is equal to the starting point of $e_1$).

We show first that $p$ does not have red edges with endpoints
labelled by $P$ and $P$. Indeed, suppose there exists such an edge
in $p$. Taking cyclic permutation of $p$ we may assume that $e_1$
is a red edge with labels $P$ and $P$. Then $e_2$ goes from a
vertex with label $P$ to a vertex with label $N$. Hence the next
red edge $e_3$ goes from $N$ to $P$, etc. This shows that every
blue edge along $p$ goes from $P$ to $N$. Hence the last edge
$e_k$ which must be blue goes from $P$ to $N$ -contradiction,
since all the labels of $e_1$ are $P$.

It follows that both colors of edges and labels of vertices in $p$
alternate. We may assume now that $p$ starts with a vertex with
label $N$ and the first edge $e_1$ is red. It follows that the end
point of $e_1$ is labelled by $N$ and all blue edges go from $N$
to $P$. Let $e_i$ be a blue edge from $v_i$ to $v_{i+1}$. Then the
variable $\mu_i$ assign to the vertex $v_i$ is principal in the
partition associated with the red edge $e_{i-1}$ , and the
variable $\mu_{i+1} = \mu_i^{\pm 1}$ associated with $v_{i+1}$ is
a non-principal side variable in the partition $f_M(\sigma)$
associated with the red edge $e_{i+1}$. Therefore, the the side
variable $\mu_{i+2}$ associated with the end vertex $v_{i+2}$ is a
principal side variable in the partition $f_M(\sigma)$ associated
with $e_{i+1}$. It follows from the definition of principal
variables that the length of $\mu_{i+2}^\alpha$ is much longer
than the length of $\mu_{i+1}^\alpha$, unless the variable $\mu_i$
is described in the Case 1). However, in the letter case the
variable $\mu_{i+2}$ cannot occur in any other partition
$f_M(\delta)$  for $\delta \neq \sigma$. This shows that there no
blue edges in $\Delta$ with endpoints labelled by such
$\mu_{i+2}$. This implies that $v_{i+2}$ has degree one in
$\Delta$ - contradiction wit the choice of $p$. This shows that
there are no vertices labelled by such variables described in Case
1).  Notice also, that the length of variables (under $\alpha$) is
preserved along blue edges: $|\mu_{i+1}^\alpha| = |(\mu_i^{\pm
1})^\alpha| = |\mu_i^\alpha|$. Therefore,
 $$  |\mu_i^\alpha| = |\mu_{i+1}^\alpha| < |\mu_{i+2}^\alpha|$$
 for every $i$.

It follows that going along $p$ the length of $ |\mu_i^\alpha|$
increases, so $p$ cannot be a loop. This implies that $C$ is a
tree.

Now we are ready to show that the component $C$ has type (3). Let
$\mu_1$ be a variable assigned to some vertex $v_1$ in $C$ with
label $N$.   If $\mu_1$ satisfies the condition (3) then we are
done. Otherwise, $\mu_1$ occurs as a label of one of
$P$-endpoints, say $v_2$  of a red edge $e_2$ in $C$ such that the
other endpoint of $e_2$, say $v_3$   is non-principal. Let $\mu_3$
be the label of $v_3$. Thus $v_1$ is connected to $v_2$ by a blue
edge and $v_2$ is connected to $v_3$ by a red edge.   If $\mu_3$
does not satisfy the condition (3) then we can repeat the process
(with $\mu_3$ in place of $\mu_1$). The graph $C$ is finite, so in
finitely many steps either we will find a variable that satisfies
(3) or we will construct a closed reduced path in $C$. Since $C$
is a tree the latter does not happen, therefore $C$ satisfies (3),
as required.

 \end{proof}

\begin{Claim} \label{21} The graph $\Delta$ is a forest, i.e., it is union of trees.
\end{Claim}
 \begin{proof}
  Let $C$ be a connected component of $\Delta$. If $C$ has type (3)
  then it is a tree, as has been shown in Claim \ref{20.}
    If $C$ of the type (2) then by Claim \ref{20.}  $C$ is an
    isolated vertex -- hence a tree.
    If $C$ is of the type (1) then $C$ is a tree because each interval corresponding to this component has exactly one
    principal variable, and the same long variable cannot be principal in two different
    intervals. Although the same argument as in (3) also
    works here.

 \end{proof}

 Now  we define the sets $\bar M_{\rm useless}, \bar M_{\rm free}$ and assign values
 to variables from $\bar M = \bar M_{\rm useless} \cup  \bar M_{free} \cup \bar M_{\rm
 veryshort}$. To do this we use the structure of connected components of
 $\Delta$. Observe first, that all occurrences of a given variable
 from ${\bar M_{\rm sides}}$ are located in the same connected
 component.

 Denote by $\bar M_{free}$ subset of $\bar M$ which consists of variables
of the following types:
 \begin{enumerate}
  \item variables which do not occur as principal in any partition
  of $(\bar\Pi _{K_1})$;
   \item \label{item: 2-2}
   one (but not the other) of the variables  $\mu$ and $\nu$ if they
   are both principal side variables of a partition of the
   type (\ref{item:3}) and such that $\nu \neq \mu^{-1}$.
 \end{enumerate}

Denote by $\bar M_{\rm useless} = \bar M_{\rm side}- \bar M_{\rm free}.$

 \begin{Claim}
For  every $\mu \in \bar M_{\rm useless}$ there exists a word
$$V_\mu \in F[X \cup \bar M_{\rm free}\cup \bar M_{\rm veryshort}]$$
such that for every map $\alpha_{\rm free} : \bar M_{\rm free}
 \rightarrow F$, and every solution
 $$\alpha_s: F[\bar M_{\rm veryshort}] \rightarrow F$$
  of the system $\Delta(\bar M_{\rm veryshort}) = 1$  the map $\alpha: F[\bar M] \rightarrow F$
  defined by
   \[ \mu^\alpha = \left\{\begin{array}{ll}
   \mu^{\alpha_{\rm free}} &\mbox{ if  $\mu \in \bar M_{\rm free}$;}\\
   \mu^{\alpha_{s}} & \mbox{ if $\mu \in \bar M_{\rm veryshort}$;}\\
     \bar V_\mu(X^\delta, \bar M_{\rm free}^{\alpha_{\rm free}},
   \bar M_{\rm veryshort}^{\alpha_s}) & \mbox{ if $\mu \in \bar M_{\rm useless}$.}
   \end{array}
 \right. \]
  is a group solution of $\bar \Pi_{K_1}$ with respect to $\beta$.

 \end{Claim}
 \begin{proof}   The claim follows from Claims \ref{20.} and
 \ref{21}. Indeed,
take as values of short variables an arbitrary solution $\alpha
_s$ of the system $\Delta (\bar M_{\rm veryshort})=1$. This system
is obviously consistent, and we fix its solution. Consider
connected components of type (1) in Claim \ref{20.}.  If $\mu$ is
a principal variable for some $\sigma$ in such a component, we
express $\mu ^{\alpha}$ in terms of values of very short variables
$\bar M_{\rm veryshort}$ and elements $t^{\psi _{K_1}},$ $t\in X$
that correspond to labels of the intervals. This expression does
not depend on $\alpha _s, \beta $ and tuples $q, p^*.$ For
connected components of $\Delta $ of types (2) and (3) we express
values $\mu ^{\alpha}$  for $\mu\in M_{\rm useless}$ in terms of
values $\nu ^{\alpha}$, $\nu\in M_{\rm free}$  and $t^{\psi
_{K_1}}$ corresponding to the labels of the intervals.
\end{proof}

We can now finish the proof of Proposition \ref{Or}. Observe, that
$M_{\rm veryshort}\subseteq \bar M_{\rm veryshort}.$ If $\lambda$
is an additional very short variable from $M^*_{\rm veryshort}$
that appears when transformation $T_1$ or $T_2$ is performed,
$\lambda ^{\alpha}$ can be expressed in terms $M_{\rm
veryshort}^{\alpha}$. Also, if a variable $\lambda$ belongs to
$\bar M_{\rm free}$ and does not belong to $M$, then there exists
a variable $\mu\in M$, such that $\mu ^{\alpha}=u^{\psi
_{K_1}}\lambda ^{\alpha},$ where $u\in F(X, C_S)$, and we can
place $\mu$ into $M_{\rm free}.$

Observe, that the argument above is based only on the tuple $p$,
it does not depend on the tuples $p^*$ and $q$. Hence the words
$V_\mu$ do not depend on $p^*$ and $q$.

 The Proposition is proved for $n\neq 0.$
If $n=0$, partitions of the intervals with  labels $z_{n-1}^{\phi
_{K_1}}$ and $z_{n}^{\phi _{K_1}}$ can have equivalent principal
right variables, but in this case the left variables will be
different and do not appear in other non-equivalent partitions.
The connected component of $\Delta$ containing these partitions
will have only four vertices  one blue edge.

In the case $n=0$ we transform equation $\Pi _{K_1}$ applying
transformation $T_1$ to the form when the intervals are labelled
by $u^{\phi _{K_1}},$ where $$ u\in\left\{z_1,\ldots ,z_m,
c_{m-1}^{z_{m-1}}, z_{m}c_{m-1}^{-z_{m-1}}\right\}.$$

If $\mu _{\rm left}$ is very short for the interval $\delta$
labelled by $(z_mc_{m-1}^{-z_{m-1}})^{\phi _{K_1}},$ we can apply
$T_2$ to
 $\delta$, and split it into intervals with labels $z_m^{\phi
_{K_1}}$ and $c_{m-1}^{-z_{m-1}^{\phi _{K_1}}}.$ Indeed, even if
we had to replace $\mu_{\rm right}$ by the product of two
variables, the first of them would be very short.

If $\mu _{\rm left}$ is not very short for the interval $\delta$
labelled by \[(z_mc_{m-1}^{-z_{m-1}})^{\phi _{K_1}}= c_mz_m^{\phi
_{K_2}}A_{m-1}^{*p_{m-1}-1},\] we do not split the interval, and
$\mu _{\rm left}$ will be considered as the principal variable for
it. If $\mu _{\rm left}$ is not very short for the interval
$\delta$ labelled by $z_m^{\phi _{K_1}}=z_m^{\phi
_{K_2}}A_{m-1}^{*p_{m-1}}$, it is a principal variable, otherwise
$\mu _{\rm right}$ is principal.

If an interval $\delta$ is labelled by $(c_{m-1}^{z_{m-1}})^{\phi
_{K_1}}=A_{m-1}^{*-p_{m-1}+1}c_m^{-z_m^{\phi
_{K_2}}}A_{m-1}^{*p_{m-1}},$ we consider $\mu _{\rm right}$
principal if $\mu _{\rm right}^{\alpha}$ ends with
$(c_m^{-z_m^{\phi _{K_2}}}A_{m-1}^{*p_{m-2}})^{\beta},$ and the
difference is not very short. If $\mu _{\rm left}^{\alpha}$ is
almost $z_m^{-\phi _k\beta}$ and $\mu _{\rm right}^{\alpha}$ is
almost $z_m^{\phi _k\beta}$, we do not call any of the side
variables principal. In all other cases $\mu _{\rm left}$ is
principal.

Definition of the principal variable in the interval with label
$z_i^{\phi _{K_1}}$,  $i=1,\dots, m-2$ is the same as in 5) for
$n\neq 0.$

A variable can be principal only in one  class of equivalent
partitions. All the rest of the proof is the same as for $n>0.$

 \end{proof}

 Now we continue the proof of Theorem A.
 Let  $L = 2K + \kappa(\Pi)4K$ and
 $$\Pi_\phi= \Pi_L\rightarrow \Pi_{L-1}\rightarrow\ldots  \rightarrow \ldots  $$
be the sequence of $\Gamma$-cut equations (\ref{eq:7.3.15}). For a
$\Gamma$-cut equation $\Pi_j$ from (\ref{eq:7.3.15}) by $M_j$ and
$\alpha_j$ we denote the corresponding set of variables and the
solution relative to $\beta$.

By Claim \ref{5.} in the sequence (\ref{eq:7.3.15}) either there
is $3K$-stabilization at $K(r+2)$ or $Comp(\Pi_{K(r+1)}) = 0$.

Case 1. Suppose there is $3K$-stabilization at $K(r+2)$ in the
sequence (\ref{eq:7.3.15}).

By Proposition \ref{Or} the set of variables $M_{K(r+1)}$ of the
cut equation $\Pi_{K(r+1)}$ can be partitioned into three subsets
 $$M_{K(r+1)} = M_{\rm veryshort} \cup M_{\rm free} \cup M_{\rm
 useless}$$
  such that there exists a finite consistent system of equations $\Delta(M_{\rm veryshort}) = 1$
   over $F$ and words $V_\mu \in F[X, M_{\rm free}, M_{\rm
   veryshort}]$, where $\mu \in M_{\rm useless}$, such that
for every solution $\delta \in {\mathcal B}$, for every map
$\alpha_{\rm free} : M_{\rm free}
 \rightarrow F$, and every solution $\alpha_{short}: F[M_{\rm veryshort}] \rightarrow F$
  of the system $\Delta(M_{\rm veryshort}) = 1$  the map $\alpha_{K(r+1)}: F[M] \rightarrow F$
  defined by
   \[ \mu^{\alpha_{K(r+1)}} = \left\{\begin{array}{ll}
   \mu^{\alpha_{\rm free}} &\mbox{ if  $\mu \in M_{\rm free}$}\\
   \mu^{\alpha_{\rm short}} & \mbox{ if $\mu \in M_{\rm veryshort}$}\\
     V_\mu(X^\delta, M_{\rm free}^{\alpha_{\rm free}},
   M_{\rm veryshort}^{\alpha_s}) & \mbox{ if $\mu \in M_{\rm useless}$}
   \end{array}
 \right. \]
  is a group solution of $\Pi_{K(r+1)}$ with respect to $\beta$.
Moreover, the words $V_\mu$  do not depend on tuples $p^*$ and
$q$.

By Claim \ref{2.0.0} if $\Pi = ({\mathcal E}, f_X, f_M)$ is
 a $\Gamma$-cut equation and  $\mu \in M$ then
there exists a word ${\mathcal M}_\mu (M_{T(\Pi)}, X)$ in the free
group $F[M_{T(\Pi)} \cup X]$ such that
$$\mu^{\alpha_{\Pi}} = {\mathcal M}_\mu \left(M_{T(\Pi)}^{{\alpha_{T(\Pi)}}},
X^{\phi_{K(r+1)}}\right)^\beta, $$  where $\alpha_\Pi$ and
$\alpha_{T(\Pi)}$ are the corresponding solutions of $\Pi$ and
$T(\Pi)$ relative to $\beta$.

Now, going along the sequence (\ref{eq:7.3.15})  from
$\Pi_{K(r+1)}$ back to the cut equation $\Pi_L$ and using
repeatedly the remark above for each $\mu \in M_L$ we obtain a
word
$${\mathcal M^\prime}_{\mu,L} (M_{K(r+1)}, X^{\phi_{K(r+1)}}) = {\mathcal M^\prime}_{\mu,L}
(M_{\rm useless},M_{\rm free}, M_{\rm veryshort},
X^{\phi_{K(r+1)}})$$
 such that
$$\mu^{\alpha_L} = {\mathcal M^\prime}_{\mu,L} (M_{K(r+1)}^{\alpha_{K(r+1)}}, X^{\phi_{K(r+1)}})^{\beta}.$$
Let $\delta = \phi_{K(r+1)} \in {\mathcal B}$ and put
\begin{multline*}
{\mathcal M}_{\mu,L}(X^{\phi_{K(r+1)}})\\ = {\mathcal M^\prime}_{\mu,L} (V_\mu(X^{\phi_{K(r+1)}},M_{\rm free}^{\alpha_{\rm free}},
M_{\rm veryshort}^{\alpha_{short}}),M_{\rm free}^{\alpha_{\rm
free}}, M_{\rm veryshort}^{\alpha_{short}}, X^{\phi_{K(r+1)}}).
\end{multline*}
 Then for every $\mu \in M_L$
$$\mu^{\alpha_L}  = {\mathcal M}_{\mu,L}(X^{\phi_{K(r+1)}})^\beta$$
If we denote by ${\mathcal M}_{L}(X)$ a tuple of words
$${\mathcal M}_{L}(X) = ({\mathcal M}_{\mu_1,L}(X), \ldots,{\mathcal
M}_{\mu_{|M_L|},L}(X)),$$
  where $\mu_1, \ldots, \mu_{|M_L|}$ is some fixed ordering of $M_L$
  then
    $$M_L^{\alpha_L} ={\mathcal M}_{L}(X^{\phi_{K(r+1)}})^\beta.$$
Observe, that  the words  ${\mathcal M}_{\mu,L}(X)$, hence
${\mathcal M}_{L}(X)$ (where $ X^{\phi_{K(r+1)}}$  is replaced by
$X$) are the same for  every $\phi_L \in {\mathcal B}_{p,q}$.

It follows from property c) of the cut equation $\Pi_\phi$ that
the solution ${\alpha_L}$ of $\Pi_\phi$ with respect to $\beta$
gives rise to a group solution of the original cut equation
$\Pi_{\mathcal L}$ with respect to $\phi_L\circ \beta$.

Now, property c) of the initial cut equation $\Pi_{\mathcal L} =
(\mathcal E,f_X,f_{M_L})$
 insures that for every $\phi_L \in {\mathcal B}_{p,q}$ the pair
$(U_{\phi_L\beta},V_{\phi_L\beta})$ defined by
  $$U_{\phi_L\beta} =  Q(M_L^{\alpha_L}) =  Q({\mathcal M}_L(X^{\phi_{K(r+1)}}))^\beta, $$
$$V_{\phi_L\beta}= P(M_L^{\alpha_L}) = P({\mathcal
M}_L(X^{\phi_{K(r+1)}}))^\beta.$$
 is a solution  of the  system$S(X) = 1 \wedge T(X,Y) = 1.$

  We claim that
 $$Y(X) = P({\mathcal M}_L(X))$$
  is a solution of the equation $T(X,Y) = 1$ in
$F_{R(S)}$.  By Theorem \ref{cy2} ${\mathcal B}_{p,q,\beta}$ is a
discriminating family of solutions for the group $F_{R(S)}$. Since
$$T(X,Y(X))^{\phi\beta} = T(X^{\phi\beta},Y(X^{\phi\beta})) = T(X^{\phi\beta}, {\mathcal
M}_L(X^{\phi\beta})) = T(U_{\phi_L\beta}, V_{\phi_L\beta}) = 1$$
for any ${\phi\beta} \in {\mathcal B}_{p,q,\beta}$ we deduce that
$T(X,Y_{p,q}(X)) = 1$ in $F_{R(S)}$.

Now we need to  show that $T(X,Y) = 1$ admits a complete $S$-lift.
 Let $W(X,Y) \neq 1$ be an inequality such that  $T(X,Y) = 1 \wedge W(X,Y) \neq 1$
  is compatible with $S(X) = 1$. In this event, one may assume
  (repeating the argument from the beginning of this section) that the set
   $$
   \Lambda = \{(U_\psi,V_\psi) \mid \psi \in {\mathcal L}_2
   \}$$
   is such that every pair
 $(U_\psi,V_\psi) \in \Lambda$ satisfies the formula $T(X,Y) = 1 \wedge W(X,Y) \neq
 1.$  In this case, $W(X,Y_{p,q}(X)) \neq
 1$ in $F_{R(S)}$, because its image in $F$ is non-trivial:
 $$W(X,Y_{p,q}(X))^{\phi\beta} = W(U_\psi,V_\psi) \neq 1.$$
 Hence $T(X,Y) =1$ admits a complete lift into  generic
point of $S(X) = 1$.

Case 2. A similar argument applies when $Comp(\Pi_{K(r+2)}) = 0$.
Indeed, in this case for every $\sigma \in {\mathcal E}_{K(r+2)}$
the word $f_{M_{K(r+1)}}(\sigma)$ has length one, so
$f_{M_{K(r+1)}}(\sigma) =\mu$ for some $\mu \in M_{K(r+2)}.$ Now
one can replace the word $V_\mu \in F[X \cup M_{\rm free} \cup
M_{\rm veryshort}]$ by the label $f_{X_{K(r+1)}}(\sigma)$ where
$f_{M_{K(r+1)}}(\sigma) =\mu$ and then repeat the argument.

 \end{proof}

\section{Implicit function theorem for NTQ systems} \label{se:7.4}

In this section we prove Theorems B, C, D from  Introduction.

We begin with the proof of Theorem B. To this end let $U(X,A)=1$
be a regular NTQ-system and $V(X,Y,A)=1$ an equation  compatible
with $U=1$. We need to show that $V(X,Y,A)=1$ admits a complete
effective $U$-lift.

 We use induction on the number $n$ of levels in the system
$U=1$. We construct a solution tree $T_{\rm sol}(V(X,Y,A)\wedge
U(X,Y))$ with parameters $X=X_1\cup\cdots\cup X_n.$ In the
terminal vertices of the tree there are generalized equations
$\Omega _{v_1},\dots ,\Omega _{v_k}$ which are equivalent to cut
equations $\Pi _{v_1},\dots ,\Pi _{v_k}$.

If $S_1(X_1,\dots ,X_n)=1$ is an empty equation, we can take
Merzljakov's words (see Introduction ) as values of variables from
$X_1$, express $Y$ as functions in $X_1$ and a solution of some
$W(Y_1,X_2,\dots ,X_n)=1$ such that for any solution of the system
\bea
S_2(X_2, \dots, X_n,A) &=& 1\\
&\vdots & \\
S_n(X_n,A) &=& 1 \eea
 equation $W=1$ has a solution.

Suppose, now that $S_1(X_1,\dots ,X_n)=1$ is a regular quadratic
equation. Let $\Gamma $ be a basic sequence of automorphisms for
the equation $S_1(X_1,\dots ,X_n,A)=1.$
 Recall that $$\phi_{j,p} = \gamma_j^{p_j} \cdots \gamma_1^{p_1} = \stackrel{\leftarrow}{\Gamma}_j^p,$$
 where $j \in {\mathbb{N}}$,  $\Gamma_j = (\gamma_1,
 \dots, \gamma_j)$ is the initial subsequence of length $j$ of the
sequence
 $\Gamma^{(\infty)}$, and $p = (p_1, \dots,p_j) \in {\mathbb{N}}^j$.
  Denote by $\psi_{j,p}$ the following solution of \newline $S_1(X_1) = 1$:
  $$\psi_{j,p}=\phi_{j,p}\alpha,$$
where $\alpha$   is a composition of a solution of $S_1=1$ in
$G_2$ and a solution from a generic family of solutions of the
system
\bea S_2(X_2, \dots, X_n,A) &=& 1\\
&\vdots & \\
S_n(X_n,A) &=& 1 \eea in $F(A).$ We can always suppose that
$\alpha$ satisfies a small cancellation condition with respect to
$\Gamma .$

  Set $$ \Phi = \left\{\phi_{j,p} \mid j \in {\mathbb{N}}, p \in {\mathbb{N}}^j \right\}$$
  and let ${\mathcal L}^{\alpha }$ be an  infinite subset of $\Phi ^{\alpha}$
  satisfying one of the cut equations above. Without loss of generality
  we can suppose it satisfies $\Pi _1$.
By Proposition \ref{Or} we can express variables from $Y$ as
functions of the set of $\Gamma $-words in $X_1$, coefficients,
variables $M_{\rm free}$ and variables $M_{\rm veryshort}$,
satisfying the system of equations $\Delta (M_{\rm veryshort})$
The system $\Delta (M_{\rm veryshort})$ can be turned into a
generalized equation with parameters $X_2\cup\cdots \cup X_n$,
such that for any solution of the system
\bea S_2(X_2, \dots, X_n,A) &=& 1\\
& \vdots & \\
S_n(X_n,A) &=& 1 \eea
 the system
$\Delta (M_{\rm veryshort})$ has a solution. Therefore, by
induction, variables\linebreak $(M_{\rm veryshort})$ can be found as
elements of $G_2$, and variables $Y$ as elements of $G_1$. Theorem
B is proved. \hfill $\Box$
\medskip

In order to prove Theorem C we need some auxiliary results.

\begin{lm}
All stabilizing automorphisms {\rm (}see \cite{Gr}{\rm )} of the
left side of the equation
\begin{equation}\label{10}
c_1^{z_1}c_2^{z_2}(c_1c_2)^{-1}=1
\end{equation}
have the form $z_1^{\phi}=c_1^kz_1(c_1^{z_1}c_2^{z_2})^n,
z_2^{\phi}=c_2^mz_2(c_1^{z_1}c_2^{z_2})^n.$ All stabilizing
automorphisms of the left side of the equation
\begin{equation}\label{11}
x^2c^{z}(a^2c)^{-1}=1\end{equation} have the form
$x^{\phi}=x^{(x^2c^z)^n}, z^{\phi}=c^kz(x^2c^z)^n$. All
stabilizing automorphisms  of the left side of the equation
\begin{equation}\label{12}
x_1^2x_2^2(a_1^2a_2^2)^{-1}=1\end{equation} have the form
$x_1^{\phi}=(x_1(x_1x_2)^m)^{(x_1^2x_2^2)^n}$,
$x_2^{\phi}=((x_1x_2)^{-m}x_2)^{(x_1^2x_2^2)^n}.$
\end{lm}
\begin{proof}  The computation of the automorphisms can be done
by software ``Magnus''. The statement of the lemma also follows
from the fact that punctured surfaces corresponding to $QH$
subgroups corresponding to these equations (see \cite{JSJ},
Section 5) do not contain two intersecting simple closed curves
that are not boundary-parallel. Therefore if G is a freely
indecomposable finitely generated fully residually free group that
has a $QH$ subgroup $Q$ corresponding to one of these equations,
then $G$ does not have two intersecting cyclic splittings with
edge groups conjugated into $Q$.
\end{proof}

If a quadratic equation $S(X) = 1$ has only commutative solutions
then the radical $R(S)$ of $S(X)$ can be described (up to a linear
change of variables) as follows (see \cite{KMNull}):
 $${\rm Rad}(S)=ncl\{[x_i,x_j], [x_i,b],  \mid \ i,j =1,\dots ,k\},$$
 where $b$ is an element (perhaps, trivial) from $F$. Observe, that if $b$ is not
 trivial then $b$ is not a proper power in $F$. This shows that $S(X) = 1$ is
 equivalent to the system
 \begin{equation}
\label{eq:U}
 U_{\rm com}(X) =  \{[x_i,x_j]=1, [x_i,b]=1,  \mid \ i,j =1,\dots ,k\}.
\end{equation}
The system $U_{\rm com}(X) = 1$ is equivalent to a single
equation, which we also denote by $U_{\rm com}(X)=1.$
 The coordinate group $H = F_{R(U_{\rm com})}$ of the system $U_{\rm com}= 1$, as well
 as of the corresponding  equation, is $F$-isomorphic to the free extension of the
 centralizer $C_F(b)$ of rank $n$. We need the following notation to deal with $H$.
 For a set $X$ and  $b \in F$    by  $A(X)$ and $A(X,b)$ we denote
  free abelian groups with  basis $X$ and $X \cup \{b\}$,  correspondingly. Now,
  $H \simeq F \ast_{b = b} A(X,b)$. In particular, in the case when $b = 1$ we
  have $H = F \ast A(X)$.

 \begin{lm}
\label{simple} Let $F = F(A)$ be a non-abelian free group and
$V(X,Y,A) = 1$, $W(X,Y,A) = 1$ be equations over $F$.  If a
formula
 $$\Phi=\forall X (U_{\rm com}(X) = 1\rightarrow \exists Y (V(X,Y,A)=1 \wedge W(X,Y,A)\not
=1))$$
 is true in $F$ then there exists a finite number of
 extensions $\phi _k$ on $H$ of $\left< b \right>$-embeddings
 $A(X,b) \rightarrow A(X,b) \ ( k \in K)$   such that:
\bi
 \item[(1)]  every formula
 $$\Phi _k = \exists Y (V(X^{\phi _k}, Y,A)=1 \wedge W(X^{\phi _k},Y,A)\not =1)$$
 holds in the coordinate group $H = F\ast_{b = b} A(X,b)$;

 \item[(2)]  for any solution  $\lambda: H \rightarrow F$  there exists a solution
 $\lambda^{\ast} : H \rightarrow F$ such that $\lambda = \phi_k
 \lambda^{\ast}$ for some $k \in K$.
 \ei
 \end{lm}

\begin {proof} We  construct a set of initial parameterized generalized
equations $${\mathcal GE}(S) = \{\Omega_1, \ldots, \Omega_r\}$$
for $V(X,Y,A) = 1$ with respect to the set of parameters $X$.
 For each $\Omega \in {\mathcal GE}(S)$, in \cite[Section 8]{JSJ},  we
constructed  the finite tree $T_{\rm sol}(\Omega)$ with respect to
parameters $X$. Observe, that non-active part
$[j_{v_0},\rho_{v_0}]$ in the root equation $\Omega \Omega_{v_0}$
of the tree $T_{\rm sol}(\Omega)$ is  partitioned into a disjoint
union of closed sections  corresponding to $X$-bases and constant
bases (this follows from the construction of the initial equations
in the set  ${\mathcal GE}(S)$). We label every closed section
$\sigma$ corresponding to a variable  $x \in X^{\pm 1}$ by $x$,
and every constant section corresponding to a constant $a$ by $a$.
Due to our construction of the tree $T_{\rm sol}(\Omega)$ moving
along a brunch $B$ from the initial vertex $v_0$ to a terminal
vertex $v$ we transfer all the bases from the non-parametric part
into parametric part until, eventually,  in  $\Omega_v$ the whole
interval consists of the parametric part. For a terminal vertex
$v$ in $T_{\rm sol}(\Omega )$ equation $\Omega _v$ is periodized
(see Section 5.4). We can consider the correspondent periodic
structure $\mathcal P$ and the subgroup $\tilde Z_2$. Denote the
cycles generating this subgroup by $z_1,\dots ,z_m$.  Let
$x_i=b^{k_i}$ and $z_i=b^{s_i}$. All $x_i$'s are cycles, therefore
the corresponding system of  equations can be written as a system
of linear equations with integer coefficients in variables
$\{k_1,\dots ,k_n\}$ and variables $\{s_1,\dots ,s_m\}$ :
\begin{equation}\label{abel}
k_i=\sum _{j=1}^m \alpha _{ij}s_j +\beta _i,\ i=1,\dots ,n.
\end{equation}

We can always suppose $m\leqslant n$ and at least for one equation
$\Omega _v$ $m=n$, because otherwise the solution set of the
irreducible system $U_{\rm com}=1$ would be represented as a union
of a finite number of proper subvarieties.

We will show now that all the tuples $(k_1,\dots ,k_n)$ that
correspond to some system (\ref{abel}) with $m<n$ (the dimension
of the subgroup $H_v$ generated by $\bar k-\bar\beta=k_1-\beta
_1,\dots ,k_n-\beta _n$ in this case is less than $n$), appear
also in the union of  systems (\ref{abel}) with $m=n$. Such
systems have form $\bar k-\bar\beta _{q} \in H_q$, $q$ runs
through some finite set $Q$, and where $H_q$ is a subgroup of
finite index in $Z^n=\left< s_1 \right>\times\cdots\times \left<
s_n \right>$. We use induction on $n$.
 If for some terminal vertex $v$, the system (\ref{abel})
has $m<n$, we can suppose without loss of generality  that the set
of tuples $H$ satisfying this system  is defined by the equations
$k_r=\dots ,k_n=0$. Consider just the case $k_n=0$. We will show
that all the tuples $\bar k_0=(k_1,\dots ,k_{n-1},0)$ appear in
the systems (\ref{abel}) constructed for the other terminal
vertices with $n=m$. First, if $N_q$ is the  index of the subgroup
$H_q$, $N_q\bar k\in H_q$ for each tuple $\bar k$. Let $N$ be the
least common multiple of $N_1,\dots ,N_Q$. If a tuple $(k_1,\dots
,k_{n-1},tN)$ for some $t$ belongs to $\bar\beta _q+H_q$ for some
$q$, then $(k_1,\dots ,k_{n-1},0)\in
 \bar\beta _q+H_q$, because $(0,\dots ,0,tN)\in H_q$. Consider the
 set $K$ of all tuples $(k_1,\dots ,k_{n-1},0)$
such that $(k_1,\dots ,k_{n-1},tN)\not\in  \bar\beta _q+H_q$ for
any $q=1,\dots ,Q$ and $t\in {\mathbb{Z}}$ . The set $\{(k_1,\dots
,k_{n-1},tN) \mid (k_1,\dots ,k_{n-1},0)\in K, t\in
{\mathbb{Z}}\}$ cannot be a discriminating set for $U_{\rm comm}=1$.
Therefore it satisfies some proper equation. Changing variables
$k_1,\dots ,k_{n-1}$ we can suppose that for an irreducible
component the equation has form
 $k_{n-1}=0$. The contradiction arises from the fact that we cannot
 obtain a discriminating set for $U_{\rm comm}=1$ which does not belong to
$\bar\beta _q+H_q$ for any $q=1,\dots ,Q.$

Embeddings $\phi _k$ are given by the systems (\ref{abel}) with
$n=m$ for generalized equations $\Omega _v$ for all terminal
vertices $v$.
\end{proof}

Below we describe two  useful  constructions.  The first one is a
{\it normalization} construction which allows one to rewrite
effectively an NTQ-system $U(X) = 1$ into a  normalized NTQ-system
$U^*=1$. Suppose we have an NTQ-system $U(X)=1$ together with a
fundamental sequence of solutions which we denote $\bar V(U)$.

 Starting from the bottom we
replace each non-regular quadratic equation $S_i = 1$ which has a
non-commutative solution by a system of equations effectively
constructed as follows.

1)  If $S_i = 1$  is in  the form
 $$c_1^{x_{i1}}c_2^{x_{i2}}=c_1c_2,$$
  where $[c_1,c_2]\not =1$,  then we replace it  by a system
  $$\{\,x_{i1}=z_1 c_1z_3,\,
x_{i2}=z_2 c_2z_3,\, [z_1,c_1]=1,\, [z_2,c_2]=1,\, [z_3,c_1c_2]=1\,\}.$$

 2) If $S_i = 1$  is in  the form
 $$x_{i1}^2c^{x_{i2}}= a^2c,$$
  where $[a,c]\not = 1$, we replace it by  a system
   $$\{x_{i1}= a^{z_1}, x_{i2}=z_2 cz_1, [z_2, c]=1,[z_1,a^2c]=1\}.$$

 3)  If $S_i = 1$  is in  the form
 $$x_{i1}^2x_{i2}^2=a_1^2a_2^2$$
 then we replace it  by the system
 $$\{\,x_{i1}=(
a_1z_1)^{z_2},\, x_{i2}=(z_1^{-1} a_2)^{z_2},\, [z_1, a_1 a_2]=1,\,
[z_2,a_1^2a_2^2]=1\,\}.$$

The normalization  construction effectively provides an NTQ-system
 $U^* = 1$ such that each solution in $\bar V(U)$ can be obtained
 from a solution of $U^*=1$.
We refer to  this system as to the normalized system of $U$
corresponding to $\bar V(U)$. Similarly, the  coordinate group of
the normalized system is called the {\it normalized} coordinate
group of $U = 1$.

\begin{lm}
\label{le:7.5embed} Let $U(X) = 1$ be an NTQ-system, and $U^* = 1
$ be  the normalized system corresponding to the fundamental
sequence $\bar V(U)$. Then the following holds: \bi \item[(1)] The
coordinate group $F_{R(U)}$ canonically  embeds into $F_{R(U^*)}$;

\item[(2)] The system $U^* = 1$ is an NTQ-system of the type \bea
S_1(X_1, X_2, \dots, X_n,A) &=& 1\\
S_2(X_2, \dots, X_n,A) &=& 1\\
 &\vdots& \\
S_n(X_n,A) &=& 1 \eea in which every $S_i = 1$ is either a regular
quadratic equation or an empty equation or a system of the type $$
U_{\rm com}(X,b) \{[x_i,x_j]=1, [x_i,b]=1  \mid  i,j =1,\dots
,k\}$$ where $b \in G_{i+1}$.

\item[(3)] Every solution $X_0$ of $U(X)=1$ that belongs to the
fundamental sequence $\bar V(U)$ can be obtained from a solution
of the system $U^*=1$. \ei
\end{lm}

\begin{proof} Statement (1) follows from the normal forms of elements in free
constructions or from the fact that applying standard
automorphisms $\phi_L$ to  a non-commuting solution  (in
particular, to a basic one) one obtains a discriminating set of
solutions (see Section 7.2). Statements (2) and (3) are obvious
from the normalization construction.\end{proof}

\begin{df} A family of solutions $\Psi$ of a regular NTQ-system $U(X,A)=1$ is called {\em generic}
if for any equation $V(X,Y,A)=1$ the following is true: if for any
solution from $\Psi$ there exists a solution of
$V(X^{\psi},Y,A)=1$, then $V=1$ admits a complete $U$-lift.

A family of solutions $\Theta$ of a regular  quadratic equation
$S(X)=1$ over a group $G$  is called {\em generic} if for any
equation $V(X,Y,A)=1$ with coefficients in $G$ the following is
true: if for any solution $\theta\in\Theta$ there exists a
solution of $V(X^{\theta},Y,A)=1$ in $G$, then $V=1$ admits a
complete $S$-lift.

A family of solutions $\Psi$ of an NTQ-system $U(X,A)=1$ is called
{\em generic} if $\Psi=\Psi _1\ldots \Psi _n$, where $\Psi _i$ is
a generic family of solutions of $S_i=1$ over $G_{i+1}$ if $S_i=1$
is a regular quadratic system, and $\Psi _i$ is a discriminating
family for $S_i=1$ if it is  a system of the type $U_{\rm com}$.

\end{df}

The second construction is a {\it correcting extension of
centralizers} of a normalized NTQ-system $U(X) = 1$ relative to an
equation $W(X,Y,A) = 1$, where $Y$ is a tuple of new variables.
Let $U(X) = 1$ be an NTQ-system in the normalized form: \bea
S_1(X_1, X_2, \ldots, X_n,A) &=& 1\\
S_2(X_2, \ldots, X_n,A) &=& 1\\
&\vdots & \\
S_n(X_n,A) &=& 1. \eea So every $S_i = 1$ is either a regular
quadratic equation or an empty equation or a system of the type
$$ U_{\rm com}(X,b) =  \{[x_i,x_j]=1, [x_i,b]=1,  \mid \ i,j =1,\dots ,k\}$$ where
$b \in G_{i+1}$.  Again, starting from the bottom we find the
first equation $S_i(X_i) = 1$ which is in the form $U_{\rm
com}(X)= 1$ and replace it with a new centralizer extending system
${\bar U}_{\rm com}(X) = 1$ as follows.

We construct $T_{\rm sol}$ for the system $W(X,Y) = 1 \wedge
U(X)=1$ with parameters $X_i,\dots ,X_n$. We obtain generalized
equations corresponding to final vertices. Each of them consists
of a periodic structure on $X_i$ and generalized equation on
$X_{i+1}\ldots X_n$. We can suppose that for the periodic
structure the set of cycles $C^{(2)}$ is empty. Some of the
generalized equations have a solution over the extension of the
group $G_i$. This extension is  given by the relations $\bar
U_{\rm com}(X_i)=1, S_{i+1}(X_{i+1},\dots , X_n)=1,\dots ,
S_n(X_n)=1$, so that there is an embedding $\phi_k : A(X,b)
\rightarrow A(X,b)$.   The others provide a proper (abelian)
equation $E_j(X_i) = 1$  on $X_i$.  The argument above shows that
replacing each  centralizer extending system $S_i(X_i) = 1$ which
is in the form $U_{\rm com}(X_i) = 1$ by a new  system of the type
${\bar U}_{\rm com}(X _i) = 1$ we eventually rewrite the system
$U(X) = 1$ into  finitely many new ones ${\bar U}_1(X) = 1,
\ldots, {\bar U}_m(X) = 1 $. We denote this set of NTQ-systems by
${\mathcal C}_W (U)$. For every NTQ-system ${\bar U}_m(X) = 1  \in
{\mathcal C}_W(U)$ the embeddings $\phi_k$ described above give
rise to embeddings ${\bar \phi}:F_{R(U)} \rightarrow F_{R({\bar
U})}$. Finally, combining normalization and correcting extension
of centralizers (relative to $W = 1$)  starting with an NTQ-system
$U = 1$ and a fundamental sequence of its solutions $\bar V(U)$ we
can obtain a finite set
 $${\mathcal NC}_W(U) =   {\mathcal C}_W(U^*) $$
 which comes equipped with a  finite set of embeddings
 $\theta_{i}: F_{R(U)}  \rightarrow
 F_{R({\bar U}_i)}$ for each ${\bar U}_i \in {\mathcal NC}_W(U)$.
These embeddings are called {\em correcting normalizing
embeddings}. The construction implies the following result.

\begin{theorem}\label{tqe}
Let $U(X,A) =1$ be an NTQ-system with a fundamental sequence of
solutions $ V_{\rm fund}(U)$. If a formula
 $$\Phi=\forall X (U(X) = 1 \rightarrow \exists Y (W(X, Y,A)=1  \wedge W_1(X, Y,A)\not =1)$$
 is true in $F$. Then for every ${\bar U}_i \in {\mathcal NC}_W(U)$
 the formula
 $$\exists Y (W(X^{\theta_{i}}, Y,A)=1 \wedge W_1(X^{\theta _{i}}, Y,A)\not =1)$$
is true in the group $F_{R({\bar U}_{i})}$ for every correcting
normalizing embedding $$\theta_i:  F_{R(U)}  \rightarrow
F_{R({\bar U}_i)}.$$

 Furthermore,  for every fundamental solution $\phi: F_{R(U)} \rightarrow F$
there exists a fundamental solution $\psi$ of one of the systems
$\bar U_i=1$, where ${\bar U}_i \in {\mathcal NC}_W(U)$ such that
$\phi =\theta _i\psi .$
\end{theorem}

\begin{cy}
Theorem C holds.
\end{cy}

Now we are ready to prove Theorem D.

\medskip
\textsc{Proof of Theorem D.}
 By   \cite[Theorem 11.1]{JSJ} for a finite system of
equations $U=1$ over $F$ one can effectively find  NTQ systems
$U_i=1,\ i=1,\dots ,k$ and homomorphisms $\theta _i:
F_{R(U)}\rightarrow F_{R(U_i)}$ such that for every solution
$\phi$ of $U=1$ there exists $i$ such that
 $\phi =\theta _i\psi$, where $\psi\in V_{\rm fund}(U_i).$  Now the result follows from
 Theorem C.
\hfill $\Box$

\section{Groups that are elementary equivalent to a free group}
\label{se:eefg}
 In this section we  prove Theorem E from the introduction.

Let $\mathcal C$ (${\mathcal C}^\ast$)  be the class  of finite
systems $U(X) = 1$ over $F$ such that every equation $T(X,Y) = 1$
compatible with $U(X) = 1$ admits $U$-lift (complete $U$-lift). We
showed in Section \ref{se:lift}, Lemma \ref{le:rat-equiv},  that
these classes are closed under rational equivalence.  Denote by
$\mathcal K$ the class of the coordinate groups $F_{R(U)}$ of
systems $U(X) = 1$ over $F$ such that every equation $T(X,Y) = 1$
over $F$ compatible with $U(X)  = 1$ admits a $U$-lift. It follows
 that every finite set of
defining relations of a group from $\mathcal K$ gives rise to a
system from $\mathcal C$ .

By Theorem B the class $\mathcal K$ contains the coordinate groups
of regular NTQ systems.

Below, in the case of a coefficient-free system $S(X)=1$ we put
$G_{cfR(S)}=F(X)/R(S)$, then $G_{R(S)}=G\ast G_{cfR(S)}.$ In this
case the group  $G_{cfR(S)}$ can be also viewed as the coordinate
group of $V(S)$. It is usually clear from the context which groups
is considered in the case of the coefficient-free system.

\begin{lm} \label{le:retract}
The class $\mathcal K$ is closed under retracts. Namely, if $H$ is
a finitely generated
 subgroup of $G$ such that there exists
a retract $\phi :G\rightarrow H$. Then:
 \begin{enumerate}
 \item  if $F \leqslant H$ then  $H = F_{R(U)}$ for some system $U = 1$ over $F$ and
  every equation compatible with $U=1$ admits a $U$-lift;
   \item  if $F \cap H = 1$ then $H = F_{R(U)}$ for some coefficient-free  system $U = 1$ over $F$ and
  every coefficient-free equation compatible with $U=1$ admits a $U$-lift into $F_{cfR(S)}.$
\end{enumerate}
\end{lm}
\begin{proof}
 We show only (1),  but a similar argument proves (2).
  Let $H = \langle F \cup X_1\rangle$ be a finitely generated subgroup of $G$ generated by
 $F$ and a finite set $X_1$.
 Then $H$ is residually free, so  $H = F_{R(U)}$ for some system
  $U(X_1) = 1$ over $F$. Since $H$ is a subgroup of $G$ it follows that
  $X_1 = P(X)$ for
  some word mapping $P$.  If $T(X_1,Y)=1$ is compatible with
$U(X_1)=1$  then $T(P(X),Y) = 1$ is  compatible with $S(X)=1$.
Therefore  $T(P(X),Y) = 1$ admits an $S$-lift, so $T(P(X),V(X)) =
1$ in $G$ for some $V(X) \in G.$  It follows that
 $$T(P(X),V(X))^\phi = T(P(X)^\phi, V(X^\phi)) =  T(P(X), V(X^\phi)) = T(X_1,V(X^\phi)) = 1$$
so $T(X_1,Y) = 1$ admits a $U$-lift.
  \end{proof}

\begin{cy} \label{freep1}
The class $\mathcal K$ is closed under  free factors. Namely, if
$G \in {\mathcal K}$ then every factor in a free decomposition of
$G$ modulo $F$ belongs to ${\mathcal K}$.
\end{cy}

\medskip
{\bf Theorem E.} {\it Let $F$ be a free non-abelian group and
$S(X) = 1$ a consistent system of equations over $F$. Then the
following conditions are equivalent:

\begin{enumerate}
\item   The system $S(X) = 1$ is rationally equivalent to a
regular NTQ system.
 \item  Every equation $T(X,Y)=1$ which is
compatible with $S(X)=1$ over $F$ admits an $S$-lift.
 \item   Every
equation $T(X,Y)=1$ which is compatible with $S(X)=1$ over $F$
admits a complete $S$-lift.
\end{enumerate}
}

\begin{proof} $(1) \Longrightarrow (3)$. It follows from Lemma \ref{le:rat-equiv}
 which states that the class ${\mathcal C}^\ast$ is closed under
 rational equivalence and the fact that ${\mathcal C}^\ast$
 contains all regular NTQ systems (Theorem B).

$(3) \Longrightarrow (2)$. Obvious.

$(2) \Longrightarrow (1)$. Suppose that every equation which is
compatible with $S=1$ over $F$ admits an $S$-lift. Consider
$G=F_{R(S)}$.

\begin{lm} \label{le:ab-subgroups}
The group $G$ does not have non-cyclic abelian subgroups.\end{lm}

\begin{proof}
Suppose  $G$ has a non-cyclic abelian subgroup, let $x,y$ be two
basis elements in  this subgroup. Consider their expressions in
generators of $G$: $x=u(X)$, $y=v(X)$. Then the system of
equations
$$S_1(X,x,y)= \left(S(X)=1\wedge x=u(X) \wedge y=v(X)\wedge [x,y]=1\right)$$ is
rationally equivalent to $S(X)=1$, therefore every system of
equations compatible with $S_1(X,x,y)=1$ admits an $S$-lift. The
formula
$$\forall X \forall x \forall y \exists u (S_1(X,x,y)=1\rightarrow (u^2=x\vee u^2=y\vee u^2=xy))$$
is true in every free group, because in a free group the images of
$x,y$ are powers of the same element. But this formula is false in
$G$. Therefore the system
 $$u^2=x\vee u^2=y\vee u^2=xy$$
  does not admit an $S$-lift. This gives a contradiction to the assumption.
\end{proof}

By Corollary \ref{freep1} we may assume that $G$ is freely
indecomposable.  There are two cases to consider,   $F \leqslant
G$ and
 $F \cap G  = 1$. Since the same argument gives a proof for both
of them we consider only one case, say $F \leqslant G$.

If $G$ does not have a non-degenerate JSJ
 $\mathbb{Z}$-decomposition \cite{JSJ} then $G$ is either a
surface group, or $G$ is an infinite cyclic group (in the case $F
\cap G = 1$). In both cases $G$ is the coordinate group of a
regular  NTQ system, as required.

Suppose now, that $G$ has a non-degenerate JSJ
$\mathbb{Z}$-decomposition of $G$, say $D$. Denote by $\langle X
\mid U \rangle$  the canonical finite presentation of $G$ as the
fundamental group of the graph of groups $D$. By Lemma
\ref{le:rat-equiv} the class ${\mathcal C}$, of systems $V = 1$
over $F$ for which every compatible equation admits an $V$-lift,
is closed under rational equivalence. Hence  $U = 1$ belongs to
${\mathcal C}$.  Since $G = F_{R(U)}$ we may assume from the
beginning that $S = U$, so $G = \langle X \mid S \rangle$  is the
canonical finite presentation of $G$ as the fundamental group of
$D$.

Let  $A_E$ be the group of automorphisms ($F$-automorphisms, in
the case $F \leqslant G$) of $G$ generated by Dehn's twists along
the edges of $D$. The group $A_E$ is abelian by Lemma 2.25
\cite{JSJ}. Recall, that two solutions $\phi_1$ and $\phi_2$ of
the equation $R(X) = 1$ are {\em $A_E$-equivalent} if there is an
automorphism $\sigma \in A_E$  such that $ \sigma \phi_1 =
\phi_2$.

Recall, that if $A$ is a group of canonical  automorphisms of $G$
then the the maximal standard quotient of $G$ with respect to $A$
is  the quotient $G/R_A$ of $G$ by the intersection $R_A$ of the
kernels of all solutions of $S(X) = 1$ which are minimal with
respect to $A$ (see \cite{JSJ}  for details).

By  \cite[Theorem 9.1]{JSJ} the maximal standard quotient
$G/R_{A_D}$ of with respect to the whole group of canonical
automorphisms $A_D$ is a proper quotient of $G$, i.e., there
exists an equation $V(X) = 1$ such that $V \not \in R(S)$ and  all
minimal solution of $S(X)=1$ with respect to the canonical group
of automorphisms $A_D$ satisfy the equation $V(X)=1.$ Now, compare
this with the following result.

\begin{lm} \label{le:max-stand-AE}
 The maximal standard quotient of $G$
 with respect to the group $A_E$ is equal to $G$, i.e.,
 the set of of minimal solutions with respect
 to $A_E$ discriminates $G$.
 \end{lm}

\begin{proof} Suppose, to the contrary, that the standard minimal
quotient $G/R_{A_E}$ of $G$ is a proper quotient of $G$, i.e.,
there exists $V \in G$ such that $V \neq 1$ and $V^\phi = 1$ for
any minimal solution of $S$ with respect to $A_E$.
 Recall that the group $A_E$ is generated by Dehn twists along the edges of $D$.
 If $c_e$ is a given generator of the cyclic subgroup associated with the edge $e$,
 then we know how the Dehn twists $\sigma = \sigma_e$ associated with $e$ acts on
 the generators from the set $X$. Namely, if $x \in X$ is a
 generator of a vertex group, then
  either $x^\sigma = x$ or $x^\sigma = c^{-1}xc$.
   Similarly, if $x \in X$ is a stable letter then
  either $x^\sigma = x$ or $x^\sigma = xc$. It follows that for $x
  \in X$ one has
  $x^{\sigma^n} = x$ or $x^{\sigma^n} = c^{-n}xc^n$
  [$x^{\sigma^n} = xc^n$] for every $n
  \in \mathbb{Z}$. Now, since the centralizer of $c_e$ in $G$ is
  cyclic (Lemma \ref{le:ab-subgroups}) the following equivalence
  holds:
   $$ \exists n \in \mathbb{Z} (x^{\sigma^n}  = z)
   \Longleftrightarrow \left \{ \begin{array}{ll}
                               \exists y ([y,c_e] = 1 \wedge y^{-1}xy = z) & \mbox{if}\; x^\sigma  = c_e^{-1}x c_e;\\
                               x= z & \mbox{if}\; x^\sigma  = x.
                               \end{array} \right.
                               $$

Similarly, since the group $A_E$ is finitely generated abelian one
can write down a formula  which describes the relation
$$\exists \alpha \in A_e (x^{\alpha}  = z)$$
One can write the  elements $c_e$ as words in
  generators $X$, say $c_e = c_e(X)$. Now the formula
$$\forall X\exists Y \exists Z \left(S(X) = 1 \rightarrow  \left(\bigwedge_{i = 1}^{m}[y_i,c_i(X)]=1 \wedge
Z = X^{\sigma_Y} \wedge V(Z)=1\right)\right)$$
 holds in the group $F$. Indeed, this
formula tells one  that each solution of $S(X)=1$ is
$A_E$-equivalent  to (a minimal)  solution that satisfies the
equation  $V(X)=1$. Since $S(X)= 1$ is in $\mathcal C$ the system
 $$\left(\bigwedge_{i = 1}^{m}[y_i,c_i(X)]=1 \wedge
Z = X^{\sigma_Y} \wedge V(Z)=1\right) $$ admits an $S$-lift. Hence
there is an automorphism $\alpha \in A_E$ such that $V(X^\alpha) =
1$ in $G$, so $V(X) = 1$ -- contradiction.
\end{proof}

\begin{lm} There exist QH subgroups  in $D$.
\end{lm}

 \begin{proof}
 By Theorem 9.1 \cite{JSJ} the maximal standard quotient $G/R_{A_D}$ of
 $G$ with respect to the whole group $A_D$ of the standard
 automorphisms of $G$ is a proper quotient of $G$. Let $E_1$ be
 the set of edges between non-QH vertex groups.
 By  \cite[Lemma 2.25]{JSJ} the
 group $A_D$   is a direct product of $A_{E_1}$ and the group generated by the
  canonical
automorphisms corresponding to $QH$ vertices and abelian
non-cyclic vertex groups. By Lemma \ref{le:ab-subgroups} there are
no abelian non-cyclic groups in $D$, so $A_D$ is a direct product
of $A_{E_1}$ and the group generated by the canonical
automorphisms of $QH$ vertices. Since the maximal standard
quotient of $G$ with respect to $A_E$ is not proper (Lemma
\ref{le:max-stand-AE}) then $A_D \neq A_E$ hence (see Section 2.20
in \cite{JSJ})  $D$ has $QH$ subgroups.
\end{proof}

 Let $K=\langle X_2\rangle$ be the fundamental group
of the graph of groups obtained from $D$ by removing all $QH$
subgroups.

\begin{lm}The natural homomorphism $G\rightarrow G/R_D$  is a
monomorphism on $K$.\end{lm}

\begin{proof} This follows from Lemma \ref{le:max-stand-AE} and
the fact that canonical automorphisms corresponding to $QH$
subgroups fix $K$.\end{proof}

\begin{lm} There is a $K$-homomorphism $\phi$ from $G$ into itself with
the non-trivial kernel.\end{lm}

\begin{proof} The generating set
$X$ of $G$ corresponding to the decomposition $D$ can be partition
 as $X=X_1\cup X_2$. Consider a formula
\begin{multline*}
\forall X_1 \forall X_2\exists Y \exists T\exists Z \left(
S(X_1,X_2)=1\right.\\  \rightarrow \left(\left.\bigwedge_{i
1}^{m}[t_i,c_i(X_2)]=1 \wedge Z=X_2^{\sigma_T}\wedge
S(Y,X_2)=1\wedge V(Y,Z)=1\right)\right).
\end{multline*}
It says that each solution of the equation $S(X_1,X_2)=1$ can be
transformed by a canonical automorphism into a solution $Y,Z$ that
satisfies $V(Y,Z)=1$. It is true in a free group, therefore the
system
$$\left(\bigwedge_{i = 1}^{m}[t_i,c_i(X_2)]=1 \wedge
Z=X_2^{\sigma_T}\wedge S(Y,X_2)=1\wedge V(Y,Z)=1\right)$$ can be
lifted in $G$. Elements $Z$ generate the same subgroup $K$ as
$X_2$, because $t_i=c_i^{n_i}$, for a fixed number $n_i$,
$i=1,\dots ,m$ in $G$. Therefore, there is a proper
$K$-homomorphism $\phi$ from $G$ into itself.  \end{proof}

For a QH subgroup $Q$ we denote by $P_Q$ the fundamental group of
the graph of groups obtained from $D$ by removing the QH-vertex
$v_Q$ and all the adjacent edges. In the following lemma, the
second statement in not needed for the proof of Theorem $E$, but
we included it for completeness.

\begin{lm}\label{le:retract1}\
\bi \item[1.] There exists a QH subgroup $Q$ in $D$ such that
$P_Q$ is a retract.

\item[2.] The maximal standard quotient $G/R_{A_Q}$ of $G$,  with
respect to the group $A_Q$ of the canonical automorphisms of $G$
corresponding to $Q$, is a proper quotient of $G$. \ei
\end{lm}

\begin{proof}
1. The image $H=\phi (G)$ cannot contain conjugates of finite
index subgroups of all the QH subgroups of $D$. Indeed, suppose it
does. Let $Q_1,\dots ,Q_s$ be QH subgroups with minimal number of
free generators. There is no homomorphism from a finitely
generated free group onto a proper finite index subgroup.
Therefore the family $Q_1,\dots ,Q_s$ has to be mapped onto the
same family of QH subgroups. Similarly, the family of all QH
subgroups would be mapped onto the conjugates of subgroups from
the same family, and different QH subgroups would be mapped onto
conjugates of different QH subgroups. In this case $H$ would be
isomorphic to $G$. This is impossible because $G$ is hopfian.
Therefore there is a QH subgroup $Q$ such that $H$ does not
intersect any conjugate $Q^g$ in a subgroup of finite index.

By construction, $G$ is the fundamental group  of the graph of
groups with  vertex $v_Q$ and vertices corresponding to connected
components $Y_1,\dots ,Y_k$ of the graph for $P_Q$. Let $P_1,\dots
,P_k$ be the fundamental groups of the graph of groups on
$Y_1,\dots ,Y_k.$ Then $P_Q=P_1\ast\cdots \ast P_k$. Let $D_Q$ be
a JSJ decomposition of $G$ modulo $K$. Then it has two vertices
$v_Q$ and the vertex with vertex group $P_Q$.



By \cite[Lemma 2.13]{JSJ}  applied to $D_Q$ and the subgroup $H$,
one of the following holds:

\begin{enumerate}

 \item  $H$ is a nontrivial  free product modulo $K$;
 \item  $H \leqslant P_Q^g$ for some $g \in G$.
\end{enumerate}

Moreover, the second statement of  this lemma is the following.
 If $H_Q=H\cap Q$ is non-trivial and has infinite index in $Q$,
 then $H_Q$ is a free product of some conjugates of $p_1^{\alpha
 _1},\dots ,p_m^{\alpha _m}, p^{\alpha}$ and a free group $F_1$ (maybe trivial) which
 does not intersect any conjugate of $\langle p_i\rangle $ for
 $i=1,\dots ,m.$

 In the  case (2) one has $H  \leqslant P_Q^g$, and, conjugating, we can suppose that $H\leqslant P_Q$.

Suppose now that the case (1) holds. For any $g$ the subgroup
$Q^g\cap H$ is either trivial or has the structure described in
the second statement of Lemma 2.13, \cite{JSJ}. Consider now the
decomposition $D_H$ induced on $H$ from $D_Q$. If the group $F_1$
is nontrivial, then $H$ is freely decomposable modulo $K$, because
the vertex group $Q_H$ in $D_H$ is a free product, and all the
edge groups belong to the other factor. If at least for one
subgroup $Q^g$, such a group $F_1$ is non-trivial, then $H$ is a
non-trivial free product and the subgroup $K$  belongs to the
other factor. Hence $H=H_1*T$, where $K\in H_1$. In this case we
consider $\phi _1=\phi\psi$, where $\psi$ is identical on $H_1$
and $\psi (x)=1$ for $x\in T$. Now each non-trivial subgroup
$H_1\cap Q^g$ is a free product of conjugates of some elements
$p_i^{\alpha _i},\ \alpha _i\in Z$, in $Q^g$.

 According to the Bass-Serre theory, for the group $G$ and its
decomposition $D_Q$ one can construct a tree such that $G$ acts on
this tree, and stabilizers correspond to vertex and edge groups of
$D_Q$. Denote this Bass-Serre tree by $T_{D_Q}$. The subgroup
$H_1$ also acts on $T_{D_Q}$. Let $T_1$ be a fundamental
transversal for this action. Either  $H_1\leqslant P_Q^g$ or $H_1$
is not conjugated into $P_Q$. The amalgamated product of the
stabilizers of the vertices of $T_1$ is a free product of
subgroups $H_1\cap P_Q^g$. Therefore $H_1$ is either such a free
product  or is obtained from such a free product by a sequence of
HNN extensions with associated subgroups belonging to distinct
factors of the free product. In both cases $H_1$ is freely
decomposable modulo $K$. Conjugating, we can suppose that one of
the factors of $\phi _1(G)$ is contained in $P_Q$. We replace now
$\phi _1$ by $\phi _2$ which is a composition of $\phi _1$ with
the homomorphism identical on the factor that is contained in
$P_Q$ and sending the other free factors into the identity. Then
$\phi _2(G)=H_2\leqslant P_Q$, where $H_2$ is freely
indecomposable modulo $K$.

A mapping $\pi$ defined on the generators $X$ of $G$ as
$$ \pi (x) =\left \{ \begin{array}{ll}
                               \phi _2 (x) & \verb if \ x\in Q;\\
                               x & \verb if \ x\not\in Q
                               \end{array} \right.$$
can be extended to a proper homomorphism $\pi$ from $G$ onto
$P_Q$. Then $\pi$ is a $P_Q$-homomorphism, and $P_Q$ is a retract.

2. Let $X=X_3\cup X_4$ be a partition of $X$ such that $X_4$ are
generators of $P_Q$. Then the following formula is true in $G$
$$\forall X_3\forall X_4\exists Y (S(X_3,X_4)=1\rightarrow
(S(Y,X_4)=1\wedge Y=r(X_4))),$$ where $Y=r(X_4)=\pi (X_3).$ This
formula is also true in $F$.

For a homomorphism $\gamma: G\rightarrow F$ there are two
possibilities:

a) $\gamma$ can be transformed by a canonical automorphism from
$A_Q$ into a homomorphism $\beta: G\rightarrow F$, such that there
exists $\alpha :G\rightarrow P_Q\ast F(Z)$ and $\psi :P_Q\ast
F(Z)\rightarrow F$ such that $\beta =\alpha\psi .$ Here $F(Z)$ is
a free group corresponding to free variables of the quadratic
equation corresponding to $Q$.

b) $\gamma$ is a solution of one of the finite number of proper
equations that correspond to the cases $\gamma (Q)$ is abelian or
$\gamma (G_e)=1$, where $e$ is an edge adjacent to $v_Q$.

Since $\ker (\alpha )=\bigcap \ker (\alpha\psi )$, where $\psi\in
Hom(P_Q\ast F(Z),F)$, the statement follows.
 \end{proof}

By Lemma \ref{le:retract} the group $P=P_Q$ belongs to ${\mathcal
K}$. If $P$ is freely undecomposable [modulo $F$] and does not
have a non-degenerate JSJ decomposition [modulo $F$] then $H$ is
either $F$ or a cyclic group, or a surface group. In this event,
$G$ is a regular NTQ (since only regular quadratic equations
belong to the class ${\mathcal C}$). If $P$ is freely decomposable
modulo $F$ or it has a non-degenerate JSJ decomposition we put
$G_0 = G$,  $Q_0 = Q$ and repeat the argument above to the group
$G_1 = P$. Thus, by induction we construct a sequence of proper
epimorphisms:
 $$G \rightarrow G_1 \rightarrow G_2 \rightarrow \ldots $$
and a sequence of $QH$ subgroups $Q_i$ of the groups $G_i$ such
that $G_i$ is the fundamental group of the graph of groups with
two vertices $Q_i$ and $G_{i+1}$ and such that $Q_i$ is defined by
a regular quadratic equation $S_i = 1$ over $G_{i+1}$ and such
that $S_i = 1$ has a solution in $G_{i+1}$. Since free groups are
equationally  Noetherian this sequence terminates in finitely many
steps either  at a surface group, or the free group $F$, or an
infinite cyclic group. This shows that the group $G$ is
$F$-isomorphic to a coordinate group of some regular NTQ system.

This proves the theorem. \end{proof}

As a corollary one can obtain the following result. To explain we
need few definitions. Let $F$ be a free group and $L_F$ be a group
theory language with constants from the group $F$, and $\Phi$ be a
set of first order sentences of the language $L_F$.  Recall, that
two groups $G$ and $H$ are $\Phi$-equivalent if they satisfy
precisely the same sentences from the set $\Phi$. In this event we
write $G\equiv_\Phi H$. In particular, $G\equiv_{\forall \exists}
H$ ($G\equiv_{ \exists \forall} H$)  means that $G$ and $H$
satisfy precisely the same $\forall \exists$-sentences ($exists
\forall$-sentences ). Notice that $G\equiv_{\forall \exists} H
\Longleftrightarrow G\equiv_{\exists \forall } H$.  We have shown
in \cite{KMIrc} that for a finitely generated group $G$
$G\equiv_{\forall \exists} H$ implies that $G$ is torsion-free
hyperbolic. Now we can prove Theorem F from the introduction:

\medskip
{\bf Theorem F.} {\em Every finitely generated group which is
$\forall \exists$-equivalent to a free non-abelian group $F$ is
isomorphic to the coordinate group of a regular NTQ system over
$F$.}

\medskip
\textsc{Proof of Theorem F.} Let $G$ be a finitely generated group
which is $\forall \exists$-equivalent to a free non-abelian group
$F$. In particular, $G$ is $\forall$-equivalent to $F$, hence by
Remeslennikov's theorem \cite{R1} the group $G$ is fully
residually free. It follows then that $G$ is the coordinate  group
of some irreducible system $S=1$ over $F$ (see \cite{BMR1}), so $G
= F_{R(S)}$. We claim that every equation compatible with $S(X)=1$
admits an $S$-lift over $F$. Indeed, if $T(X,Y) = 1$ is compatible
with $S(X) = 1$ over $F$ then the formula $$\forall X\exists Y
(S(X)=1\rightarrow T(X,Y)=1)$$ is true in $F$, hence  in $G$.
Therefore, the equation $T(X^\mu,Y)=1$ has a solution in $G$ for
any specialization of variables from $X$ in $G$, in particular,
for the canonical generators $X$ of $G$. This shows that every
equation compatible with $S= 1$ admits $S$-lift. By Theorem E, the
group $G$ is isomorphic to the coordinate group of a regular NTQ
system, as required. \hfill$\Box$

\end{document}